\documentclass[12pt,a4paper,twoside]{article}
\usepackage[LGR,T1]{fontenc} %"LGR" permet d'utiliser des caractères \text... avec "\texorpdfstring"
\usepackage[utf8]{inputenc}
\usepackage[english]{babel} %for hyphenation
\usepackage{titlesec}
\usepackage{marginnote}
\usepackage[hyperfootnotes=false,unicode,psdextra]{hyperref}
\usepackage{multicol,amsmath,amsfonts}
\usepackage{indentfirst} %\frenchspacing
\usepackage{color}

\usepackage{lineno}%[mathlines]%\usepackage{lineno,fnlineno}
\pagewiselinenumbers%On lance la numérotation des lignes
\newcommand{\llabel}[1]{\hypertarget{llineno:#1}{\linelabel{#1}}}
\newcommand{\lref}[1]{\hyperlink{llineno:#1}{\ref*{#1}}}
\newcommand{\lpageref}[1]{\hyperlink{llineno:#1}{\pageref*{#1}}}

\makeatletter

\def\ps@headings{\let\@mkboth\@gobbletwo
\def\@oddhead{\ifnum\value{page}=1\firstheadline
	\else\rightheadline\rm\thepage\fi}
\def\@oddfoot{}
\def\@evenhead{\ifnum\value{page}=1\firstheadline
	\else\rm\thepage\leftheadline\fi}
\def\@evenfoot{}
}
\newtoks\literat
\def\[#1 #2\par{\literat={#2\unskip.}%
\hbox{\vtop{\hsize=.1\hsize\noindent [#1]\hfill}%
\vtop{\hsize=.9\hsize
\noindent\the\literat}}\par
\vskip.3\baselineskip}
\newdimen\authorswidth
\newdimen\maxauthorswidth
\newskip\vspacebetweenbibitems
\makeatletter
\newcount\@refno
\def\nextref#1#2#3#4{\advance\@refno\@ne
\if@filesw \immediate\write\@auxout
	{\string\bibcite{#1}{\number\@refno}}\fi\ignorespaces
	\[{\number\@refno} #2, {#3}, #4\par
	\global\setbox5=\hbox{#2}
	\ifnum\wd5>\maxauthorswidth\authorswidth=\maxauthorswidth
	\else\authorswidth=\wd5\fi
	\def\sameauthor{\rule{\authorswidth}{0.5truept}}}

\def\trivlist{\bigbreak\vskip-\parskip
 \@trivlist \labelwidth\z@ \leftmargin\z@
 \itemindent\z@ \def\makelabel##1{##1}}
\def\@thmcounter#1{\noexpand\arabic{#1}}
\def\@thmcountersep{}
\def\@begintheorem#1#2{\it \trivlist \item[\hskip
\labelsep{\bf #1\ #2.\quad}]}
\def\@opargbegintheorem#1#2#3{\it \trivlist
	\item[\hskip \labelsep{\bf #1\ #2.\quad{\rm #3}}]}

\makeatother
\newcommand{\hlabel}{\phantomsection\label} %On remplace "\label" hors (sous-)sections par "\hlabel"
\titleformat{\part}[display]{\slshape\bfseries\boldmath}{}{0pt} {\centering\thepart}
\titlespacing{\part}{0pt}{10pt}{10pt}

\makeatletter
\renewcommand{\section}{
	\@startsection{section}{1}{\z@}{0.8truecm}{\medskipamount}
	{\centering\bfseries\boldmath}}
\makeatother
\newtheorem{Theoreme}{Theorem}[section]%\newtheorem{Theoreme}{Théorème}[section]
\newtheorem{Proposition}[Theoreme]{Proposition}
\newtheorem{Definition}[Theoreme]{Definition}%\newtheorem{Definition}[Theoreme]{Définition}
\newtheorem{Lemme}[Theoreme]{Lemma}%\newtheorem{Lemme}[Theoreme]{Lemme}
\newtheorem{Remarque}[Theoreme]{Remark}%\newtheorem{Remarque}[Theoreme]{Remarque}

\newenvironment{Dem}[1]
	{\begin{trivlist}\item[\hskip\labelsep{\bf#1.\quad}]}{\end{trivlist}}

\newtheorem{Conventions}{Conventions}[part]
\newtheorem{Hypothese}{Hypothesis\!\!} 
\newtheorem{résultats}[Hypothese]{Results\!\!}

\def\thebibl#1#2#3
{\renewcommand{\thepart}{}
\part*{{\normalsize\bf #2}}
\list
	{[\arabic{enumi}]}{\settowidth\labelwidth{[#1]}\leftmargin\labelwidth
	\advance\leftmargin\labelsep
	\setlength{\itemsep}{0cm}
	\usecounter{enumi}\setcounter{enumi}{#3}}
	\def\newblock{\hskip .11em plus .33em minus -.07em}
	\sloppy
	\sfcode`\.=1000\relax}

\newenvironment{index des notations}[1]
 {\renewcommand{\thepart}{}\begin{multicols}{4}[\part*{\slshape Index of notation}\smallskip]
 \begin{trivlist}%\addcontentsline{toc}{part}{Index of notations}
 }
 {\end{trivlist}\end{multicols}}
\def\labind#1{\index{#1@\protect#1}\hlabel{#1}}

\newskip\cqfdpostskipamount\cqfdpostskipamount=8pt plus 3pt minus 1.5pt
\def\cqfdpostskip{\vskip\cqfdpostskipamount}

\def\cqfd{\nolinebreak\quad\nolinebreak\hfill
	\vbox{\hrule
	\hbox to 6pt{\vrule height 5,2pt \hfil \vrule}
	\hrule}
	\cqfdpostskip}
\def\cqfdpartiel{\nolinebreak\quad\nolinebreak\hfill
	\hbox{\ldots}
	\cqfdpostskip}
\def\cqfr{\nolinebreak\quad\nolinebreak\hfill
	\rule{4pt}{4pt} \cqfdpostskip}
\newcommand{\cala}{\mathcal{A}}
\newcommand{\calb}{\mathcal{B}}
\newcommand{\calc}{\mathcal{C}}
\newcommand{\cald}{\mathcal{D}}
\newcommand{\cale}{\mathcal{E}}
\newcommand{\calf}{\mathcal{F}}
\newcommand{\calh}{\mathcal{H}}
\newcommand{\call}{\mathcal{L}}
\newcommand{\calo}{\mathcal{O}}
\newcommand{\calr}{\mathcal{R}}
\newcommand{\cals}{\mathcal{S}}
\newcommand{\calv}{\mathcal{V}}
\newcommand{\calw}{\mathcal{W}}
\def\scalo{{\scriptstyle \calo}}

\newcommand{\CC}{\mathbb{C}}
\newcommand{\GG}{\mathbb{G}}
\newcommand{\MM}{\mathbb{M}}
\newcommand{\NN}{\mathbb{N}}
\newcommand{\Proj}{\mathbb{P}}
\newcommand{\RR}{\mathbb{R}}
\newcommand{\ZZ}{\mathbb{Z}}

\let\goth\mathfrak
\newcommand{\gotha}{\mathfrak{a}}
\newcommand{\gothb}{\mathfrak{b}}
\newcommand{\gothc}{\mathfrak{c}}
\newcommand{\gothg}{\mathfrak{g}}
\newcommand{\gothh}{\mathfrak{h}}
\newcommand{\gothj}{\mathfrak{j}}
\newcommand{\gothk}{\mathfrak{k}}
\newcommand{\gothm}{\mathfrak{m}}
\newcommand{\gothn}{\mathfrak{n}}
\newcommand{\gothq}{\mathfrak{q}}
\newcommand{\gotht}{\mathfrak{t}}
\newcommand{\gothu}{\mathfrak{u}}
\newcommand{\gothw}{\mathfrak{w}}
\def\Sum{\displaystyle\sum}
\def\Sumpetit{\sum\limits}
\def\Prod{\displaystyle\prod}
\def\Prodpetit{\prod\limits}
\def\Int{\displaystyle\int}
\def\Lim{\displaystyle\lim}
\def\Opluspetit{\mathop{\oplus}\limits}

\def\Re{\mathop{\mathrm{Re}}\nolimits}
\def\Im{\mathop{\mathrm{Im}}\nolimits}
\def\sg{\mathop{\mathrm{sg}}\nolimits}
\def\pr{\mathop{\mathrm{pr}}\nolimits}
\def\id{\mathop{\mathrm{id}}\nolimits}
\def\abs#1{{\left\vert#1\right\vert}}
\def\rg{\mathop{\mathrm{rk \,}}\nolimits}%\def\rg{\mathop{\mathrm{rg \,}}\nolimits}
\def\tr{\mathop{\mathrm{tr}}\nolimits}
\def\Ker{\mathop{\mathrm{Ker}}\nolimits}
\def\Sp{\mathop{\mathrm{Sp}}\nolimits}

\def\moins0{\mathchoice{\setminus\!\{0\}}{\setminus\!\{0\}}
	{\setminus\{0\}}{\setminus\{0\}}}

\def\egdef{:=}

\def\ddotsc{\smash{\raisebox{-1ex}{$\ddots$}}}
\makeatletter
\def\diagd#1#2{
\mathchoice
{\begin{array}{cc} {#1}&\\&{#2} \end{array}}
{\fontsize{\f@size}{\ssf@size}\selectfont
\begin{array}{@{\,}c@{}c@{\,}} {#1}&\\&{#2} \end{array}}
{\fontsize{\sf@size}{\ssf@size}\selectfont
\begin{array}{@{\,}c@{}c@{\,}} {#1}&\\&{#2} \end{array}}
{\fontsize{\ssf@size}{\ssf@size}\selectfont
\begin{array}{@{\,}c@{}c@{\,}} {#1}&\\&{#2} \end{array}}
}
\def\diagt#1#2#3{
\mathchoice
{\begin{array}{ccc} {#1}&&\\&{#2}&\\&&{#3} \end{array}}
{\fontsize{\f@size}{\ssf@size}\selectfont
\begin{array}{@{\,}c@{}c@{}c@{\,}} {#1}&&\\&{#2}&\\&&{#3} \end{array}}
{\fontsize{\sf@size}{\ssf@size}\selectfont
\begin{array}{@{\,}c@{}c@{}c@{\,}} {#1}&&\\&{#2}&\\&&{#3} \end{array}}
{\fontsize{\ssf@size}{\ssf@size}\selectfont
\begin{array}{@{\,}c@{}c@{}c@{\,}} {#1}&&\\&{#2}&\\&&{#3} \end{array}}
}
\def\diagq#1#2#3#4{
\mathchoice
{\begin{array}{cccc}
{#1}&&&\\&{#2}&&\\&&{#3}&\\&&&{#4} \end{array}}
{\fontsize{\f@size}{\ssf@size}\selectfont
\begin{array}{@{\,}c@{}c@{}c@{}c@{\,}}
{#1}&&&\\&{#2}&&\\&&{#3}&\\&&&{#4} \end{array}}
{\fontsize{\sf@size}{\ssf@size}\selectfont
\begin{array}{@{\,}c@{}c@{}c@{}c@{\,}}
{#1}&&&\\&{#2}&&\\&&{#3}&\\&&&{#4} \end{array}}
{\fontsize{\ssf@size}{\ssf@size}\selectfont
\begin{array}{@{\,}c@{}c@{}c@{}c@{\,}}
{#1}&&&\\&{#2}&&\\&&{#3}&\\&&&{#4} \end{array}}
}
\def\diags#1#2#3#4#5#6{
\mathchoice
{\begin{array}{cccccc}
{#1}&&&&&\\&{#2}&&&&\\&&{#3}&&&\\&&&{#4}&&\\&&&&{#5}&\\&&&&&{#6} \end{array}}
{\fontsize{\f@size}{\ssf@size}\selectfont
\begin{array}{@{\,}c@{}c@{}c@{}c@{}c@{}c@{\,}}
{#1}&&&&&\\&{#2}&&&&\\&&{#3}&&&\\&&&{#4}&&\\&&&&{#5}&\\&&&&&{#6} \end{array}}
{\fontsize{\sf@size}{\ssf@size}\selectfont
\begin{array}{@{\,}c@{}c@{}c@{}c@{}c@{}c@{\,}}
{#1}&&&&&\\&{#2}&&&&\\&&{#3}&&&\\&&&{#4}&&\\&&&&{#5}&\\&&&&&{#6} \end{array}}
{\fontsize{\ssf@size}{\ssf@size}\selectfont
\begin{array}{@{\,}c@{}c@{}c@{}c@{}c@{}c@{\,}}
{#1}&&&&&\\&{#2}&&&&\\&&{#3}&&&\\&&&{#4}&&\\&&&&{#5}&\\&&&&&{#6} \end{array}}
}
\makeatother

\makeatletter
\def\Rot#1{
\mathchoice
{\begin{array}{cc} 0 & -{#1} \\ {#1} & 0 \end{array}}
{\fontsize{\f@size}{\ssf@size}\selectfont
\begin{array}{@{\,}c@{\;}c@{\,}} 0 & -{#1} \\ {#1} & 0 \end{array}}
{\fontsize{\sf@size}{\ssf@size}\selectfont
\begin{array}{@{\,}c@{\;}c@{\,}} 0 & -{#1} \\ {#1} & 0 \end{array}}
{\fontsize{\ssf@size}{\ssf@size}\selectfont
\begin{array}{@{\,}c@{\;}c@{\,}} 0 & -{#1} \\ {#1} & 0 \end{array}}
}
\def\Rotneg#1{
\mathchoice
{\begin{array}{cc} 0 & {#1} \\ -{#1} & 0 \end{array}}
{\fontsize{\f@size}{\ssf@size}\selectfont
\begin{array}{@{\,}c@{\;}c@{\,}} 0 & {#1} \\ -{#1} & 0 \end{array}}
{\fontsize{\sf@size}{\ssf@size}\selectfont
\begin{array}{@{\,}c@{\;}c@{\,}} 0 & {#1} \\ -{#1} & 0 \end{array}}
{\fontsize{\ssf@size}{\ssf@size}\selectfont
\begin{array}{@{\,}c@{\;}c@{\,}} 0 & {#1} \\ -{#1} & 0 \end{array}}
}
\makeatother

\def\sym#1{\mathfrak{S}_{#1}}
\def\interieur{\mathop{\mathrm{int}}\nolimits}
\def\centregr#1{\mathrm{Z}(#1)}
\def\centrealg#1{\mathrm{Z}(#1)}
\def\derivegr#1{\mathrm{D}(#1)}
\def\derivealg#1{\mathrm{D}(#1)}
\def\centragr#1#2{\mathrm{C}_{#1}(#2)}
\def\centraalg#1#2{\mathrm{C}_{#1}(#2)}
\def\normagr#1#2{\mathrm{N}_{#1}(#2)}
\def\normaalg#1#2{\mathrm{N}_{#1}(#2)}
\newcommand{\me}{\mathrm{e}} %base des logarithmes népériens
\newcommand{\mi}{\mathrm{i}} %racine carrée particulière de -1

\def\Log{\mathop{\mathrm{Log}}\nolimits}
\def\Supp{\mathop{\mathrm{Supp}}\nolimits}
\def\diff{\mathrm{d}}
\def\norm#1{{\left\Vert#1\right\Vert}}

\def\restriction#1#2{\mathchoice
	{\setbox1\hbox{${\displaystyle #1}_{\scriptstyle #2}$}
	\restrictionaux{#1}{#2}}
	{\setbox1\hbox{${\textstyle #1}_{\scriptstyle #2}$}
	\restrictionaux{#1}{#2}}
	{\setbox1\hbox{${\scriptstyle #1}_{\scriptscriptstyle #2}$}
	\restrictionaux{#1}{#2}}
	{\setbox1\hbox{${\scriptscriptstyle #1}_{\scriptscriptstyle #2}$}
	\restrictionaux{#1}{#2}}}
\def\restrictionaux#1#2{{#1\,\smash{\vrule height .8\ht1 depth
.85\dp1}}_{\,#2}}

\def\fonctioncar#1{\mbox{\rm 1\hspace{-.25em}l}_{#1}}
\def\ad{\mathop{\mathrm{ad}}\nolimits}
\def\Ad{\mathop{\mathrm{Ad}}\nolimits}
\def\Ind{\mathop{\mathrm{Ind}}\nolimits}
\def\Car{\mathop{\mathrm{Car}}\nolimits}

\def\ssreg#1{#1_{\scriptscriptstyle\mathit{ss\,reg}}}
\def\ssreggr#1{#1_{\!\scriptscriptstyle\mathit{ss\,reg}}}

\def\chiinfty#1#2{\chi_{#1}^{U \smash{{#2}_\CC}}}
\def\F+{{\calf^+}}
\def\Fh+{{\calf^+_{\!\gothh}}}
\def\Fp+{{{{\calf}'}^+}}
\def\Fep+{{{{\calf}_{\!e}'}^+}}
\def\Fphp+{{{\calf_{\!\gothh'}'}^+}}
\def\XInd{X^{\textit{Ind}}}
\def\Xirr{X^{\textit{irr}}}
\def\Xirrp{X^{\textit{irr,+}}}
\def\Xfin{X^{\textit{final,+}}}

\let\wt\widetilde
\def\wtlie#1#2{\wt{{#1}^{\vrule width 0ex height .4ex \smash{*}}}_{\!\!\!#2}}
\def\lambdatilde{{\mathchoice
	{\wt{\lambda}}
	{\wt{\lambda}}
	{\tilde{\lambda}}
	{\tilde{\lambda}}}}
\def\classe#1{\mathop{#1}\limits^.}
\def\mutilde{{\mathchoice
	{\wt{\mu}}
	{\wt{\mu}}
	{\tilde{\mu}}
	{\tilde{\mu}}}}

\def\gstregtilde{\wtlie{\gothg}{\textit{reg}}}
\def\gstregGtilde{\wtlie{\gothg}{\textit{reg},G}}
\def\gestregtilde{\wtlie{\gothg(e)}{\textit{reg}}}
\def\gejstregtilde{\wtlie{\gothg(e_j)}{\textit{reg}}}
\def\mpestregtilde{\wtlie{\gothm'(e)}{\textit{reg}}}

\def\gstssG{\gothg^*_{\textit{ss},G}}%\def\gstssG{\gothg^*_{\textit{ss,G}}
\def\gstfondtilde{\wtlie{\gothg}{\textit{fund}}}%\def\gstfondtilde{\wtlie{\gothg}{\textit{fond}}}
\def\gstfondGtilde{\wtlie{\gothg}{\textit{fund},G}}%\def\gstfondGtilde{\wtlie{\gothg}{\textit{fond},G}}
\def\mpstfondMptilde{\wtlie{\gothm'}{\textit{fund},M'}}%\def\mpstfondMptilde{\wtlie{\gothm'}{\textit{fond},M'}}

\def\gstssI{\gothg^*_{\textit{ssI}}}
\def\gstItilde{\wtlie{\gothg}{\textit{I}}}
\def\gestItilde{\wtlie{\gothg(e)}{\textit{I}}}
\def\gezerostItilde{\wtlie{\gothg(e_0)}{\textit{I}}}
\def\gstssIG{\gothg^*_{\textit{ssI},G}}
\def\gstIGtilde{\wtlie{\gothg}{\textit{I},G}}
\def\mstIMtilde{\wtlie{\gothm}{\textit{I},M}}

\def\gstssInc{\gothg^*_{\textit{ssInc}}}
\def\gstInctilde{\wtlie{\gothg}{\textit{Inc}}}
\def\gestInctilde{\wtlie{\gothg(e)}{\textit{Inc}}}
\def\gstssIncG{\gothg^*_{\textit{ssInc},G}}
\def\gstIncGtilde{\wtlie{\gothg}{\textit{Inc},G}}

\def\gstreg{\gothg^*_{\textit{reg}}}
\def\glzerostreg{\gothg(l_0)^*_{\textit{reg}}}
\def\glambdastreg{\gothg(\lambda)^*_{\textit{reg}}}

\def\gstss{\gothg^*_{\textit{ss}}}
\def\gestss{\gothg(e)^*_{\textit{ss}}}
\def\gstssreg{\ssreg{\gothg^*}}
\def\gestssreg{\ssreg{\gothg(e)^*}}
\def\mpstssreg{\ssreg{{\gothm'}^*}}

\def\gssreg{\ssreg{\gothg}}
\def\gessreg{\ssreg{\gothg(e)}}
\def\mpessreg{\ssreg{\gothm'(e)}}
\def\Gssreg{\ssreggr{G}}
\def\Gzerossreg{\ssreggr{(G_0)}}
\def\MpGzerossreg{\ssreggr{(M'G_0)}}
\def\Mpssreg{\ssreggr{M'}}

\makeatletter
\def\cases#1{\left\{\vcenter{\normalbaselines\m@th
	\ialign{$##\hfil$&\quad##\hfil\crcr#1\crcr}}\,\right.}
\def\@oddhead{\ifnum\value{page}=137\firstheadline
	\else\rightheadline\rm\thepage\fi}
\makeatother
\begin{document}
\linenumbers

\setcounter{page}{137}

\newcommand{\rightheadline}{\hfill\sc Ducloux\hfill}
\newcommand{\leftheadline}{\hfill\sc Ducloux\hfill}
\newcommand{\firstheadline}{\begin{minipage}[t]{9cm}
\emph{Translation from}
Journal of Lie Theory\hfill\\
{\bf 12}\ (2002)\ 137-190\hfill\\
\copyright\ 2002\ Heldermann Verlag\hfill
\end{minipage}\hfill}

\renewcommand\thefootnote{}
\footnote{\vskip-1mm\normalsize ISSN: 0949-5932}

%===============================================================================
%
%Firstpage
%
%==============================================================================

\vbox{\vskip5mm}

{\font\bfone=ecbx1000 scaled\magstep2
\bfone
\centerline{The Orbit Method and Character Formulas}
\vskip1.5mm
\centerline{for Tempered representations of a Nonconnected}
\vskip1.5mm
\centerline{Reductive Real Algebraic Group$^{(\star)}$}
}
%
%\newlength{\largeurnote}\setlength{\largeurnote}{\linewidth}\addtolength{\largeurnote}{-7mm}%
\footnote{%
%Les lignes rectifiées sont signalées par "\llabel{}"
\vskip-3mm\noindent
\parbox[t]{\linewidth}{
\reversemarginpar\marginnote{$(\star)\!\!$}
In this translation written in December 2025, the following lines contain corrections:\\
simplification of notation pages \lpageref{gstssG} and \pageref{RepIndex} (“$\gothg^*_{\textit{ssfund},G}\!$”$\,\mapsto\,$“$\gstssG\!$”),
line \lref{page217} p.$\!\!$~\lpageref{page217} (“p.$\!\!$~221”$\,\mapsto\,$“p.$\!\!$~217”),
lines \lref{OrbitesCoadjointes1}--\lref{OrbitesCoadjointes2} p.$\!\!$~\lpageref{OrbitesCoadjointes1} (brief ideas rewritten),
line \lref{chapeau} p.$\!\!$~\lpageref{chapeau} (“$\beta$”$\,\mapsto\,$“$\widehat{\beta}$” twice),
lines \lref{PbFp} and \lref{PbFpSuite} p.$\!\!$~\lpageref{PbFp} and line \lref{PbFpDém} p.$\!\!$~\lpageref{PbFpDém} (we do not ask $\Fp+\!$ to be in $C(\gothg(\lambda),\gothh)_{\textit{reg}}$),
line \lref{theta} p.$\!\!$~\lpageref{theta} (“$\theta_q$”$\,\mapsto\,$“$\theta_n$”),
line \lref{RacinesCompactes} p.$\!\!$~\lpageref{RacinesCompactes} (“$R^+(\gothm_\CC,\gotht_\CC)$”$\,\mapsto\,$“$R(\gothm_\CC,\gotht_\CC)$”),
line \lref{a->c} p.$\!\!$~\lpageref{a->c} (“(a)”$\,\mapsto\,$“(c)”),
lines \lref{Mackey1} and \lref{Mackey2} p.$\!\!$~\lpageref{Mackey1} (“is a multiple of”$\,\mapsto\,$“has a common multiple with”),
line \lref{d->e} p.$\!\!$~\lpageref{d->e} (“(d)”$\,\mapsto\,$“(e)”).%
}%
}

\centerline{\bf Jean-Yves Ducloux}
\vskip5mm\rm

\centerline{\rm Communicated by J. Faraut}
\vskip7mm

%===============================================================================
%
%Abstract
%
%==============================================================================
\renewenvironment{abstract}{\narrower\footnotesize\bf
 \noindent\abstractname.\quad\footnotesize\rm}{\par\bigskip}

\begin{abstract}
%\nonfrenchspacing %\begin{english}
Let $G$ be a possibly disconnected reductive real Lie group.
In this paper, I parametrize the set of irreductible tempered characters of $G$.
I~then describe these characters using certain “Kirillov's formulas,”
based on the descent method near each elliptic element in $G$.

If $G$ is linear and connected, the parameters I use are “final basic”
parameters in the sense of Knapp and Zuckerman (see \cite[p.$\!\!$~453] {KZ82}).
%\frenchspacing %\end{english}
\end{abstract}
\vskip3mm

%===============================================================================
%
\centerline{\slshape\bfseries Table of contents}
%
%==============================================================================
\def\ldotfill{\leaders\hbox spread6pt{\hss .\hss}\hfill}

\vskip3mm\noindent
{{\sl Introduction and general notation}
\ldotfill{\sl \pageref{RepIntro}}}%
\medskip

\noindent
{{\sl I. The “linear form” parameters }
\ldotfill{\sl \pageref{RepI}}}%

\hspace{-2em}
{1.\quad The parameters $\,\lambdatilde \in \gstregGtilde$
\ldotfill\pageref{Rep1}}%

\hspace{-2em}
{2.\quad The measures $\beta_{G \cdot \lambdatilde}$
\ldotfill\pageref{Rep2}}%

\hspace{-2em}
{3.\quad Fixed points of an elliptic element
\ldotfill\pageref{Rep3}}%
\medskip

\noindent
{{\sl II. The “projective representation” parameters}
\ldotfill{\sl \pageref{RepII}}}%

\hspace{-2em}
{4.\quad About special-metalinear and metaplectic groups
\ldotfill\pageref{Rep4}}%

\hspace{-2em}
{5.\quad The parameters $\,\tau \in \XInd_G(\lambdatilde)$
\ldotfill\pageref{Rep5}}%

\hspace{-2em}
{6.\quad Integrability condition
\ldotfill\pageref{Rep6}}%
\medskip

\noindent
{{\sl III. Construction of representations}
\ldotfill{\sl \pageref{RepIII}}}%

\hspace{-2em}
{7.\quad The case where $G$ is connected
\ldotfill\pageref{Rep7}}%

\hspace{-2em}
{8.\quad The representations $T_{\lambdatilde,{\gotha^*}^+\!,\tau_+}^G\!$
\ldotfill\pageref{Rep8}}%

\hspace{-2em}
{9.\quad The injection
$\,G \cdot (\lambdatilde,\tau)\mapsto T_{\lambdatilde,\tau}^G$
from $G\,\backslash\,\XInd_G$ to $\widehat{G}$
\ldotfill\pageref{Rep9}}%
\medskip

\noindent
{{\sl IV. Characters of representations}
\ldotfill{\sl \pageref{RepIV}}}%

\hspace{-2em}
{10.\quad Character restriction formula
\ldotfill\pageref{Rep10}}%

\hspace{-2em}
{11.\quad Transition from $\,T_{\lambdatilde,\tau_{M'}}^{M'}\!$ to
$\,T_{\lambdatilde,{\gotha^*}^+\!,\tau_+}^G\!$
\ldotfill\pageref{Rep11}}%

\hspace{-2em}
{12.\quad Translation in the sense of Zuckerman
\ldotfill\pageref{Rep12}}%

\hspace{-2em}
{13.\quad Transition from \smash{$\,T_{\lambda_r,\sigma_r}^{M'}\!$ to
$\,T_{\lambdatilde,\tau_{M'}}^{M'}\!$}
\ldotfill\pageref{Rep13}}%
\medskip

\noindent
{{\sl Index of notation}
\ldotfill{\sl \pageref{RepIndex}}}%
\medskip

\noindent
{{\sl References}
\ldotfill{\sl \pageref{RepRef}}}%
%\vfill\eject
\bigskip

%===============================================================================
%
\part*{Introduction and general notation}\label{RepIntro}
%
%==============================================================================

In this paper, I will describe the tempered dual of a real reductive Lie group using the orbit method, based on the following four works:

\noindent\quad
-- the description of the tempered dual of a connected real reductive Lie group by Adams, Barbasch, and Vogan in \cite[Ch. 11]{ABV92};

\noindent\quad
-- the parametrization by Duflo of the equivalence classes of representations of a real reductive Lie group that are irreducible tempered with a regular infinitesimal character, obtained via “Mackey theory” in \cite[III]{Df82a};

\noindent\quad
-- the “Kirillov's character formula” of Bouaziz in \cite{Bo87}, which expresses the character of these representations in terms of the Fourier transforms of regular semisimple coadjoint orbits;

\noindent\quad
-- Rossmann's formulas in \cite{Ro82} (where “principal” corresponds to “regular” and “regular” to “regular semisimple”), which relate the Fourier transforms of regular orbits to those of regular semisimple orbits.
\medskip

A theorem of Duflo in \cite[p.$\!\!$~189] {Df82b} reduces the classification of unitary duals of real algebraic linear groups to the classification of unitary duals of “almost algebraic reductive groups with finite kernels.”
It is therefore these reductive groups that interest me.

%------------------------------------------------------------------------------
\begin{Hypothese}
%------------------------------------------------------------------------------

Let $G$\labind{G} be a real Lie group with Lie algebra $\gothg$\labind{g}, a linear algebraic group $\GG$\labind{GG} defined over $\RR$, and a Lie group homomorphism from $G$ to $\GG(\RR)$ with a finite central kernel whose image is open in the usual topology.
Thus, every element of $G$ (respectively $\gothg$) admits a real Jordan decomposition into elliptic, positively hyperbolic, and unipotent components (respectively infinitesimally elliptic, hyperbolic, and nilpotent components), as described in \cite[top of p.$\!\!$~36 and Lem. 31 p.$\!\!$~38] {DV93}.
Assume that $\gothg$ is reductive with center consisting of semisimple elements.%
\end{Hypothese}
\medskip

Duflo parameterized in \cite[Lem. 8 p.$\!\!$~173] {Df82a} a subset of the unitary dual $\widehat{G}$\labind{G^} of $G$ (specified above), in terms of regular semisimple coadjoint orbits, by reducing via induction to the discrete series representations of certain subgroups of~$G$.
When $G$ is connected, this subset of $\widehat{G}$ was described by Harish-Chandra (see \cite[Th. 1 p.$\!\!$~198] {Ha76}).
Here, I parameterize, in Theorem \ref{injection dans le dual}~(b) p.~\!\pageref{injection dans le dual}, in terms of \emph{generalized regular coadjoint orbits}, a larger subset of $\widehat{G}$, namely the tempered dual of $G$, obtained by replacing, at the initial step of the induction, the discrete series representations with limits of discrete series representations.
When~$G$ is connected, my statement reproduces part of Theorem 11.14 de \cite[p.$\!\!$~131] {ABV92} (see also \cite[Th. 14.76 p.$\!\!$~598] {Kn86}) with the following addition: the equality of two representations corresponds exactly to the conjugacy, under $G$, of the associated pairs of \emph{good parameters}.
\medskip

The expression “good parameters” refers to the constraints imposed on these parameters.
The following is the heuristic principle (see \cite[Th. 1 p.$\!\!$~193]{Df82a} and \cite[Th. 19 p.$\!\!$~211]{Df82b}).
An irreducible unitary representation $T$ of $G$ will, in favorable cases, be parameterized by the orbit under $G$ of an ordered pair~$(\lambdatilde,\tau)$.
The term $\lambdatilde$ must be a kind of linear form whose semisimple component~$l$ corresponds to the character $\chiinfty{\mi\,l}{\gothg}$ through which the set $(U \gothg_\CC)^G$ of $G$-invariant elements of $U\gothg_\CC$ acts on the space of $C^\infty$ vectors of the representation space of $T$ (see \ref{caractères canoniques} (a)).
In this paper, this first parameter is related to a limit of discrete series, of which $T$ will, roughly speaking, be an induced representation.
For this reason, it will be written in the form $\lambdatilde = (\lambda,\F+)$, where $\F+\!$ is a Weyl chamber for imaginary roots (see \cite[Th. 5.7 p.$\!\!$~305] {Zu77}).
Its “semisimple component” is the linear form $l\egdef\lambda$, in a sense consistent with the definitions relevant to the case of coadjoint orbits (see \ref{bijection entre orbites}~(a)).
Let $G(\lambdatilde)$ denote the stabilizer of $\lambdatilde$ in $G$.
The term $\tau$ is to be a projective representation of $G(\lambdatilde)/G(\lambdatilde)_0$, lifted to a unitary representation of a degree $2$ covering of $G(\lambdatilde)$, such that $\,\diff_1\tau=\mi\,l\id\,$ and the central character of $\tau$ extends that of $T$ (see \ref{caractères canoniques} (b)).
In this paper, this second parameter is a unitary representation (possibly reducible) of a suitable analogue for $\lambdatilde$ of \emph{Duflo’s double cover} of the stabilizer $G(f)$ in $G$ of an element $f$ of $\gothg^*\!$.
In fact, the constructions developed here will use a pair of parameters $((\lambdatilde,{\gotha^*}^+\!),\tau_+)$ instead of $(\lambdatilde,\tau)$, where ${\gotha^*}^+$ is a certain Weyl chamber for restricted roots, and $\tau_+$ is an irreducible unitary representation of a covering of $\normagr{G(\lambdatilde)}{{\gotha^*}^+}$.
This new pair of parameters still satisfies the conditions described above.
The pair $(\lambdatilde,\tau)$ will turn out to be “more canonical” once the character of $T$ has been calculated, but it will play only a secondary role.
The $G$-orbits of $(\lambdatilde,{\gotha^*}^+)$ and $((\lambdatilde,{\gotha^*}^+),\tau_+)$ are respectively identified with those of $\lambdatilde$ and $(\lambdatilde,\tau)$.
In the statement of results below, I will use the set $\XInd_G$ of pairs $(\lambdatilde,\tau)$.
It will be defined more precisely in Part~\ref{RepII}, in connection with the notion of “final character” (see \ref{paramètres adaptes} (c)).
\smallskip

The main result of this paper is Theorem \ref{formule du caractère} p.~\!\pageref{formule du caractère} that describes the character of tempered representations $T$ of $G$ in terms of Fourier transforms of canonical measures on certain regular coadjoint orbits associated with a parameter $\lambdatilde$ attached to $T$.
These character formulas show the need to adapt the “nilpotent orbits” perspective of the orbit method, disregarding the “coadjoint orbit” aspect of $G \cdot \lambdatilde$.
For example, the coadjoint orbits associated with the two limits of the discrete series representations of $SL(2,\RR)$, which correspond to the two open nilpotent half-cones, have no fixed point under the action of the “rotation” $\Ad^* \!\left( {\scriptstyle \Rot{1}} \right)$ (see Remark \ref{insuffisance de la méthode des orbites} (3));
they must instead be replaced by the sets of half-lines contained in the interior of their respective convex hulls, up to multiplication by $\mi $.
Remark \ref{autre insuffisance de la méthode des orbites} will also show that this new viewpoint clarifies why certain “admissible” nilpotent orbits in the sense of Duflo do not correspond to any representation.
\medskip

My two theorems can be summarized as follows (see the index of notation).%

%------------------------------------------------------------------------------
\begin{résultats}
%------------------------------------------------------------------------------

There exists a “canonical” bijection
$\;G \cdot (\lambdatilde,\tau)\mapsto T_{\lambdatilde,\tau}^G\;$ from
$\,G\,\backslash\,\XInd_G\,$ onto the tempered dual of $G$.
It is determined by character formulas of the form:\\[1mm]
\hspace*{\fill}%
$k_e(X)\; \tr T_{\lambdatilde,\tau}^G (e \exp X)
\;\;=\;\;
\Sum_{\substack
{\classe{\lambdatilde'} \,\in\,
G(e)\backslash \,G \cdot \lambdatilde \, \cap \, \gstregtilde(e)\,\\
\textrm{such that } \,\lambdatilde'\![e] \,\in\, \gestregtilde}}\!\!\!
c_{\widehat{e'},\lambdatilde}\;
\tr\tau(\widehat{e'})\;
\widehat{\beta}_{G(e)\cdot \lambdatilde'\![e]} (X)$
\hspace*{\fill}\\[1mm]
where
$\;(\lambdatilde,\tau)\in \XInd_G\!$,
$e$ ranges over the set of elliptic elements in $G$,\\
$\,X\,$ is an element of a certain neighborhood ${\cal V}_{\!e}$ of $0$ in $\mathfrak{g}(e)$ such that $\;e \exp X$ is regular semisimple in $G$ (in which case $X$ is regular semisimple in $\mathfrak{g}(e)\!$),\\
$\;\widehat{e'}$ is an element of the double cover
$\,G(\lambdatilde)^{\mathfrak{g}/\!\mathfrak{g}(\lambda)(\mi \rho_\F+)}\,$
(with $\lambdatilde = (\lambda,\F+)$) that projects to $\,e' = g^{-1}eg\,$ for a $\,g \in G\,$ such that $\;\lambdatilde' = g\lambdatilde$,\\
and, $\,k_e(X)\,$ and $\,c_{\widehat{e'},\lambdatilde}$ are certain nonzero complex numbers
(see $\!(\boldsymbol{F})\!$ p.$\!\!$~\pageref{F}).
\end{résultats}
\medskip

The ideas that enabled me to obtain these results are as follows.
Duflo's constructions, carried out step by step, yield a bijection from a set of orbits onto the tempered dual of $G$, by means of \emph{Lemma \ref{résultat inattendu} (b)} p.~\!\pageref{résultat inattendu}.
To move from the connected case to the disconnected case, I will need a generalization of a result of Vogan used by Duflo in a homological construction.
This can be found on page 555 of the book \cite{KV95} by Knapp and Vogan.
Given the connected case (see \cite[p.$\!\!$~64] {Ro80}), it is then natural to seek to recover the characters of the irreducible tempered representations by taking a limit from those associated with the regular semisimple orbits.
In the disconnected case, the calculation of the character of the representations will be broken down into two steps, the first of which is rather cumbersome.
\medskip

\noindent\quad
\emph{Step 1.}
We induce a limit of discrete series representation of a subgroup of $G$ to the neutral component of the Levi component $M'$ (which is generally outside the Harish-Chandra class) of a certain “parabolic subgroup” $M'U$ of $G$. The analogue of this induced representation in the connected case was simply a limit of discrete series representation of the Levi component $M$ of a cuspidal parabolic subgroup $MAN$ of $G$.
The method here, as in the connected case, is to apply Zuckerman's translation functor (adapted to the case of the disconnected group $M'$) to reduce to representations of $M'$ whose character is known.
One difficulty is that the move from $M'_0$ to $M'$involves a homological step.
On the advice of Bouaziz, I drew inspiration from points (i) and (ii) on page 550 of his paper \cite{Bo84} to calculate the action of a group $M'(\lambdatilde)^{\gothm'\!/\gothh}$ on a weight space of the homology of a $M'_0$-module translated in the sense of Zuckerman.
\medskip

\noindent\quad
\emph{Step 2.}
We then induce from $M'U$ to $G$.
The groundwork will have been laid in Part~\ref{RepI} with a review and extension of Rossmann's results for calculating the limits of Fourier transforms of canonical measures on regular semisimple orbits.
In particular, in Section~\ref{Rep3}, I will have highlighted the effect of the passage from $G$ to $G(e)$ in terms of the parameters $\lambdatilde$.
The rest of this step consists in applying the methods and results of Bouaziz.
\bigskip

The following definitions and notations will allow me, first, to specify my general conventions, and second, to introduce the conventions related to the reductive groups that I will use most frequently.

%------------------------------------------------------------------------------
\begin{Conventions}\rm
\hlabel{conventions générales}
%------------------------------------------------------------------------------

{\bf (a)}
Let $\abs{A}$ denote the cardinality of a set $A$ and let $\dot{a}$ denote the class of an element $a$ in $A$ modulo an equivalence relation $\sim$ on $A$, denote by $\,\restriction{f}{B'}: B' \to C\,$ the restriction of a map $\,f: B \to C\,$ to a subset $B'$ of $B$,
by $\,\overline{\!F}$ the closure of a subset $F$ in a topological space $E$ (in Sections \ref{Rep2} and \ref{Rep12}), by $l^\CC$ the complexification of a real linear map $l$.
Let $\overline{w}$ and $\overline{W}$ denote the conjugate of a vector $w$ and of a vector subspace $W$ in the complexification $V_\CC$ of a real vector space $V$ (in Sections \ref{Rep1}, \ref{Rep3}, \ref{Rep4}, \ref{Rep5}, \ref{Rep8} and \ref{Rep9}).
Let $u_Y$ and $u_{X/Y}$ denote the endomorphisms induced on $Y$ and $X/Y$ by an endomorphism $u$ of a vector space $X$ leaving invariant a vector subspace $Y$ of $X$.
\smallskip

{\bf (b)}
Let $A$ be a real Lie group (we will always assume that its topology is Hausdorff with a countable base) and let $\gotha$ be its Lie algebra.
%When $A$ is a topological group, $\Bbb A$ is a linear algebraic group defined over~$\RR$ (and therefore a metric space that is the countable union of compact sets), and $A \to {\Bbb A}(\RR)$ is a Lie group morphism whose image is open under the usual topology and with a finite central kernel that is discrete, the topology of $A$ is Hausdorff with a countable base.
Let $1$ denote the~neutral element of $A$,
and $A_0$ the neutral component of $A$.
Denote by~$\exp_A$ (or $\exp$) the exponential map from $\gotha$ to $A$,
by $\Ad^*$ and $\ad^*$ the coadjoint represen\-tations of $A$ and of $\gotha$,
by $\interieur\gotha$\labind{int(a)} the subgroup of the linear group of $\gotha$ generated by the elements $\exp(\ad X)$ as $X$ ranges over $\gotha$,
by $\centregr{A}$ and $\derivegr{A}$ the center~and the derived subgroup of $A$,
by $\centrealg{\gotha}$ and $\derivealg{\gotha}$ the center and the derived Lie algebra of $\gotha$,
by $\centrealg{U \gotha_\CC}$ the center of the enveloping algebra $U \gotha_\CC$ of~$\gotha_\CC$,
and by $\widehat{A}$ the unitary dual of $A$.
To every Lie subalgebra $\gothb$ of $\gotha$, we associate the Lie~subgroups\\
\hspace*{\fill}%
$\centragr{A}{\gothb}
= \{ x \in A \mid \restriction{(\Ad x)}{\gothb} = \id \}\;$
and
$\;\normagr{A}{\gothb}
= \{ x \in A \mid \Ad x \cdot \gothb \subseteq \gothb \}$
\hspace*{\fill}\\
of $A$, with Lie algebras\\
\hspace*{\fill}%
$\centraalg{\gotha}{\gothb}
= \{ X \in \gotha \mid \restriction{(\ad X)}{\gothb} = 0 \}\;$
and
$\;\normaalg{\gotha}{\gothb}
= \{ X \in \gotha \mid \ad X \cdot \gothb \subseteq \gothb \}$.
\hspace*{\fill}%
\smallskip

{\bf (c)}
Let $M$ be a smooth manifold of dimension $m$ that is Hausdorff and second-countable.
A~$C^\infty$ density on $M$ is a complex Radon measure $\rho$ on $M$ that, in any chart of $M$ centered at a point $x$, takes the form $\,a\,\diff x_1 \!\cdots \diff x_m$, where $\,a\,$ is a $C^\infty$ function defined on an open neighborhood of $0$ in $\RR^m\!$.
It is “identified” with the family of maps $\rho(x)$ with $x \in M$, which send an element of $\bigwedge^m T_x M\,\moins0$ whose image under the chart is $\,\omega_0\,$ to $\,a(0)\,\abs{(\diff x_1 \wedge \cdots \wedge \diff x_m)(\omega_0)}$.
A~generalized function on $M$ is a continuous linear form on the space of compactly supported $C^\infty$ densities on $M$, equipped with the Schwartz topology.
\smallskip

{\bf (d)}
Let $A$ be a real Lie group and $\,\diff_A$ be a left Haar measure on $A$.
A~continuous representation $T$ of $A$ in a complex Hilbert space is said to be trace class if the operators $\;T(\varphi \, \diff_A)$ with $\,\varphi \in C^\infty_c(A)\,$ are trace class;
in this case
$\;\tr T : \varphi \, \diff_A \mapsto \tr T(\varphi \, \diff_A)\;$
is a generalized function on $A$.
For any continuous unitary representation $\pi$ of a closed subgroup $B$ of $A$ equipped with a left Haar measure $\,\diff_B$, let $\Ind_B^A \pi$ denote the unitary representation of $A$ “induced” from $\pi$ as in \cite[p.$\!\!$~99] {B.72}.
\smallskip

{\bf (e)}
Let $V$ be a finite-dimensional real vector space.
The Fourier transform of a tempered distribution $\mu$ on $V^*$ is the generalized function $\widehat{\mu}$ on $V$ defined by the equality
$\;\; \widehat{\mu}(v) = \Int_{V^*} \me ^{\mi \,l(v)} \diff \mu(l)\;\;$
of generalized functions in $v \in V$.
\end{Conventions}

%------------------------------------------------------------------------------
\begin{Conventions}\rm
\hlabel{conventions fréquentes}
%------------------------------------------------------------------------------

{\bf (a)}
Let $\Car \gothg$\labind{Carg} denote the set of Cartan subalgebras of $\gothg$.
Let $\gothh \in \Car \gothg$.
Fix a system of positive roots $R^+(\gothg_\CC,\gothh_\CC)$\labind{R+(gC,hC)}, arbitrarily (unless otherwise indicated), in the set $R(\gothg_\CC,\gothh_\CC)$\labind{R(gC,hC)} of roots of $\gothh_\CC$ in $\gothg_\CC$.

Denote by $W(G,\gothh)$\labind{W(G,h)} the finite group $\,\normagr{G}{\gothh}/\centragr{G}{\gothh}$,
by $\gotht$\labind{t} and $\gotha$\labind{a} the infinitesimally elliptic and hyperbolic components of $\gothh$,
$\,\gothh_{(\RR)} \!= \mi \, \gotht \oplus \gotha$\labind{hr},
$T_0=\exp \gotht$\labind{T0} and $A=\exp \gotha$\labind{A}.
Let $\rho_{\gothg,\gothh}$\labind{ρ g,h} denote the half sum of the elements of $R^+(\gothg_\CC,\gothh_\CC)$,
and let $H_{\alpha}$\labind{Hα} denote the coroot of $\alpha$ in the root system $R(\gothg_\CC,\gothh_\CC)$ (such that $\alpha(H_\alpha)=2$).

A root $\,\alpha \in R(\gothg_\CC,\gothh_\CC)\,$ is said to be complex (respectively real, imaginary, or compact) when its conjugate
$\;\overline{\alpha}: X \in \gothh_\CC\mapsto \overline{\alpha(\overline{X})}\;$
satisfies $\,\overline{\alpha} \notin \{\alpha,-\alpha\}\,$ (respectively $\,\overline{\alpha}=\alpha$, $\,\overline{\alpha}=-\alpha$, or
$\;(\CC H_{\alpha} \oplus \gothg_\CC^\alpha \oplus \gothg_\CC^{-\alpha})\cap \gothg \simeq {\goth su}(2)$).
\smallskip

{\bf (b)}
Fix a nondegenerate $G$-invariant bilinear form
${\scriptstyle \langle} ~,\!~ {\scriptstyle \rangle}$\labind{<,>ev}
on $\gothg$ whose complexification (again denoted by ${\scriptstyle \langle} ~,\!~ {\scriptstyle \rangle}$) restricts to a scalar product on each $\gothh_{(\RR)}$, $\gothh \in \Car \gothg$.
Let
$\,{\scriptstyle \perp}^{\!{\scriptscriptstyle \langle} \!~,\!~ {\scriptscriptstyle \rangle}}\,$
denote the relation of ${\scriptstyle \langle} ~,\!~ {\scriptstyle \rangle}$-orthogonality.
Let $x \in G$ and $X \in \gothg$ be semisimple.
Let
$G(x)$\labind{G(x)} (resp. $G(X)$\labind{G(X)})
and
$\gothg(x)$\labind{g(x)} (resp. $\gothg(X)$\labind{g(X)})
denote the centralizers of $x$ (resp. $X$) in $G$ and $\gothg$.
The $G(x)$-module $\gothg(x)^*$ (resp. the $G(X)$-module $\gothg(X)^*$) is canonically isomorphic, via restriction, to the set $\gothg^*(x)$\labind{g*(x)} (resp. $\gothg^*(X)$\labind{g*(X)}) of elements of $\gothg^*$ fixed under $\Ad ^*x$ (resp. annihilated by $\,\ad^*X$), i.e., vanishing everywhere on
$\,\gothg(x)^{{\scriptstyle \perp}^{\!{\scriptscriptstyle \langle} \!~,\!~ {\scriptscriptstyle \rangle}}}\!$
(resp. $\,\gothg(X)^{{\scriptstyle \perp}^{\!{\scriptscriptstyle \langle} \!~,\!~ {\scriptscriptstyle \rangle}}}\!$).
\smallskip

{\bf (c)}
Let $f \in \gothg^* \!$.
Denote by $G(f)$\labind{G(f)} the stabilizer of $f$ in $G$ and by $\gothg(f)$\labind{g(f)} the Lie algebra of $G(f)$.
Whenever possible, $\mu$\labind{μ}, $\nu$\labind{ν}, $\lambda = \mu\!+\!\nu$\labind{λ}, and $\xi$\labind{ξ} will be used to denote the “infinitesimally elliptic, hyperbolic, semisimple, and nilpotent components of $f$” (concepts from $\gothg$ using ${\scriptstyle \langle} ~,\!~ {\scriptstyle \rangle}$).
We identify $\gothg(\lambda)^*$ with the set of elements of $\gothg^*$ that are vanishing everywhere on
$\,\gothg(\lambda)^{{\scriptstyle \perp}^{\!{\scriptscriptstyle \langle} \!~,\!~ {\scriptscriptstyle \rangle}}}\!$.
\smallskip

{\bf (d)}
Let $\gothm$ be the centralizer in $\gothg$ of a semisimple element $X$ of $\gothg$.
Given Haar measures $\diff_{\gothg}$ and $\diff_{\gothm}$ on $\gothg$ and $\gothm$, we set for all $f \in \gothm^*$:\\
\hspace*{\fill}%
$\bigl|\Pi_{\gothg,\gothm}\bigr|(f)
= \frac{1}{k!}\;\bigl|
\bigl( \restriction{B_f}{[\gothg,X]^2} \bigr)^{\!k}
(\omega)\bigr|
\,{\scriptstyle \times}
\left( \diff_{[\gothg,X]} (0) (\omega)\right)^{\!-1}\;$\labind{Π g,m}
(see (b) and Conventions \ref{conventions générales} (c)),%
\hspace*{\fill}\\
with
$k = \frac {1}{2} \dim [\gothg,X]$,
$B_f = f([\cdot,\cdot])$,
$\omega \in \bigwedge ^{2k}[\gothg,X]\moins0$,
and
$\,\displaystyle \diff_{[\gothg,X]} = \diff_{\gothg} / \diff_{\gothm}$.
\end{Conventions}
\medskip

The other notations used in the statements can be found in the index of notations at the end of this paper.

%===============================================================================
%
\part{.\quad The “linear form” parameters}\label{RepI}
%
%==============================================================================

The orbit method, first proposed by Kirillov in the case of simply connected nilpotent Lie groups, consists of parameterizing the irreducible unitary representations of a Lie group using the orbits of its coadjoint representation.
When the representation is trace class, its character in some neighborhood of the neutral element must be related to the Fourier transform (assuming it exists) of the canonical measure on the associated orbit.

A slight modification of this point of view will be necessary here.
Each irreducible tempered representation of $G$ will be associated with an orbit $\wt{\Omega}$ of $G$ in a certain set $\gstregtilde$ other than $\gothg^* \!$, and to this orbit $\wt{\Omega}$ will be attached a finite sum of canonical measures on coadjoint orbits.
In preparation for an upcoming paper, I also introduce below parts $\gstItilde$ and $\gstInctilde$ of $\gstregtilde$.
\medskip

In this section, I will make use of certain results of Rossmann.

%-------------------------------------------------------------------------------
\section{The parameters \texorpdfstring{$\,\lambdatilde \in \gstregGtilde$}{\~\lambda \in \~g*\_\{reg,G\}}}\label{Rep1}
%------------------------------------------------------------------------------

%\noindent
In this section, we construct a set $\gstregtilde$ equipped with an action of $G$ in which the set of $\lambda \!\in \gothg^*\!$ regular semisimple embeds in a $G$-equivariant manner.

%------------------------------------------------------------------------------
\begin{Definition}\rm
%------------------------------------------------------------------------------

Let $\gstreg$\labind{g*reg} denote the set of $f \!\in \gothg^*\!$ regular (i.e., the set of $f \!\in \gothg^*\!$ for which the dimension of $\gothg(f)$ is equal to the rank of $\gothg$),
$\gstss$\labind{g*ss} the set of $f \!\in \gothg^*\!$ semisimple and
$\,\gstssreg = \gstss \cap \gstreg$\labind{g*ssreg}.

Also let $\gstssI$\labind{g*ssI} (respectively: $\gstssInc$\labind{g*ssInc}) denote the set of $\lambda \! \in \gstss$ such that $\gothg(\lambda)$ has a Cartan subalgebra $\gothh$ for which the roots of $(\gothg(\lambda)_\CC,\gothh_\CC)$ are imaginary (respectively: noncompact imaginary).
\smallskip

Let $\lambda \! \in \gstss$ and $\gothh \in \Car\gothg(\lambda)$.
Let $C(\gothg(\lambda),\gothh)$\labind{Cgh} denote the set of (open) chambers in $\,\gothh_{(\RR)}^{~~*}$ for the imaginary roots of $(\gothg(\lambda)_\CC,\gothh_\CC)$ and let
$C(\gothg(\lambda),\gothh)_{reg}$\labind{Cghreg}
denote the set of $\F+\! \!\in C(\gothg(\lambda),\gothh)$ such that the imaginary roots of $(\gothg(\lambda)_\CC,\gothh_\CC)$ that are simple relative to $\F+\!$ are noncompact.

Set
\hspace*{\fill}%
$\gstregtilde
\,=\, \Bigl\{
(\lambda,\F+)\,;
\lambda \!\in\! \gstss, \,
\gothh \!\in\! \Car\gothg(\lambda)\textrm{ and }
\F+\! \!\in C(\gothg(\lambda),\gothh)_{reg}
\Bigr\}$\labind{g*reg tilde},
\hspace*{\fill}\\
\hspace*{\fill}%
and
$\quad\gstregGtilde
\,=\, \Bigl\{
(\lambda,\F+)\in \gstregtilde \mid
\forall Z \in \Ker \exp_{T_0}\;\,
\me ^{(\mi \,\lambda+\rho_{\gothg,\gothh})(Z)}=1
\Bigr\}\;\;$\labind{g*reg,G tilde}
(see \ref{conventions fréquentes} (a))
\hspace*{\fill}\\
independently of the choice of a system of positive roots $R^+(\gothg_\CC,\gothh_\CC)$ associated with the element $\gothh$ of $\Car\gothg(\lambda)$ attached to $\F+\!$.

Also set (in the following “fundamental” means “without real root”):\\[0.5mm]
\hspace*{\fill}%
$\gstfondtilde
= \Bigl\{
(\lambda,\F+)\,;
\lambda \!\in\! \gstss, \,
\gothh \!\in\! \Car\gothg(\lambda)\textrm{ fundamental and }
\F+\! \!\in C(\gothg(\lambda),\gothh)_{reg}
\Bigr\}$\labind{g*fond tilde},
\hspace*{\fill}\\
\hspace*{\fill}%
and
$\quad\gstssG
\,=\, \Bigl\{ \lambda \!\in\! \gstss \mid
\forall Z \in \Ker \exp_{T_0}\;
\me ^{(\mi \,\lambda+\rho_{\gothg,\gothh})(Z)}=1
\Bigr\}$\labind{g*ss,G}\llabel{gstssG}
\hspace*{\fill}\\
independently of the choice of a fundamental $\gothh \in \Car\gothg(\lambda)$
and of a system of positive roots $R^+(\gothg_\CC,\gothh_\CC)$,
$\,\gstssIG = \gstssI \cap \gstssG$\labind{g*ssI,G} and
$\,\gstssIncG = \gstssInc \cap \gstssG$\labind{g*ssInc,G}.%
\smallskip

Finally, let $\gstItilde$\labind{g*I tilde},
$\gstInctilde$\labind{g*Inc tilde}, $\gstfondGtilde$\labind{g*fond,G tilde},
$\gstIGtilde\!$\labind{g*I,G tilde}, and $\gstIncGtilde\!$\labind{g*Inc,G tilde}
denote the inverse images of
$\gstssI$, $\gstssInc$, $\gstssG$, $\gstssIG\!$, and $\gstssIncG\!$
under the first projection from $\gstfondtilde\!$ to $\gstss$.%
\end{Definition}
%\smallskip

%------------------------------------------------------------------------------
\begin{Remarque}\rm
%------------------------------------------------------------------------------

{\bf (1)}
By applying \cite[Th. p.$\!\!$~217] {Ro82} to $\gothg(\lambda)$ for three Cartan subalgebras (one equal to $\gothh$, one without imaginary roots, and one fundamental) and taking into account \cite[Lem. {\sc a} $(a)\!\Rightarrow\! (c)$ p.$\!\!$~220 and Suppl. {\sc c} p.$\!\!$~218] {Ro82},we see that the images of the canonical maps from $\gstregtilde$ and $\gstfondtilde$ to $\gothg^*$ consist of those $\lambda \! \in \gstss$ such that $\gothg(\lambda)$ has a Cartan subalgebra without imaginary roots.%
\medskip

{\bf (2)}
According to \cite[Th. 6.74 p.$\!\!$~341 and Th. 6.88 p.$\!\!$~344] {Kn96}, each complex semisimple Lie algebra has, up to isomorphism, a unique real form $\gothg_0$ that has a Cartan subalgebra $\gothh_0$ for which the root system $\,R({\gothg_0}_\CC,{\gothh_0}_\CC)\,$ has a basis consisting of noncompact imaginary roots.
Given (1) above (or \cite[Pb. 18 p.$\!\!$~369 and p.$\!\!$~557] {Kn96}), the image of the canonical map from $\gstItilde$ to $\gothg^*$ consists of those $\lambda \! \in \gstss$ such that each of the simple ideals of $\gothg(\lambda)$ is isomorphic to one of the following Lie algebras:\\
$\;\goth{su}(p,p)$ with $p \geq 1$ and $\goth{su}(p,p-1)$ with $p \geq 2$,
$\goth{so}(p,p-1)$ with $p \geq 3$,
$\goth{sp}(2n,\RR)$ (denoted $\goth{sp}(n,\RR)$ in \cite{Kn96})
with $n \geq 3$,
$\goth{so}(p,p)$ with $p$ even $\geq 4$
and $\goth{so}(p,p-2)$ with $p$ even $\geq 6$,
$E \, I\!I$, $E \, V$, $E \, V\!I\!I\!I$, $F \, I$, $G$.
\medskip

{\bf (3)}
The sum of two noncompact imaginary roots is never a noncompact root according to \cite[(6.99) p.$\!\!$~352] {Kn96}. Thus, $\gstssInc$ consists of those semisimple $\lambda \! \in \gothg^*$ such that the simple ideals of the reductive Lie algebra $\gothg(\lambda)$ are isomorphic to $\,\goth{sl}(2,\RR)$.
\cqfr
\end{Remarque}
%\smallskip

Starting from certain enriched parameters $(\lambdatilde,{\gotha^*}^+)$, I will canonically construct two systems of positive roots for $(\gothg_\CC,\gothh_\CC)$, which will have the same noncomplex roots.
The system of positive roots $R^+(\gothg_\CC,\gothh_\CC)$ is attached to a certain regular linear form $\lambda_+$.
It will allow, by limiting processes (Section~\ref{Rep11}) and translation (Section~\ref{Rep13}), to construct representations in Part~\ref{RepIII} by reducing to the case of regular infinitesimal character.
The system of positive roots $R^+_{\lambdatilde,{\gotha^*}^+}\!$ has, as a subset of complex roots with negative conjugate, a part that is invariant under the action of the group $G(\lambdatilde)$ defined below.
It will be used in Definition \ref{paramètres adaptes} (c) and in Lemma \ref{résultat inattendu} (b).
We will see in Lemma \ref{nouveaux paramètres} another parametrization in which the parameter $(\lambdatilde,{\gotha^*}^+)$ will be replaced by~$\lambdatilde$, and the role of $(\smash{R^+_{\lambdatilde},{\gotha^*}^+},\lambda_+)$ will be played by a certain pair $(\smash{R^+_{\lambdatilde}},\lambda_{\textit{can}})$ obtained from $\lambdatilde$.

%------------------------------------------------------------------------------
\begin{Definition}\rm
\hlabel{choix de racines positives}
%------------------------------------------------------------------------------

Let $\lambda \! \in \gstss$, $\gothh \in \Car\gothg(\lambda)$ and $\,\F+\! \!\in C(\gothg(\lambda),\gothh)$ be given.
Let $\mu$ and $\nu$ (resp. $\gotht$ and $\gotha$) denote the infinitesimally elliptic and hyperbolic components of $\lambda$ (resp. $\gothh$).
Set $\lambdatilde = (\lambda,\F+)$\labind{λtilde}.
\smallskip

{\bf (a)}
Let $G(\lambdatilde)$\labind{G(λtilde)} denote the normalizer of $\F+\!$ in $\,G(\lambda)$ and let $\gothg(\lambdatilde)$\labind{g(λtilde)} denote its Lie algebra.
Associate with $\lambdatilde$ the half sum $\rho_\F+\! \in \mi\gotht^*$\labind{ρ F+} of the imaginary $\,\alpha \in R(\gothg(\lambda)_\CC,\gothh_\CC)\,$ such that $\,\F+(H_{\alpha})\subseteq \RR_+ \!\moins0$.
(Thus $\,\gothh$ is a Cartan subalgebra of $\gothg(\lambda)(\mi \rho_\F+)$ without imaginary roots, and therefore $G \cdot (\lambda,\mi \rho_\F+)$ determines $G \cdot \lambdatilde$.)

Fix, for the rest of this definition a chamber ${\gotha^*}^+\!$ of $(\gothg(\lambda)(\mi \rho_\F+),\gotha)$.%
\smallskip

{\bf (b)}
Introduce the following system of positive roots:\\[1mm]
\hspace*{\fill}%
$R^+(\gothg(\lambda)(\mi \rho_\F+)_\CC,\gothh_\CC)
= \{ \, \alpha \! \in \! R(\gothg(\lambda)(\mi \rho_\F+)_\CC,\gothh_\CC)\mid
{\gotha^*}^+ (H_{\alpha})\subseteq \RR_+ \!\moins0 \, \}$.
\hspace*{\fill}%
\smallskip

Fix $\epsilon > 0$ small enough so that, upon setting
$\;\nu_+
= \nu + \epsilon \, \rho_{\gothg(\lambda)(\mi \rho_\F+),\gothh} \in \gotha^*$\labind{ν+}
we have:
$\;\,a \cdot \nu_+ \not= \nu_+$ for any automorphism $a$ of the root system $R(\gothg_\CC,\gothh_\CC)$ such that $a \cdot \nu \not= \nu$.
Choose $\epsilon$ so that the previous property remains valid when we replace $\epsilon$ by $t\epsilon$, $\,t \!\in\! \left]0,1\right]$.
We therefore have:
$\;\gothg(\nu_+) = \gothg(\nu)(\rho_{\gothg(\lambda)(\mi \rho_\F+),\gothh})$.
\smallskip

We will use the following systems of positive roots:\\
\hspace*{\fill}%
$R^+(\gothg(\nu_+)_\CC,\gothh_\CC)
= \left\{\, \alpha \in R(\gothg(\nu_+)_\CC,\gothh_\CC)\;\,\Big|\;\,
\mi \mu(H_{\alpha}) > 0
\textrm{ or }
\Big\{\begin{smallmatrix}
\mi \mu(H_{\alpha}) = 0 \\ \textrm{and}\hfill \\ \rho_\F+\!(H_{\alpha}) > 0
\end{smallmatrix}
\,\right\}$
\hspace*{\fill}\\
and
$\quad R^+(\gothg_\CC,\gothh_\CC)
= \{ \alpha \in R(\gothg_\CC,\gothh_\CC)\mid
\nu_+ (H_{\alpha})\! > \! 0 \}
\, \cup \, R^+(\gothg(\nu_+)_\CC,\gothh_\CC)$\labind{R+(gC,hC)noncan}.

Set
$\,\mu_+ = \mu - 2\mi \rho_{\gothg(\nu_+),\gothh} \in \gotht^*$\labind{μ+}
and
$\;\lambda_+ = \mu_+\!+\nu_+$\labind{λ+}
\\
(therefore
$\,R^+(\gothg(\nu_+)_\CC,\gothh_\CC)
= \{\, \alpha \in R(\gothg(\nu_+)_\CC,\gothh_\CC)\mid
\mi \mu_+(H_{\alpha}) > 0 \}$).
\smallskip

In some cases, we will write
$\lambda_{\gothg,\lambdatilde,{\gotha^*}^+\!,\epsilon}$,
$\mu_{\gothg,\lambdatilde,{\gotha^*}^+}\!$ and
$\nu_{\gothg,\lambdatilde,{\gotha^*}^+\!,\epsilon}$
for $\,\lambda_+\!$, $\mu_+\!$ and~$\nu_+\!$.%
\smallskip

{\bf (c)}
Introduce the set\\
\hspace*{\fill}%
$R^+_\lambdatilde
= \left\{ \alpha \in R(\gothg_\CC,\gothh_\CC)\;\Big|\;
\nu (H_{\alpha})\! > \! 0
\textrm{ or }
\Big\{ \begin{smallmatrix}
\nu(H_{\alpha}) = 0 \\ \textrm{and}\hfill \\ \mi \mu(H_{\alpha}) > 0
\end{smallmatrix}
\textrm{ or }
\Big\{\begin{smallmatrix}
\lambda(H_{\alpha}) = 0 \\ \textrm{and}\hfill \\ \rho_\F+\!(H_{\alpha}) > 0
\end{smallmatrix}
\right\}$\labind{R+lambdatilde},
\hspace*{\fill}\\
the half sum $\rho_{\textit{can}}$ of the elements of $R^+_\lambdatilde \,\cap\, R(\gothg(\nu)_\CC,\gothh_\CC)$, and the system of positive roots
$\;R^+_{\lambdatilde,{\gotha^*}^+}
= R^+_\lambdatilde \,\cup\,
R^+(\gothg(\lambda)(\mi \rho_\F+)_\CC,\gothh_\CC)\,$\labind{R+(gC,hC)can}
of $R(\gothg_\CC,\gothh_\CC)$.

Set
$\,\mu_{\textit{can}} = \mu - 2\mi \rho_{\textit{can}} \in \gotht^*$
and
$\;\lambda_{\textit{can}} = \mu_{\textit{can}}\!+\nu$\labind{λcan}\\
(so that, using \cite[Cor. 4.69 p.$\!\!$~271] {KV95}, we find that
$\,\gothg(\lambda_{\textit{can}}) = \gothg(\lambda)(\mi \rho_\F+)\,$
and
$\;R^+_\lambdatilde \,\cap\, R(\gothg(\nu)_\CC,\gothh_\CC)
= \{\, \alpha \in R(\gothg(\nu)_\CC,\gothh_\CC)\mid
\mi \mu_{\textit{can}}(H_{\alpha}) > 0 \}$).
\end{Definition}
\medskip

Here is the key lemma that will allow us to reduce $\gstregtilde$ to $\gstssreg\!$.

%------------------------------------------------------------------------------
\begin{Lemme}
\hlabel{lemme clef}
%------------------------------------------------------------------------------

We place ourselves in the situation of the previous definition.\\
Let $\,t \!\in\! \left]0,1\right]$.
Set
$\;\nu_t = \nu + t\epsilon \, \rho_{\gothg(\lambda)(\mi \rho_\F+),\gothh}$,
$\mu_t = \mu - 2\mi t\rho_{\gothg(\nu_+),\gothh}$
and
$\lambda_t = \mu_t+\!\nu_t$.

We have: $\;\,\gothg(\lambdatilde) = \gothh$, $\;\lambda_t \in \gstssreg\,$
and $\;\normagr{G(\lambdatilde)}{{\gotha^*}^+} = G(\lambda_t)$.

Therefore $\,G(\lambda_+) = G(\lambdatilde)\,$ when $\,{\gotha^*}^+ = \gotha^*$ (for example, when $\,\lambdatilde \in \gstItilde$).%
\end{Lemme}

\begin{Dem}{Proof of the lemma}

We have
$\;\,G(\lambdatilde)\subseteq \normagr{G}{\gothh}$,
therefore
$\;\gothg(\lambdatilde) = \gothh$.

By construction of $\epsilon$, we have
$\;\gothg(\nu_t)
= \gothg(\nu)(\rho_{\gothg(\lambda)(\mi \rho_\F+),\gothh})
=\gothg(\nu_+)$.
Furthermore,
$\;\mi \mu_t(H_{\alpha})
= \mi \mu(H_{\alpha}) + 2t\rho_{\gothg(\nu_+),\gothh}(H_{\alpha})
> 0\;$
holds for all $\,\alpha \in R^+(\gothg(\nu_+)_\CC,\gothh_\CC)$.
Therefore $\;\lambda_t \in \gstssreg$.
A~continuity argument in $t$ implies that:\\[1mm]
\hspace*{\fill}%
$R^+(\gothg_\CC,\gothh_\CC)
= \left\{\, \alpha \in R(\gothg_\CC,\gothh_\CC)\;\,\Big|\;\,
\nu_t(H_{\alpha}) > 0
\textrm{ or }
\Big\{\begin{smallmatrix}
\nu_t(H_{\alpha}) = 0 \\ \textrm{and}\hfill \\ \mi \mu_t(H_{\alpha}) > 0
\end{smallmatrix}
\,\right\}$.%
\hspace*{\fill}%

The inclusion $\,\normagr{G(\lambdatilde)}{{\gotha^*}^+} \subseteq G(\lambda_t)\,$ is immediate.
Let $g \in G(\lambda_t)$.\\
We have $\,g \cdot \gothh = g \cdot \gothg(\lambda_t) = \gothh\,$ and $\restriction{(\Ad g)}{\gothh}$ is an automorphism of the root system $R(\gothg_\CC,\gothh_\CC)$ such that $\,g \nu_t = \nu_t$,
therefore $\,g (\nu_t - \nu) = (\nu_t - \nu)$ then $\,g \nu_+ = \nu_+$.
Since $g$ stabilizes $R(\gothg(\nu_+)_\CC,\gothh_\CC)\cap R^+(\gothg_\CC,\gothh_\CC)$, we have:
$\,g \rho_{\gothg(\nu_+),\gothh} = \rho_{\gothg(\nu_+),\gothh}$
then $\,g \mu = \mu$.
Next, since $g$ stabilizes the set of imaginary roots of $R(\gothg(\lambda)_\CC,\gothh_\CC)\cap R^+(\gothg_\CC,\gothh_\CC)$, we have:
$\,g \F+ =\F+$.
Finally, since $g$ stabilizes the set of restrictions $\restriction{\alpha}{\gotha}$ for
$\alpha \in R(\gothg(\lambda)(\mi \rho_\F+)_\CC,\gothh_\CC)\cap R^+(\gothg_\CC,\gothh_\CC)$, we have:
$\,g \cdot {\gotha^*}^+ = {\gotha^*}^+$.
\cqfd
\end{Dem}
\medskip

The following lemma will explain why, in certain formulas that I will write later in another paper, and in which the absolute values of Pfaffians $\,\bigl|\Pi_{\gothg,\centraalg{\gothg(\mu)}{\gotha}}\bigr|$ appear as factors, only the elements of $\gstssI$ are involved.

To make these Pfaffians more amenable to calculation, I first specify a formula from \cite[p.$\!\!$~301] {DV88}. The bilinear form ${\scriptstyle \langle} ~,\!~ {\scriptstyle \rangle}$ from Conventions \ref{conventions fréquentes} (b) determines a Haar measure $\diff_V$ on any vector subspace $V$ of $\gothg$ on which its restriction is nondegenerate, by the condition
$\;\diff_V (0) (v_1 \wedge \cdots \wedge v_n)
= |\! \det ( \langle v_i,v_j \rangle )_{_{1 \leq i,j \leq n}}|^{\frac{1}{2}}\;$
for any basis $(v_1,\dots,v_n)$ of~$V$.
Let $\,\gothh \in \Car \gothg$.
Fix a regular $\lambda_0 = \mu_0 + \nu_0\in \gothh^*$ in $\gothg^*\!$.
Associate it with the system of positive roots\\
\hspace*{\fill}%
$\;R_0^+(\gothg_\CC,\gothh_\CC)
= \left\{\, \alpha \in R(\gothg_\CC,\gothh_\CC)\;\,\Big|\;\,
\nu_0(H_{\alpha}) > 0
\textrm{ or }
\Big\{\begin{smallmatrix}
\nu_0(H_{\alpha}) = 0 \\ \textrm{and}\hfill \\ \mi \mu_0(H_{\alpha}) > 0
\end{smallmatrix}
\,\right\}$.
\hspace*{\fill}\\
Let $\diff_{\gothg}$ and $\diff_{\gothh}$ be the Haar measures on $\gothg$ and $\gothh$ obtained from ${\scriptstyle \langle} ~,\!~ {\scriptstyle \rangle}$.
For all $\lambda \in \gothh^*$ and all $\,\omega_0 \in \bigwedge ^{2k}[\gothg,\gothh]\backslash \{0\}\,$ element of the orientation $\scalo (B_{\lambda_0})_{\gothg/\gothh}$ (see \ref{géométrie métaplectique} (b) and \ref{paramètres de Duflo} (a)) of $(\gothg/\gothh,B_{\lambda_0})$, we have\\[1mm]
\hspace*{\fill}%
{\bf (\boldmath$*$)}
$\;\frac{1}{k!}
\left( \restriction{B_{\lambda}}{[\gothg,\gothh]^2} \right)^{\!k}
\!(\omega_0)
\,{\scriptstyle\times}\,
\left( \diff_{[\gothg,\gothh]} (0) (\omega_0)\right)^{\!\!-1} \!\!
= \,
\mi^{n_0}
\!\!\!
\Prod_{\alpha \in R_0^+(\gothg_\CC,\gothh_\CC)} \!\!
\langle \lambda,\alpha \rangle$
\hspace*{\fill}\\
where
$\,k = \frac {1}{2} \dim [\gothg,\gothh]$,
$\,\displaystyle \diff_{[\gothg,\gothh]}
= \diff_{\gothg} / \diff_{\gothh}$,
$\,n_0
=\abs{\{\alpha \in R_0^+(\gothg_\CC,\gothh_\CC)\mid
\overline{\alpha} \notin R_0^+(\gothg_\CC,\gothh_\CC)\}}$,
denoting again by
${\scriptstyle \langle} ~,\!~ {\scriptstyle \rangle}$\labind{<>dual}
the bilinear form on $\gothh_\CC^{~*}$ dual to the restriction to
$\gothh_\CC$ of~${\scriptstyle \langle} ~,\!~ {\scriptstyle \rangle}$.

%------------------------------------------------------------------------------
\begin{Lemme}
\hlabel{description d'ensembles de formes linéaires}
%------------------------------------------------------------------------------

Let $\,f \!\in \gothg^*\,$ be of infinitesimally elliptic, hyperbolic, semisimple, and nilpotent components $\mu$, $\nu$, $\lambda=\mu\!+\!\nu$, and $\xi$.
Fix $\gothh \in \Car\gothg$ such that $\lambda \in \gothh^* \!$.
Let $\gotha$ denote the hyperbolic component of $\gothh$.
\smallskip

{\bf (a)}
We have $\;f \!\in \gstreg$ if and only if $\;\,\xi \!\in \glambdastreg$.
\smallskip

{\bf (b)}
We have $\;\lambda \in \gstssI\,$ and $\,\gothh \in \Car\gothg(\lambda)\,$ is fundamental,\\
if and only if
$\;\;\bigl|\Pi_{\gothg,\centraalg{\gothg(\mu)}{\gotha}}\bigr|(\lambda)
\not= 0\;\;$
(see Conventions \ref{conventions fréquentes} (d)).
\end{Lemme}

\begin{Dem}{Proof of the lemma}

{\bf (a)}
It is clear that the equality $\;\dim \gothg (f) = \rg \gothg\;$ is equivalent to $\;\dim \gothg (\lambda) (\xi) = \rg \gothg (\lambda)$.
\smallskip

{\bf (b)}
For all $X \in \gothh$ and $\lambda' \in \gothh^*\!$, the complementary vector subspaces $\,[\gothg,\gothg(X)] = [\gothg,X]\,$ and $\,[\gothg(X),\gothh] \subseteq \gothg(X)\,$ of $[\gothg,\gothh]$ being orthogonal for both $B_{\lambda'}$ and ${\scriptstyle \langle} ~,\!~ {\scriptstyle \rangle}$, we have:
$\;\,\bigl|\Pi_{\gothg,\gothh}\bigr| (\lambda')
= \bigl|\Pi_{\gothg,\gothg(X)}\bigr| (\lambda')
\,{\scriptstyle \times}\,
\bigl|\Pi_{\gothg(X),\gothh}\bigr|(\lambda')$.
By choosing $X \in \gothh$ such that $\gothg(X) = \centraalg{\gothg(\mu)}{\gotha}$, it follows from ($*$) just before this lemma that:\\
\hspace*{\fill}%
$\bigl|\Pi_{\gothg,\centraalg{\gothg(\mu)}{\gotha}}\bigr| (\lambda)\;\,
= \!\!\Prod_{\substack
{\alpha \in R(\gothg_\CC,\gothh_\CC)\\
\alpha \notin R(\centraalg{\gothg(\mu)}{\gotha}_\CC,\gothh_\CC)}} \!\!\!
\abs{\langle \lambda,\alpha \rangle}^{\frac{1}{2}}$.
\hspace*{\fill}%
\\[-2mm]
This gives the result.
\cqfd
\end{Dem}

%-------------------------------------------------------------------------------
\section{The measures \texorpdfstring{$\beta_{G \cdot \lambdatilde}$}{\textbeta\_\{G \cdot \~\lambda\}}}\label{Rep2}
%------------------------------------------------------------------------------

%\noindent
We will now associate to each orbit $\,\wt{\Omega} \!\in\! G \backslash \gstregtilde$ a Radon measure $\beta_{\wt{\Omega}}$ on~$\gothg^*\!$ such that $\beta_{G\cdot\lambdatilde}$ coincides with the canonical measure on $G \cdot \lambda$ when $\,\lambda \in \gstssreg$, $\gothh=\gothg(\lambda)$ and $\lambdatilde = (\lambda,\gothh_{(\RR)}^{~~*})$.

%------------------------------------------------------------------------------
\begin{Definition}\rm
%------------------------------------------------------------------------------

Let $\lambda \!\in \gstss$, $\gothh \in \Car\gothg(\lambda)$
and $\,\F+\! \!\in C(\gothg(\lambda),\gothh)$.

{\bf (a)}
Let $\Fh+$ denote $\gotha^* + {\gotht^*}^+ \!$\labind{Fh+},
where $\,\gotht$ and $\gotha$ are the infinitesimally elliptic and hyperbolic components of $\gothh$ and $\F+\!$ is written $\,\F+\! = \gotha^* + \mi {\gotht^*}^+\!$ with ${\gotht^*}^+ \subseteq \gotht^*\!$.%
\smallskip

{\bf (b)}
Let $\,\Supp_{\gothg^*}(G \cdot (\lambda,\F+))\,$\labind{SuppgstOmegatilde}
denote the union in $\gstreg$ of $\,G \cdot (\omega+\{\lambda\})$, where $\omega$ ranges over the set of regular nilpotent orbits of $G(\lambda)$ in $\gothg(\lambda)^*$ included in the closure of $G(\lambda)\cdot \Fh+$ (see Lemma \ref{description d'ensembles de formes linéaires} (a) and Remark \ref{explication de la notation support} (1)).
\end{Definition}
\medskip

The following proposition is given in order to facilitate the use of the character formula \ref{formule du caractère}.
It will also, in some cases, make it possible to pass directly from the characters of the representations to the associated parameters (see \ref{recupération des paramètres}).%

%------------------------------------------------------------------------------
\begin{Proposition}
\hlabel{bijection entre orbites}
%------------------------------------------------------------------------------

{\bf (a)}
The set $\,G \,\backslash \gstfondtilde\,$ is in bijection with $\,G \,\backslash \gstreg\,$ via the map $R_G$\labind{RG} that sends $\,G \cdot \lambdatilde\,$ to $\,\Supp_{\gothg^*}(G \cdot \lambdatilde)$.
\smallskip

{\bf (b)}
Let $\lambda \in \gstss$, $\gothh \in \Car\gothg(\lambda)$ and $\,\F+\! \in C(\gothg(\lambda),\gothh)$.\\
The part $\,\Supp_{\gothg^*}(G \cdot (\lambda,\F+))\,$ of $\gstreg$ is the union of $R_G(G \cdot (\lambda,\Fp+))$ with $\gothh' \in \Car\gothg(\lambda)$ fundamental and $\Fp+\! \in C(\gothg(\lambda),\gothh')_{reg}$ satisfying $\,G(\lambda)\cdot \rho_\F+ \cap \overline{\Fp+} \not= \emptyset$.%
\vspace*{-1mm}
\end{Proposition}

\begin{Dem}{Proof of the proposition}

{\bf (a)}
First, let $\;(\lambda,\F+)\in \gstfondtilde$.
According to \ \cite[Th. p.$\!\!$~217 and Suppl. {\sc a}(b) p.$\!\!$~218] {Ro82} applied to $G_0(\lambda)$, there exists a unique regular nilpotent orbit $\omega_0$ of $G_0(\lambda)$ in $\;\overline{G_0(\lambda)\cdot \Fh+}$.
The set $\,\omega \egdef G(\lambda) \cdot \omega_0\;$ is therefore a regular nilpotent orbit of $G(\lambda)$ in $\;\overline{G(\lambda)\cdot \Fh+}$.
Any regular nilpotent orbit $\omega'$ of $G(\lambda)$ in $\;\overline{G(\lambda)\cdot \Fh+}$ intersects $\,u \cdot \overline{G_0(\lambda)\cdot \Fh+}$ for a certain $u \in G(\lambda)$, and therefore satisfies $\,\omega' \supseteq u \cdot \omega_0\,$ then $\omega' = \omega$.\\
This allows us to define $R_G$.
\smallskip

We now abandon the notation used at the beginning of this proof.

Let $\,f \in \gstreg$ be of semisimple and nilpotent components $\lambda$ and $\xi$.

Fix a fundamental $\gothh \in \Car\gothg(\lambda)$.
By Lemma \ref{description d'ensembles de formes linéaires} (a) and \cite[Lem. {\sc a} p.$\!\!$~220] {Ro82}, we can go through a finite sequence of inverse Cayley transformations (see \cite[p.$\!\!$~218] {Ro82}) from $\gothh$ to a Cartan subalgebra of $\gothg(\lambda)$ without imaginary roots.
%By \cite[Suppl. {\sc a}(a) and {\sc b} and {\sc c} p.$\!\!$~218, and Th. p.$\!\!$~221]{Ro82},
By \cite[Suppl. {\sc a}(a) and {\sc b} and {\sc c} p.$\!\!$~218, and Th. p.$\!\!$~217\llabel{page217}]{Ro82},
there exists $\F+\! \!\in C(\gothg(\lambda),\gothh)_{reg}$ unique up to the action of $\normagr{G_0(\lambda)}{\gothh}$ for which $\,\omega_0 \egdef G_0(\lambda)\cdot \xi\,$ is the regular nilpotent orbit of $G_0(\lambda)$ in $\;\overline{G_0(\lambda)\cdot \Fh+}$.
Thus, $G \cdot f$ is the image of $\,G \cdot (\lambda,\F+)\,$ by $R_G$.

Let $\,G \cdot (\lambda',\Fp+)$ be a preimage of $G \cdot f$ by $R_G$.
We have $\Fp+\!\! \in C(\gothg(\lambda'),\gothh')_{reg}$ for some fundamental $\gothh' \in \Car\gothg(\lambda')$.
Let $\omega'_0$ denote the regular nilpotent orbit of $G_0(\lambda')$ in $\;\overline{G_0(\lambda') \cdot \Fphp+}$.
Given the beginning of this proof, there exists $g \in G\,$ such that $g\xi \in \omega'_0$ and $g\lambda = \lambda'$.
Since the Cartan subalgebras $\gothh$ and $g^{-1}\gothh'$ of $\gothg(\lambda)$ are fundamental, there exists $x \in G(\lambda)_0$ such that $g^{-1}\gothh' = x\gothh$.
Therefore, $\omega_0$ which is equal to $g^{-1} \omega'_0$ is also the regular nilpotent orbit of $G_0(\lambda)$ in $\;\overline{G_0(\lambda)\cdot x^{-1}g^{-1}\Fp+\!}$.
Hence:
$\;\,x^{-1}g^{-1} \cdot (\lambda',\Fp+)\in
\normagr{G_0(\lambda)}{\gothh} \cdot (\lambda,\F+)$.
It follows that $G \cdot (\lambda,\F+)$ is the unique preimage of $G \cdot f$ by $R_G$.

{\bf (b)}
Using an argument from the beginning of the proof of (a), we see that $\,\Supp_{\gothg^*}(G \cdot (\lambda,\F+))\,$ is the union of the $G$-orbits that intersect $\,\Supp_{\gothg^*}(G_0 \cdot (\lambda,\F+))$.
We are therefore reduced to proving (b) for $G_0$ with $\lambda = 0$.
Let us consider this case.
Fix a fundamental Cartan subalgebra $\gothh'$ of $\gothg$ that contains the infinitesimally elliptic linear form $-\mi \,\rho_\F+$.
Set $\,l_0 = -\mi \,\rho_\F+$ in~$\Fh+\!$.
We will briefly review the ideas of \cite[Proof of Suppl. {\sc c} p.$\!\!$~228] {Ro82}.

From \cite[Th. 23 p.$\!\!$~50] {Va77} and \cite[$(L1)\!\!\Leftrightarrow\!(L2)\!\!\Leftrightarrow\!\!(L3)$ p.$\!\!$~216]{Ro82}, we~have:\\[0.5mm]
\hspace*{\fill}%
$\Sum_{\Omega \,\in\, G_0 \backslash \Supp_{\gothg^*}(G_0 \cdot (0,\F+))} \!\!\!
\beta_{\Omega} (\varphi)
\;=\;
\lim_{t \to 0^+}\,
\Sum_{\Omega' \in L_t}
\beta_{\Omega'} (\varphi)\;\;$
for all $\,\varphi \in C^\infty_c(\gothg^*)$,
\hspace*{\fill}\\[0.5mm]
where \emph{we use \ref{mesures de Liouville} (c)} and set
$\;L_t = \Lim_{\substack {l' \in \Fh+ \cap \, \gstreg \\ l' \to tl_0}}\! G_0 \cdot l'
\;\subseteq\; G_0 \backslash \gstreg\;$ when $t>0$.

By \cite[lines 2--9 p.$\!\!$~217] {Ro82}, the canonical map from $G_0(l_0)\backslash \glzerostreg$ to $G_0 \backslash \gstreg$ restricts to a bijection from
$\;\Lim_{\substack {l'' \in \F+[l_0]_\gothh \cap \, \glzerostreg \\ l'' \to tl_0}}\! G_0(l_0)\cdot l''\;$
onto $L_t$ for all $t>0$, where the limit is taken in $G_0(l_0)\backslash \glzerostreg$ and $\F+[l_0]$ is the element of $C(\gothg(l_0),\gothh)$ containing $\F+\!$.

The roots of $(\gothg(l_0)_\CC,\gothh_\CC)$ are nonimaginary and those of $(\gothg(l_0)_\CC,\gothh'_\CC)$ are nonreal.
From \cite[Suppl. {\sc c} p.$\!\!$~218] {Ro82}, we deduce that:\\
\hspace*{\fill}%
$\Lim_{\substack {l'' \in \F+[l_0]_\gothh \cap \, \glzerostreg \\
l'' \to tl_0}}\!
G_0(l_0)\cdot l''
\;=
\!\bigcup_{\dot{\Fp+} \,\in\, W(G_0(l_0),\gothh')\backslash \cale}\;
\Lim_{\substack {l'' \in \Fp+[l_0]_{\gothh'} \cap \, \glzerostreg \\
l'' \to tl_0}}\!
G_0(l_0)\cdot l''$
\hspace*{\fill}\\
where
$\,\cale \egdef \big\{\, \Fp+ \! \in C(\gothg,\gothh')\mid l_0 \in \overline{\Fphp+} \,\Big\} \,$ and $\Fp+[l_0]$ is the element of $C(\gothg(l_0),\gothh')$ containing $\Fp+\!$.
Furthermore, the previous union is a union of pairwise disjoint sets, by \cite[Suppl. {\sc b} p.$\!\!$~218] {Ro82}.

The Radon measures $\beta_{\Omega}$ with $\,\Omega \in G_0 \backslash \gstreg\,$ (see \ref{mesures de Liouville} (c)) are linearly independent.
Using \cite[Th. p.$\!\!$~217] {Ro82}, we find that $\Supp_{\gothg^*}(G_0 \cdot (0,\F+))$ is the disjoint union of $R_{G_0}(G_0 \cdot (0,\Fp+))$ where $\dot{\Fp+} \!\in\, W(G_0(l_0),\gothh')\backslash C(\gothg,\gothh')_{reg}$ and $\rho_\F+ \in \overline{\Fp+}$.
This yields the result.
\cqfd
\end{Dem}

%------------------------------------------------------------------------------
\begin{Definition}\rm
\hlabel{mesures de Liouville}
%------------------------------------------------------------------------------

{\bf (a)}
Let $\,\cald_{\gothg}$\labind{Dg} denote the function on $\gothg$ whose value at $\,X \in \gothg$ is the coefficient of $T^r$ in $\det(T\id - \ad X)$, where $\,r$ is the rank of $\gothg$.
It is polynomial and invariant under $\Ad G$.
\smallskip

{\bf (b)}
Let $\,\gssreg\!$\labind{g ssreg} denote the set of semisimple $X \in \gothg$ such that $\gothg(X)$ is abelian.
Therefore $\,\gssreg\!$ is the dense open set (with negligible complement) of $\gothg$ formed by the points where $\,\cald_{\gothg}$ does not vanish (see \cite[(2) and Lem. 1 p.$\!\!$~9] {Va77}).
\smallskip

{\bf (c)}
Let $\Omega \in G\,\backslash \gothg^*\!$.
Let $\beta_{\Omega}$\labind{β_Ω} denote the pushforward measure of the Liouville measure on the symplectic manifold $\Omega$ by the canonical injection $\Omega \hookrightarrow \gothg^*$ (see \cite[2.6 p.$\!\!$~20] {B.72}).
This is a Radon measure on $\gothg^*\!$ (see \cite[Th. 2 p.$\!\!$~509] {Ra72}).
\smallskip

{\bf (d)}
Let $\;\wt{\Omega}_0 \in G_0\,\backslash \gstregtilde$.
Set
$\;\; \beta_{\wt{\Omega}_0}
= \Sumpetit_{\Omega_0 \,\in\,
G_0 \backslash \Supp_{\gothg^*}(\wt{\Omega}_0)} \beta_{\Omega_0}$.

{\bf (e)}
Let $\;\wt{\Omega} \in G\,\backslash \gstregtilde$.
Set
$\;\; \beta_{\wt{\Omega}}
= \Sumpetit_{\wt{\Omega}_0 \,\in\, G_0
\backslash\wt{\Omega}}\beta_{\wt{\Omega}_0}$\labind{Omegatilde}
(see end of \ref{limites de mesures de Liouville} (a)).%
\vspace*{-1mm}
\end{Definition}

%------------------------------------------------------------------------------
\begin{Remarque}\rm
\hlabel{explication de la notation support}
%------------------------------------------------------------------------------

{\bf (1)}
Let $\,\wt{\Omega} \in G\,\backslash \gstregtilde$.
From the beginning of the proof of Proposition \ref{bijection entre orbites} (b), the set $\Supp_{\gothg^*}(\wt{\Omega})$ is the union of the orbits $\Omega \in G\,\backslash \gstreg$ such that $\beta_{\wt{\Omega}}(\Omega)\not= 0$.
\smallskip

It also follows from the above definition that
$\;\; \beta_{\wt{\Omega}}
= \Sumpetit_{\smash{\Omega_0 \,\in\,G_0\,\backslash \Supp_{\gothg^*}(\wt{\Omega}})}
m_{\Omega_0} \, \beta_{\Omega_0}\;$
where $\,m_{\Omega_0}$ is the (nonzero) number of $\wt{\Omega}_0 \in G_0 \backslash \wt{\Omega}$ satisfying $\Omega_0 \subseteq \Supp_{\gothg^*}(\wt{\Omega}_0)$.
\\
In particular:
$\;\; \beta_{\wt{\Omega}} = \beta_{R_G(\wt{\Omega})}\;$ when $\;\wt{\Omega} \in G\,\backslash \gstfondtilde$.
\smallskip

{\bf (2)}
The example where $G$ is equal to the semidirect product of $\ZZ/ 2\ZZ$ by $SL(2,\RR)^2$ where the nontrivial element of $\ZZ/ 2\ZZ$ acts on $SL(2,\RR)^2$ by permuting the coordinates, and where $\wt{\Omega}$ is of the form $G \cdot (0,\F+)$ for a chamber $\F+$ associated to $\gothh \in \Car\gothg$ neither fundamental nor split, shows that the coefficients $m_{\Omega_0}$ of (1) above can take values other than $1$.
\smallskip

{\bf (3)}
From (1) and \ref{bijection entre orbites} (a), the map $\wt{\Omega} \mapsto \beta_{\wt{\Omega}}$ is injective on $G\,\backslash
\gstfondtilde$.

When $G = SL(3,\RR)$, the two elements $\wt{\Omega}$ of $G\,\backslash \gstregtilde$ of the form $G \cdot (0,\F+)$ (one of which lies in $G\,\backslash \gstfondtilde$) correspond to the same measure $\beta_{\wt{\Omega}}$, equal to the Liouville measure of the unique regular nilpotent orbit of $G$ in~$\gothg^*$.
\cqfr
\end{Remarque}
\medskip

Here is a well-known result.
Since (b) will not be used, I will give only an abbreviated proof of it, suggested by Charbonnel.

%------------------------------------------------------------------------------
\begin{Proposition}
\hlabel{mesures de Liouville tempérées}
%------------------------------------------------------------------------------

{\bf (a)}
The function $\,\abs{\cald_{\gothg}}^{-1/2}$ is locally integrable on $\gothg$.
\smallskip

{\bf (b)}
The Radon measures $\beta_{\Omega}$ on $\gothg^*$ with $\,\Omega \!\in\! G\,\backslash\, \gothg^*$ are tempered.
\smallskip

{\bf (c)}
Let $\Omega \!\in\! G\,\backslash\, \gothg^*\!$.
The generalized function $\,\widehat{\beta}_{\Omega}$ is (modulo equality almost everywhere) a locally integrable function on $\gothg$ and analytic on $\gssreg$.%
\end{Proposition}

\begin{Dem}{Proof of the proposition}

{\bf (a)}
This first assertion is proven in \cite[Prop. 15 p.$\!\!$~66]{Va77}.
\smallskip

{\bf (b)}
Consider a real linear algebraic group~$H$.
Let $\gothh$ denote its Lie algebra.
Let us consider a coadjoint orbit $\Omega$ of $H$ whose Liouville measure $\beta_{\Omega}$ is a Radon measure on $\gothh^*$.
We follow the arguments of \cite[3.2 p.$\!\!$~220] {Ch96}.
Let $W_0$ be an open set of $\Proj(\gothh^* \times \RR)$ corresponding to the nonvanishing of a homogeneous coordinate, and let $U$ be a relatively compact open set in $W_0$ for the usual topology.
We construct a certain regular function $\,q :U \to \RR_+ \!\moins0$.
\llabel{OrbitesCoadjointes1}
According to \cite[p.$\!\!$~217] {Ch96}, the meromorphic function $\Phi$ on $\CC$ such that
$\,\langle\Phi(s),\phi\rangle = \int_{\Omega \cap U} q(l)^s\,\phi(l)\, \diff \beta_{\Omega}(l)\,$
for $\,\phi \in C^\infty_c(\gothh^* \cap U)$,
when $\Re(s)$ is large enough,
is continuous for the topology of $\cals(\gothh^*)$.
Hence the measure $\beta_{\Omega}$ on $\gothh^*$ is tempered by [Cha 96, line 3 p.$\!\!$~222].%
\llabel{OrbitesCoadjointes2}%
\smallskip

{\bf (c)}
See \cite[Prop.\,13 p.$\!\!$~65,\,Th.\,17 p.$\!\!$~66,\,Th.\,28 p.$\!\!$~95,\,bottom\,of p.$\!\!$~105] {Va77}.\cqfd
\end{Dem}
\medskip

The following proposition will be of interest when attempting to reduce, by limiting processes or by translation, to regular semisimple linear forms.

%------------------------------------------------------------------------------
\begin{Proposition}
\hlabel{limites de mesures de Liouville}
%------------------------------------------------------------------------------

Let $\lambda \in \gstss$, $\,\gothh \in \Car\gothg(\lambda)$, and $\,\F+\! \in C(\gothg(\lambda),\gothh)$.
Set $\lambdatilde = (\lambda,\F+)$.
Let $\gotha$ denote the hyperbolic component of $\gothh$, and fix a chamber ${\gotha^*}^+\!$ of $(\gothg(\lambda)(\mi \rho_\F+),\gotha)$.
\smallskip

{\bf (a)}
The positive measures $\,\beta_{\wt{\Omega}}$ with $\,\wt{\Omega} \!\in\! G\,\backslash \gstregtilde$ are tempered Radon measures.
Let $\;\varphi\in \cals(\gothg^*)$.
We have\\[-1.5mm]
\hspace*{\fill}%
$\Lim_{\substack
{\lambda' \in \Fh+ \cap \, \gstreg \\ \lambda' \to \lambda}}
\bigl(
{\scriptstyle \abs{W(G,\gothh)\, (\lambda')}}\;\;
\beta_{G\cdot \lambda'}(\varphi)\bigr)
= \cases{
{\scriptstyle \abs{W(G,\gothh)\, (\lambda_+)}}\;\;
\beta_{G \cdot \lambdatilde}(\varphi)
& if $\;\lambdatilde \in \gstregtilde$
\cr \hfil 0
& otherwise.
}$
\hspace*{\fill}\\[0.5mm]
In particular
$\;\; \beta_{G \cdot \lambdatilde}(\varphi)
= \Lim_{t \to 0^+}\; \beta_{G \cdot \lambda_t}(\varphi)\;\,$
when $\;\lambdatilde \in \gstregtilde\;$ (see \ref{lemme clef}).
\smallskip

{\bf (b)}
For all $\,X \in \gssreg$, taking into account \ref{mesures de Liouville tempérées} (c) we have\\
\hspace*{\fill}%
$\Lim_{\substack
{\lambda' \in \Fh+ \cap \, \gstreg \\ \lambda' \to \lambda}}
\bigl( {\scriptstyle \abs{W(G,\gothh)\, (\lambda')}}\;\;
\widehat{\beta}_{G\cdot \lambda'}(X)\bigr)
= \cases{
{\scriptstyle \abs{W(G,\gothh)\, (\lambda_+)}}\;\;
\widehat{\beta}_{G \cdot \lambdatilde}(X)
& if $\;\lambdatilde \in \gstregtilde$
\cr \hfil 0
& otherwise.
}$
\hspace*{\fill}\\
In particular
%$\;\; \beta_{G \cdot \lambdatilde}(X) = \Lim_{t \to 0^+}\; \beta_{G \cdot \lambda_t}(X)\;\,$
$\;\; \widehat{\beta}_{G \cdot \lambdatilde}(X) = \Lim_{t \to 0^+}\; \widehat{\beta}_{G \cdot \lambda_t}(X)\;\,$\llabel{chapeau}
when $\;\lambdatilde \in \gstregtilde\;$ (see \ref{lemme clef}).

More precisely, given $\gothj \in \Car\gothg$, $y \in \interieur(\gothg_\CC)$ for which $\,\gothj_\CC= y \gothh_\CC$, a system of positive roots $R^+(\gothg_\CC,\gothj_\CC)$ of $R(\gothg_\CC,\gothj_\CC)$, a connected component $\Gamma$ of $\gothj \cap \gssreg$, and a connected component $\,{\gothh^*}^+$ of $\,\Fh+ \cap \gstreg\,$ to which $\lambda$ is close to, there exists a family $\bigl(c_w\bigr)_{w \in W(\gothg_\CC,\gothj_\CC)}\!$ of complex numbers such that, for all $X \!\in\! \Gamma\!$, we have
\smallskip

$\Sum_{w \in W(\gothg_\CC,\gothj_\CC)}
c_w \;\; \me ^{\mi \,wy\lambda'\,(X)}
= \Prod_{\alpha \in R^+(\gothg_\CC,\gothj_\CC)}
\!\! \alpha(X)\;\;
{\scriptstyle \times}\;\;
\widehat{\beta}_{G_0\cdot \lambda'}(X)\;\,$
when $\,\lambda' \in {\gothh^*}^+$\\
and
\vskip-2mm

$\Sum_{w \in W(\gothg_\CC,\gothj_\CC)}
c_w \;\, \me ^{\mi \,wy\lambda\,(X)}
= \cases{
\Prod_{\alpha \in R^+(\gothg_\CC,\gothj_\CC)}
\!\! \alpha(X)\;\;
{\scriptstyle \times}\;\,
\widehat{\beta}_{G_0 \cdot \lambdatilde}(X)
& if $\;\lambdatilde \in \gstregtilde$
\cr \; \hfil 0
& otherwise.
}$
\end{Proposition}

\begin{Dem}{Proof of the proposition}

{\bf (a)}
Let $\,\varphi \in C^\infty_c(\gothg^*)$ be nonnegative.
According to
\cite[$(L1)\!\Leftrightarrow\!(L2)\!\Leftrightarrow\!(L3)$ p.$\!\!$~216, lines 9--12 p.$\!\!$~217, and Th. p.$\!\!$~217] {Ro82}, we have\\
\hspace*{\fill}%
$\Lim_{\substack
{\lambda' \in \Fh+ \cap \, \gstreg \\ \lambda' \to \lambda}}
\beta_{G_0\cdot\lambda'}(\varphi)
= \cases{
\beta_{G_0 \cdot \lambdatilde}(\varphi)
& if $\;\lambdatilde \in \gstregtilde$
\cr \hfil 0
& otherwise.
}$
\hspace*{\fill}\\
Fix a norm $\norm{\cdot}$ on $\gothg^*\!$.
According to \cite[(i) p.$\!\!$~40] {Va77}, there exists $N \!\in \NN$ such that the set formed by $\;\Int_{\gothg^*} (1 + \norm{l})^{-N} \diff \beta_{\Omega_0}(l)\;$ with $\,\Omega_0 \!\in\! G_0 \backslash \gstssreg\,$ is bounded.
It follows that~$\,\beta_{G_0 \cdot \lambdatilde}\,$ is a tempered Radon measure when $\,\lambdatilde \in \gstregtilde$, and that the limiting process remains valid when the condition “$\varphi \in C^\infty_c(\gothg^*)$ nonnegative” is replaced by “$\varphi \in \cals(\gothg^*)$.”

We then obtain, for all $\varphi \in \cals(\gothg^*)$:\\
\hspace*{\fill}%
$\Lim_{\substack
{\lambda' \in \Fh+ \cap \, \gstreg \\ \lambda' \to \lambda}}\;
{\textstyle \Sum_{\dot{x} \in G/G_0} \!\! \beta_{G_0\cdot x\lambda'}
(\varphi)} = \cases{
\Sumpetit_{\dot{x} \in G/G_0} \!\! \beta_{G_0 \cdot x\lambdatilde} (\varphi)
& if $\;\lambdatilde \in \gstregtilde$
\cr \hfil 0
& otherwise.
}$
\hspace*{\fill}\\
We first transform the sum on the left for $\lambda' \in \Fh+ \cap \gstreg$:\\
\hspace*{\fill}%
$\Sum_{\dot{x} \in G/G_0} \!\! \beta_{G_0\cdot x\lambda'}\textstyle
\,=\, {\scriptstyle \abs{G(\lambda') G_0 / G_0}}\;\;
\beta_{G\cdot\lambda'}
\,=\, \frac{\abs{W(G,\gothh)(\lambda')}}
{\abs{\centragr{G}{\gothh} / \centragr{G_0}{\gothh}}}\;\;
\beta_{G\cdot\lambda'}$.
\hspace*{\fill}\\
Furthermore, the chambers of $(\gothg(\lambda)(\mi \rho_\F+),\gotha)$ are conjugate under $\normagr{G(\lambda)(\mi \rho_\F+)_0}{\gotha}$ and \emph{a fortiori} under $G_0(\lambdatilde)$, which proves that $\;G(\lambdatilde) = G_0(\lambdatilde)\;\normagr{G(\lambdatilde)}{{\gotha^*}^+}$.\\
When $\;\lambdatilde \in \gstregtilde$, we therefore have, using Lemma \ref{lemme clef}:\\
\hspace*{\fill}%
$\Sum_{\dot{x} \in G/G_0} \!\! \beta_{G_0\cdot x\lambdatilde}\textstyle
\,=\, {\scriptstyle \abs{G(\lambdatilde) G_0 / G_0}}\;\;
\beta_{G\cdot\lambdatilde}
\,=\, \frac{\abs{W(G,\gothh)(\lambda_+)}}
{\abs{\centragr{G}{\gothh} / \centragr{G_0}{\gothh}}}\;\;
\beta_{G\cdot\lambdatilde}$.
\hspace*{\fill}\\
This proves (a).
\smallskip

{\bf (b)}
The penultimate equality in this statement is taken from \cite[Th. 4 p.$\!\!$~108]{Va77}.
It shows that $\;\widehat{\beta}_{G_0 \cdot \lambda'}\,$ converges pointwise on $\,\gssreg\,$ to a continuous $G_0$-invariant function when $\lambda' \to \lambda$ with $\lambda' \in {\gothh^*}^+\!$.
Let $\psi \in C^\infty_c(\gssreg)$.
From the proof in (a), we have\\[1mm]
\hspace*{\fill}%
$\Lim_{\substack
{\lambda' \in {\gothh^*}^+ \\ \lambda' \to \lambda}}
\Int_{\gssreg}
\widehat{\beta}_{G_0 \cdot \lambda'}(X)\;
\psi(X)\,\diff_{\gothg}(X)
= \cases{
\widehat{\beta}_{G_0 \cdot \lambdatilde}(\psi\,\diff_{\gothg})
& if $\;\lambdatilde \in \gstregtilde$
\cr \hfil 0
& otherwise.
}$%
\hspace*{\fill}%

In light of the upper bounds $\;\abs{c_w \; \me ^{\mi \,wy\lambda'\,(X)}} \,\leq\, \abs{c_w}\;$ from \cite[Th. 7 p.$\!\!$~111] {Va77}, we see that the $G_0$-invariant functions
$\,\abs{\cald_{\gothg}}^{1/2}\! {\scriptstyle \times}\, \widehat{\beta}_{G_0 \cdot \lambda'}$
on $\gssreg$ with $\lambda' \in {\gothh^*}^+$ are uniformly bounded.
Proposition \ref{mesures de Liouville tempérées} (a) therefore allows us to apply Lebesgue's dominated convergence theorem in the above limit.
This yields the final equality in the statement.

The passage from $G_0$ to $G$ is carried out as in the proof of (a).
\cqfd
\end{Dem}

%-------------------------------------------------------------------------------
\section{Fixed points of an elliptic element}\label{Rep3}
%------------------------------------------------------------------------------

%\noindent
Let us consider an elliptic element $e$ of $G$.

The Lie algebra $\gothg(e)$ is reductive because $\gothg(e)_\CC$ has as real form $\gothc(e)$, where~$\gothc$ is the Lie algebra of a maximal compact subgroup of $\GG(\CC)$ that contains the projection of $e$ into $\GG(\RR)$.
Furthermore, any Cartan subalgebra $\gothh_e$ of $\gothg(e)$ intersects $\,\gssreg$, which means that the centralizer $\gothh$ of $\gothh_e$ in $\gothg$ is a Cartan subalgebra of $\gothg$ (see \cite[Lem. 1.4.1 p.$\!\!$~6] {Bo87}).
Therefore $G(e)$ satisfies the hypothesis stated in the introduction concerning $G$ (see \cite[bottom of~p.$\!\!$~38] {DV93}).%

%------------------------------------------------------------------------------
\begin{Definition}\rm
%------------------------------------------------------------------------------

Let $\,\gstreg(e)$\labind{g*reg(e)}
(respectively $\,\gstItilde(e)$\labind{g*I tilde(e)},
$\gstInctilde(e)$\labind{g*Inc tilde(e)},
$\gstfondtilde(e)$\labind{g*fond tilde(e)}, and
$\gstregtilde(e)$\labind{g*reg tilde(e)})
denote the set of elements of $\gstreg$
(respectively $\gstItilde$, $\gstInctilde$, $\gstfondtilde$, and
$\gstregtilde$)
that are fixed under the action of $\,e$.
\end{Definition}

Among the results of the following lemma, only point~(a) is essential.
Point~(b) will be used in Remark \ref{recupération des paramètres} and point~(c) in a later paper.%

%------------------------------------------------------------------------------
\begin{Lemme}
\hlabel{descente pour les formes linéaires}
%------------------------------------------------------------------------------

Let $\lambda \in \gstss$, $\,\gothh \in \Car\gothg(\lambda)\,$ with infinitesimally elliptic and hyperbolic components $\gotht$ and $\gotha$, and let $\,\F+\! \!\in C(\gothg(\lambda),\gothh)\,$ such that $e$ fixes $\;\lambdatilde \egdef (\lambda,\F+)$.
Assume there exists a chamber ${\gotha^*}^+\!$ of $(\gothg(\lambda)(\mi \rho_\F+),\gotha)$ stable under $e$ (we will refer to this condition by saying that “$\lambdatilde[e]$ exists”, see (a) below).
\smallskip

{\bf (a)}
We have $\,\lambda \in \gestss$, $\gothh(e)\in \Car\gothg(e)(\lambda)$, and the element $\rho_\F+\!$ of $\gothh(e)_{(\RR)}^{~~*}$ is regular for the imaginary roots of $(\gothg(e)(\lambda)_\CC,\gothh(e)_\CC)$.
By definition, $\F+[e]$\labind{F+[e]} denote the element of $C(\gothg(e)(\lambda),\gothh(e))$ containing $\rho_\F+\!$, and $\lambdatilde[e] = (\lambda,\F+[e])$\labind{λtilde[e]}.
\smallskip

{\bf (b)}
There exists $e_0 \in e\exp\gotht(e)$ such that:
if $\Fp+ \!\in C(\gothg(\lambda),\gothh)$ and $\rho_\Fp+ \!\in \F+[e_0]$, then $\Fp+ \!= \F+\!$.
When $\lambdatilde \in \gstItilde$ we can choose $e_0$ such that $\,\lambdatilde[e_0] \in \gezerostItilde$.%
\smallskip

{\bf (c)}
Assume that $\,\lambdatilde \in \gstItilde$.
The cardinality of the set of $\,\lambdatilde' \!=\! (\lambda',\Fp+)\!$\llabel{PbFp} with $\Fp+ \!\in C(\gothg(\lambda),\gothh)$ satisfying $\,e \lambdatilde' = \lambdatilde'\,$ and $\lambdatilde'[e] \!=\! \lambdatilde[e]$\llabel{PbFpSuite}
is equal to
$\frac{|W(\gothg(\lambda)_\CC,\gothh_\CC)(\Ad^*\! e)|}{|W(\gothg(e)(\lambda)_\CC,\gothh(e)_\CC)|}\!$,
where $\,W(\gothg(\lambda)_\CC,\gothh_\CC)(\Ad^*\! e)$ denotes the centralizer of $(\Ad^* e^\CC)_{\gothh_\CC^*}$ in $\,W(\gothg(\lambda)_\CC,\gothh_\CC)$.%
\end{Lemme}

\begin{Dem}{Proof of the lemma}

{\bf (a)}
By Lemma \ref{lemme clef}, we have:
$\;\lambda_+ \in \gestssreg$, then $\,\gothh(e)\!=\! \gothg(e)(\lambda_+)\in \Car\gothg(e)(\lambda)$.
Let $\,\alpha' \in R(\gothg(e)(\lambda)_\CC,\gothh(e)_\CC)\,$ be imaginary.
There exists $\,\alpha \in R(\gothg(\lambda)_\CC,\gothh_\CC)\,$ such that $\,\alpha' = \restriction{\alpha}{\gothh(e)_\CC}$.
We have:
$\,\langle \nu_+,\alpha \rangle = \langle \nu_+,\alpha' \rangle = 0\,$
and therefore $\alpha \in R(\gothg(\lambda)_\CC,\gothh_\CC) \cap R(\gothg(\nu_+)_\CC,\gothh_\CC)$.
Furthermore $\langle \rho_\F+,\alpha' \rangle \!=\! \langle \rho_\F+,\alpha \rangle$,
which gives $\;\langle \rho_\F+,\alpha' \rangle \not= 0\;$ (see \ref{choix de racines positives} (b)).
\smallskip

{\bf (b)}
Let $B$ be the base of $R(\centraalg{\gothg}{\gotha}(\lambda)_\CC,\gothh_\CC)$ associated with $\F+\!$.
For each $\alpha \!\in\! B$, let $m_{\dot{\alpha}}$ denote the cardinality of the orbit ${\dot{\alpha}}$ of $\alpha$ under the action of the subgroup $\langle e \rangle$ of $G$ generated by $e$ on $B$, and choose $\,t_{\dot{\alpha}} \in \RR\,$ such that
$\;((\Ad e^\CC)^{\,m_{\dot{\alpha}}})_{\gothg_\CC^\alpha} = \me ^{\mi \,t_{\dot{\alpha}}} \id$.
Take $\,e_0 = e \exp (X_0)$, where $X_0 \in \gothh_\CC\cap \derivealg{\centraalg{\gothg}{\gotha}\,(\lambda)_\CC}\,$ is such that $\,\alpha(X_0) = -\mi t_{\dot{\alpha}} / m_{\dot{\alpha}}\,$ for all $\alpha \!\in\! B$.
It is clear that $\lambdatilde[e_0]$ exists.\\
Each $\restriction{\alpha}{\gothh(e)_\CC}$ with $\,\alpha \in B\,$ is a weight of $\gothh(e)_\CC$ in $\gothg(e_0)_\CC$ whose eigenspace contains
$(1 + e_0 + \cdots + {e_0}^{\! m_{\dot{\alpha}} - 1})\cdot \gothg_\CC^{\alpha}$.
Let $\Fp+\! \in C(\gothg(\lambda),\gothh)$ be such that $\rho_\Fp+\! \in \F+[e_0]$, hence $e_0\Fp+\!= \Fp+\!$.
For all $\alpha \in B$, we have:
$\langle \rho_\F+,\restriction{\alpha}{\gothh(e)_\CC} \rangle \!= \langle \rho_\F+,\alpha \rangle > 0$
and
$\langle \rho_\Fp+,\alpha \rangle = \langle \rho_\Fp+,\restriction{\alpha}{\gothh(e)_\CC} \rangle$,
so the root $\restriction{\alpha}{\gothh(e)_\CC}\!$ of $R(\centraalg{\gothg(e_0)}{\gotha(e)}(\lambda)_\CC,\gothh(e)_\CC)$ is positive with respect to $\F+[e_0]$,
therefore $\langle \rho_\Fp+,\alpha \rangle > 0$.
It follows that $\Fp+\! = \F+\!$.

Now suppose that $\lambdatilde \in \gstItilde$.
Every root of $(\gothg(e_0)(\lambda)_\CC,\gothh(e)_\CC)$ is a restriction of a root of $(\gothg(\lambda)_\CC,\gothh_\CC)$.
The roots $\,\alpha' = \restriction{\alpha}{\gothh(e)_\CC}$ with $\,\dot{\alpha} \in \langle e \rangle \backslash B\,$ therefore describe as $\dot{\alpha}$ varies (injectively because $\dim \gothg(e_0)_\CC^{\alpha'} = 1$), the set of simple roots with respect to $\F+[e_0]$ in $R(\gothg(e_0)(\lambda)_\CC,\gothh(e)_\CC)$.
These roots $\alpha'$ are noncompact because $\alpha([\, \cdot \, ,\overline{\vrule width 0ex height 1ex \, \cdot \,}])$ is nonnegative on the vector subspace $\gothg(e_0)_\CC^{\alpha'}$ of the orthogonal direct sum $\,\gothg_\CC^{\alpha} \oplus e_0 \gothg_\CC^{\alpha} \oplus \cdots \oplus {e_0}^{\! m_{\dot{\alpha}} - 1} \gothg_\CC^{\alpha}$.
\smallskip

{\bf (c)}
The map
$\,\pi_e : W(\gothg(\lambda)_\CC,\gothh_\CC) \cdot \rho_\F+ \cap \gothh(e)_\CC^*
\to C(\gothg(e)(\lambda),\gothh(e))\,$
that sends a linear form to the chamber that contains it, commutes with the action of the group $\;W(\gothg(e)(\lambda)_\CC,\gothh(e)_\CC)$.
The groups $\;W(\gothg(\lambda)_\CC,\gothh_\CC)(\Ad^* e)\;$ and $\;W(\gothg(e)(\lambda)_\CC,\gothh(e)_\CC)\;$ act simply transitively on the domain and codomain of $\pi_e$ respectively.
The fibers of $\pi_e$ therefore have cardinality
$\,\frac{|W(\gothg(\lambda)_\CC,\gothh_\CC)(\Ad^* e)|}{|W(\gothg(\lambda)(e)_\CC,\gothh(e)_\CC)|}$.
The set of
%$\lambdatilde' \!=\! (\lambda',\Fp+)\in \gstItilde(e)$ such that $\lambdatilde'\![e] = \lambdatilde[e]$,
$\lambdatilde' \!=\! (\lambda',\Fp+)$ with $\Fp+ \in C(\gothg(\lambda),\gothh)$ such that $\,e \lambdatilde' = \lambdatilde'\,$ and $\lambdatilde'[e] \!=\! \lambdatilde[e]$,\llabel{PbFpDém}
whose elements satisfy $\lambda' = \lambda$ and
$\gothg(\lambdatilde')
= \centraalg{\gothg}{\gothg(e)(\lambdatilde'\![e])}
= \gothh$,
is in bijection with $\pi_e^{-1}(\{\F+[e]\})$ via the map $(\lambda',\Fp+)\mapsto \pi_1^{-1}(\Fp+)$.
%This allows us to calculate its cardinality.
This provides its cardinal.
\cqfd
\end{Dem}

%------------------------------------------------------------------------------
\begin{Remarque}\rm
\hlabel{insuffisance de la méthode des orbites}
%------------------------------------------------------------------------------

{\bf (1)}
It can happen that $\,\lambdatilde \!\in\! \gstItilde(e)\,$ satisfies $\,\lambdatilde[e] \!\notin\! \gestItilde$.
(However, the condition $\,\lambdatilde \in \gstInctilde(e)\,$ implies $\,\lambdatilde[e] \in \gestInctilde$.)
For example:
$\,G = Sp(4,\RR)$, $\,\lambdatilde = (0,\F+)\,$ such that $\gothh \egdef \gothg(\lambdatilde)$ is equal to its infinitesimally elliptic component, and the system of positive roots of $\,(\gothg_\CC,\gothh_\CC)\,$ associated with $\F+$ is written $\,\{ \, \alpha_1,\alpha_2,\alpha_1+\alpha_2,2\alpha_1+\alpha_2 \, \}$ with the only compact root $\alpha_1+\alpha_2$.
Taking $\,e = \exp E$, where $\,E \in \gothh\,$ is determined by $\;\alpha_1(E) = -\alpha_2(E) = \mi \pi$, the Lie algebra $\gothg(e)$ is isomorphic to $\gothu(2)$ and therefore $\,C(\gothg(e),\gothh(e))_{reg} = \emptyset$.
\smallskip

{\bf (2)}
We can also find $\,\lambdatilde_e \in \gestregtilde\,$ that is not of the form $\lambdatilde[e]$ for a $\,\lambdatilde \in \gstregtilde(e)$.
For example:
$\,G = SU(2)$, $\,e = \left( {\scriptstyle \Rot{1}} \right)\,$ and $\,\lambdatilde_e = (0,\mi \,\goth{so}(2)^*)$.%
\smallskip

{\bf (3)}
For certain $\,\lambdatilde \!\in\! \gstfondtilde(e)$ there is no $f \!\in\! \gstreg(e)$ such that $G \cdot f = R_G(G \cdot \lambdatilde)$.
This is the case, with $\,G = SL(2,\RR)\,$ and $\,e = \left( {\scriptstyle \Rot{1}} \right)$, for elements of $\,\gstfondtilde(e)\,$ of the form $\,(0,\F+)\,$ such that $\,\F+ \!\subseteq \mi \,\goth{so}(2)^* \!$.
\cqfr
\end{Remarque}

%===============================================================================
%
\vspace*{-10pt}
\part{.\quad The “projective representation” parameters}\label{RepII}
%
%==============================================================================

The parameters in this section are based on those introduced by Duflo, which were suitable for the case of regular semisimple coadjoint orbits of $G$.

%-------------------------------------------------------------------------------
\section{About special-metalinear and metaplectic groups}\label{Rep4}
%------------------------------------------------------------------------------

%\noindent
My reference regarding the metaplectic group and the orientation functions of Duflo and Vergne attached to it is their paper \cite{DV93}.
Like them, I choose to avoid using the character of the metaplectic representation, following the viewpoint of \cite[Conj. p.$\!\!$~291]{Ve94} (except for the sign of the symplectic form on coadjoint orbits).

%------------------------------------------------------------------------------
\begin{Definition}\rm
%------------------------------------------------------------------------------

Consider a finite-dimensional real vector space $V$.
\smallskip

{\bf (a)}
Denote by $\,DL(V)\to SL(V)$\labind{DL(V)} “the” double cover of $SL(V)$, unique up to isomorphism of coverings, which is connected when $\dim V \geq 2$.
(Item (b) below will provide a canonical description of such a cover, as a set of ordered pairs consisting of an element of $SL(V)$ and an orientation of a certain vector subspace~of~$V\!$.)
\smallskip

{\bf (b)}
Let $\hat{a} \in DL(V)$ lie over an elliptic $a \in SL(V)$.

When $\dim V \geq 2$, set\\[0.5mm]
\hspace*{\fill}%
$\scalo (\hat{a})_{(1-a)\cdot V}
= (-1)^{\frac{\alpha_1 +\cdots+ \alpha_p}{2\pi}} \,
\sg (\sin (\frac{\beta_1}{2})\ldots\sin (\frac{\beta_q}{2}))\;
{\scriptstyle \times}\;
\RR_+ \!\moins0 \,
(w_1 \wedge \cdots \wedge w_{2q})$\labind{Oa},
\hspace*{\fill}%
\\[0.5mm]
where “$\sg$\labind{sg}” denotes the sign function on $\RR\moins0$, independently of the choice of an infinitesimally elliptic $A \in \goth{sl} (V)$ such that $\,\hat{a} = \exp_{\!DL(V)} \! A$, and of a basis $(v_1,\dots,v_{2p},w_1,\dots,w_{2q})$ of $A \cdot V$ in which the matrix of the restriction of $A$ has the form\\[-1mm]
\hspace*{\fill}%
{\bf (\boldmath$**$)}
$\left(\!
{\scriptscriptstyle\diags
{\bigg(\Rot{\alpha_1}\bigg)\!}
{\ddotsc}
{\bigg(\Rot{\alpha_p}\bigg)}
{\bigg(\Rot{\beta_1}\bigg)\!}
{\ddotsc}
{\bigg(\Rot{\beta_q}\bigg)}}
\!\right)$
\hspace*{\fill}%
\\[2mm]
with
$\;\alpha_1,\dots,\alpha_p \in 2\pi\ZZ\moins0\,$
and
$\,\beta_1,\dots,\beta_q \in \RR\setminus\! 2\pi\ZZ$.
\smallskip

When $\dim V \leq 1$, $\;\scalo (\hat{a})_{\{0\}}$ is equal to $\RR_+ \!\moins0$ or $\RR_- \!\moins0$ depending on whether $\hat{a}$ is trivial or not.
\smallskip

{\bf (c)}
Let $A \in \goth{sl} (V)$ be infinitesimally elliptic.
Set\\[0.5mm]
\hspace*{\fill}%
$\scalo (A)_{A \cdot V}
= \sg (\beta_1\ldots\beta_q)\;
{\scriptstyle \times}\;
\RR_+ \!\moins0 \,
(w_1 \wedge \cdots \wedge w_{2q})$\labind{OA},
\hspace*{\fill}%
\\[0.5mm]
independently of the choice of a basis $(w_1,\dots,w_{2q})$ of $A \cdot V$ in which the matrix of the restriction of $A$ has the form\\[3mm]
\hspace*{\fill}%
$\left(\!
{\scriptstyle\diagt
{\Big(\Rot{\beta_1}\Big)}
{\ddotsc}
{\vrule width 0ex depth 1.2ex
\Big(\Rot{\beta_q}\Big)}}
\!\right)$
\hspace*{\fill}%
\\
with $\;\beta_1,\dots,\beta_q \in \RR\moins0$.
\end{Definition}

%------------------------------------------------------------------------------
\begin{Definition}\rm
\hlabel{géométrie métaplectique}
%------------------------------------------------------------------------------

Consider a finite-dimensional real vector space $V$ equipped with a nondegenerate alternating bilinear form $B$.
We also denote by $B$ the complex bilinear extension of this bilinear form to $V_\CC$.
\smallskip

{\bf (a)}
Denote by $Mp(V)$\labind{Mp(V)} the inverse image of $Sp(V)$ in $DL(V)$.
\smallskip

{\bf (b)}
Denote by $\scalo (B)_V$\labind{OB} the orientation of $V$ (which means some $\RR_+ \!\moins0 \,\omega_0$ with $\,\omega_0 \!\in\! \bigwedge^{\scriptscriptstyle\mathrm{max}}V\moins0\big\}$) on which $B^{\frac{1}{2}\dim V}\!$ is positive.
A “symplectic basis” of $(V,B)$ is any basis $(P_1,\dots,P_n,Q_1,\dots,Q_n)$ of $V$ such that $\;B(P_i,Q_j) = \delta_{i,j}\;$ (Kronecker symbol) \ and $\;B(P_i,P_j) = B(Q_i,Q_j) = 0\;$ for $1 \leq i,j \leq n$.%
\smallskip

{\bf (c)}
A “Lagrangian” of $(V_\CC,B)$ is any complex vector subspace of $V_\CC$ equal to its orthogonal with respect to $B$.
A~“positive Lagrangian” of $(V_\CC,B)$ is any Lagrangian of $(V_\CC,B)$ on which the Hermitian form $\;(v,w)\mapsto \mi \,B(v,\overline{w})\;$ is nonnegative.

Let $\,\call$ be a Lagrangian of $(V_\CC,B)$.
Let $Mp(V)_\call$\labind{Mp(V)L} denote the normalizer of $\,\call$ in $Mp(V)$.
To each $\,\hat{x} \!\in\! Mp(V)_\call\,$ lying over an element $x$ of $Sp(V)$, associate the number $\,n_\call(x)\,$\labind{nL} of eigenvalues counted with multiplicity in $\left] 1,+\infty \right[$ of the restriction of $x^\CC$ to $\,\call$, and the number $\,q_\call(x)\,$\labind{qL} of (strictly) negative coefficients in the matrix of the Hermitian form $\;(v,w)\mapsto \mi \,B(v,\overline{w})\;$ on $(1-x^\CC)\cdot \call$ with respect to an orthogonal basis.
\smallskip

{\bf (d)}
Let $\,\call$ be a Lagrangian of $(V_\CC,B)$.
For every $\,\hat{x} \in Mp(V)_\call\,$ whose elliptic component is $\hat{x}_e$ (in the sense of \cite[Lem. 31 p.$\!\!$~38]{DV93}) lying over an element $x$ of $Sp(V)$ whose elliptic component is $x_e$, set\\[1mm]
\hspace*{\fill}%
$\rho_\call(\hat{x})
\,=\, (-1)^{q_\call(x_e)}\;
\frac
{\scalo (\hat{x}_e)_{(1-x_e)\cdot V}}
{\scalo (B)_{(1-x_e)\cdot V}}\;
\Prod_{1 \leq k \leq n} (\sqrt{r_k} \, \me ^{\mi \,\theta_k/2})$\labind{ρ L},
\hspace*{\fill}%
\\
where $\;r_1 \me ^{\mi \,\theta_1},\dots,r_n \me ^{\mi \,\theta_n}$ are the eigenvalues counted with multiplicity of the restriction of $x^\CC$ to $\,\call$, with $\;r_1,\dots,r_n \in \RR_+ \!\moins0\,$ and $\;\theta_1,\dots,\theta_n \in \left] -2\pi,0 \right]$.%
\smallskip

{\bf (e)}
Let $\hat{x} \in Mp(V)$ be of elliptic component $\hat{x}_e$ lying over a semisimple element $x$ of $Sp(V)$ with elliptic component $x_e$.
Set\\[1mm]
\hspace*{\fill}%
$\delta(\hat{x})
\,=\; \frac
{\scalo (\hat{x}_e)_{(1-x_e)\cdot V}}
{\scalo (B)_{(1-x_e)\cdot V}}\;
\Prod_{1 \leq k \leq n} \me ^{\mi \,\theta_k/2}$\labind{δ},
\hspace*{\fill}%
\\
independently (given the proof of (c) of the proposition below) of the choices of an infinitesimally elliptic $\,E \in \goth{sp} (V)$ such that $\,x_e = \exp_{\!Sp(V)} \! E$, and of a symplectic basis $(P_1,\dots,P_n,Q_1,\dots,Q_n)$ of $(V,B)$ for which the matrix of $E$ relative to $(P_1,Q_1,\dots,P_n,Q_n)$ has the form\\
\hspace*{\fill}%
$\left(\!
{\scriptstyle\diagt
{\Big(\Rotneg{\theta_1}\Big)}
{\ddotsc}
%{\Big(\Rotneg{\theta_q}\Big)}}
{\Big(\Rotneg{\theta_n}\Big)}}
\!\!\right)$\llabel{theta}
\hspace*{\fill}%
\\
with $\;\theta_1,\dots,\theta_n \in \left] -2\pi,0 \right]$.
\smallskip

{\bf (f)}
Let $\hat{x} \in Mp(V)$ be of elliptic component $\hat{x}_e$ lying over a semisimple element $x$ of $Sp(V)$ with elliptic component $x_e$.
We set\\[1mm]
\hspace*{\fill}%
$\Phi(\hat{x})
\,=\; \frac
{\scalo (\hat{x}_e)_{(1-x_e)\cdot V}}
{\scalo (B)_{(1-x_e)\cdot V}}\;\,
\mi^{-\frac{1}{2}\dim(1-x_e)\cdot V}\;\,
|\det \;(1-x)_{_{\scriptstyle (1-x)\cdot V}}|^{-1/2}$\labind{Φ}.
\hspace*{\fill}%
\par
\end{Definition}
\medskip

We will now see, as Vergne suggests in \cite[Prop. p.$\!\!$~289]{Ve94} (with the Weil representation contragredient to the one used in \cite{DHV84}, derived from that of \cite{DHV84} by replacing, at will $B$ with $-B$ \,or\, the central character of the Heisenberg group with its conjugate), that the functions introduced in points (d), (e), (f) are the functions $\rho_\call$ and $\delta$ of \cite{Df82a} and the function $\Phi$ of \cite{DHV84}.%

%------------------------------------------------------------------------------
\begin{Proposition}
\hlabel{morphisme rho et fonction delta}
%------------------------------------------------------------------------------

Consider a finite-dimensional real vector space $V$ endowed with a nondegenerate alternating bilinear form $B$.
\smallskip

{\bf (a)}
The cover $\,Mp(V)\to Sp(V)$ is connected when $V \not= \{0\}$.
\smallskip

{\bf (b)}
Let $\,\call$ be a Lagrangian of $(V_\CC,B)$.
The function $\rho_\call$ is a Lie group homomorphism from $Mp(V)_\call$ to $\CC\moins0$ such that $\;\diff_1 \rho_\call = \frac{1}{2}\,\tr (\,\cdot\,^\CC)_\call$.
For every $\hat{x} \in Mp(V)_\call$ lying over a semisimple element $x$ of $Sp(V)$, we have\\[1mm]
\hspace*{\fill}%
$\Phi(\hat{x})
= (-1)^{n_\call(x)\,+\, q_\call(x)}\;
\rho_\call(\hat{x})\;\,
(\det \;(1-x^\CC)_{_{\scriptstyle (1-x^\CC)\cdot \call}})^{-1}$.
\hspace*{\fill}%
\smallskip

{\bf (c)}
Let $\hat{x} \in Mp(V)$ lie over an $x \in Sp(V)$ with elliptic component~$x_e$.
There exists a positive Lagrangian of $(V_\CC,B)$ stable under $x^\CC$.
For every Lagrangian $\,\call$ of $(V_\CC,B)$ stable under $x^\CC$, we have\\
\hspace*{\fill}%
$\rho_\call(\hat{x})\;
\abs{\rho_\call(\hat{x})}^{-1}
\,=\; \delta(\hat{x})\;
\Prod_{z \in \Sp(x_e)}
z^{q_z}$,
\hspace*{\fill}\\
where $\;\Sp(x_e)$ is the spectrum of $x_e$ and, for each $z \in \Sp(x_e)$, $(p_z,q_z)$ is the signature $(n_+,n_-)$ of the Hermitian form $\,(v,w)\!\mapsto\! \mi \,B(v,\overline{w})\,$ on the eigenspace of $({x_e\!}^\CC)_\call$ associated with $z$.%
\vspace*{-1mm}
\end{Proposition}

\begin{Dem}{Proof of the proposition}

{\bf (a)}
Assume $V$ is nonzero.
Fix a symplectic basis $(P_1,\dots,P_n,Q_1,\dots,Q_n)$ of $(V,B)$.
Let $T_{Sp}$ denote the maximal torus of $Sp(V)$ consisting of endomorphisms of $V$ whose matrices in the basis $(P_1,Q_1,\dots,P_n,Q_n)$ are block diagonal with blocks in $SO(2)$, and let $\gotht$ denote its Lie algebra.
The torus $\,T_{Mp} \egdef \exp_{\!Mp(V)} \! \gotht\,$ is a maximal torus of $DL(V)$.

Choose an element $A$ of $\Ker \exp_{T_{Sp}}\!$ such that $\;\scalo (\exp_{\!DL(V)} \! A)_{\{0\}} = \RR_- \! \moins0$ (it exists).
The path $\,t \in [0,1] \mapsto \exp_{T_{Mp}} \! (t A)\,$ connects, in $Mp(V)$ the two points of the kernel of the canonical group homomorphism from $Mp(V)$ to $Sp(V)$.
\smallskip

{\bf (b)}
Let $\hat{x} \in Mp(V)_\call$ lie over a semisimple element $x$ of $Sp(V)$ with elliptic component $x_e$.

Let $\wt{\rho}_\call$ denote the Lie group homomorphism from $Mp(V)_\call$ to $\CC \moins0$ defined in \cite[p.$\!\!$~107]{Df84}.
According to \cite[(10) p.$\!\!$~108]{Df84}, we have\\
\hspace*{\fill}%
$\abs{\wt{\rho}_\call(\hat{x})}
= \abs{\det(x^\CC)_\call}^{1/2}
= \abs{\rho_\call(\hat{x})}\quad$
and
$\quad\diff_1 \wt{\rho}_\call
= \frac{1}{2}\,\tr (\,\cdot\,^\CC)_\call$.
\hspace*{\fill}\\
These equalities, in view of \cite[end of Lem. 30 p.$\!\!$~ 37]{DV93}, allow us to prove that the maps $\wt{\rho}_\call$ and $\rho_\call$ commute with taking elliptic, positively hyperbolic, and unipotent components on the groups $Mp(V)_\call$ and $GL(\CC)$, and that they agree on each of the subsets of $Mp(V)_\call$ formed by its positively hyperbolic elements or its unipotent elements.

Let $\wt{\Phi}$ denote the function (independent of $\call$) from the set of semisimple elements of $Mp(V)$ to $\CC \moins0$ defined in \cite[p.$\!\!$~102]{DHV84}.
It is written as\\[1mm]
\hspace*{\fill}%
$\wt{\Phi}(\hat{x})
= (-1)^{n_\call(x)\,+\, q_\call(x)}\;
\wt{\rho}_\call(\hat{x})\;\,
(\det \;(1-x^\CC)_{_{\scriptstyle (1-x^\CC)\cdot \call}})^{-1}$.
\hspace*{\fill}\\
Furthermore, the map $\,\dot{v} \mapsto B(v,.)\,$ from $\,(1-x^\CC)\cdot V_\CC\,/\, (1-x^\CC)\cdot \call\,$ to $\,((1-x^\CC)\cdot \call)^*\,$ is a linear bijection that commutes with the action of $\,1 \!-\! x$.
It follows that\\[1mm]
\hspace*{\fill}%
$|\wt{\Phi}(\hat{x})|
= |\det \;(1-x)_{_{\scriptstyle (1-x)\cdot V}}|^{-1/2}
= \abs{\Phi(\hat{x})}$.%
\hspace*{\fill}%

To show that $\,\wt{\Phi}=\Phi\,$ (respectively $\,\wt{\rho}_\call=\rho_\call$), given \cite[(9) p.$\!\!$~108]{Df84} it remains to prove that $\,\wt{\Phi}/|\wt{\Phi}|\,$ and $\,\Phi/\abs{\Phi}\,$ (respectively $\,\wt{\rho}_\call\,$ and $\,\rho_\call$) coincide at one of the two points of $Mp(V)$ lying over $x$ (respectively $x_e$).
\smallskip

Let $V^z$ (respectively $V_\CC^z$) denote the eigenspace of $x$ (respectively $x^\CC$) associated with an eigenvalue $z \!\in\! \RR \moins0$ (respectively $z \!\in\! \CC \moins0$), and let $\,{\scriptstyle \perp}^{\!\scriptscriptstyle B}\,$ denote the relation of $B$-orthogonality.
Now consider an eigenvalue $z \in \CC \moins0$ of~$x^\CC$ such that $\,\Im z \geq 0\,$ and $\,\abs{z} \geq 1$.
We will associate to it a certain symplectic basis $\calb_z$ of $\;(( V_\CC^z + V_\CC^{z^{-1}} \!\!+ V_\CC^{\bar{z}} + V_\CC^{{\bar{z}}^{-1}} )\cap V,B)$.

\quad
If $\,z \in \RR\,$ and $\,\abs{z} \geq 1$:\
when $\,\abs{z} = 1$, fix a symplectic basis $\;\calb_z = (P_1^0,\dots,P_{n^0}^0,Q_1^0,\dots,Q_{n^0}^0)\;$ of $V^z$;
when $\,\abs{z} > 1$, fix a basis $(P_1^0,\dots,P_{n^0}^0)$ of $V^z$, and deduce from it a unique symplectic basis
$\,\calb_z = (P_1^0,\dots,P_{n^0}^0,Q_1^0,\dots,Q_{n^0}^0)\,$
of $\,V^z \oplus V^{z^{-1}}$ such that
$\,Q_1^0,\dots,Q_{n^0}^0 \in V^{z^{-1}}\!$,
because $V^{z^{-1}}$ identifies with $(V^z)^*$ by means of $B$ since $\,V^z \,{\scriptstyle \perp}^{\!\scriptscriptstyle B}\, V^z\,$ and $\,V^{z^{-1}} {\scriptstyle \perp}^{\!\scriptscriptstyle B}\, V^{z^{-1}}$.

\quad
If $\,\Im z > 0\,$ and $\,\abs{z} > 1$:\
fix a basis
$\;(P_1 + \mi P_2,\dots,P_{2n-1} + \mi P_{2n})\;$
of $V_\CC^z$ such that $\,P_1,\dots,P_{2n} \in V$;
we deduce from it a unique symplectic basis
$\;(P_1 + \mi P_2,\dots,P_{2n-1} + \mi P_{2n},
Q_1 -\mi Q_2,\dots,Q_{2n-1} -\mi Q_{2n})\;$
of $\,V_\CC^z \oplus V_\CC^{z^{-1}}\!$ endowed with $B/2$ such that
$\,Q_1 \!-\! \mi Q_2,\dots,Q_{2n-1} \!-\! \mi Q_{2n} \in V_\CC^{z^{-1}}\,$
and $\,Q_1,\dots,Q_{2n} \in V\!$,
because $V^{z^{-1}}$ identifies with $(V^z)^*$ by means of $B$ since
$\;V_\CC^z \,{\scriptstyle \perp}^{\!\scriptscriptstyle B}\, V_\CC^z\;$
and $\;V_\CC^{z^{-1}} {\scriptstyle \perp}^{\!\scriptscriptstyle B}\, V_\CC^{z^{-1}}$;
we~then \ obtain \ the \ symplectic \ basis $\;\calb_z \egdef (P_1,\dots,P_{2n},Q_1,\dots,Q_{2n})\;$ of
$\,( V_\CC^z \oplus V_\CC^{z^{-1}}\! \oplus V_\CC^{\bar{z}} \oplus V_\CC^{{\bar{z}}^{-1}} )\cap V\!$,
because $\;V_\CC^z \,{\scriptstyle \perp}^{\!\scriptscriptstyle B} \,V_\CC^{\bar{z}}$,
$\;V_\CC^{z^{-1}} \,{\scriptstyle \perp}^{\!\scriptscriptstyle B}\, V_\CC^{\bar{z}^{-1}}\,$
and $\;V_\CC^z \,{\scriptstyle \perp}^{\!\scriptscriptstyle B}\, V_\CC^{\bar{z}^{-1}}\!$.

\quad
If $\,\Im z > 0\,$ and $\,\abs{z} = 1$:\
fix a basis
$(P'_1 \!+\! \mi Q'_1,\dots,P'_p \!+\! \mi Q'_p,P''_1 \!-\! \mi Q''_1,$ $\!\!\dots,P''_q \!-\!\mi Q''_q)$
of $V_\CC^z$ such that $\,P'_1,Q'_1,\dots,P'_p,Q'_p,P''_1,Q''_1,\dots,P''_q,Q''_q \in V\!$,
in which the nondegenerate Hermitian form $\;(v,w)\mapsto iB/2 \, (v,\overline{w})\;$ has matrix
$\left(
{\scriptstyle \begin{smallmatrix} I_p & 0 \\ 0 & -I_p \end{smallmatrix}}
\right)$;
it~yields the symplectic basis
$\calb_z \egdef (P'_1,\dots,P'_p,P''_1,\dots, $ $\!\!P''_q,Q'_1,\dots,Q'_p,Q''_1,\dots,Q''_q)$
of $\;(V_\CC^z \oplus V_\CC^{\bar{z}})\cap V\!$,
because
$\;V_\CC^z \,{\scriptstyle \perp}^{\!\scriptscriptstyle B}\, V_\CC^z\;$
and
$\;V_\CC^{\bar{z}} \,{\scriptstyle \perp}^{\!\scriptscriptstyle B}\, V_\CC^{\bar{z}}$.
\smallskip

Note that the vector subspace $\,\call_x^+$ (respectively $\,\call_{x,e}^+$) of $V_\CC$ spanned by the various vectors
$P_1^0,...,P_{n^0}^0$,
and,
$P_1,...,P_{2n}$,
and,
$P'_1 + \mi Q'_1,...,P'_p + \mi Q'_p,P''_1 + \mi Q''_1,...,P''_q + \mi Q''_q$
(respectively
$P_1^0 + \mi Q_1^0,...,P_{n^0}^0 + \mi Q_{n^0}^0$,
and,
$P_1 + \mi Q_1,...,P_{2n} + \mi Q_{2n}$,
and,
$P'_1 + \mi Q'_1,...,P'_p + \mi Q'_p,P''_1 + \mi Q''_1,...,P''_q+\mi Q''_q\,$)
as $z$ varies, is a positive Lagrangian of $(V_\CC,B)$ stable under~$x^\CC$ (respectively ${x_e\!}^\CC$).
\smallskip

Let $E$ (respectively $H$) denote the element of ${\goth sp} (V)$ that stabilizes the vector spaces
$\;( V_\CC^z + V_\CC^{z^{-1}}\!\! + V_\CC^{\bar{z}} + V_\CC^{\bar{z}^{-1}} )\cap V\;$
with $z = r \me ^{\mi \,\theta}$, $r \geq 1$, $-2\pi < \theta \leq 0$,
and whose restriction to such a space has, in the basis $\calb_z$, a block diagonal matrix whose blocks are the matrices\\[1mm]
$\big(\!{\scriptstyle \Rotneg{\theta}}\big)$
(resp.
$\big(\!{\scriptstyle \begin{smallmatrix} \ln r & 0 \\ 0 & -\ln r \end{smallmatrix}}\!\big)$)
with respect to $(P_k^0,Q_k^0)$,\\[1mm]
or,
$\left(\!\!{\scriptstyle\diagd{\Rotneg{\theta}}{\Rot{\theta}}}\!\!\right)$
(resp.
$\left( {\scriptstyle \diagq{\ln r}{\ln r}{-\ln r}{-\ln r}} \right)$)
with respect to $(P_{2k-1},P_{2k},Q_{2k-1},-Q_{2k})$,\\[1mm]
or,
$\big(\!{\scriptstyle \Rotneg{\theta}}\big)$
(resp. $\big( \begin{smallmatrix} 0 & 0 \\ 0 & 0 \end{smallmatrix} \big)$)
with respect to $(P'_k,Q'_k)$ and
$\big({\scriptstyle \Rot{\theta}}\!\big)$
(resp. $\big( \begin{smallmatrix} 0 & 0 \\ 0 & 0 \end{smallmatrix} \big)$)
with respect to $(P''_l,Q''_l)$.

The elements $\,\hat{x}_{E,H} \!=\! \exp_{\!Mp(V)} \! E \, \exp_{\!Mp(V)} \! H$ and $\,\hat{x}_E \!=\! \exp_{\!Mp(V)} \! E$ of $Mp(V)$ lie respectively over $x$ and $x_e$.
By \cite[(15) p.$\!\!$~109]{Df84}, we have:\\
\hspace*{\fill}%
${\displaystyle\frac{\wt{\rho}_{\call_x^+}}{|\wt{\rho}_{\call_x^+}|}}
(\hat{x}_{E,H})
= \wt{\rho}_{\call_{x,e}^+}(\exp_{\!Mp(V)} \! E)
=\, \me ^{\frac{1}{2} \tr (E^\CC)_{\call_{x,e}^+}}$.
\hspace*{\fill}\\[1mm]
Hence we have the following equalities, which yield the result:\\[1mm]
\hspace*{\fill}%
${\displaystyle\frac{\wt{\Phi}}{|\wt{\Phi}|}}(\hat{x}_{E,H})
= (-1)^{n_{\call_x^+}(x)}\; \me ^{\frac{1}{2}\,\tr (E^\CC)_{\call_{x,e}^+}}\;
\biggl(
\frac{\det \;(1-x^\CC)_{(1-x^\CC)\cdot \call_x^+}}
{\abs{\det \;(1-x^\CC)_{(1-x^\CC)\cdot \smash{\call_x^+}}}}
\biggr)^{\!\!-1} \!
= \cdots
= {\displaystyle\frac{\Phi}{\abs{\Phi}}}(\hat{x}_{E,H})$
\hspace*{\fill}\\
\hspace*{\fill}%
and
$\quad \wt{\rho}_\call(\hat{x}_E)
\!= (-1)^{q_\call(x_e)}\;
{\displaystyle\frac{\Phi}{\abs{\Phi}}}(\hat{x}_E)\;
\frac{\det \;(1-\exp E^\CC)_{(1-\exp E^\CC)\cdot \call}}
{\abs{\det \;(1-\exp E^\CC)_{(1-\exp E^\CC)\cdot \call}}}
= \cdots
= \rho_\call(\hat{x}_E).$
\hspace*{\fill}%
\medskip

{\bf (c)}
The existence of a positive Lagrangian of $(V_\CC,B)$ stable under $x^\CC$ is proven in \cite[Cor. p.$\!\!$~82]{B.72}.
Assume $V$ is nonzero and denote by $\hat{x}_e$ the elliptic component of $\hat{x}$.
By definition $\delta(\hat{x})$ is equal to $\rho_{\CC(P_1 \!+ \mi Q_1) +\cdots+ \CC(P_n \!+ \mi Q_n)} \! (\hat{x}_e)$, with the notation of \ref{géométrie métaplectique} (e).
The equality \cite[(15) p.$\!\!$~109]{Df84} proves that this complex number is independent of the choice of $(P_1,\dots,P_n,Q_1,\dots,Q_n)$.
\smallskip

Consider a Lagrangian $\call$ of $(V_\CC,B)$ stable under $x^\CC$ and a complement $\cals$ of $\,\call \cap \overline{\call}\,$ in $\call$ that is stable under ${x_e\!}^\CC$.
Let $W$ denote the symplectic subspace of $(V,B)$ equal to the orthogonal complement of $\,\cals + \overline{\cals}\,$ in $V$.
Fix a basis $\,(P_1^0,\dots,P_{n^0}^0,P_1^\RR,\dots,P_{2m}^\RR)\,$ of the Lagrangian $\,\call \cap \overline{\call}\,$ of $(W_\CC,B)$, consisting of vectors of $V$, such that $\,P_1^0,\dots,P_{n^0}^0\,$ are eigenvectors of $x_e$ with real eigenvalues and $\,P_1^\RR+ \mi P_2^\RR,\dots,P_{2m-1}^\RR+ \mi P_{2m}^\RR\,$ are eigenvectors of ${x_e\!}^\CC$ with nonreal eigenvalues.
It extends to a symplectic basis $(P_1^0,\dots,P_{n^0}^0,P_1^\RR,\dots,P_{2m}^\RR,Q_1^0,\dots,Q_{n^0}^0,2Q_1^\RR,\dots,2Q_{2m}^\RR)$ of $(W,B)$ such that $Q_1^0,\dots,Q_{n^0}^0$ and $Q_1^\RR+ \mi Q_2^\RR,\dots,Q_{2m-1}^\RR+ \mi Q_{2m}^\RR$ are eigenvectors of~${x_e\!}^\CC$.
Indeed, we can construct $\,(Q_1^\RR,\dots,Q_{2m}^\RR)\,$ by induction on $m = \frac{1}{2} \dim (x_e^2-1)\cdot W$, associating to the vector $\,P \egdef P_{2m-1}^\RR+ \mi P_{2m}^\RR\,$ of $\,(x_e^2-1)\cdot (\call \cap \overline{\call})\,$ an eigenvector $\,Q = Q_{2m-1}^\RR- \mi Q_{2m}^\RR\,$ of ${x_e\!}^\CC$ with $\,Q_{2m-1}^\RR,Q_{2m}^\RR\in (x_e^2-1)\cdot W\,$ such that $B(P,Q) = 1$ and $B(Q,\overline{Q}) = 0$ (the latter equality is achieved by replacing $Q$ with $\,Q - \frac{1}{2} B(Q,\overline{Q})\, \overline{P}$).
Also fix a basis $\,(P'_1 + \mi Q'_1,$ $\dots,P'_p + \mi Q'_p,P''_1 - \mi Q''_1,\dots,P''_q -\mi Q''_q)\,$ of $\,\cals$ consisting of eigenvectors of ${x_e\!}^\CC$ with $\,P'_1,Q'_1,\dots,P'_p,Q'_p,P''_1,Q''_1,\dots,P''_q,Q''_q \in V$, in which the nondegenerate ${x_e\!}^\CC$-invariant Hermitian form $\,(v,w)\mapsto iB/2 \, (v,\overline{w})\,$ has matrix~%
$\left(\begin{smallmatrix} I_p & 0 \\ 0 & -I_p \end{smallmatrix}\right)$.
\smallskip

The basis of $\call$ consisting of the eigenvectors $\,P_1^0,...,P_{n^0}^0,P_1^\RR+ \mi P_2^\RR,$ $P_1^\RR- \mi P_2^\RR,$ $...,P_{2m-1}^\RR+ \mi P_{2m}^\RR,P_{2m-1}^\RR- \mi P_{2m}^\RR$ and $P'_1 + \mi Q'_1,...,P'_p + \mi Q'_p,P''_1 - \mi Q''_1,...,P''_q -~\mi Q''_q\,$ of~${x_e}^\CC$ yields an expression for $\,\rho_\call(\hat{x}_e)$.
We calculate $\delta(\hat{x})$ using the symplectic basis $\;(P_1^0,\dots,P_{n^0}^0,P_1^+,\dots,P_{2m}^+,P'_1,\dots,P'_p,P''_1,\dots,P''_q,Q_1^0,\dots,Q_{n^0}^0,Q_1^+,\dots,Q_{2m}^+,$ $Q'_1,\dots,Q'_p,Q''_1,\dots,Q''_q)\;$ of $(V,B)$,\\
where $\;P_{2k-1}^+ + \mi Q_{2k-1}^+ \egdef (P_{2k-1}^\RR+ \mi P_{2k}^\RR) + \mi (Q_{2k-1}^\RR+ \mi Q_{2k}^\RR)\;$\\
and $\;P_{2k}^+ - \mi Q_{2k}^+ \egdef -\mi (P_{2k-1}^\RR+ \mi P_{2k}^\RR) - (Q_{2k-1}^\RR+ \mi Q_{2k}^\RR)\;$ for $1 \leq k \leq m$.
\smallskip

From this we can easily deduce the result.
\cqfd
\end{Dem}

%-------------------------------------------------------------------------------
\section{The parameters \texorpdfstring{$\,\tau \in \XInd_G(\lambdatilde)$}{\tau \in X\^{}Ind\_G(\~\lambda)}}\label{Rep5}
%------------------------------------------------------------------------------

%\noindent
In the following two sections, fix $\,\lambdatilde = (\lambda,\F+)\in \gstregtilde$.
Set $\,\gothh = \gothg(\lambdatilde)$.
Let $\gotha$ denote the hyperbolic component of $\gothh$ and fix a chamber ${\gotha^*}^+\!$ of $(\gothg(\lambda)(\mi \rho_\F+),\gotha)$.
Let $\mu$ and $\nu$ denote the infinitesimally elliptic and hyperbolic components of $\lambda$.
\smallskip

In order to provide a canonical parametrization of a subset of $\widehat{G}$ that allows us to describe the Plancherel formula, Duflo introduced the objects in (a) of the definition below.
When $\lambda\in\gstssreg$, these objects will be compatible with mine, because we can take $t \!\in\! \left[0,1\right]$ by setting $\lambda_0 \!=\! \lambda$ in Lemma \ref{fonctions canoniques sur le revêtement double}, and we will have $\,\XInd_G(\lambdatilde) = \Xirr_G(\lambda)$ in Definition \ref{paramètres adaptes} (c).

%------------------------------------------------------------------------------
\begin{Definition}\rm
\hlabel{paramètres de Duflo}
%------------------------------------------------------------------------------

{\bf (a)}
Let $f \in \gothg^*\!$.
We denote by\\
$B_f$\labind{Bf} the nondegenerate alternating bilinear form $\,(\dot{X},\dot{Y})\mapsto f([X,Y])\,$ on $\,\gothg / \gothg(f)$,\\
$G(f)^{\gothg/\gothg(f)}\,$\labind{G(f)g/g(f)} the Lie subgroup of $\,G(f)\times DL(\gothg / \gothg(f))\,$ consisting of the ordered pairs $(x,\hat{a})$ such that $(\Ad x)_{\gothg / \gothg(f)}$ is the image of $\hat{a}$ in $SL(\gothg / \gothg(f))$ (see \ref{morphisme rho et fonction delta} (a)),\\
$\{1,\iota\}$\labind{ι} and $G(f)^{\gothg/\gothg(f)}_0\!$\labind{G(f)g/g(f)0} the inverse images of $\{1\}$ and $G(f)_0$ under the canonical surjective Lie group homomorphism from $G(f)^{\gothg/\gothg(f)}\!$ to $G(f)$,\\
$\Xirr_G(f)$\labind{Xirr} the finite set of isomorphism classes of irreducible unitary representations $\tau$ of $G(f)^{\gothg/\gothg(f)}$ such that
$\tau(\iota) = -\id$ and $\,\tau(\exp X) = \me ^{\mi \,f(X)}\id\,$
for $X \!\in\! \gothg(f)$ (in this case $\tau$ is finite-dimensional).
%\smallskip

{\bf (b)}
Denote by $G(\lambdatilde)^{\gothg/\!\gothg(\lambda)(\mi \rho_\F+)}\!$\labind{G(λtilde)g/glambdairho} the inverse image of $G(\lambdatilde)$ in $G(\lambda_{\textit{can}})^{\gothg/\gothg(\lambda_{\textit{can}})}\!\!$.
\\
We will also set $\,G(\lambdatilde)^{\gothg/\gothh}_0 \!= G(\lambda_+)^{\gothg/\gothh}_0$\labind{G(λtilde)g/h0}
(see $\;G(\lambdatilde)_0 = G(\lambda_+)_0$ and $\gothg(\lambda_+) = \gothh$).

When $\,{\gotha^*}^+ = \gotha^*$, we have $\;\gothg(\lambda)(\mi \rho_\F+) = \gothh\,$ and $\,G(\lambdatilde) = G(\lambda_+)$, and in this case we will denote $G(\lambdatilde)^{\gothg/\!\gothg(\lambda)(\mi \rho_\F+)}\!$ and $\,G(\lambda_+)^{\gothg/\gothh}\!$ by $\,G(\lambdatilde)^{\gothg/\gothh}\!$\labind{G(λtilde)g/h}.
\end{Definition}
\medskip

In what follows, we will identify complex functions on $\,G(\lambda_+)\,$ with complex functions on the cover $\,G(\lambda_+)^{\gothg/\gothh}\,$ that are constant on the fibers.

%------------------------------------------------------------------------------
\begin{Lemme}
\hlabel{fonctions canoniques sur le revêtement double}
%------------------------------------------------------------------------------

Fix a real number $\,t \!\in\! \left]0,1\right]$.
We use the positive root systems from Definition \ref{choix de racines positives} (b) and the notation $\lambda_t$ from Lemma \ref{lemme clef}.\\
We set
$\;\,\call_{\lambda_t}
= \gothh_\CC\oplus
\Sumpetit_{\alpha \in R^+(\gothg_\CC,\gothh_\CC)} \gothg_\CC^{\alpha}\;$
and, $\,\gothm' = \gothg(\nu_+)\,$ in anticipation of Part~\ref{RepIII}

{\bf (a)}
We have $G(\lambda_t)^{\gothg/\gothh} \!=\! G(\lambda_+)^{\gothg/\gothh}\!$ and $\call_{\lambda_t}/\gothh_\CC$ is a Lagrangian of $(\gothg_\CC/\gothh_\CC,\!B_{\lambda_t})$ stable under $G({\lambda_t})$. The complex function $\,\delta_{\lambda_t}^{\gothg/\gothh}$ on $\,G(\lambda_+)^{\gothg/\gothh}$ and the Lie group homomorphism $\;\rho_{\lambda_t}^{\gothg/\gothh}$ from $\,G(\lambda_+)^{\gothg/\gothh}\,$ to $\,\CC \moins0$ deduced from it (see \ref{géométrie métaplectique} (e) and \ref{morphisme rho et fonction delta} (b)) are independent of $t$, and therefore equal to $\,\delta_{\lambda_+}^{\gothg/\gothh}\!$\labind{δ_λ+} and $\,\rho_{\lambda_+}^{\gothg/\gothh}\!$\labind{ρ λ+} respectively.\\
We have $\;\,\diff_1 \rho_{\lambda_+}^{\gothg/\gothh} = \rho_{\gothg,\gothh}$.
\smallskip

{\bf (b)}
Let $\hat{e} \in G(\lambda_+)^{\gothg/\gothh}$ lie over an elliptic element $e$ of $G(\lambda_+)$.
We have\\[1mm]
\hspace*{\fill}%
$\rho_{\lambda_+}^{\gothg/\gothh}(\hat{e})
\;=\;
\delta_{\lambda_+}^{\gothg/\gothh}(\hat{e})\;{\scriptstyle \times}\;
\det\left(\Ad e^\CC\right)\!\!\!_{\Sumpetit_{\scriptstyle\alpha' \in R^+_K(\gothg_\CC,\gothh_\CC)}
\!\!\gothg_\CC^{\alpha'}}
\;\,{\scriptstyle \times}\,
\Prod_{\calo \in \langle e \rangle \backslash
\wt{R}^+_\CC(\gothm'_\CC,\gothh_\CC)}
\!\!\! \! (-1)^{m_\calo-1} \, u_\calo$
\hspace*{\fill}%
\\
where $R^+_K(\gothg_\CC,\gothh_\CC)$ is the set of compact $\alpha' \!\in\! R^+(\gothg_\CC,\gothh_\CC)$, $\wt{R}^+_\CC(\gothm'_\CC,\gothh_\CC)$ is the set of classes of complex $\beta \!\in\! R^+(\gothm'_\CC,\gothh_\CC)$ modulo the identification of $\beta$ with $-\overline{\!\beta}$, and, for each orbit $\calo$ of an element $\{\beta,\!-\overline{\!\beta}\}$ of $\wt{R}^+_\CC(\gothm'_\CC,\gothh_\CC)$ under the action of the subgroup $\langle e \rangle$ of $G$ generated by $e$, we denote by $m_\calo$ the cardinality of $\calo$ and by~$u_\calo$ the unique eigenvalue of $(\Ad e^\CC)^{m_\calo}$ associated with an eigenvector that can be written $X_\beta - \overline{\!X_{-\beta}\!}$ for at least one $X_\beta \in \gothg_\CC^{\beta}$ and one $X_{-\beta} \in \gothg_\CC^{-\beta}$ satisfying~$[X_\beta,X_{-\beta}] = H_\beta$.%
%\vspace*{-2mm}
\end{Lemme}

\begin{Dem}{Proof of the lemma}

{\bf (a)}
The first equality follows from Lemma \ref{lemme clef}.
We can see immediately that $\call_{\lambda_t} / \gothh_\CC$ is a Lagrangian of $(\gothg_\CC/\gothh_\CC,B_{\lambda_t})$ stable under~$G({\lambda_t})$.%

Let $\varphi_t$ denote the Hermitian form $\,(v,w)\mapsto \mi \,B_{\lambda_t}(v,\overline{w})$ on $(\gothg/\gothh)_\CC$.
We will use the notation from statement (b).
The vector space $\,\call_{\lambda_t} / \gothh_\CC$ is the direct sum of the projections of the following vector spaces, which are pairwise orthogonal for $\varphi_t$ with the signature of the restricted Hermitian form $(n_+,n_-)$ specified in parentheses:
the $\,\gothg_\CC^{\alpha'}$ with $\alpha' \!\in\! R^+_K(\gothg_\CC,\gothh_\CC)$
(signature $(0,1)$),
the $\,\gothg_\CC^{\alpha''}$ with $\alpha'' \!\in\! R^+(\gothg_\CC,\gothh_\CC)$ imaginary noncompact
(signature $(1,0)$),
the $\,\gothg_\CC^{\beta}\oplus\gothg_\CC^{-\,\overline{\!\beta}}$ with $\,\{\beta,\!-\,\overline{\!\beta}\} \in \wt{R}^+_\CC(\gothm'_\CC,\gothh_\CC)\,$
(signature $(1,1)$),
and the $\,\gothg_\CC^{\gamma}$ with $\gamma \!\in\! R^+(\gothg_\CC,\gothh_\CC)$ outside $R^+(\gothm'_\CC,\gothh_\CC)$
(signature $(0,0)$).

Let $e$ be an elliptic element of $G(\lambda_+)$.
By the above, the Hermitian form $\varphi_t$ is nondegenerate on $(\call_{\lambda_t} \cap \gothm'_\CC) / \gothh_\CC$ and vanishes everywhere on its orthogonal.
Furthermore, $(\call_{\lambda_t} \cap \gothm'_\CC) / \gothh_\CC$ is an orthogonal direct sum, with respect to $\varphi_t$, of the image and kernel of $\,1 - \Ad e^\CC$.
Therefore, the restriction of $\varphi_t$ to each of them is nondegenerate.
A~continuity argument then proves that the signatures of these restrictions, and, \emph{a fortiori} $q_{\call_{\lambda_t}}(e)$, are independent of $t$.
Likewise, the orientation $\,\scalo (B_{\lambda_t})_{(1-\Ad e)\cdot (\gothg/\gothh)}\,$ is independent of $t$.
This can be clarified as in the proof of \ref{espace symplectique canonique}.
We then use Definition \ref{géométrie métaplectique} (d) and Proposition \ref{morphisme rho et fonction delta} (c).
\vskip-0.5mm

The formula for $\,\diff_1 \rho_{\lambda_+}^{\gothg/\gothh}\,$ is deduced from Proposition \ref{morphisme rho et fonction delta} (b).
\smallskip

{\bf (b)}
It is a matter of applying Proposition \ref{morphisme rho et fonction delta} (c) taking into account the proof in (a).
We use the notation $\varphi_t$ from that proof.

Let $\calo$ be the orbit under $\langle e \rangle$ of a $\{\beta,\!-\,\overline{\!\beta}\} \!\in\! \wt{R}^+_\CC(\gothm'_\CC,\gothh_\CC)$.
Fix $(X_\beta,X_{-\beta})\!\in\! \gothg_\CC^{\beta} \,{\scriptstyle \times}\, \gothg_\CC^{-\beta}$ such that $\;[X_\beta,X_{-\beta}] = H_\beta$.
The condition “$X_\beta - \overline{\!X_{-\beta}\!}\,$ is an eigenvector of $\,(\Ad e^\CC)^{\,m_\calo}$” is satisfied (for a unique eigenvalue of $\,(\Ad e^\CC)^{\,m_\calo}$) when $\,e^{\,m_\calo}\,\beta = \beta$, and is realized (again with uniqueness of the eigenvalue) by replacing $X_\beta$ and $X_{-\beta}$ with certain multiples of them when $\,e^{\,m_\calo}\,\beta = -\,\overline{\!\beta}$.
Under this condition, the eigenvectors
$\,\Sumpetit_{0 \leq k \leq m_\calo-1}\!\! \zeta^{-k} (\Ad e^\CC)^k \, (X_\beta \!- \overline{\!X_{-\beta}\!})\,$
of $\Ad e^\CC\!$, where $\zeta$ ranges over the set of $m_\calo^{\textrm{th}}\!$ roots of $u_\calo$,
have projections in $\,\call_{\lambda_t} / \gothh_\CC$ that are pairwise orthogonal for $\varphi_t$ with (strictly) negative “squares.”
This leads to the stated formula.
\cqfd
\end{Dem}
\vspace*{-5mm}

%------------------------------------------------------------------------------
\begin{Remarque}\rm
\hlabel{sous-algèbre de Borel canonique}
%------------------------------------------------------------------------------

Let $\,f \in \Supp_{\gothg(\lambda)^*}(G(\lambda)_0 \cdot \lambdatilde)$.

Let $\xi$ denote the nilpotent component of $f$.
There exists a unique Borel subalgebra $\gothb$ of $\,\gothg(\lambda)_\CC\,$ on which $\xi$ is zero.
Indeed, by Lemma \ref{description d'ensembles de formes linéaires} (a) and \cite[Cor. 5.3 p.$\!\!$~997 and Cor. 5.6 p.$\!\!$~1001] {Ko59}, the element of $\derivealg{\gothg(\lambda)}$ with which $\xi\,$ is identified using the bilinear form ${\scriptstyle \langle} ~,\!~ {\scriptstyle \rangle}$ belongs to a unique maximal subalgebra consisting of nilpotent elements of $\gothg(\lambda)_\CC$, whose orthogonal complement with respect to ${\scriptstyle \langle} ~,\!~ {\scriptstyle \rangle}$ is suitable.
The uniqueness of $\gothb$ shows that $\;\overline{\gothb} = \gothb$.

Let $\,\gothn_{\lambda}\,$ denote the sum of the eigenspaces for the action of $\,\centrealg{\gothg(\lambda)}_\CC\,$ on~$\gothg_\CC$, whose weight $\,\alpha \in \centrealg{\gothg(\lambda)}_\CC^*\,$ (always a restriction of a weight of $\gothh_\CC$) satisfies:
$\,\langle \nu,\alpha \rangle >0\,$ or $\,(\langle \nu,\alpha \rangle =0\;$ and $\;\mi\langle \mu,\alpha \rangle >0)$.
Set $\,\call_f = \gothb \oplus \gothn_{\lambda}$.
Note that $\call_f$ is a Borel subalgebra of $\gothg_\CC$ by decomposing $\call_f$ into eigenspaces relative to a Cartan subalgebra of $\gothg_\CC$ included in $\gothb$.

The inclusions $\;[\gothb,\gothb] \subseteq \gothb \subseteq \Ker \xi\;$ and
$\;[\call_f,\gothn_{\lambda}] \subseteq \gothn_{\lambda}
\subseteq (\gothg(\lambda)_\CC)^{{\scriptstyle \perp}^
{\!{\scriptscriptstyle \langle} \!~,\!~ {\scriptscriptstyle \rangle}}}\,$
provide the equalities $\;f([\gothb,\gothb]) = 0\;$ and $\;f([\call_f,\call_f]) = 0$.
Therefore $\,\gothg(f)_\CC\subseteq \gothb\,$ and $\,\call_f / \gothg(f)_\CC$ is a Lagrangian of $(\gothg_\CC/ \gothg(f)_\CC,B_f)$, clearly stable under $G(f)$.
By~\ref{morphisme rho et fonction delta}~(b)~the Lie group homomorphism $\rho_f^{\gothg / \gothg(f)}$ from $G(f)^{\gothg/\gothg(f)}$ to $\,\CC\moins0$ deduced from it satisfies
$\;\,\diff_1 \rho_f^{\gothg / \gothg(f)} = \frac{1}{2} \tr ({\ad \cdot\,}^\CC)_{\call_f / \gothg(f)_\CC}$.
\cqfr
\end{Remarque}
\medskip

I will now adapt Duflo's parameters to the situation that interests me.

%------------------------------------------------------------------------------
\begin{Lemme}
\hlabel{espace symplectique canonique}
%------------------------------------------------------------------------------

{\bf (a)}
The vector space
$\;\call_{\lambdatilde,{\gotha^*}^+}
= \gothh_\CC\oplus \Sumpetit_{\alpha \in \smash{R^+_{\lambdatilde,{\gotha^*}^+}}}\!\gothg_\CC^{\alpha}\;$
(see \ref{choix de racines positives} (c)) yields the Lagrangian $\,\call_{\lambdatilde,{\gotha^*}^+}/\gothh_\CC$ of $(\gothg_\CC/ \gothh_\CC,B_{\lambda_+})$.
Let $\rho_{\lambdatilde,{\gotha^*}^+}^{\gothg/\gothh}$\labind{ρ λtilde,a*+ g,h} denote the complex character of $G(\lambda_+)^{\gothg/\gothh}$ associated with it (see \ref{morphisme rho et fonction delta} (b)).
The restriction of $\rho_{\lambdatilde,{\gotha^*}^+}^{\gothg/\gothh} (\delta_{\lambda_+}^{\gothg/\gothh})^{-1}\!$ to $\centragr{G_0}{\gothh}$ is independent of ${\gotha^*}^+\!$.
\smallskip

{\bf (b)}
The restriction of $B_{\lambda_+}\!$ to the vector subspace $\gothg(\lambda)(\mi \rho_\F+)/\gothh$ of $\gothg/\gothh$ is nondegenerate.
The orthogonal complement of $\gothg(\lambda)(\mi \rho_\F+)/\gothh$ in $\gothg/\gothh$ endowed with $B_{\lambda_+}\!$ is the projection in $\gothg/\gothh$ of the intersection of $\gothg$ with the sum of $\gothg_\CC^{\alpha}$ as $\alpha$ ranges over $R(\gothg_\CC,\gothh_\CC) \setminus R(\gothg(\lambda)(\mi \rho_\F+)_\CC,\gothh_\CC)$.
We will identify it with $\gothg/\gothg(\lambda)(\mi \rho_\F+)$.

Let $e \in G(\lambda_+)$ be elliptic.
The orientations $\,\scalo (B_{\lambda_+})_{(1-\Ad e)(\gothg/\gothg(\lambda)(\mi\rho_\F+))}\,$
and $\,\scalo (B_{\lambda_{\textit{can}}})_{(1-\Ad e)(\gothg/\gothg(\lambda)(\mi\rho_\F+))}\,$
are equal.
\smallskip

{\bf (c)}
The vector space
$\;\call_{\lambdatilde}
= \gothg(\lambda)(\mi \rho_\F+)_\CC\oplus \Sumpetit_{\alpha \in R^+_\lambdatilde}\!\gothg_\CC^{\alpha}\;$
(see \ref{choix de racines positives} (c)) yields the Lagrangian $\call_{\lambdatilde}/\gothg(\lambda)(\mi \rho_\F+)_\CC\!$ of $\gothg_\CC/ \gothg(\lambda)(\mi \rho_\F+)_\CC$, both for $B_{\lambda_{\textit{can}}}$ and for~$B_{\lambda_+}\!$.
Let $\,\rho_\lambdatilde^{\gothg/\!\gothg(\lambda)(\mi \rho_\F+)}\!$\labind{ρ λtilde g,h}
denote the complex character of $G(\lambdatilde)^{\gothg/\!\gothg(\lambda)(\mi \rho_\F+)}\!$ associated with
$(\gothg_\CC/ \gothg(\lambda)(\mi \rho_\F+)_\CC,B_{\lambda_{\textit{can}}})$
(see \ref{morphisme rho et fonction delta} (b)).
It extends the complex character of the inverse image of $G(\lambda_+)$ in $G(\lambdatilde)^{\gothg/\!\gothg(\lambda)(\mi \rho_\F+)}\!$ associated with $(\gothg_\CC/ \gothg(\lambda)(\mi \rho_\F+)_\CC,B_{\lambda_+})$.%
\end{Lemme}
\vspace*{-5mm}

\begin{Dem}{Proof of the lemma}

{\bf (a)}
Let $\;e \in \centragr{G_0}{\gothh}$ be elliptic.
The independence of
$\big( \rho_{\lambdatilde,{\gotha^*}^+}^{\gothg/\gothh} (\delta_{\lambda_+}^{\gothg/\gothh})^{-1} \big) (e)$
with respect to ${\gotha^*}^+$ is obtained by resuming the proof of Lemma \ref{fonctions canoniques sur le revêtement double}~(b).
\smallskip

{\bf (b)}
At the outset, we use the fact that $\gothg_\CC\!/\gothh_\CC\!$ is the $B_{\lambda_+}\!$-orthogonal direct sum of the projections of the vector subspaces $\gothg_\CC^{\alpha} \oplus \gothg_\CC^{-\alpha}$ of $\gothg_\CC$ with $\alpha \in R^+(\gothg_\CC,\gothh_\CC)$.

Set $V = \gothg/\gothg(\lambda)(\mi\rho_\F+)$ and $V(e) = \gothg(e)/\gothg(e)(\lambda)(\mi\rho_\F+)$.
The vector space $V$ is the direct sum of the vector subspaces $V(e)$ and $(1-\Ad e)(V)$, which are simultaneously orthogonal for $B_{\lambda_{\textit{can}}}$ and for $B_{\lambda_+}\!$.
For each $\alpha \in R^+(\gothg_\CC,\gothh_\CC)$, fix
$X_{\alpha} \in \gothg_\CC^{\alpha}$ and
$X_{-\alpha} \in \gothg_\CC^{-\alpha}$ such that
$[X_{\alpha},X_{-\alpha}] = H_{\alpha}$.
Note that the orientations $\scalo (B_{\lambda_+})_V$ and $\scalo (B_{\lambda_{\textit{can}}})_V$ are both oriented by the exterior product of the following vectors:
the $X_{\alpha} \wedge X_{-\alpha}$ with real $\alpha \in R^+_\lambdatilde$,
the $\mi (X_{\alpha} \wedge X_{-\alpha})$ with imaginary~$\alpha \in R^+_\lambdatilde$,
and the $X_{\alpha} \wedge X_{-\alpha} \wedge \overline{X_{\alpha}\!}\, \wedge \overline{X_{-\alpha}\!}\,$
where $\dot{\alpha}$ ranges over the classes of complex
$\beta \in R(\gothg_\CC,\gothh_\CC)\!\setminus\!R(\gothg(\lambda)(\mi\rho_\F+)_\CC,\gothh_\CC)$
modulo the action of the group generated by conjugation and sign change.
Likewise, we can prove that the orientations $\scalo (B_{\lambda_+})_{V(e)}$ and $\scalo (B_{\lambda_{\textit{can}}})_{V(e)}$ are equal, using the root system $R(\gothg(e)_\CC,\gothh(e)_\CC)$ and the ordered pair $(\lambda,\rho_\F+)$.
Therefore
$\;\scalo (B_{\lambda_+})_{(1-\Ad e)(V)} = \scalo (B_{\lambda_{\textit{can}}})_{(1-\Ad e)(V)}$.%
\smallskip

{\bf (c)}
Let $e \in G(\lambda_+)$ be elliptic and let $\widehat{e} \in G(\lambdatilde)^{\gothg/\!\gothg(\lambda)(\mi \rho_\F+)}\!$
lie over $e$.
Let $\varphi_{\textit{can}}$ and $\varphi_+$ denote the Hermitian forms
$(v,w)\mapsto \mi \,B_{\lambda_{\textit{can}}}(v,\overline{w})$
and $(v,w)\mapsto \mi \,B_{\lambda_+}(v,\overline{w})$ on $\gothg_\CC/ \gothg(\lambda)(\mi \rho_\F+)_\CC$.
The vector space $\call \egdef \call_{\lambdatilde}/\gothg(\lambda)(\mi \rho_\F+)_\CC\!$
is the direct sum of the vector subspaces
$\;\call(e)\egdef (\call_{\lambdatilde} \cap \gothg(e)_\CC) /\gothg(e)(\lambda)(\mi \rho_\F+)_\CC\;$
and $\;(1-\Ad e)(\call)$, which are simultaneously orthogonal for $\varphi_{\textit{can}}$ and for $\varphi_+$.

We~follow the beginning of the proof of Lemma \ref{fonctions canoniques sur le revêtement double} (b).
The number of (strictly) negative coefficients in the matrix of either of the Hermitian forms $\varphi_{\textit{can}}$ and $\varphi_+$ on $\call$, with respect to an orthogonal basis, is equal to:
half the number of compact roots in
$R(\gothg_\CC,\gothh_\CC) \!\setminus\!R(\gothg(\lambda)(\mi\rho_\F+)_\CC,\gothh_\CC)$
plus one quarter of the number of complex roots in
$R(\gothg(\nu)_\CC,\gothh_\CC) \!\setminus\!R(\gothg(\lambda)(\mi\rho_\F+)_\CC,\gothh_\CC)$.
An analogous property holds for $\call(e)$.
By (b) and Definition \ref{géométrie métaplectique} (d), the functions $\rho$ associated with $B_{\lambda_{\textit{can}}}$ and $B_{\lambda_+}\!$ therefore take the same values at $\widehat{e}$.
\cqfd
\end{Dem}

%------------------------------------------------------------------------------
\begin{Definition}\rm
\hlabel{paramètres adaptes}
%------------------------------------------------------------------------------

{\bf (a)}
Let $\alpha \!\in\! R(\gothg_\CC,\gothh_\CC)$ be real.
Set
$\,n_\alpha = \frac{1}{2} \Sumpetit_{\beta} \beta(H_{\alpha})\,$
where the sum ranges over
$\beta \in R(\gothg_\CC,\gothh_\CC)$
such that
$\beta + \overline{\beta} \in \RR_+ \alpha$.
To two vectors
$X_{\alpha} \!\in \gothg_\CC^{\alpha}\cap\gothg$
and $X_{-\alpha} \!\in \gothg_\CC^{-\alpha}\cap\gothg$
such that $[X_{\alpha},X_{-\alpha}] = H_{\alpha}$,
we associate the elliptic element
$\,\gamma_\alpha = \exp(\pi(X_{\alpha} - X_{-\alpha}))\,$
of $\centragr{G_0}{\gothh}$.
(We know that $\,n_\alpha \in \NN\,$ and that the set $\{\gamma_\alpha,\gamma_\alpha^{-1}\}$ is independent of the choice of $(X_{\alpha},X_{-\alpha})$, see \cite[p.$\!\!$~327] {DV88}.)%
\smallskip

{\bf (b)}
Let $\,\Xirrp_G(\lambdatilde,{\gotha^*}^+)\,$\labind{Xirr+} denote the finite set of isomorphism classes of irreducible unitary representations $\tau_+$\labind{τ+} of $G(\lambda_+)^{\gothg/\gothh}$
such that $\;\tau_+(\iota) = -\id\;$ and
$\;\tau_+(\exp X) = \me ^{\mi \,\lambda(X)}\id\;$ for $X \in \gothh$
(such a $\tau_+$ is of finite-dimensional),\\
and let $\Xfin_G(\lambdatilde,{\gotha^*}^+)$\labind{Xfinal} denote the set of
$\tau_+ \in \Xirrp_G(\lambdatilde,{\gotha^*}^+)$ “final” in the sense that for every real root
$\,\alpha \in R(\gothg(\lambda)_\CC,\gothh_\CC)$
the endomorphism $(\delta_{\lambda_+}^{\gothg/\gothh} \tau_+) (\gamma_\alpha)$
does not admit the eigenvalue $(-1)^{n_\alpha}$.
(When $G$ is connected, the proof of \ref{résultat inattendu} (b) will show that this condition can be replaced by
$(\delta_{\lambda_+}^{\gothg/\gothh} \tau_+) (\gamma_\alpha)
\not= (-1)^{n_\alpha} \id$.)
\\
Let $\,\Xfin_G$\labind{Xf} denote the set of $(\lambdatilde',{{\gotha'}^*}^+\!,\tau_+')$
with $\;\lambdatilde' = (\lambda',\Fp+)\in \gstregtilde$,
${{\gotha'}^*}^+\!$ is a chamber of $(\gothg(\lambda')(\mi \rho_\Fp+),\gotha')$
where $\gotha'$ is the hyperbolic component of $\gothh' \egdef \gothg(\lambdatilde')$,
and $\;\tau_+' \in \Xfin_G(\lambdatilde',{{\gotha'}^*}^+)$.
\smallskip

{\bf (c)}
Let $X_G(\lambdatilde)$ (resp. $\Xirr_G(\lambdatilde)$\labind{Xirr tilde})
denote the set of isomorphism classes of unitary (resp. irreducible unitary) representations
$\tau$\labind{τ} of $G(\lambdatilde)^{\gothg/\!\gothg(\lambda)(\mi \rho_\F+)}\!$ such that
$\,\tau(\iota)\!=\! -\id\,$ and
$\;\tau(\exp X) = \me ^{\mi \,\lambda(X)}\id\;$ for $X \in \gothh$,\\
and let $\;\XInd_G(\lambdatilde)$\labind{XInd} denote the set (independent of the choice of ${\gotha^*}^+\!$) consisting of the isomorphism classes of unitary representations $\tau$ of
$G(\lambdatilde)^{\gothg/\!\gothg(\lambda)(\mi \rho_\F+)}$
such that\\
\hspace*{\fill}%
$\displaystyle
\frac{\rho_\lambdatilde^{\gothg/\!\gothg(\lambda)(\mi \rho_\F+)}}
{|\rho_\lambdatilde^{\gothg/\!\gothg(\lambda)(\mi \rho_\F+)}|}\;
\tau
\,=\, \Ind_{G(\lambda_+)} ^{G(\lambdatilde)}
\Big(\frac{\rho_{\lambdatilde,{\gotha^*}^+}^{\gothg/\gothh}}
{|\rho_{\lambdatilde,{\gotha^*}^+}^{\gothg/\gothh}|}\;
\tau_+\Big)$
\hspace*{\fill}\\
for some $\,\tau_+ \in \Xfin_G(\lambdatilde,{\gotha^*}^+)\,$
(hence $\,\XInd_G(\lambdatilde)\subseteq X_G(\lambdatilde)$).
\\[0.5mm]
We also define:
$\;\XInd_G
= \Bigl\{
(\lambdatilde',\tau')\,;\,
\lambdatilde' \in \gstregtilde
\,\textrm{ and }\,
\tau' \in \XInd_G(\lambdatilde')
\Bigr\}$\labind{XI}.
\end{Definition}

%-------------------------------------------------------------------------------
\section{Integrability condition}\label{Rep6}
%------------------------------------------------------------------------------

%\noindent
We retain the notation from the previous section.

The following lemma is an adaptation of \cite[Rem. 2 p.$\!\!$~154] {Df82a}.
It will not be used in the remainder of this paper.

%------------------------------------------------------------------------------
\begin{Lemme}
\hlabel{condition d'intégrabilité}
%------------------------------------------------------------------------------

The following properties are equivalent:
\smallskip

{\bf (i)}
$\;\; \Xirrp_G(\lambdatilde,{\gotha^*}^+)\not= \emptyset$;\smallskip

{\bf (ii)}
there exists a (unique) unitary character $\,\chi_{\lambdatilde}^G\,$\labind{χ λ+G}
of $\,G(\lambdatilde)^{\gothg/\gothh}_0$ such that:\\
\hspace*{\fill}%
$\chi_{\lambdatilde}^G (\iota)\, = \, -1\;$
and
$\;\diff_1 \chi_{\lambdatilde}^G \, = \, \mi \restriction{\lambda}{\gothh}$;%
\hspace*{\fill}%

{\bf (iii)}
$\;\; \lambdatilde \in \gstregGtilde$.
\end{Lemme}

\begin{Dem}{Proof of the lemma}

{\bf (i)\boldmath$\Leftrightarrow\!\!$ (ii)} \quad
The implication “(i) $\!\Rightarrow\!$ (ii)” is clear.

Assume (ii).
The restriction to $\,G(\lambdatilde)^{\gothg/\gothh}_0$ of the unitary representation
$\;\Ind_{G(\lambda_+)^{\gothg/\gothh}_0}^{G(\lambda_+)^{\gothg/\gothh}} \chi_{\lambdatilde}^G\;$
of $\,G(\lambda_+)^{\gothg/\gothh}$ is nonzero (since it is induced from a nonzero space) and a multiple of
$\,\chi_{\lambdatilde}^G$.
Theorem 8.5.2 of \cite[p.$\!\!$~153]{Di64} therefore shows the existence of an element of $\,\Xirrp_G(\lambdatilde,{\gotha^*}^+)$.
\smallskip

{\bf (ii)\boldmath$\Leftrightarrow\!\!$ (iii)} \quad
According to Lemma \ref{fonctions canoniques sur le revêtement double} (a), by multiplying $\,\chi_{\lambdatilde}^G\,$ by $\;\rho_{\lambda_+}^{\gothg/\gothh}\,$ condition (ii) is equivalent to the existence of a complex character of the Lie group $\,G(\lambda_+)_0\,$ with differential $\,\mi \lambda + \rho_{\gothg,\gothh}$.
This property can be written as:~%
$\;\lambdatilde \in \gstregGtilde$.
\cqfd
\end{Dem}

%------------------------------------------------------------------------------
\begin{Remarque}\rm
\hlabel{autre insuffisance de la méthode des orbites}
%------------------------------------------------------------------------------

Let $\,f \in \Supp_{\gothg(\lambda)^*}(G(\lambda)_0 \cdot \lambdatilde)$.
In particular, we have $\,f \in \gstreg$.

The real Lie group $G(f)_0$ is abelian (see \cite[Th. p.$\!\!$~17]{B.72}).
The stabilizer of the nilpotent component of $f$ in the adjoint group of $G(\lambda)_0$ is a unipotent linear algebraic subgroup by \cite[Th. 5.9 (b) p.$\!\!$~138] {Sp66}.
The maximal torus of $G(f)_0$ is therefore equal to that of $\centregr{G(\lambda)_0}$.
Given Remark \ref{sous-algèbre de Borel canonique}, the arguments of the proof of the previous lemma yield:\\[1mm]
\hspace*{\fill}%
$\Xirr_G(f)\not= \emptyset \iff
\forall Z \in \Ker \exp_{\centregr{G(\lambda)_0}}\;
\me ^{(\mi \,\lambda+\rho_{\gothg,\gothh})(Z)}=1$.
\hspace*{\fill}\\
Therefore, the condition “$\Xfin_G(\lambdatilde,{\gotha^*}^+)\not= \emptyset$” implies the condition “$\Xirr_G(f)\not= \emptyset$.”%

The converse is false, as shown by the case where $\,G = PSL(2,\RR)$, $\lambda = 0$ and $\gothh = \goth{so}(2)$, therefore $\Xfin_G(\lambdatilde,\{0\}) = \emptyset$, $\Xirr_G(f)\not= \emptyset$, and also $\Xirr_G(-f)\not= \emptyset$.
This reflects the fact that $PSL(2,\RR)$ has no limit of discrete series representation outside its discrete series representations.
To generalize the orbit method proposed by Duflo, it was natural to imagine that the regular coadjoint orbits of $G$ that could be the orbit $G \!\cdot\! l$ associated with a representation $T_{l,\tau}^G$ of zero infinitesimal character would be $G \!\cdot\! f$ and $G \!\cdot\! (-f)$ (see \ref{caractères canoniques} (a), replacing $(\lambdatilde,{\gotha^*}^+)$ and $\lambda$ with~$l$).
This generalization would have created difficulties, because $\,G \!\cdot\! f \not= G \!\cdot\! (-f)$ and the unitary dual of $G$ has only one element with zero infinitesimal character.
\cqfr
\end{Remarque}

%===============================================================================
%
\vspace*{-40pt}
\part{.\quad Construction of representations}\label{RepIII}
%
%==============================================================================

This section closely follows the paper \cite{Df82a} by Duflo, pp. 160--180.
\medskip

Given
$\,\lambdatilde = (\lambda,\F+)\in \gstregtilde$.\hlabel{λtilde,τ}
Set $\gothh = \gothg(\lambdatilde)$.
Let $\mu$ and $\nu$ (resp. $\gotht$ and $\gotha$) denote the infinitesimally elliptic and hyperbolic components of $\lambda$ (resp. of $\gothh$).
We also fix a chamber ${\gotha^*}^+\!$ of $(\gothg(\lambda)(\mi \rho_\F+),\gotha)$
and $\,\tau_+ \in \Xfin_G(\lambdatilde,{\gotha^*}^+)$.
This provides $\;R^+(\gothg_\CC,\gothh_\CC)\,$ and $\,\lambda_+ = \mu_+\!+\nu_+$
(see \ref{choix de racines positives} (b)).
Set $\,\mutilde = (\mu,\F+ \!\cap \mi\gotht^*)$\labind{μtilde}.%
\vspace*{-2mm}

%-------------------------------------------------------------------------------
\section{The case where \texorpdfstring{$G$}{G} is connected}\label{Rep7}
%------------------------------------------------------------------------------

%\noindent
Throughout this section, we assume $G$ is connected.
\smallskip

The inclusion $\,\Ad G \subseteq \interieur\gothg_\CC\,$ ensures that $\,G(\lambda_+)\,$ is equal to $\,\centragr{G}{\gothh}$, and that $\,\delta_{\lambda_+}^{\gothg/\gothh}\,$ is a unitary character of
$\,G(\lambda_+)^{\gothg/\gothh}\,$ (see \ref{morphisme rho et fonction delta}, and also \cite[top of p.$\!\!$~118] {Df84}).
Indeed, for every connected component ${\gotha^*}^+_0$ of the set of elements of $\gotha^*$ that do not vanish on any $H_{\alpha}$ with nonimaginary $\alpha \in R(\gothg_\CC,\gothh_\CC)$, the vector space\\[1mm]
\hspace*{\fill}%
$\call^+
\;\egdef\; \gothh_\CC\oplus\!
\Sum_{\substack
{\alpha \in R^+(\gothg_\CC,\gothh_\CC)\\
\textrm{$\alpha$ noncompact imaginary}}} \!\!\!
\gothg_\CC^{\alpha}
\oplus
\Sum_{\substack
{\alpha \in R^+(\gothg_\CC,\gothh_\CC)\\ \textrm{$\alpha$ compact}}} \!
\gothg_\CC^{-\alpha}
\oplus
\Sum_{\substack
{\alpha \in R(\gothg_\CC,\gothh_\CC)\\
{\gotha^*}^+_0(H_{\alpha})\,\subseteq\, \RR_+ \!\moins0}}
\! \gothg_\CC^{\alpha}$
\hspace*{\fill}\\
yields the positive Lagrangian $\call^+ / \gothh_\CC$ of $(\gothg_\CC/ \gothh_\CC,B_{\lambda_+})$ stable under $G(\lambda_+)$.
\smallskip

Let
$M$\labind{M} denote the intersection of the kernels of the positive real characters of $\centragr{G}{\gotha}$ and let $\gothm$\labind{m} be its Lie algebra (thus, we have $\,\mutilde \in \mstIMtilde$ by \ref{construction de représentations, cas connexe} (a) and $\Xirr_M(\mutilde) = \Xfin_M(\mutilde,\{0\})$),\\
$R^+(\gothm_\CC,\gotht_\CC) =\! R(\gothm_\CC,\gotht_\CC)\cap R^+(\gothg_\CC,\gothh_\CC)$
(Definition \ref{choix de racines positives} (b) relative to $\gothm$ and $(\mutilde,\{0\})$),\\
$\mu_{+,\gothm} = \mu_{\gothm,\mutilde,\{0\}}$\labind{μ+m}
the element $\,\mu - 2\mi \rho_{\gothm,\gotht}\,$ of $\gotht^*$
(thus, we have $M(\mutilde)\!=\! \centragr{M}{\gotht}$ and $\delta_{\mu_{+,\gothm}}^{\gothm / \gotht}$ is a unitary character of $M(\mutilde)^{\gothm / \gotht}$),\\
${\gothn}_M = \!\!\! \Sum_{\alpha \in R^+(\gothm_\CC,\gotht_\CC)}\! \gothg_\CC^{\alpha}\;$
and
$\;\gothb_M = \gotht_\CC\oplus {\gothn}_M$,\\
$q = \abs{\{ \alpha \in R^+(\gothm_\CC,\gotht_\CC)\mid \textrm{$\alpha$ compact} \}}$\labind{q},\\
$\gothk_M
= \gotht \oplus
\bigl( \Sum_{\substack
{\alpha \in R(\gothm_\CC,\gotht_\CC)\\ \textrm{$\alpha$ compact}}}
\! \gothg_\CC^{\alpha}
\; \bigr)
\cap \gothg\;$\llabel{RacinesCompactes}
and
$\,K_{M_0} = \exp \gothk_M$.
\vspace*{-2mm}

%------------------------------------------------------------------------------
\begin{Definition}\rm
\hlabel{limites de séries discrètes}
%------------------------------------------------------------------------------
\vskip-1mm

{\bf (a)}
We denote by $T_{\mutilde}^{M_0}$ the “limit of discrete series representations of $M_0$” belonging to $\widehat{M_0}$, whose space of $K_{M_0}$-finite vectors is isomorphic to the image of the unitary character of $T_0$ (see \ref{conventions fréquentes} (a)) with differential $\,\mi \mu \!-\! \rho_{\gothm,\gotht}\,$ under the “cohomological induction” functor
$\;\calr_{M_0}^q
\egdef
\left({}^{u\!}\calr_{\gothb_M,T_0}^{\gothm_\CC,K_{M_0}}\right)^{\!q}
(\, {\scriptscriptstyle\bullet} \otimes \bigwedge^{\scriptscriptstyle\mathrm{max}}{\gothn}_M )\;$
described in \cite[(11.73) p.$\!\!$~677] {KV95}.
(It is denoted $\pi_{\mi \mu,\gothb_M}$ in \cite[bottom of p.$\!\!$~734] {KV95}, taking into account \cite[Prop. 11.180 p.$\!\!$~733] {KV95}.)
\smallskip

{\bf (b)}
We set
$\;T_{\mutilde,\sigma}^M
= \, \Ind_{\centragr{M}{\gothm} . M_0}^M
(\delta_{\mu_{+,\gothm}}^{\gothm / \gotht} \sigma
\otimes T_{\mutilde}^{M_0} )\;$
for every $\,\sigma \in \Xirr_M(\mutilde)$,\\
where we can define a representation
$\delta_{\mu_{+,\gothm}}^{\gothm / \gotht} \sigma \otimes T_{\mutilde}^{M_0}$
of $\centragr{M}{\gothm} M_0$ by the equality\\
\hspace*{\fill}%
$\left( \delta_{\mu_{+,\gothm}}^{\gothm / \gotht} \sigma
\otimes T_{\mutilde}^{M_0} \right) (xy)
\,=\,
(\delta_{\mu_{+,\gothm}}^{\gothm / \gotht} \sigma)(x)
\otimes T_{\mutilde}^{M_0}(y)\;$
for $\,x \in \centragr{M}{\gothm}\,$ and $\,y \in M_0$
\hspace*{\fill}\\
(in which $\,(\delta_{\mu_{+,\gothm}}^{\gothm / \gotht} \sigma)(x) = \sigma(x,1)$),
because the group $\centregr{M_0}$ operates in “the space” of $T_{\mutilde}^{M_0}$ via the unitary character
$\,\restriction{(\delta_{\mu_{+,\gothm}}^{\gothm / \gotht} \sigma)} {\centregr{M_0}}$
by \cite[(11.184c) p.$\!\!$~734]{KV95}.
\vspace*{-1mm}
\end{Definition}

%------------------------------------------------------------------------------
\begin{Remarque}\rm
\hlabel{cohomologie pour $M$}
%------------------------------------------------------------------------------

Fix a maximal compact subgroup $K$ of $G$ whose Cartan involution normalizes $\gothh$.
Therefore $\,K_M \egdef K \cap M\,$ is a maximal compact subgroup of $M$ containing $\centragr{M}{\gothm} K_{M_0}$.
Every $\,\pi \in (\centragr{M}{\gothm} M_0)\widehat{~}\,$ is related to certain representations
$\,\pi_1 \in (\centragr{M}{\gothm})\widehat{~}\,$ and
$\,\pi_2 \in (M_0)\widehat{~}\,$ by the equalities
$\;\pi(xy) = \pi_1(x)\otimes \pi_2(y)\;$ for $x \in \centragr{M}{\gothm}$
and $y \in M_0$ (see \cite[Prop. 13.1.8 p.$\!\!$~251]{Di64}).
Let $\,\sigma \in \Xirr_M(\mutilde)$.
By \cite[(11.187) p.$\!\!$~735 and proof of Prop. 11.192 (a) p.$\!\!$~737] {KV95}
and \ref{fonctions canoniques sur le revêtement double} (b) to pass from $M_0$ to
$\centragr{M}{\gothm} M_0$,
and \cite[Prop. 11.57 p.$\!\!$~672 and end of (5.8) p.$\!\!$~332] {KV95} to pass from the “parabolic subgroup” $\centragr{M}{\gothm} M_0$ of $M$ to $M$, the $(\gothm_\CC,K_M)$-module associated with $T_{\mutilde,\sigma}^M$ is isomorphic to
$\,\calr_M^q((\rho_{\mu_{+,\gothm}}^{\gothm / \gotht}\!)^{-1} \sigma)$, where
$\calr_M^q$ is the cohomological induction functor relative to $\gothb_M$.
\cqfr
\end{Remarque}

%------------------------------------------------------------------------------
\begin{Lemme}
\hlabel{construction de représentations, cas connexe}
%------------------------------------------------------------------------------

{\bf (a)}
We have $\;\,G(\lambda_+) = \centragr{M}{\gothm} \exp \gothh\;$
and $\;M(\mutilde) = \centragr{M}{\gothm} \exp \gotht$.\\
There therefore exists a unique $\,\tau_M \in \Xirr_M(\mutilde)\,$\labind{τ_M}
such that
$\;\,\delta_{\mu_{+,\gothm}}^{\gothm / \gotht} \tau_M
= \restriction{(\delta_{\lambda_+}^{\gothg/\gothh} \tau_+)}{M(\mutilde)}$.\\
(This notation does not take into account the fact that $\tau_M$ depends on ${\gotha^*}^+\!$.)
\smallskip

{\bf (b)}
The representation class
$\;\Ind_{M . A . N}^G
( T_{\mutilde,\tau_M}^M \!\otimes \me ^{\mi \, \nu \circ \Log} \otimes
\fonctioncar{N} )\;$
constructed from the choice of a subgroup $\,N\,$ of $\,G\,$ which is the unipotent radical of a parabolic subgroup of $G$ with split component $A$, belongs to $\widehat{G}$ and is independent of $N\!$.
\end{Lemme}

\begin{Dem}{Proof of the lemma}

{\bf (a)}
We have
$\;G(\lambda_+) = \centragr{G}{\gothh}
\subseteq \centragr{G}{\gotha} = MA$
and
$\;M(\mutilde) = \centragr{M}{\gotht}$.
Given \cite[Th. 17 p.$\!\!$~199] {Va77}, it follows that
$\;M(\mutilde) = \centragr{M}{\gothm} \exp \gotht$,
and then
$\;G(\lambda_+) = \centragr{M}{\gothm} \exp \gothh$.
\smallskip

{\bf (b)}
Let $N$ be the unipotent radical of a parabolic subgroup of $G$ with split component $A$.
By \cite[Th. 18 p.$\!\!$~289] {Va77}, there exists a connected component ${\gotha^*}^+_0$ of the set of elements of $\gotha^*$ that do not vanish on any $H_{\alpha}$ with nonimaginary $\,\alpha \in R(\gothg_\CC,\gothh_\CC)$, such that\\
\hspace*{\fill}%
$N = \exp \gothn \;\,$
with
$\;\,\gothn
= \bigl(
\Sum_{\alpha \in R(\gothg_\CC,\gothh_\CC)\textrm{ and }
{\gotha^*}^+_0(H_{\alpha})\,\subseteq\, \RR_+ \!\moins0}
\! \gothg_\CC^{\alpha}
\; \bigr)\cap\, \gothg$.
\hspace*{\fill}\\
Let $\nu_0 \in {\gotha^*}^+_0$ and
$\,\lambda_0 \in (\gotht^* + \nu_0)\cap \gstreg$.
We use Proposition \ref{caractères d'induites} (a) with $M'$  built from  $\lambda_0$ instead of $\lambda$ 
(in which case “$(M',U)$” becomes “$(\centragr{M}{\gothm} \, M_0 \, A,N)$”)
and $\pi$ equal to
$\,(\delta_{\mu_{+,\gothm}}^{\gothm / \gotht} \tau_M \otimes T_{\mutilde}^{M_0})
\otimes \me ^{\mi \, \nu \circ \Log}$.
Since the class of
$\,\Ind_{M . A . N}^G
( T_{\mutilde,\tau_M}^M \!\otimes \me ^{\mi \, \nu \circ \Log} \otimes
\fonctioncar{N} )\,$
is determined by its character (see \cite[p.$\!\!$~64] {Co88}), it does not depend on $N$.
\smallskip

We choose $z=1$ and $\Lambda = (\Lambda^{can},R^+_{\mi \, \RR},R^+_\RR)$
in \cite[line 7 bottom of p.$\!\!$~121] {ABV92}, where
$\;\Lambda^{can} \otimes \rho(R^+(\gothg_\CC,\gothh_\CC))
= \rho_{\lambda_+}^{\gothg/\gothh} \, \tau_+\;$
(see \cite[lines 9--10 bottom of p.$\!\!$~129] {ABV92})
and the systems of positive roots $R^+_{\mi \, \RR}$ and $R^+_\RR$ are included in $R^+(\gothg_\CC,\gothh_\CC)$.\\
By Lemma \ref{fonctions canoniques sur le revêtement double} (b), the representation $\wt{\Lambda}$ of \cite[p.$\!\!$~123] {ABV92} is $\,\wt{\Lambda} = \delta_{\lambda_+}^{\gothg/\gothh} \tau_+$.
Given Proposition \ref{caractères d'induites} (a), we have
$\;\pi(\Lambda)
= \Ind_{M . A . N}^G
( T_{\mutilde,\tau_M}^M \!\otimes \me ^{\mi \, \nu \circ \Log} \otimes
\fonctioncar{N} )\;$
in \cite[(11.2)(e) p.$\!\!$~122] {ABV92}.
The ideas in the proof of Lemma 5 of \cite[p.$\!\!$~335] {DV88} allow us to show that $\Lambda$ is “final” in the sense of \cite[Def. 11.13 p.$\!\!$~130] {ABV92}.
Finally, we apply \cite[Th. 11.14 (a) p.$\!\!$~131] {ABV92}.
\cqfd
\end{Dem}

%------------------------------------------------------------------------------
\begin{Remarque}\rm
%------------------------------------------------------------------------------

We place ourselves in the case where $G = SL(3,\RR)$ and $\lambda$ is the element of $\gstssIncG$ which is the image of $\left( {\scriptscriptstyle \diagt{1}{1}{-2}} \right)$ by the isomorphism of $G$-modules that maps $A \in \gothg$ to $\tr (A\,.)\in \gothg^*\!$.
The group $MA$ which is here equal to $G(\lambda)$, is the image of $GL(2,\RR)$ by the embedding of Lie groups that maps $\,x \in GL(2,\RR)\,$ to
$\;\left( {\scriptstyle \diagd{x}{(\det x)^{-1}}} \!\right)\in G$.
But the representation
$\,\Ind_{M . A . N}^G
( T_{\mutilde,\tau_M}^M \!\otimes \me ^{\mi \, \nu \circ \Log} \otimes
\fonctioncar{N} )\,$
is not given by “nondegenerate data” in the sense of \cite[p.$\!\!$~473] {KZ82}, because the Weyl group of a fundamental Cartan subalgebra of $\goth{gl}(2,\RR)$ is represented in $GL(2,\RR)$.
\cqfr
\end{Remarque}

%------------------------------------------------------------------------------
\begin{Definition}\rm
\hlabel{définition de représentations, cas connexe}
%------------------------------------------------------------------------------

Define
$\;\,T_{\lambdatilde,{\gotha^*}^+\!,\tau_+}^G\!
\!= \Ind_{M . A . N}^G
( T_{\mutilde,\tau_M}^M \!\otimes \me ^{\mi \, \nu \circ \Log} \otimes \fonctioncar{N} )
\,\in\, \widehat{G}$\labind{Tλtilde,a*,τ,connexe},
independently of the choice of $N$ as in the lemma above.

Thus, when $M$ is connected, $\Xirr_M(\mutilde)$ has as its only element $\chi_{\mutilde}^M\!$, and the representation classes $\,T_{\mutilde,\chi_{\mutilde}^M}^M$ and $\,T_{\mutilde,\{0\},\chi_{\mutilde}^M}^M$ are both equal to $\,T_{\mutilde}^{M_0} \!$.
\end{Definition}

%-------------------------------------------------------------------------------
\section{The representations \texorpdfstring{$T_{\lambdatilde,{\gotha^*}^+\!,\tau_+}^G\!$}{T\_\{\~\lambda,a*+,\tau+\}\^{}G}}\label{Rep8}
%------------------------------------------------------------------------------

%\noindent
Given $(\lambdatilde,{\gotha^*}^+\!,\tau_+)$ as on page \pageref{λtilde,τ}, we will apply the construction of the connected case to $\,M'_0 \egdef G(\nu_+)_0$.
A~detour through homology will, thanks to a result by Vogan, allow us to normalize some intertwining operators in order to extend $T_{\lambdatilde}^{M'_0}$ into a representation of a certain group $M'$.
\smallskip

Let $M'$ denote $G(\lambda_+) G(\nu_+)_0\,$\labind{M'} and
$\,\gothm' = \gothg(\nu_+)\,$\labind{m'} its Lie algebra\\
(thus $\,\lambdatilde \in \mpstfondMptilde$ by \ref{construction équivalente, cas connexe} (a),
$\gothm'(\lambda)(\mi \rho_\F+) = \gothh\,$ and
$\,\Xirr_{M'}(\lambdatilde) = \Xfin_{M'}(\lambdatilde,\gotha^*)$),\\
$R^+(\gothm'_\CC,\gothh_\CC)
= R(\gothm'_\CC,\gothh_\CC)
\cap R^+(\gothg_\CC,\gothh_\CC)\;$
(see \ref{choix de racines positives} (b) relative to $\gothm'$ and $(\lambdatilde,\gotha^*)$),\\
$\lambda_{+,\gothm'} = \lambda_{\gothm'\!,\lambdatilde,\gotha^*\!,1}\;$\labind{λ+m'}
the element $\mu_+ + \nu$ of $\gothh^*$ where $\mu_+$ still denotes
$\mu_{\gothg,\lambdatilde,{\gotha^*}^+}\!$\\
(thus $\,B_{\lambda_{+,\gothm'}}\!$ is a restriction of $B_{\lambda_+}\!$, $M'(\lambdatilde) = G(\lambda_+)$
and
$\,\rho_\lambdatilde^{\gothm'\!/\gothh}
= \rho_{\lambda_{+,\gothm'}}^{\gothm'\!/\gothh}$),\\
$\gothn_{M'}
= \! \Sum_{\alpha \in R^+(\gothm'_\CC,\gothh_\CC)} \!
\gothg_\CC^{\alpha}\;\,$\labind{nM'}
and
$\;\gothb_{M'} = \gothh_\CC\oplus \gothn_{M'}$\labind{bM'},\\
$q'
= \abs{\{ \alpha \in R^+(\gothm'_\CC,\gothh_\CC)\mid
\textrm{$\alpha$ compact} \}}
\,+\,
\frac{1}{2} \, \abs{\{ \alpha \in R^+(\gothm'_\CC,\gothh_\CC)\mid
\textrm{$\alpha$ complex} \}}$\labind{q'},\\
$M_{\nu_+}\!$ the intersection of the kernels of the positive real characters of $\centragr{M'_0}{\gotha}$,\\
$T_{\mutilde}^{M_{\nu_+}} \!\!= T_{\mutilde,\chi_{\mutilde}^M}^{M_{\nu_+}}$
(see \ref{limites de séries discrètes} (b) and $\!M_{\nu_+}(\mutilde) \!=\! \exp \gotht$),
$T_{\lambdatilde}^{M'_0}
\!= T_{\lambdatilde,\gotha^*,\chi_{\lambdatilde}^{M'}}^{M'_0}\!$\labind{Tλztilde,M'0}
(see $M'_0(\lambdatilde) \!=\! \exp \gothh$),\\
$\gothu
=\! \bigl( \!\Sum_{\alpha \in R(\gothg_\CC,\gothh_\CC)\textrm{ and }
\nu_+(H_{\alpha})>0} \!\!\! \gothg_\CC^{\alpha} \bigr)
\cap \gothg$\labind{u},
$U \!=\! \exp \gothu$\labind{U}
(hence $\interieur{M'} . U \!\subseteq\! U\!$ and $M' \cap U \!=\! \{1\}$).%

%------------------------------------------------------------------------------
\begin{Proposition}
\hlabel{construction de représentations, cas de $M'$}
%------------------------------------------------------------------------------

Denote by $\calh$\labind{H} “the space” of $\,T_{\lambdatilde}^{M'_0}$ and by
$\,\calh^{\infty}$ the set of its $C^{\infty}\!$ vectors.
\smallskip

{\bf (a)}
The eigenspace of weight $-(\mi \lambda + \rho_{\gothm'\!,\gothh})$ under the action of $\gothh_\CC$ on  $H_{q'} (\gothn_{M'},\calh^{\infty})^*$ resulting from the action of $\gothh_\CC$ on $\bigwedge \!\gothn_{M'} \otimes_\CC\calh^{\infty}\!$, is one-dimensional.
Denote it by
$\,(H_{q'} ( \, \gothn_{M'},\calh^{\infty} )^*)_{-(\mi \,\lambda + \rho_{\gothm'\!,\gothh})}$.
\smallskip

{\bf (b)}
There exists a unique continuous unitary representation $S$\labind{S} of
$M'(\lambdatilde)^{\gothm'\!/\gothh}$ in $\calh$ satisfying the following conditions (i) and (ii):

\hspace{1cm} (i)
$\;S(\hat{x})\, T_{\lambdatilde}^{M'_0}(y)\,
S(\hat{x})^{-1}
=\, T_{\lambdatilde}^{M'_0}(xyx^{-1})\;\;$
for $\,\hat{x} \in M'(\lambdatilde)^{\gothm'\!/\gothh}$
lying over $\,x \in M'(\lambdatilde)$ and $\,y \in M'_0$;

\hspace{1cm} (ii)
the action of $M'(\lambdatilde)^{\gothm'\!/\gothh}$ on
$\,(H_{q'} ( \, \gothn_{M'},\calh^{\infty} )^*)_{-(\mi \,\lambda + \rho_{\gothm'\!,\gothh})}$
resulting from the action of $M'(\lambdatilde)^{\gothm'\!/\gothh}$
on $\,\bigwedge \gothn_{M'} \otimes_\CC\calh^{\infty}$ obtained from $S$, is $\,(\rho_{\lambda_{+,\gothm'}}^{\gothm'\!/\gothh})^{-1} \id$.
\end{Proposition}

\begin{Dem}{Proof of the proposition}

{\bf (a)}
Assume $G$ is connected and use the notation from Section~\ref{Rep7}

Fix a maximal compact subgroup $K$ of $G$ whose Cartan involution normalizes $\gothh$.
Therefore $\,K_{M'_0} \egdef K \cap M'_0\,$
and $\,K_{M_{\nu_+}} \egdef K \cap M_{\nu_+}\,$
are maximal compact subgroups of $M'_0$ and $M_{\nu_+}\!$.
By Remark \ref{cohomologie pour $M$}, the
$(\gothm_\CC,K_{M_{\nu_+}})$-module associated with
$T_{\mutilde}^{M_{\nu_+}}$ is isomorphic to
$\,\calr_{M_{\nu_+}}^q
( (\rho_{\mu_{+,\gothm}}^{\gothm / \gotht}\!)^{-1} \chi_{\mutilde}^M)$,
where $\,\calr_{M_{\nu_+}}^q$ is the cohomological induction functor relative to $\gothb_M$.

Given Definition \ref{géométrie métaplectique} (d), and Lemma \ref{fonctions canoniques sur le revêtement double} (b) with its notation, for every $\,e \in \centragr{M}{\gothm} T_0\,$ (see \ref{conventions fréquentes} (a)) we have\\[0.5mm]
\hspace*{\fill}%
$\displaystyle \frac
{(\rho_{\lambda_{+,\gothm'}}^{\gothm'\!/\gothh}
\delta_{\lambda_{+,\gothm'}}^{\gothm'\!/\gothh})(e)}
{(\rho_{\mu_{+,\gothm}}^{\gothm / \gotht}
\delta_{\mu_{+,\gothm}}^{\gothm / \gotht})(e)}
= \det (\Ad e^\CC)_{\gothn_{M'}\!/\gothn_M}
\,{\scriptstyle\times}\,
\frac
{((\rho_{\lambda_{+,\gothm'}}^{\gothm'\!/\gothh})^{-1}
\delta_{\lambda_{+,\gothm'}}^{\gothm'\!/\gothh})(e)}
{((\rho_{\mu_{+,\gothm}}^{\gothm / \gotht}\!)^{-1}
\delta_{\mu_{+,\gothm}}^{\gothm / \gotht})(e)}
=\, u_1 \,\cdots\, u_b$
\hspace*{\fill}\\
where
$\,\{\beta_1,\!-\overline{\beta_1}\},\dots,\{\beta_b,\!-\overline{\beta_b}\}$
are the elements of $\wt{R}^+_\CC(\gothm'_\CC,\gothh_\CC)$
and $u_1,\dots,u_b$ are the ratios of the homotheties
$(\Ad e^\CC)_{\gothg_\CC^{\beta_1}},\dots,(\Ad e^\CC)_{\gothg_\CC^{\beta_b}}$.

We apply \cite[Th. 11.225 p.$\!\!$~759] {KV95} (where (11.222) need to change “$z_1z_2$” to “$\xi_{\delta_c(\gothu)}(t)^{-1}z_1z_2$” on line 19 p.$\!\!$~753)
and \cite[Cor. 8.28 p.$\!\!$~566] {KV95}
to the group~$M'_0$ and to the parabolic subalgebra $\gothb_{M'}$ of $\gothm'_\CC$.
The $(\gothm'_\CC,K_{M'_0})$-module $\,\calh^f$ associated with $T_{\lambdatilde}^{M'_0}$ is therefore irreducible and isomorphic to
$\,\calr_{M'_0}^{q'} ( (\rho_\lambdatilde^{\gothm'\!/\gothh})^{-1} \chi_{\lambdatilde}^{M'})$,
where $\,\calr_{M'_0}^{q'}$ is the cohomological induction functor relative to $\gothb_{M'}$.
We obtain\\[-1mm]
\hspace*{\fill}%
$\dim \, \mathrm{Hom}_{\gothh_\CC,T_0}
\bigl( H_{q'} (\gothn_{M'},\calh^f),
\rho_\lambdatilde^{\gothm'\!/\gothh} \chi_{\lambdatilde}^{M'} \bigr)
= 1$
\hspace*{\fill}\\
by taking $\,X = \calr_{M'_0}^{q'}(Z)\,$ with
$\,Z = (\rho_\lambdatilde^{\gothm'\!/\gothh})^{-1}
\chi_{\lambdatilde}^{M'}$
in \cite[Prop. 8.11 p.$\!\!$~555] {KV95}, and taking into account Schur's Lemma (see \cite[Lem. 3.3.2 p.$\!\!$~80] {Wa88}).
Since the canonical homomorphism of $(\gothh_\CC,T_0)$-modules from $H_{q'} (\gothn_{M'},\calh^f )$ to $H_{q'} (\gothn_{M'},\calh^{\infty} )$ is bijective according to \cite[Lem. 4 p.$\!\!$~165] {Df82a}, this yields the result.%
\smallskip

{\bf (b)}
By \ref{injection dans le dual} (a), we have
$\;\interieur{x} \cdot T_{\lambdatilde}^{M'_0} \!= T_{\lambdatilde}^{M'_0}$
for all $\,x \in M'(\lambdatilde)$.
We then repeat the proof of \cite[Lem. 6 p.$\!\!$~169] {Df82a}, almost verbatim.
(The uniqueness of $S$ follows, of course, from \ref{construction de représentations, cas connexe} (b) and Schur's Lemma.)
\cqfd
\end{Dem}

%------------------------------------------------------------------------------
\begin{Lemme}
\hlabel{construction équivalente, cas connexe}
%------------------------------------------------------------------------------

We retain the notation from the previous proposition.
\smallskip

{\bf (a)}
There exists a unique $\,\tau_{M'} \in \Xirr_{M'}(\lambdatilde)\,$\labind{τ_M'}
such that
$\,\delta_{\lambda_{+,\gothm'}}^{\gothm'\!/\gothh} \tau_{M'}\!
= \delta_{\lambda_+}^{\gothg/\gothh} \tau_+$.\\
Denote by $\,T_{\lambdatilde,\tau_{M'}}^{M'} \!\!\!=\! \tau_{M'} \otimes S T_{\lambdatilde}^{M'_0}$
the representation of $M'$ defined (unambiguously)~by\\
\hspace*{\fill}%
$\bigl( \tau_{M'} \otimes S T_{\lambdatilde}^{M'_0} \bigr)\!(xy)
= \tau_{M'}(\hat{x})\otimes S(\hat{x})\, T_{\lambdatilde}^{M'_0}\!(y)$
\hspace*{\fill}\\
for $\,\hat{x} \in M'(\lambdatilde)^{\gothm'\!/\gothh}$ lying over $\,x \in M'(\lambdatilde)$ and $\,y \in M'_0$.
\smallskip

{\bf (b)}
When $G$ is connected, we have:
$\,T_{\lambdatilde,{\gotha^*}^+\!,\tau_+}^G\!\!
=\, \Ind_{M' \!. U}^G
( (\tau_{M'} \otimes S \, T_{\lambdatilde}^{M'_0})\otimes \fonctioncar{U} )$.%
\end{Lemme}

\begin{Dem}{Proof of the lemma}

{\bf (a)}
Lemma \ref{fonctions canoniques sur le revêtement double} (b) allows us to define $\tau_{M'}$ using the equality
%By Lemma \ref{fonctions canoniques sur le revêtement double} (b), we can define $\tau_{M'}$ using\\
$\displaystyle
\frac{\rho_{\lambda_{+,\gothm'}}^{\gothm'\!/\gothh}}
{|\rho_{\lambda_{+,\gothm'}}^{\gothm'\!/\gothh}|} \,
(\delta_{\lambda_{+,\gothm'}}^{\gothm'\!/\gothh})^{-1}
=\, \frac{\rho_{\lambda_+}^{\gothg/\gothh}}
{|\rho_{\lambda_+}^{\gothg/\gothh}|} \,
(\delta_{\lambda_+}^{\gothg/\gothh})^{-1}$,
where $\rho_{\lambda_{+,\gothm'}}^{\gothm'\!/\gothh}\!$ and
$\rho_{\lambda_+}^{\gothg/\gothh}\!$ are~Lie~group homomorphisms.
Furthermore, we have
$S(\exp X) = \me^{-\mi\lambda(X)}T_{\lambdatilde}^{M'_0}(\exp X)$ for $X \!\in\! \gothh$.%
\smallskip

{\bf (b)}
Assume $G$ is connected and use the notation from Section ~\ref{Rep7}\\
We therefore have
$\,M' = \centragr{M}{\gothm} M'_0$,
$\,M'(\lambdatilde) = \centragr{M}{\gothm} \exp \gothh\,$
and $\delta_{\lambda_{+,\gothm'}}^{\gothm'\!/\gothh}$ is a unitary character of
$M'(\lambdatilde)^{\gothm'\!/\gothh}\!$.

Choose a connected component ${\gotha^*}^+_0$ of the set of elements of~$\gotha^*$ that do no vanish on any $H_{\alpha}$ with nonimaginary $\,\alpha \in R(\gothg_\CC,\gothh_\CC)$, whose closure contains~$\nu_+$.
Set
$\;\;\gothn
\,=\, \bigl(
\Sum_{\alpha \in R(\gothg_\CC,\gothh_\CC)\textrm{ and }
{\gotha^*}^+_0(H_{\alpha})\,\subseteq\, \RR_+ \!\moins0}
\!\! \gothg_\CC^{\alpha}
\; \bigr)\,\cap\, \gothg$,
$\;\,N = \exp \gothn\;$ and $\;N_{\nu_+} \!= N \cap M'_0$.\\
The group $N$ (resp. $N_{\nu_+}\!$) is the unipotent radical of a parabolic subgroup of $G$ (resp. $M'_0$) with split component $A$.
Furthermore, the product of $G$ restricts to a diffeomorphism from $\,N_{\nu_+} \times U\,$ onto $N$.
\nopagebreak

Denote by $\calv$ and $\calw$ “the spaces” of $\tau_{M'}$ and $\,T_{\mutilde}^{M_0}\!$.
The space of $\,T_{\lambdatilde}^{M'_0}$ is realized in the form
$\,\calh
\egdef
\Ind_{M_0 . A . N_{\nu_+}}^{M'_0}
( \calw \otimes \me ^{\mi \, \nu \circ \Log} \otimes
\fonctioncar{N_{\nu_+}} )$.
Identify $\calv \otimes_\CC\calh$ with a set of functions on $M'_0$ with values in $\calv \otimes_\CC\calw \!$.
The map $\,\Phi \mapsto \restriction{\Phi}{M'_0}$ from the space of
$\,\Ind_{\centragr{M}{\gothm} . M_0 . A . N_{\nu_+}}^
{\centragr{M}{\gothm}\,M'_0}
(\delta_{\lambda_{+,\gothm'}}^{\gothm'\!/\gothh} \! \tau_{M'}
\otimes T_{\mutilde}^{M_0}
\otimes \me ^{\mi \, \nu \circ \Log} \otimes \fonctioncar{N_{\nu_+}} )\,$
to $\calv \otimes_\CC\calh$ equipped with $1 \otimes T_{\lambdatilde}^{M'_0}\!$
is a unitary isomorphism of $M'_0$-modules.
The action of a $\,x \in \centragr{M}{\gothm}$ on $\calv \otimes_\CC\calh$ that follows from this is given by
$\tau_{M'}(\hat{x})\otimes S_0(\hat{x})$, where
$\hat{x} \in M'(\lambdatilde)^{\gothm'\!/\gothh}$ lies over $x$ and
$\;S_0(\hat{x})\cdot\varphi
= \delta_{\lambda_{+,\gothm'}}^{\gothm'\!/\gothh}\!(\hat{x})
\,{\scriptstyle \times}\,
\varphi(x^{-1} \,\cdot\, x)\;$
for $\,\varphi \in \calh$.

Consider a closed subgroup $\Gamma$ of $\centragr{M}{\gothm}$.
Set $\,M_\Gamma = \Gamma M_{\nu_+}\!$ and $\,M'_\Gamma = \Gamma M'_0$.
Fix an irreducible component $\tau_{M'\!,\Gamma}$ (in $\Xirr_{M'_\Gamma}(\lambdatilde)$) of the restriction of $\tau_{M'}$ to $M'_\Gamma(\lambdatilde)^{\gothm'\!/\gothh}$.
As above, we obtain a unitary isomorphism of $M'_0$-modules, that allows us to transport the action of $\Gamma$, replacing $\centragr{M}{\gothm}$ by $\Gamma$,
$\tau_{M'}$ by $\tau_{M'\!,\Gamma}$
and $\calv$ by the space $\calv_{\Gamma}$ of $\tau_{M'\!,\Gamma}$.
Determine an element $\tau_{M,\Gamma}$ of $\Xirr_{M_\Gamma}(\mutilde)$ by the equality
$\;\delta_{\mu_{+,\gothm}}^{\gothm / \gotht} \tau_{M,\Gamma}
= \restriction{(\delta_{\lambda_{+,\gothm'}}^{\gothm'\!/\gothh} \!
\tau_{M'\!,\Gamma})}{\Gamma T_0}\;$
(see \ref{conventions fréquentes} (a)).

We will see that $S_0 = S$ (see \ref{construction de représentations, cas de $M'$} (b)).
The asserted equality will then follow from the induction in stages theorem and the fact that
$\;\Ind_{H_1}^{G_1} (\pi \circ \restriction{p}{H_1})
\simeq (\Ind_{H_2}^{G_2} \pi)\circ p\;$
when $p: G_1 \to G_2$ is a surjective homomorphism of real Lie groups,
$H_1$ is the inverse image by $p$ of a closed subgroup $H_2$ of $G_2$,
and $\pi$ is a continuous unitary representation of $H_2$.

Resume the arguments and notation from the proof of \ref{construction de représentations, cas de $M'$} (a).
The groups $\,K_{M'_\Gamma} \!\egdef K \cap M'_\Gamma\,$
and $\,K_{M_\Gamma} \!\egdef K \cap M_\Gamma$
are maximal compact subgroups of $M'_\Gamma$ and $M_\Gamma$.
Note that the $(\gothm_\CC,K_{M_\Gamma})$-module associated to the representation
$\,\Ind_{\Gamma . M_0}^{M_\Gamma} (\delta_{\mu_{+,\gothm}}^{\gothm / \gotht}
\tau_{M,\Gamma} \otimes T_{\mutilde}^{M_0})\,$
of the group $M_\Gamma$ is isomorphic to
$\,\calr_{M_\Gamma}^q ((\rho_{\mu_{+,\gothm}}^{\gothm / \gotht}\!)^{-1} \tau_{M,\Gamma})$,
where $\,\calr_{M_\Gamma}^q$ is the cohomological induction functor relative to $\gothb_M$.
Next, using the induction in stages theorem, we see that the
$(\gothm'_\CC,K_{M'_\Gamma})$-module associated to
$\,\Ind_{\Gamma . M_0 . A . N_{\nu_+}}^{M'_\Gamma}
( \delta_{\lambda_{+,\gothm'}}^{\gothm'\!/\gothh} \! \tau_{M'}
\otimes T_{\mutilde}^{M_0}
\otimes \me ^{\mi \, \nu \circ \Log} \otimes \fonctioncar{N_{\nu_+}} )\,$
is isomorphic to
$\,\calr_{M'_\Gamma}^{q'}
( (\rho_\lambdatilde^{\gothm'\!/\gothh})^{-1} \tau_{M'\!,\Gamma})$,
where $\,\calr_{M'_\Gamma}^{q'}$ is the cohomological induction functor relative to $\gothb_{M'}$.
Finally, we find that\\
\hspace*{\fill}%
$\dim \, \mathrm{Hom}_{\gothh_\CC,\Gamma T_0} \,
\bigl( \calv_{\Gamma} \otimes_\CC H_{q'} (\gothn_{M'},\calh^f),
\rho_\lambdatilde^{\gothm'\!/\gothh} \tau_{M'\!,\Gamma} \bigr)
= 1$,
\hspace*{\fill}\\
where each $y \in \Gamma\,$ projection of a $\,\hat{y} \in M'(\lambdatilde)^{\gothm'\!/\gothh}$
acts on $\calv_{\Gamma} \otimes_\CC H_{q'} (\gothn_{M'},\calh^f)$
by the tensor product of $\,\tau_{M'\!,\Gamma}(\hat{y})\,$
and the endomorphism obtained from $\,S_0(\hat{y})$.

Consider $\,\hat{x} \in M'(\lambdatilde)^{\gothm'\!/\gothh}$ lying over some $\,x \in \centragr{M}{\gothm}$.
We choose for $\Gamma$ the closed subgroup of $\centragr{M}{\gothm}$ generated by $x$.
Therefore $\dim \calv_{\Gamma} = 1$ and the preceding calculation proves that the action of $\hat{x}$ on
$\,(H_{q'} ( \, \gothn_{M'},\calh^{\infty} )^*)_{-(\mi \,\lambda + \rho_{\gothm'\!,\gothh})}$
arising from $S_0(\hat{x})$ is indeed the homothety with ratio
$\rho_\lambdatilde^{\gothm'\!/\gothh} (\hat{x})^{-1} \!$.
\cqfd
\end{Dem}

%------------------------------------------------------------------------------
\begin{Definition}\rm
%------------------------------------------------------------------------------

We retain the notation from Proposition \ref{construction de représentations, cas de $M'$}
and Lemma \ref{construction équivalente, cas connexe}.
Define
$\;\,T_{\lambdatilde,{\gotha^*}^+\!,\tau_+}^G\!
= \Ind_{M' .\, U}^G
( ( \tau_{M'} \otimes S \, T_{\lambdatilde}^{M'_0} )\otimes
\fonctioncar{U} )$\labind{Tλtilde,a*,τ}.

To simplify the notation
$T_{\lambdatilde,{\gotha^*}^+\!,\tau_+}^G\!$\labind{Tλzr,M'0}\labind{Tλzr,sigmar,M'},
I will omit ${\gotha^*}^+\!$ when $\,{\gotha^*}^+ = \gotha^*$
(compatible with \ref{notation pour les représentations})
and $\tau_+$ when $\Xirrp_G(\lambdatilde,{\gotha^*}^+)$ is a singleton,
and I will replace $\lambdatilde$ with $\lambda$ when $\,\lambdatilde = (\lambda,\gothh_{(\RR)}^{~~*})$
(conventions applied for
$\,T_{\mutilde}^{M_0}\!$, $T_{\mutilde,\sigma}^M$,
$\,T_{\mutilde}^{M_{\nu_+}}\!$
and $\,T_{\lambdatilde}^{M'_0}$).
\end{Definition}

%------------------------------------------------------------------------------
\begin{Proposition}
\hlabel{caractères canoniques}
%------------------------------------------------------------------------------

{\bf (a)}
For each semisimple $l \in \gothg_\CC^*$ denote by $\chiinfty{l}{\gothg}$\labind{χ lUgC} the character of
$\centrealg{U \gothg_\CC}$ canonically associated with the orbit of $l$ under the action of~$\interieur\gothg_\CC$ (see \cite[(4.114) p.$\!\!$~297] {KV95}).
The set $(U \gothg_\CC)^G\!$ of $G$-invariant elements of $U \gothg_\CC$~acts on the space of $C^{\infty}\!$ vectors of $T_{\lambdatilde,{\gotha^*}^+\!,\tau_+}^G\!\!$ via restriction of the character
$\chiinfty{\mi \lambda}{\gothg}\!$ of~$\centrealg{U \gothg_\CC}\!$.%
\smallskip

{\bf (b)}
The group homomorphism $z \!\mapsto\! (z,1)$ from $\centregr{G}\!$ to $G(\lambda_+)^{\gothg/\gothh}\!$ is injective.
The central character of $T_{\lambdatilde,{\gotha^*}^+\!,\tau_+}^G\!$ is “the restriction” of that of $\tau_+$.
\end{Proposition}

\begin{Dem}{Proof of the proposition}

{\bf (a)}
Corollary 5.25 (b) of \cite[p.$\!\!$~344] {KV95} allows us to calculate the infinitesimal character of
$\;\tau_{M'} \otimes S \, T_{\lambdatilde}^{M'_0}\,$
(equal to that of $T_{\lambdatilde}^{M'_0}$).
To then move on to $G$, the proof of \cite[Prop. 11.43 p.$\!\!$~665] {KV95} can be immediately adapted.
\smallskip

{\bf (b)}
Let $z\in\centregr{G}$.
The equality $S(z,1)=\id$ yields the result.
\cqfd
\end{Dem}

%-------------------------------------------------------------------------------
\section{The injection \texorpdfstring{$\,G \cdot (\lambdatilde,\tau)\mapsto T_{\lambdatilde,\tau}^G$}{G \cdot (\~\lambda,\tau)\mapsto T\_\{\~\lambda,\tau\}\^{}G}
from \texorpdfstring{$G\,\backslash\,\XInd_G$}{G \textbackslash X\^{}Ind\_G} to \texorpdfstring{$\widehat{G}$}{\^G}}\label{Rep9}
%------------------------------------------------------------------------------

%\noindent
Theorem \ref{injection dans le dual} will generalize a large part of the one written by Duflo in \cite[Lem. 8 p.$\!\!$~173] {Df82a}.
My proof is similar to his, except for one unexpected result:
Lemma \ref{résultat inattendu} (b).
Recall that the equivalence class of an irreducible unitary representation of trace class of a real Lie group is determined by its character (see \cite[p.$\!\!$~64] {Co88}).
\smallskip

Given $(\lambdatilde,\!{\gotha^*}^+\!,\tau_+\!)$ as on page \pageref{λtilde,τ}, construct $\tau$ from $\tau_+\!$ as in Definition~\ref{paramètres adaptes}~(c).%

%------------------------------------------------------------------------------
\begin{Lemme}
\hlabel{nouveaux paramètres}
%------------------------------------------------------------------------------

The parameters $\,\lambdatilde\,$ and $\,\tau\,$ determine \smash{$\,T_{\lambdatilde,{\gotha^*}^+\!,\tau_+}^G\!$}.
\end{Lemme}

\begin{Dem}{Proof of the lemma}

This result will follow from Theorem \ref{formule du caractère}, by first removing all references to
$T_{\lambdatilde,\tau}^G\!$ (for consistency of reasoning).
\cqfd
\vspace*{-5mm}
\end{Dem}

%------------------------------------------------------------------------------
\begin{Definition}\rm
\hlabel{notation pour les représentations}
%------------------------------------------------------------------------------

{\bf (a)}
Define
$\;T_{\lambdatilde,\tau}^G\! = T_{\lambdatilde,{\gotha^*}^+\!,\tau_+}^G\!$\labind{Tλtilde,τ}.
In the notation $T_{\lambdatilde,\tau}^G\!$, I will omit~$\tau$ when
$\XInd_G(\lambdatilde)$ is a singleton, and I will replace $\lambdatilde$ with
$\lambda$ when $\,\lambdatilde = (\lambda,\gothh_{(\RR)}^{~~*})$.%
\smallskip

{\bf (b)}
Let $\,a$ be an automorphism of the Lie group $G$.
It induces canonical isomorphisms, still denoted by $\,a$,
of vector spaces from $\gothg / \gothg(\lambda)(\mi \rho_\F+)$ to
$\gothg / \gothg(a\lambda)(\mi \rho_{a\F+})$,
of Lie algebras from
$\,\goth{sl} (\gothg / \gothg(\lambda)(\mi \rho_\F+))$ to
$\,\goth{sl} (\gothg / \gothg(a\lambda)(\mi \rho_{a\F+}))$ by conjugacy,
of Lie groups from
$\,DL (\gothg / \gothg(\lambda)(\mi \rho_\F+))$ to
$\,DL (\gothg / \gothg(a\lambda)(\mi \rho_{a\F+}))$ by integration,
and finally of Lie groups from
$\,G(\lambdatilde)^{\gothg / \gothg(\lambda)(\mi \rho_\F+)}$ to
$\,G(a\lambdatilde)^{\gothg / \gothg(a\lambda)(\mi \rho_{a\F+})}$.%
\end{Definition}

%------------------------------------------------------------------------------
\begin{Lemme}
\hlabel{résultat inattendu}
%------------------------------------------------------------------------------

Let $\,\tau_0 \in (G_0(\lambda_+)^{\gothg/\gothh})\widehat{~~}$ occurring in
$\restriction{\tau_+}{G_0(\lambda_+)^{\gothg/\gothh}}\!$.
\smallskip

{\bf (a)}
We have:
$\;\tau_0 \in \Xfin_{G_0}(\lambdatilde,{\gotha^*}^+)\,$
and
$\,G_0(\lambda_+) = \centragr{G_0}{\gothh} = \normagr{G_0(\lambdatilde)}{{\gotha^*}^+}$.
\vskip-1mm

{\bf (b)}
Define:
$\;\dot{\tau_0}
= \frac{\rho_{\lambdatilde,{\gotha^*}^+}^{\gothg/\gothh}}
{|\rho_{\lambdatilde,{\gotha^*}^+}^{\gothg/\gothh}|}\,
\tau_0\,$
(see \ref{paramètres adaptes} (c)).\llabel{a->c}
The action of the group $G_0(\lambdatilde)$ on the representation class $\dot{\tau_0}$ of $G_0(\lambda_+)$ is trivial.%
\vspace*{-3mm}
\end{Lemme}

\begin{Dem}{Proof of the lemma}

{\bf (a)}
Definition \ref{paramètres adaptes} (b) ensures that $\tau_0$ is final.
The equalities $\,G_0(\lambda_+) = \centragr{G_0}{\gothh} = \normagr{G_0(\lambdatilde)}{{\gotha^*}^+}\,$ follows from Lemma \ref{lemme clef}.
\smallskip

{\bf (b)}
We prove that $\,\tr\dot{\tau_0}(x\gamma x^{-1}) = \tr\dot{\tau_0}(\gamma)\,$
for $x \in G_0(\lambdatilde)$ and $\gamma \in G_0(\lambda_+)$.%
\smallskip

Let $G_1$ be the subgroup of $G_0$ consisting of the elements that fix $\mu$ and~$\mi \rho_\F+$.
Denote by $\gothg_1$ its Lie algebra.
Since the linear forms $\mu$ and $\mi \rho_\F+$ are infinitesimally elliptic,
$\gothg_1$ is reductive, $\gothh \in \Car \gothg_1$, and $G_1$ is connected.
The root system $\,R({\gothg_1}_\CC,\gothh_\CC)\,$ has no imaginary root.
Denote by $\,W(\gothg_1(\nu),\gotha)\,$ the Weyl group of the root system formed by the restricted roots of $(\gothg_1(\nu),\gotha)$.
By \cite[Th. 5.17 p.$\!\!$~125, (5.5) p.$\!\!$~126 and Lem. 5.16 p.$\!\!$~124] {Kn86}
applied to $G_1$ on the one hand, and \cite[Prop. 4.12 p.$\!\!$~81] {Kn86} on the other hand,
the canonical map from $\,G_0(\lambdatilde)/G_0(\lambda_+)\,$ to $W(\gothg_1(\nu),\gotha)$ is bijective.
Consider $\beta \in R({\gothg_1(\nu)}_\CC,\gothh_\CC)$, $X \in \gothg_1(\nu)^{\restriction{\beta}{\gotha}}$ and $Y \in \gothg_1(\nu)^{-\restriction{\beta}{\gotha}}$ such that $[X,Y]$ is the coroot of
$\restriction{\beta}{\gotha}$ in $\gotha$.
Set~$\,g = \exp(\frac{\pi}{2}(X-Y))$.
Therefore $g \in G_0(\lambdatilde)$ projects into $W(\gothg_1(\nu),\gotha)$
on the root reflection $s_{(\restriction{\beta}{\gotha})}$ given the proof of \cite[Prop. 5.15 (c) p.$\!\!$~123] {Kn86}.
It now suffices to show the equality
$\tr\dot{\tau_0}(g\gamma g^{-1}) = \tr\dot{\tau_0}(\gamma)$
for $\gamma \in G_0(\lambda_+)$.%

A calculation in $SL(2,\CC)$ shows that for every real root $\alpha \in R(\gothg_\CC,\gothh_\CC)$, we have
$\;\Ad \gamma_\alpha = \exp(\mi\pi \ad H_\alpha)\;$ and
$\;\Ad \gamma_\alpha^2 = \id$,  with $\gamma_\alpha$ from \ref{paramètres adaptes}~(a).
Therefore $\gamma_\alpha^2 \in \centregr{G_0}$.
By \cite[Lem. 12.30 (c) p.$\!\!$~469] {Kn86} and (a), every element
$\gamma$ of $G_0(\lambda_+)$ can be written as
$\;\gamma = h\,z\,\gamma_{\alpha_1}\!\dots\gamma_{\alpha_k}\;$ with
$h \in \exp\gothh$, $z \in \centregr{G_0}$ and
real roots $\alpha_1,\dots,\alpha_k \in R(\gothg_\CC,\gothh_\CC)$.
Given Proposition \ref{morphisme rho et fonction delta} (b), the differential at $1$ of
$\vrule width 0ex height 3.2ex depth 2.2ex
\smash{\frac{\rho_{\lambdatilde,{\gotha^*}^+}^{\gothg/\gothh}}
{|\rho_{\lambdatilde,{\gotha^*}^+}^{\gothg/\gothh}|}}$
is the half sum of $\alpha \in R(\gothg(\nu)_\CC,\gothh_\CC)$ such that:
$\mi \mu(H_{\alpha}) > 0$ or,
$\mi \mu(H_{\alpha}) = 0$ and $\rho_\F+\!(H_{\alpha}) > 0$.
Therefore $G_0(\lambdatilde)$ leaves invariant the unitary characters by which $\exp\gothh$
and $\centregr{G_0}$ act in the space of $\dot{\tau_0}$.
In the rest of this proof, we will therefore take $\gamma$ of the form
$\,\gamma = \gamma_{\alpha_1}\cdots\gamma_{\alpha_k}\,$
with real roots $\alpha_1,\dots,\alpha_k \in R(\gothg_\CC,\gothh_\CC)$.
\vskip-1mm%\smallskip

For every real root $\alpha \in R(\gothg_\CC,\gothh_\CC)$, the element $\,\gamma_\alpha g \gamma_\alpha^{-1}\,$ (resp. $\gamma_\alpha g^{-1} \gamma_\alpha^{-1}$) is equal to $g$ (resp. $g^{-1}$) if $\beta(H_\alpha)$ is even and to $g^{-1}$ (resp. $g$) otherwise.
Conse\-quently, we have:\\
$\;g (\gamma g^{-1} \gamma^{-1})
= \cases{\; 1 & if $\,\beta(H_{\alpha_1}+\cdots+H_{\alpha_k})\,$ is even \cr\; g^2 & otherwise.}$%
%\smallskip

First, consider the case where $\beta$ is real and $\;g \gamma g^{-1} \gamma^{-1} \not= 1$.
We make the choices of $g$ and $\gamma_\beta$ compatible  by choosing
$X\!=\!X_\beta$ and $Y\!=\!X_{-\beta}$, so $\,\gamma_\beta = g^2\,$ (see \ref{paramètres adaptes} (a)).
Denote by $\zeta$ the ratio of the homothety $\dot{\tau_0}(\gamma_\beta^2)$.
We have
$\,\gamma_\beta \gamma \gamma_\beta^{-1} = \gamma_\beta^2 \gamma\,$
and
$\,\gamma_\beta (g \gamma g^{-1})\gamma_\beta^{-1}
= \gamma_\beta^2 (g \gamma g^{-1})$.
Therefore $\;\tr\dot{\tau_0}(\gamma) = 0 = \tr\dot{\tau_0}(g\gamma
g^{-1})\;$ when
$\zeta \not= 1$.
We will reason differently when $\zeta = 1$, using
$\;\dot{\tau_0}(g \gamma g^{-1})\dot{\tau_0}(\gamma)^{-1}
= \dot{\tau_0}(\gamma_\beta)$.
The element $\gamma_\beta+\gamma_\beta^{-1}$ of the group algebra of $G_0(\lambda_+)$ is central.
Therefore $\dot{\tau_0}(\gamma_\beta)$ is an homothety when $\zeta\not= -1$,
equal to the 
$\;(\dot{\tau_0}(\gamma_\beta)+\dot{\tau_0}(\gamma_\beta^{-1}))
(\dot{\tau_0}(\gamma_\beta^2)^{-1}+\id)^{-1}$.
In particular, since $\tau_0$ is final, we have
$\;(\delta_{\lambda_+}^{\gothg/\gothh} \tau_0) (\gamma_\beta) = (-1)^{n_\beta+1} \id\;$
when $\zeta = 1$.
Let $F$ denote the subgroup of $G$ generated by the $\gamma_\alpha$ with $\alpha \in R(\gothg_\CC,\gothh_\CC)$ real.
According to \cite[(42) p.$\!\!$~335] {DV88}, to be corrected, there exists a $G(\lambdatilde)$-invariant character $\chi_1$ of $F$ that sends $\gamma_\alpha$ to $(-1)^{n_\alpha+1}$ for every real root $\alpha \in R(\gothg_\CC,\gothh_\CC)$.
Furthermore, the restriction $\chi_2$ of
$\,\vrule width 0ex height 3.3ex depth 2.3ex
\smash{\frac{\rho_{\lambdatilde,{\gotha^*}^+}^{\gothg/\gothh}}
{|\rho_{\lambdatilde,{\gotha^*}^+}^{\gothg/\gothh}|}}
\big(\delta_{\lambda_+}^{\gothg/\gothh}\big)^{\!-1}$
to $G_0(\lambda_+)$ is $G(\lambdatilde)$-invariant (see \ref{espace symplectique canonique}~(a)).
Thus, $\chi \egdef \chi_1 \restriction{\chi_2}{F}$ is a $G(\lambdatilde)$-invariant unitary character of
$F$, such that when $\zeta = 1$ we have:
$\;\dot{\tau_0}(\gamma_\beta)
= \chi(\gamma_\beta)\id\;$
then
$\;\dot{\tau_0}(\gamma_\beta)
= \chi(g\gamma g^{-1})\,\chi(\gamma)^{-1} \id
= \id$.%
\vskip-1mm%\smallskip

Finally, assume that $\beta$ is nonreal and that $g \gamma g^{-1} \gamma^{-1} \not= 1$.
Denote by $R_\beta$ the system of roots
$(\RR\beta \oplus \RR\overline{\beta})\cap R(\gothg_1(\nu)_\CC,\gothh_\CC)$
whose roots are neither imaginary nor real.
Examining the Satake diagram of the semisimple Lie subalgebra of $\gothg_1(\nu)$ of complexication
$(\CC H_\beta \oplus \CC H_{\overline{\beta}})\oplus \Sumpetit_{\alpha \in R_\beta} \gothg_\CC^{\alpha}$,
we see that there is an injective Lie algebra homomorphism
$\;\eta : \gothg_0 \to \gothg_1(\nu)$, with
$\,\gothg_0 = \restriction{\goth{sl}(2,\CC)}{\RR}\,$
or
$\,\gothg_0 = \goth{su}(1,2)$,
which sends, for every $\,(x,y)\in \RR^2\,$ and according to the value of $\gothg_0$,\\
the matrix
$\left(\begin{smallmatrix}
x+\mi y & 0 \\ 0 & -(x+\mi y)
\end{smallmatrix}\right)$
or the matrix
$\left(\begin{smallmatrix}
\mi y & x & 0 \\ x & \mi y & 0 \\ 0 & 0 & -2\mi y
\end{smallmatrix}\right)\,$
to
$\,x\,(H_\beta+H_{\overline{\beta}}) + \mi y\,(H_\beta-H_{\overline{\beta}})$.
Choose
$X \in (\gothg_\CC^\beta \oplus \gothg_\CC^{\overline{\beta}})\cap \gothg$
and
$Y \in (\gothg_\CC^{-\beta} \oplus \gothg_\CC^{-\overline{\beta}})\cap \gothg$,
equals according to $\gothg_0$ to:\\
\hspace*{\fill}%
$X
= \eta\big(\!\!\left(\begin{smallmatrix}
0 & 1 \\ 0 & 0
\end{smallmatrix}\right)\!\!\big)\;$
and
$\;Y
= \eta\big(\!\!\left(\begin{smallmatrix}
0 & 0 \\ 1 & 0
\end{smallmatrix}\right)\!\!\big)$,
or,
$\;X
= \eta\big(\!\!\left(\begin{smallmatrix}
0 & 0 & 1 \\ 0 & 0 & 1 \\ 1 & -1 & 0
\end{smallmatrix}\right)\!\!\big)\;$
and
$\;Y
= \eta\big(\!\!\left(\begin{smallmatrix}
0 & 0 & 1 \\ 0 & 0 & -1 \\ 1 & 1 & 0
\end{smallmatrix}\right)\!\!\big)$.
\hspace*{\fill}%
\\
A~calculation in one of the simply connected groups $SU(2)$ or
$\{1\} \times SU(2)$ then shows that $g^2$ belongs to
$\,\exp(\gothh\cap\derivealg{\gothg_1(\nu)})$.
The definition of $\gothg_1$ ensures that $W(\gothg_1(\nu)_\CC,\gothh_\CC)$ fixes
$\diff_1\dot{\tau_0}$.
Therefore $\dot{\tau_0}$ is trivial on
$\,\exp(\gothh\cap\derivealg{\gothg_1(\nu)})$.
Thus, we obtain:
$\,\dot{\tau_0}(g\gamma g^{-1}\gamma^{-1}) = \dot{\tau_0}(g^2) = \id$.
\cqfd
\end{Dem}
\medskip

In the following proposition, I state the results of Mackey's little group theory that will be useful in proving the next theorem.
I will not review the concepts of “type I group” and “Hilbert integral” (see \cite{Di69}).
For this proposition, we consider a (Hausdorff) locally compact group $A$, whose topology have a countable base, and a normal closed subgroup $B$ of $A$ which is of type $I$.
The group $A$ acts canonically on $\widehat{B}$.

A few words about terminology.
A measurable unitary cocycle of $A$ is a measurable map $c$ from
$A \times A$ to the set of complex numbers of absolute value~$1$ that satisfies
$\,c(1,1) = 1\,$
and
$\,c(xy,z)\,c(x,y) = c(x,yz)\,c(y,z)\,$
for $x,y,z \in A$.
In this case, given a separable complex Hilbert space $V$, we will call $c$-projective representation of $A$ in $V$ a map $\wt{\pi}$ from $A$ to the unitary group of $V$, for which each of the functions
$\,x \in A \mapsto \langle \wt{\pi}(x)\cdot v,w \rangle \in \CC\,$
with $v,w \in V$ is measurable, and such that
$\,\wt{\pi}(1) = \id\,$
and
$\,\wt{\pi}(xy) = c(x,y)\,{\scriptstyle \times}\,\wt{\pi}(x)\,\wt{\pi}(y)\,$
for $x,y \in A$.

%------------------------------------------------------------------------------
\begin{Proposition}
\hspace{-15pt}{\bf{(Mackey)}}\quad %\proposition[Mackey]
\hlabel{théorie de Mackey}
%------------------------------------------------------------------------------
%
{\bf (a)}
Let $\pi \in \widehat{B}$.
Denote by $A_\pi$ the stabilizer of $\pi$ in $A$ and by $\pr$ the canonical projection from $A_\pi$ onto $A_\pi/B$.
The group $A_\pi$ is closed in $A$.
There exists a measurable unitary cocycle $c$ of $A_\pi/B$ and an isomorphism class $\wt{\pi}$ of
$c{\scriptstyle \,\circ\,}(\pr{\scriptstyle \times}\pr)$-projective representation of $A_\pi$ that extends $\pi$ in the same Hilbert space.
Furthermore, there exists a nonzero $\sigma$-finite measure $m$ on $\widehat{B}$, unique up to equivalence, such that the orbit of $\pi$ under $A$ in $\widehat{B}$ has an $m$-negligible complement and $a \cdot m$ is equivalent to $m$ pour for all $a \in A$.%
\smallskip

{\bf (b)}
Retain the notation from (a).
We can define the mapping $\;\wt{\eta} \mapsto
\Ind_{A_\pi}^A (\wt{\eta} \circ \pr \otimes \wt{\pi})\;$
from the set of isomorphism classes of irreducible $c^{-1}$-projective representations of $A_\pi/B$, to the set of isomorphism classes of unitary representations of $A$.
It is injective and its image\llabel{Mackey1}
consists of the elements of $\widehat{A}$ whose restriction to $B$ has a common multiple with $\,\int_{\widehat{B\;}} \rho\,\diff m (\rho)$.%
\llabel{Mackey2}
\smallskip

{\bf (c)}
When the measurable quotient space $A \backslash \widehat{B}$ is countably separated (i.e., we can find a sequence $(E_n)_{n \in \NN}$ of measurable subsets of $E = A \backslash \widehat{B}$ such that for $x \not= y$ in $E$, there exists $m \in \NN$ satisfying $x\in E_m$ and $y\notin E_m$), every element of $\widehat{A}$ is obtained as in (b) for a unique orbit of $A$ in $\widehat{B}$.%
\end{Proposition}

\begin{Dem}{Proof of the proposition}

{\bf (a)}
Mackey's “metrically smooth of type $I$” hypothesis translates to “of type $I$” by
\cite[Prop. 4.6.1 p.$\!\!$~95 and Th. 9.1 p.$\!\!$~168]{Di64}.
The group $A_\pi$ is closed in $A$ and $m$ exists by \cite[Th. 7.5 p.$\!\!$~295] {Ma58}.
The existence of $c$ and of $\wt{\pi}$ is given by \cite[Th. 8.2 p.$\!\!$~298] {Ma58}.%
\smallskip

{\bf (b)}
This result is cited in \cite[Th. 8.1 p.$\!\!$~297 and 8.3 p.$\!\!$~300]{Ma58}.
It refers to the representation class of $B$ associated with $m$, constructed at the end of its Section 7 p.$\!\!$~296, which refers to lines 4--6 p.$\!\!$~273.
\smallskip

{\bf (c)}
Follows from \cite[Th. 9.1 p.$\!\!$~302 and end of Section 7 p.$\!\!$~296] {Ma58}.
\cqfd
\end{Dem}
\medskip

The representations classified a little further on are the “tempered” irreducible unitary representations of $G$.
This notion deserves a definition.

%------------------------------------------------------------------------------
\begin{Definition}\rm
\hlabel{représentations tempérées}
%------------------------------------------------------------------------------

{\bf (a)}
Let $\pi_0 \in \widehat{G_0}$.
We say that $\pi_0$ is tempered if it is equivalent to a subrepresentation of a representation
$\Ind_{\tilde{P}}^{G_0} ( \tilde{\sigma} \otimes \tilde{\eta} \otimes
\fonctioncar{\tilde{N}})$,
where $\tilde{P}$ is a parabolic subgroup of $G_0$ of Langlands decomposition
$\tilde{P}\!=\!\tilde{M}\tilde{A}\tilde{N}$,
$\tilde{\sigma}$ is a limit of discrete series representation of $\tilde{M}$,
and $\tilde{\eta}$ is a unitary character of $\tilde{A}$.
It can be shown, as in \cite[Th. 8.53 p.$\!\!$~260 and Th. 12.23~p.$\!\!$~456]{Kn86} that this is equivalent to the following condition:
$\;\abs{D_{G_0}}^{1/2} \tr\pi_0\,$ is bounded on $\Gzerossreg$.%
%\smallskip

{\bf (b)}
Let $\pi \in \widehat{G}$.
By \ref{théorie de Mackey} end of (b) (c), there exist $n \in \NN$ and
$\pi_1,\dots,\pi_n \in \widehat{G_0}$ such that
$\,\restriction{\pi}{G_0}=\Opluspetit_{1 \leq k \leq n}\pi_k$
(see \cite[Prop. 2.1 p.$\!\!$~425] {Di69}, \cite[Prop. 4.6.1 p.$\!\!$~95 and Th. 9.1 p.$\!\!$~168] {Di64},
and 2.1.21 of Zimmer’s book “\href{https://doi.org/10.1007/978-1-4684-9488-4}{Ergodic theory and semi-simple groups}” for $(G/G_0)\backslash\widehat{G_0}$).
We say that $\pi$ is tempered if $\,\pi_1,\dots,\pi_n$ are tempered.%
\end{Definition}

%------------------------------------------------------------------------------
\begin{Theoreme}
\hlabel{injection dans le dual}
%------------------------------------------------------------------------------

{\bf (a)}
We have:
$\;a \cdot T_{\lambdatilde,\tau}^G =\, T_{a\cdot\lambdatilde,a\cdot\tau}^G\;$
for every automorphism $\,a$ of the Lie group $G$.
\smallskip

{\bf (b)}
The set $G\,\backslash\,\XInd_G$ (resp. $G\,\backslash\,\Xfin_G\!$)
injects into $\widehat{G}$ by the map that sends the orbit of $(\lambdatilde',\tau')$
(resp. $(\lambdatilde',{{\gotha'}^*}^+\!,\tau_+')$) to $T_{\lambdatilde',\tau'}^G\!$
(resp. $T_{\lambdatilde',{{\gotha'}^*}^+\!,\tau_+'}^G\!$).
\vskip-0.5mm

The image of this map consists of the representation classes in $\widehat{G}$ that are tempered.
(By Lemma \ref{condition d'intégrabilité}, the $\lambdatilde' \in \gstregtilde$ involved here belong to $\gstregGtilde$).%
\end{Theoreme}

\begin{Dem}{Proof of the theorem}

We will use the notation from the beginning of Part~\ref{RepIII} and Sections \ref{Rep7} and \ref{Rep8}%
%\smallskip

{\bf (a)}
Given a property of induction cited in the proof of Lemma \ref{construction équivalente, cas connexe} (b), we are reduced to proving that
$a \cdot T_{\mutilde}^{M_0} \!= T_{a \cdot \mutilde}^{a \cdot M_0}$
when $G$ is connected.
This equality comes from some formulas for the characters of
$a \cdot T_{\mutilde}^{M_0}$ and $T_{a \cdot \mutilde}^{a \cdot M_0}$
obtained (without Theorem \ref{injection dans le dual}) in the last part of the proof of Theorem \ref{formule du caractère}.
\smallskip

{\bf (b)}
First assume that $G$ is connected.
To begin with, the limits of discrete series representations of $G$ are tempered
according to \cite[Th. 6.8.1 p.$\!\!$~202, Prop. p.$\!\!$~142, Prop. p.$\!\!$~139] {Wa88} (see also \cite[Cor. 12.27 p.$\!\!$~461]{Kn86} for a study of the character property given in Definition \ref{représentations tempérées}~(a)).
It follows that the representation class $T_{\lambdatilde,{\gotha^*}^+\!,\tau_+}^G\!\!$ of Definition \ref{définition de représentations, cas connexe} is tempered by adapting the proof of the implication $(3)\!\Rightarrow\! (1)$ of \cite[Prop. IV 3.7]{BW80}, and applying the implication $(1)\!\Rightarrow\! (3)$ of the same proposition.\\
Consider a tempered irreducible unitary representation of $G$ in a Hilbert space $V\!$.
Let $K$ be a maximal compact subgroup of $G$.
 By Theorem 11.14 (b) of \cite[p.$\!\!$~131] {ABV92} with $z=1$, there exists a “final limit character$\Lambda_V$” such that the $(\gothg_\CC,K)$-module underlying $V\!$ is isomorphic to that of the representation $\pi(\Lambda_V)$ which is defined in \cite[p.$\!\!$~122]{ABV92} modulo \cite[5) p.$\!\!$~742] {KV95}.
 I use the dictionary proposed in the proof of Lemma \ref{construction de représentations, cas connexe} (b).
 The equality (11.196) of \cite[p.$\!\!$~740] {KV95} shows the representation $\pi(\Lambda_V)$ as a quotient of a representation
 $\,\Ind_{\wt{M} . \wt{A} . \wt{N}}^G(
\smash{T_{\lambdatilde',{{\gotha'}^*}^+\!,\tau_+'}^{\wt{M}}}\! \!\otimes
\wt{\eta} \otimes \fonctioncar{\wt{N}})\,$
induced from a parabolic subgroup $\wt{M}\wt{A}\wt{N}$ of $G$, for some
$(\lambdatilde',{{\gotha'}^*}^+\!,\tau_+')\in \smash{\Xfin_{\wt{M}}}\!$ and some real character
$\wt{\eta}$ of $\wt{A}$.
We therefore have $\,\wt{M}\wt{A}\wt{N}=G\,$ according to \cite[Lem. IV 4.9] {BW80} and \cite[(11.197) p.$\!\!$~741] {KV95}, and the group $G$ acts in $V$ by a representation whose isomorphism class is $\,T_{\lambdatilde',{{\gotha'}^*}^+\!,\tau_+'}^G\!$.\\
Consider the case where $(\lambdatilde',{{\gotha'}^*}^+\!,\tau_+')$ satisfies
$\;T_{\lambdatilde',{{\gotha'}^*}^+\!,\tau_+'}^G\!\! = T_{\lambdatilde,{\gotha^*}^+\!,\tau_+}^G\!$.
Introduce the ordered triples $(\Lambda^{can},R^+_{\mi \, \RR},R^+_\RR)$ and
$({\Lambda'}^{can},{R'}^+_{\mi \, \RR},{R'}^+_\RR)$ attached to $\tau_+$ and
$\tau_+'$ as in the proof of \ref{construction de représentations, cas connexe} (b),
and set $\gothh' = \gothg(\lambdatilde')$.
The differentials of the central characters of $\Lambda^{can}$ and ${\Lambda'}^{can}$ are equal to
$\restriction{\mi\lambda}{\gothh}$ and $\restriction{\mi\lambda'}{\gothh'}$.
By \cite[Th. 11.14 (c) p.$\!\!$~131 and Def. 11.6 p.$\!\!$~124]{ABV92}, there exists $g \in G$ such that, denoting by $\zeta$ the map that sends $x \in G(\lambda_+)$ to the determinant of the restriction of $\Ad x^\CC$ to the sum of $\gothg_\CC^{\alpha}$ with $\alpha \in R^+_\RR$ and
$g \alpha \notin {R'}^+_\RR$, we have:
$\,\lambdatilde' = g \lambdatilde\,$
and
$\,{\Lambda'}^{can} =\, g \cdot \big(\Lambda^{can}\otimes \frac{\zeta}{\abs{\zeta}}\big)$.
Using Lemma \ref{espace symplectique canonique} (a) and a calculation made in the proof of \ref{construction de représentations, cas connexe} (b), we obtain:
$\,\vrule width 0ex height 3.2ex depth 2.1ex
\smash{\frac{\rho_{\lambdatilde,{{\gotha'}^*}^+}^{\gothg/\gothh'}}
{|\rho_{\lambdatilde,{{\gotha'}^*}^+}^{\gothg/\gothh'}|}}\,\tau_+'
= g \cdot \big(\smash{\frac{\rho_{\lambdatilde,{\gotha^*}^+}^{\gothg/\gothh}}
{|\rho_{\lambdatilde,{\gotha^*}^+}^{\gothg/\gothh}|}}\,\tau_+\big)$.
Since the chambers of $\,(\gothg(\lambda)(\mi \rho_\F+),\gotha)\,$ are conjugate under
$\,\normagr{G(\lambda)(\mi \rho_\F+)_0}{\gotha}$,
there exists $u \in G(\lambdatilde)$ such that
$\,(\lambdatilde',{{\gotha'}^*}^+) = gu \cdot (\lambdatilde,{\gotha^*}^+)$.
\emph{A~fortiori}, we have $\tau_+' = gu \cdot \tau_+$ by Lemma \ref{résultat inattendu} (b).%
\smallskip

We no longer assume $G$ connected.
We adapt Duflo's proof in \cite[p.$\!\!$~176 to p.$\!\!$~179] {Df82a}.
Note that the closed normal subgroups $B_1 = G_0(\lambda_+)^{\gothg/\gothh}$ and $B_2 = G_0$, respectively of $A_1 = G(\lambda_+)^{\gothg/\gothh}$ and $A_2 = G$, are of type $I$ by \cite[Prop. 2.1 p.$\!\!$~425] {Di69}.
In particular, the measurable spaces $\widehat{B_1}$ and $\widehat{B_2}$ are standard by \cite[Prop. 4.6.1 p.$\!\!$~95 and Th. 9.1 p.$\!\!$~168]{Di64}.
Their quotients under the left actions of the finite groups respectively equal to $A_1/B_1$ and $A_2/B_2$ are therefore countably separated.
%\smallskip

We will introduce below a certain $\tau_0 \in (G_0(\lambda_+)^{\gothg/\gothh})\widehat{~~}$.
Let $G(\lambda_+)_{\tau_0}^{\gothg/\gothh}$ denote the stabilizer of $\tau_0$ in $G(\lambda_+)^{\gothg/\gothh}\!$,
let $G(\lambda_+)_{\tau_0}$ denote the image of $G(\lambda_+)_{\tau_0}^{\gothg/\gothh}$ in $G(\lambda_+)$
and let $M'(\lambdatilde)_{\tau_0}^{\gothm'\!/\gothh}$ denote the inverse image of $G(\lambda_+)_{\tau_0}$ in $M'(\lambdatilde)^{\gothm'\!/\gothh}\!$.
We will use the same notation $\pr$ to denote the canonical projections from $G(\lambda_+)_{\tau_0}^{\gothg/\gothh}$, $G(\lambda_+)_{\tau_0} G_0$ and $M'(\lambdatilde)_{\tau_0}^{\gothm'\!/\gothh}$ onto $G(\lambda_+)_{\tau_0} / G_0(\lambda_+)$.
By Proposition \ref{théorie de Mackey} (a) and (c), there exists an irreducible subrepresentation $\tau_0$ of
$\restriction{\tau_+}{G_0(\lambda_+)^{\gothg/\gothh}}\!$,
a measurable unitary cocycle~$c$ of $G(\lambda_+)_{\tau_0} / G_0(\lambda_+)$,
an isomorphism class $\wt{\tau_0}$ of $c{\scriptstyle \,\circ\,}(\pr{\scriptstyle \times}\pr)$-projective representation of $G(\lambda_+)_{\tau_0}^{\gothg/\gothh}$ that extends $\tau_0$ in its space,
and an isomorphism class $\wt{\eta}_0$ of irreducible $c^{-1}$-projective representations of $G(\lambda_+)_{\tau_0} / G_0(\lambda_+)$, such that:\\
\hspace*{\fill}%
$\tau_+
= \Ind_{G(\lambda_+)_{\tau_0}^{\gothg/\gothh}}^{G(\lambda_+)^{\gothg/\gothh}}
(\wt{\eta}_0 \circ \pr \otimes \wt{\tau_0})$.
\hspace*{\fill}%
\smallskip

We now construct three isomorphisms of unitary representations.\\[1mm]
\hspace*{3mm}
The first, of $M'(\lambdatilde)^{\gothm'\!/\gothh}\!$-modules, goes from the space of
$\,\Ind_{G(\lambda_+)_{\tau_0}^{\gothg/\gothh}}^{G(\lambda_+)^{\gothg/\gothh}}
(\wt{\eta}_0 \circ \pr \otimes \wt{\tau_0})\,$
equipped with $\tau_{M'}$ (see Lemma \ref{construction équivalente, cas connexe} (a)) to that of
$\,\Ind_{M'(\lambdatilde)_{\tau_0}^{\gothm'\!/\gothh}}^{M'(\lambdatilde)^{\gothm'\!/\gothh}}
(\wt{\eta}_0 \circ \pr \otimes (\wt{\tau_0})_{M'}\!)$,
where $(\wt{\tau_0})_{M'}$ is the
$c{\scriptstyle \,\circ\,}(\pr{\scriptstyle \times}\pr)$-projective representation of
$M'(\lambdatilde)_{\tau_0}^{\gothm'\!/\gothh}$ such that
$\,\delta_{\lambda_{+,\gothm'}}^{\gothm'\!/\gothh} (\wt{\tau_0})_{M'}
= \delta_{\lambda_+}^{\gothg/\gothh} \wt{\tau_0}\,$
on $G(\lambda_+)_{\tau_0}$.
It is written $\;\iota_1 \colon \varphi \mapsto \varphi_{M'}\;$ with:\\
\hspace*{\fill}%
\smash{$(\delta_{\lambda_{+,\gothm'}}^{\gothm'\!/\gothh})^{-1} \,\varphi_{M'}
= (\delta_{\lambda_+}^{\gothg/\gothh})^{-1} \,\varphi\;$}
(see the proof of \ref{construction équivalente, cas connexe} (a)).
\hspace*{\fill}%
\\[1mm]
\hspace*{3mm}
The second, for $M'\!$, goes from the space of
$\,\bigl( \Ind_{M'(\lambdatilde)_{\tau_0}^{\gothm'\!/\gothh}}^{M'(\lambdatilde)^{\gothm'\!/\gothh}}
(\wt{\eta}_0 \circ \pr \otimes (\wt{\tau_0})_{M'}\!)\bigr)
\otimes S T_{\lambdatilde}^{M'_0}$
to that of
$\Ind_{G(\lambda_+)_{\tau_0} M'_0}^{M'}\!
\bigl( (\wt{\eta}_0 \circ \pr \otimes
(\wt{\tau_0})_{M'}\!)\otimes S T_{\lambdatilde}^{M'_0}
\bigr)$,
where the representations in tensor products of $\,M'\!=M'(\lambdatilde)M'_0\,$ and $\,G(\lambda_+)_{\tau_0} M'_0\,$ referred to are constructed as $\,T_{\lambdatilde,\tau_{M'}}^{M'}\!$ in \ref{construction équivalente, cas connexe} (a).
It sends an element of the form $\varphi \otimes v$ to the element $\psi$ satisfying\\
\hspace*{\fill}%
$\psi(xy) \,=\, \varphi(\hat{x})\otimes T_{\lambdatilde}^{M'_0}(y)^{-1} S(\hat{x})^{-1} \!\cdot v$
\hspace*{\fill}\\
for $\hat{x} \!\in\! M'(\lambdatilde)^{\gothm'\!/\gothh}$
lying over $x \in G(\lambda_+)$ and $y \in M'_0$.\\[1mm]
\hspace*{3mm}
Denote by $(M'_{G_0},(\tau_0)_{M'})$ the “pair $(M',\tau_{M'})$” obtained from $(G_0,\lambdatilde,{\gotha^*}^+\!,\tau_0)$ instead of $(G,\lambdatilde,{\gotha^*}^+\!,\tau_+)$.
There exists a unique
$c{\scriptstyle \,\circ\,}(\pr{\scriptstyle \times}\pr)$-projective representation
$\wt{T_0}$ of $G(\lambda_+)_{\tau_0} G_0$ that extends
$T_{\lambdatilde,{\gotha^*}^+\!,\tau_0}^{G_0}\!$ so that its action on an element $\phi$ of the space of
$\,\Ind_{M'_{G_0} . U}^{G_0} (((\tau_0)_{M'} \otimes S T_{\lambdatilde}^{M'_0})\otimes \fonctioncar{U})\,$
satisfies\\
\hspace*{\fill}%
$(\wt{T_0}(xg)\!\cdot\! \phi)(h)
= \abs{\det(\Ad x)_{\gothu}}^{1/2}
\,{\scriptstyle \times}\,
((\wt{\tau_0})_{M'}(\hat{x})\otimes S(\hat{x}))\,\cdot\,
\phi(g^{-1}x^{-1}hx)$
\hspace*{\fill}\\
for $\hat{x} \in M'(\lambdatilde)_{\tau_0}^{\gothm'\!/\gothh}$ lying over
$x \in G(\lambda_+)_{\tau_0}$ and $g,h \in G_0$.\\[2mm]
\hspace*{3mm}
The third isomorphism, of $G(\lambda_+)_{\tau_0} G_0$-modules, is provided by the restriction map to $G_0$, from the space of
$\Ind_{(G(\lambda_+)_{\tau_0} M'_0) . U}^{G(\lambda_+)_{\tau_0} G_0}
\bigl( ((\wt{\eta}_0 \circ \pr \otimes
(\wt{\tau_0})_{M'}\!)\otimes S T_{\lambdatilde}^{M'_0})
\otimes \fonctioncar{U} \bigr)$
to that of $\wt{\eta}_0 \circ \pr \otimes \wt{T_0}$
(identified with a space of functions on $G_0$).\\[1.5mm]
\hspace*{3mm}
To summarize, we find that:\\[1.5mm]
\hspace*{\fill}%
$\;T_{\lambdatilde,{\gotha^*}^+\!,\tau_+}^G\!
= \Ind_{G(\lambda_+)_{\tau_0} G_0}^G
(\wt{\eta}_0 \circ \pr \otimes \wt{T_0})$.
\hspace*{\fill}%
\smallskip

It remains to apply Mackey's theory.
According to the connected case study above, we have
$\,T_{\lambdatilde,{\gotha^*}^+\!,\tau_0}^{G_0}\! \in \widehat{G_0}$
and the stabilizer of $T_{\lambdatilde,{\gotha^*}^+\!,\tau_0}^{G_0}\!$ in $G$ is equal to $G(\lambda_+)_{\tau_0} G_0$.
Proposition \ref{théorie de Mackey} (b) provides the bijection
$\wt{\sigma} \mapsto
\Ind_{G(\lambda_+)_{\tau_0} G_0}^G
(\wt{\sigma} \circ \pr \otimes \wt{T_0})$
from the set of isomorphism classes of irreducible $c^{-1}$-projective representations of
$G(\lambda_+)_{\tau_0} / G_0(\lambda_+)$ onto the set of elements of $\widehat{G}$ whose restriction to $G_0$ is equal to the Hilbert direct sum of the
$\,g \cdot T_{\lambdatilde,{\gotha^*}^+\!,\tau_0}^{G_0}$
with $\,\dot{g} \in G/G(\lambda_+)_{\tau_0} G_0$.
Thus, given the connected case, we have: the representation class
$T_{\lambdatilde,{\gotha^*}^+\!,\tau_+}^G\!$ is tempered irreducible,
every $\,\pi \in \widehat{G}\,$ which is tempered is obtained by Proposition \ref{théorie de Mackey} (c),
and by an easy calculation, any
$(\lambdatilde',{{\gotha'}^*}^+\!,\tau_+')\in \Xfin_G\!$ such that
$\;T_{\lambdatilde',{{\gotha'}^*}^+\!,\tau_+'}^G\!\!
=\, T_{\lambdatilde,{\gotha^*}^+\!,\tau_+}^G$
is conjugate to $(\lambdatilde,{\gotha^*}^+\!,\tau_+)$ under $G$.
\emph{A~fortiori} for such a triple
$(\lambdatilde',{{\gotha'}^*}^+\!,\tau_+')$ we~have
$(\lambdatilde',\tau')\in G \cdot (\lambdatilde,\tau)$,
where $\tau'$ is deduced from $\tau_+'$ as in Definition \ref{paramètres adaptes} (c).
\cqfd%
\end{Dem}%

%------------------------------------------------------------------------------
\begin{Remarque}\rm
\vspace*{-7mm}
\hlabel{recupération des paramètres}
%------------------------------------------------------------------------------

In this remark we assume that $\,\lambdatilde \in \gstItilde\,$ and set $\tau = \tau_+$.\\
We will explicitly find the orbit of $(\lambdatilde,\tau)$ under $G$ starting from $\,\tr T_{\lambdatilde,\tau}^G$.
\smallskip

First, I retrieve $G \cdot \lambdatilde$.
Theorem \ref{formule du caractère} provides the equality\\[0.5mm]
\hspace*{\fill}%
$\,\bigl( \tr T_{\lambdatilde,\tau}^G \bigr)_{\!1}
= \dim \tau \,{\scriptstyle \times}\,
\restriction{~\widehat{\beta}_{G \cdot \lambdatilde}}{\calv_1}$
\hspace*{\fill}%
\\[0.5mm]
where $\calv_1$ is an open neighborhood of $0$ in $\gothg$, which therefore intersects every connected component of $\gssreg$.
By Proposition \ref{mesures de Liouville tempérées}~(c) and the end of Remark \ref{explication de la notation support} (1), we deduce that $T_{\lambdatilde,\tau}^G$ determines $\beta_{R_G(G \cdot \lambdatilde)}$.
The tempered Radon measures $\beta_{\Omega}$ on $\gothg^*$ with $\Omega \in G \backslash \gstreg$ are linearly independent, because each $\beta_{\Omega}$ is nonzero and concentrated on $\Omega$.
Given Proposition \ref{bijection entre orbites} (a), the representation class $T_{\lambdatilde,\tau}^G$ determines the orbit $G \!\cdot\! \lambdatilde$.
\smallskip

I will now reconstruct $\tr \tau$ from $\,T_{\lambdatilde,\tau}^G\,$ and $\lambdatilde$.
We already know that $\;\tau(\exp X) = \me ^{\mi \,\lambda(X)}\id\;$ for $X \in \gothh$.
Let $\hat{e} \in G(\lambda_+)^{\gothg/\gothh}$ lie over an elliptic $e \in~G(\lambda_+)$.
Associate it to a $e_0 \in e\exp\gotht(e)$ as in Lemma \ref{descente pour les formes linéaires} (b).
Theorem \ref{formule du caractère} expresses
$\,\bigl( \tr T_{\lambdatilde,\tau}^G \bigr)_{\!e_0}\,$
as a linear combination of the locally integrable functions
$\restriction {~\widehat{\beta}_{G(e_0)\cdot \lambdatilde'_{e_0}}}{\calv_{e_0}} \!$
with
$\,\smash{\classe{\lambdatilde'_{e_0}} \!\!\in G(e_0)\backslash \gezerostItilde}$
that are linearly independent.
The coefficient of $\bigl( \tr T_{\lambdatilde,\tau}^G \bigr)_{\!e_0}\!$ following
$\restriction{~\widehat{\beta}_{G(e_0)\cdot \lambdatilde[e_0]}}{\calv_{e_0}}$
provides $\tr\tau(\hat{e})$.
\cqfr
\end{Remarque}

%===============================================================================
%
\vspace*{-30pt}
\part{.\quad Characters of representations}\label{RepIV}
%
%==============================================================================

Again, consider an element $\,\lambdatilde \!=\! (\lambda,\F+)$ of $\gstregtilde$
and set $\gothh = \gothg(\lambdatilde)$.
Denote by $\mu$ and $\nu$ (resp. $\gotht$ and $\gotha$)
the infinitesimally elliptic and hyperbolic components of $\lambda$ (resp. $\gothh$).
Let ${\gotha^*}^+$ be a chamber of $(\gothg(\lambda)(\mi \rho_\F+),\gotha)$.
I will use the other notations from Part~\ref{RepIII} (including $M$, $\gothm$, $\tau_M$ for Section~\ref{Rep7}, and, $M'$, $\gothm'$, $U$, $\tau_{M'}$ for Section~\ref{Rep8}) without recalling what they represent.
\smallskip

In this section, we consider an elliptic element $e$ of $G$ (not to be confused with the base $\me$ of natural logarithms).

%-------------------------------------------------------------------------------
\section{Character restriction formula}\label{Rep10}
%------------------------------------------------------------------------------

%\noindent
The formula I will now write generalizes that of Bouaziz in \cite[Th. 5.5.3 p.$\!\!$~52]{Bo87} (see also \cite[Th. (7) p.$\!\!$~106]{DHV84}), while reducing to it, and that of Rossmann \cite[p.$\!\!$~64]{Ro80} relating to the case $e=1$.

%------------------------------------------------------------------------------
\begin{Definition}\rm
\hlabel{notations pour les restrictions}
%------------------------------------------------------------------------------

{\bf (a)}
Denote by $D_G$\labind{DG} the function on $G$ whose restriction to each connected component $G^+$ of $G$ sends $x \in G^+\!$ to the coefficient of $T^{r^{\!+}}\!$ in $\det(T \id + \id - \Ad x)$, where $\,r^+\!$ is the common rank of the reductive Lie algebras~$\gothg(x^+)$ as $\,x^+\!$ ranges over the set of semisimple elements of $G^+$ (see \cite[Lem. 1.4.1 p.$\!\!$~6] {Bo87}).
It is analytic and invariant under the action of the group $\interieur{G} \egdef \{\interieur{g} \,;\, g \in G\}$\labind{int(G)}, where $\interieur{g}$\labind{int(g)} denote the map $\,x \in G \mapsto gxg^{-1} \in~G$.%

\smallskip

{\bf (b)}
Denote by $\Gssreg$\labind{G ssreg} the set of semisimple $x \in G$ such that $\gothg(x)$ is abelian.
Therefore $\Gssreg$ is the dense open set of $G$ with negligible complement, formed by the points where $D_G$ does not vanish (see \cite[1.3 p.$\!\!$~5] {Bo87}).
\smallskip

{\bf (c)}
Write
$\;d_e = \frac {1}{2} \dim (1-\Ad e)(\gothg/\gothj)$\labind{de}
and
$D_e = \det (1-\Ad e)_{_{\scriptstyle (1-\Ad e)(\gothg/\gothj)}} > 0$\labind{De},
\\[0.5mm]
independently of the choice of a Cartan subalgebra $\gothj_e$ of $\gothg(e)\,$ to which we associate the Cartan subalgebra $\gothj \egdef \centraalg{\gothg}{\gothj_e}$ of $\Car \gothg$.%
\smallskip

{\bf (d)}
Set\\
\hspace*{\fill}%
$\calv_e\,
=\, \{\,X\in\gothg(e)\mid \abs{\Im z} < \varepsilon_e \,\textrm{ for any eigenvalue }z\textrm{ of }\ad_{\gothg} X\,\}$\labind{Ve}
\hspace*{\fill}\\
\hspace*{\fill}%
and
$\;\,\displaystyle k_e(X)
= \left( \!
\det { \left(
\frac{\me ^{\ad X/2} - \me ^{- \ad X/2}}{\ad X}
\right)}_{\!\!\gothg(e)}
\frac{\det (1 - \Ad (e \exp X))_{\gothg/\gothg(e)}}
{\det (1 - \Ad e)_{\gothg/\gothg(e)}}
\right)^{\!\!1/2} \!\! > 0$\labind{ke}
\hspace*{\fill}\\[1mm]
for $\,X \in \calv_e$,
where
$\;\,\varepsilon_e
= \inf \{\,\theta \in \, \left] 0,2\pi \right] \mid \me ^{\mi \,\theta} \textrm{ is an eigenvalue of } \Ad_G e\,\} \;\leq\; 2\pi$.%
\smallskip

{\bf (e)}
To every $\interieur{G}$-invariant generalized function $\Theta$ on $G$, we associate (see \cite[(6) p.$\!\!$~98] {DHV84}) the $\,\Ad G(e)$-invariant generalized function $\Theta_e$\labind{Θ_e} on $\calv_e$,
determined by the equality
$\;\Theta_e(X) =\, k_e(X)\;{\scriptstyle \times}\; \restriction{\Theta}{e \exp \calv_e} \, ( e \exp X)\;$
of generalized functions in $\,X \in \calv_e$.
\end{Definition}
\medskip

The equality
$\;\displaystyle
G = \!
\bigcup_{\substack
{e_0 \in G \\ \textrm{$e_0$ elliptique}}} \!\!\!
\interieur{G} \cdot (e_0 \exp \calv_{e_0})\;$
of \cite[Lem. 8.1.1 p.$\!\!$~72] {Bo87} (see \cite[Lem. 40 p.$\!\!$~41] {DV93}) shows that in the situation of (e)\llabel{d->e}, the set of $\Theta_{e_0}$ with $e_0 \!\in\! G$ elliptic determine $\Theta$.
\smallskip

According to \cite[3.1 top of p.$\!\!$~21] {Bo87}, any topologically irreducible unitary representation $\,\pi\,$ of $G$ is trace class, ant its character is a locally integrable function$\,\tr \pi\,$ on $G$
invariant under $\interieur{G}$ whose restriction to $\Gssreg$ is analytic.

%------------------------------------------------------------------------------
\begin{Theoreme}
\hlabel{formule du caractère}
%------------------------------------------------------------------------------

Given $\;\tau \in \XInd_G(\lambdatilde)\;$ and $\;\tau_+ \in \Xfin_G(\lambdatilde,{\gotha^*}^+)\;$ associated as in \ref{paramètres adaptes} (c), we have\\[0.5mm]
$\bigl( \tr T_{\lambdatilde,\tau}^G \bigr)_{\!e}
=\; D_e^{-\frac{1}{2}}\!\!\!
\Sum_{\scriptstyle \classe{\lambdatilde'_e} \,\in\,
G(e)\backslash\gestregtilde}
\!\!
{\scriptstyle \abs{\{ {{\gotha'_e}^*}^+ \}}^{-1}}\!
\Biggl(
\Sum_{\substack
{\lambdatilde' \in\, G \cdot \lambdatilde \, \cap \, \gstregtilde(e)\\
\textrm{such that } \lambdatilde'\![e] \textrm{ exists}\\
\hfill \textrm{and } \,\lambdatilde'\![e] \,=\, \lambdatilde'_e}}\!\!\!\!
\pm_{\widehat{e'}}\;\, \mi^{-d_{e',\lambdatilde}}\;
\tr\tau(\widehat{e'})\!
\Biggr)\,
\restriction{\widehat{\beta}_{G(e)\cdot \lambdatilde'_e}}{\calv_e}
\\
= \bigl( \tr T_{\lambdatilde,{\gotha^*}^+\!,\tau_+}^G \bigr)_{\!e}
\!=\, \mi^{-d_e} D_e^{-\frac{1}{2}}\!\!\!\!\!
\Sum_{\scriptstyle \classe{\lambdatilde'_e} \,\in\, G(e)\backslash
\gestregtilde}
\!\!\!\!\!\!\!\!
{\scriptstyle \abs{\{ {{\gotha'_e}^*}^+ \}}^{-1}}
\Biggl(\!\!\!\!\!\!\!
\Sum_{\substack{\lambda_+' \in\, G \cdot \lambda_+ \cap\, \gothg(e)^*\\
\textrm{(hence } \lambdatilde'\![e] \textrm{ exists, see \ref{lemme clef})}\\
\textrm{such that } \,\lambdatilde'\![e] \,=\, \lambdatilde'_e}}\!\!\!\!\!\!\!\!\!\!\!\!\!\!
\pm_{\widehat{e'_+}}
\tr\tau_+(\widehat{e'_+})\!
\Biggr)\,
\restriction{\widehat{\beta}_{G(e)\cdot \lambdatilde'_e}}{\calv_e}$\\[0.5mm]
where the sum over $\lambdatilde'$ is understood by choosing
$\,g \in G\,$ such that $\,\lambdatilde' = g\lambdatilde\,$ and then fixing
$\,\widehat{e'} \in G(\lambdatilde)^{\gothg/\!\gothg(\lambda)(\mi \rho_\F+)}\!$
lying over $\,e' = g^{-1}eg$\\
(resp.: the sum over $\lambda'_+$ is understood by choosing $\,g \in G\,$ such that $\,\lambda'_+ = g\lambda_+$ then setting $\,\lambdatilde' = g\lambdatilde\,$ and fixing $\,\widehat{e'_+} \in G(\lambda_+)^{\gothg/\gothh}$ above $\,e' = g^{-1}eg$),\\
$\pm_{\widehat{e'}}$ is the sign such that
$\;\scalo (\widehat{e'})_{(1-\Ad e')(\gothg/\gothg(\lambda)(\mi \rho_\F+))}
=\, \pm_{\widehat{e'}}\;
\scalo (B_{\lambda_{\textit{can}}})_{(1-\Ad e')(\gothg/\gothg(\lambda)(\mi\rho_\F+))}$\\
(resp.: $\pm_{\widehat{e'_+}}$ is the sign such that
$\;\scalo (\widehat{e'_+})_{(1-\Ad e')(\gothg/\gothh)}
=\, \pm_{\widehat{e'_+}}\;\,
\scalo (B_{\lambda_+})_{(1-\Ad e')(\gothg/\gothh)}$),\\
$d_{e'}^{\gothg(\lambda)(\mi \rho_{\F+})}\!$ is the coefficient $d_{e'}$ relative to $\gothg(\lambda)(\mi \rho_{\F+})$, and $\;d_{e',\lambdatilde} = d_{e'} - d_{e'}^{\gothg(\lambda)(\mi\rho_{\F+})}\!$,\\
$\gothh'_e = \gothg(e)(\lambdatilde'\![e])$ and $\gotha'_e$ is the infinitesimally hyperbolic component of $\gothh'_e$,\\
($\lambda',\Fep+)=\lambdatilde'\![e]$
and $\,\{ {{\gotha'_e}^*}^+ \}$ is the set of chambers of $(\gothg(e)(\lambda')(\mi \rho_{\Fep+}),\gotha'_e)$.
\\[1mm]
The summations relating to $G(e)\backslash \gestregtilde$ involve a finite number of nonzero terms.
The other summations relate to finite sets.
The condition “$\lambdatilde'\![e]$ exists” (see \ref{descente pour les formes linéaires}) can be omitted, as it is satisfied when $\,\tr\tau(\widehat{e'})\not= 0\,$
(calculation of $\tr\tau$ from \ref{paramètres adaptes}~(c)).%
\smallskip

In particular, for almost all $\,X \in \calv_{1}\,$ we have
\\[1.5mm]
\hspace*{\fill}%
$\smash{\det\left( \frac{\me ^{\ad X/2} - \me ^{- \ad X/2}}{\ad X} \right)_{\!\!\gothg}^{\!1/2}}
\,{\scriptstyle \times}\; \tr T_{\lambdatilde,{\gotha^*}^+\!,\tau_+}^G(\exp X)
\;=\; \dim\tau_+ \;{\scriptstyle \times}\;
\widehat{\beta}_{G \cdot \lambdatilde}(X)$.
\hspace*{\fill}
\end{Theoreme}

%------------------------------------------------------------------------------
\begin{Remarque}\rm
%------------------------------------------------------------------------------

Here are two examples and an important clarification.
\smallskip

{\bf (1)}
First example:
$\,G = GL(2,\RR)$ and $\,e = \left({\scriptstyle\Rot{1}}\right)$,
$\,\lambdatilde = (0,\F+)\,$ with
$\F+ \subseteq \mi \,\goth{so}(2)^* \oplus \centrealg{\gothg}^*$
and $\tau = \chi_{\lambdatilde}^G$.
So $\lambdatilde[e]$ exists in $\gestregtilde$.
The sum over $G \cdot \lambdatilde \cap \gstregtilde(e)$ here involves two nonzero, opposite complex numbers.
\smallskip

{\bf (2)}
Second example:
$G$ is the semidirect product of $\ZZ/ 2\ZZ$ by $SL(2,\RR)^2$ where $\dot{1} \in \ZZ/ 2\ZZ$ operates on $SL(2,\RR)^2$ by permutation of coordinates and
$\,e = (\left({\scriptstyle\Rot{1}}\right),\left({\scriptstyle\Rot{1}}\right))$,
$\lambdatilde = (0,\F+)\,$ with $\,\F+ \subseteq (\mi \,\goth{so}(2)^*)^2$
not stable under $\ZZ/ 2\ZZ$ and $\tau = \chi_{\lambdatilde}^G$.
Therefore $\lambdatilde[e]$ exists in $\gestregtilde$.
The two elements $\lambdatilde'$ of $G \cdot \lambdatilde \cap \gstregtilde(e)$ are conjugate under $G(e)$, and satisfy: $\lambdatilde'\![e]$ exists and $\lambdatilde'\![e]=\lambdatilde[e]$.
The corresponding summed complex numbers are nonzero and, of course, equal.%
\smallskip

{\bf (3)}
It may happen that nonlinearly independent measures $\,\widehat{\beta}_{G(e)\cdot \lambdatilde'_e}$
occur with nonzero coefficients in the character formula.\\
Indeed, consider the following case:
$G$ is the canonical semidirect product of the symmetric group $\,\sym{4}\,$ by $\,SL(3,\RR)^4$,
$\lambda=0$, $\gothh$ is the product of two copies of a split Cartan subalgebra of $\,\goth{sl}(3,\RR)\,$ with two copies of a fundamental Cartan subalgebra of $\goth{sl}(3,\RR)$,
${\gotha^*}^+$ is stable under the permutation of the first two components.
Define $\,e = (1\,2)\,$ and $\,g = (1\,3)(2\,4)\,$ in $\sym{4}$.
Therefore, we have:
$\;G \cdot \lambdatilde \cap \gstregtilde(e) = G(e)\cdot \lambdatilde \,\cup\, G(e)\cdot g\lambdatilde$,
$\,\lambdatilde[e]\,$ and $\,(g\lambdatilde)[e]\,$ exist in $\gestregtilde$,
$\;G(e)\cdot \lambdatilde[e] \not= G(e)\cdot (g\lambdatilde)[e]\;$
and
$\;\widehat{\beta}_{G(e)\cdot \lambdatilde[e]} = \widehat{\beta}_{G(e)\cdot (g\lambdatilde)[e]}\;$
(see \ref{explication de la notation support} (3)).
Furthermore, there exists $\,\tau_+ \in \Xfin_G(\lambdatilde,{\gotha^*}^+)$ for which the coefficients of
$\restriction{~\widehat{\beta}_{G(e)\cdot \lambdatilde[e]}}{\calv_e}$
and
$\restriction{~\widehat{\beta}_{G(e)\cdot (g\lambdatilde)[e]}}{\calv_e}$
as sums over $\lambdatilde'$ are equal and nonzero.
\cqfr
\end{Remarque}

%------------------------------------------------------------------------------
\begin{Remarque}\rm
%------------------------------------------------------------------------------

Assume that $G$ is connected and use the notations of Section~\ref{Rep7}
Various authors have, in certain situations, used character formulas equivalent to those in the previous theorem, but relating to other parametrizations.

To link these formulas, fix a system of positive roots $R_0^+(\gothg_\CC,\gothh_\CC)$ of
$\,R(\gothg_\CC,\gothh_\CC)$ that contains the conjugates of its nonimaginary roots.
Let $\,\rho_{\gothg,\gothh,0}$ and $\,\rho_{\gothm,\gotht,0}$ denote the half sums of positive roots associated with $\,R_0^+(\gothg_\CC,\gothh_\CC)\,$ and with
$\,R_0^+(\gothg_\CC,\gothh_\CC) \cap R(\gothm_\CC,\gotht_\CC)$.%
\smallskip

{\bf (a)}
The map from $C(\gothg(\lambda),\gothh)$ to the set of chambers of $(\gothm_\CC,\gotht_\CC)$ in $\,\mi \,\gotht^*\,$ whose closure contains $\mi \mu$, which to an element $\Fp+$ of $C(\gothg(\lambda),\gothh)$
associates the chamber of $(\gothm_\CC,\gotht_\CC)$ containing $\mi \mu_{\gothm,(\lambda,\Fp+),\{0\}}$, is bijective.
Denote by $\calc^+$ the image of $\F+$ under this map.
Let\\[1mm]
\hspace*{\fill}%
$\sg_0 (\F+,\mu)
= \Prod_{\substack
{\alpha \, \in \, R_0^+(\gothg_\CC,\gothh_\CC)\\
\textrm{$\alpha$ compact}}}\!
\sg(\calc^+(H_\alpha))\;\;
{\scriptstyle \times}
\Prod_{\substack
{\alpha \, \in \, R_0^+(\gothg_\CC,\gothh_\CC)\\
\textrm{$\alpha$ noncompact imaginary}}}\!\!\!
\sg(-\calc^+(H_\alpha))$.
\hspace*{\fill}%
\smallskip

{\bf (b)}
Let $\,T = \centragr{M}{\gotht} = \centragr{M}{\gothm} \exp \gotht$.
The map from $\Xirr_M(\mutilde)$ to the set of unitary characters of $T$ of differential $\,\mi \,\mu - \rho_{\gothm,\gotht,0}\,$ which to $\sigma \in \Xirr_M(\mutilde)$ associates the unitary character $\eta$ of $T$ such that
$\;\eta (x\exp X) = \sigma(x,1)\, \me ^{(\mi \,\mu - \rho_{\gothm,\gotht,0})(X)}\;$
for $\,x \in \centragr{M}{\gothm}\,$ and $\,X \in \gotht$, is bijective.
\smallskip

{\bf (c)}
Let $\sigma \in \Xirr_M(\mutilde)$.
Denote by $\eta$ its image under the map of (b).
The function $\Theta_\eta^{(\calc^+)}\!$ of \cite[Th. 23 p.$\!\!$~260] {Va77}
relative to the group $M$ is given by:\\[0.5mm]
\hspace*{\fill}%
$\Theta_\eta^{(\calc^+)}\!
= \,\sg_0 (\F+,\mu)\,{\scriptstyle \times}\, \tr T_{\mutilde,\sigma}^M$.
\hspace*{\fill}\\[1mm]
In particular, when $M_0$ is semisimple and “acceptable,” the function $\Theta_{\mi \,\mu,\calc^+}$ of \cite[p.$\!\!$~305] {Ha65} relative to the group $M_0$ is given by:\\[1mm]
\hspace*{\fill}%
$\Theta_{\mi \,\mu,\calc^+}
= \sg_0 (\F+,\mu)\,{\scriptstyle \times}\, \tr T_{\mutilde}^{M_0}$.%
\hspace*{\fill}%
\vskip1mm

Furthermore, when $G$ is linear, $\,\tr T_{\mutilde,\sigma}^M$ is denoted
$\,\Theta^M \bigl( \mi \mu,\calc^+\!,
\restriction{\sigma(\,\cdot\,,1)}{\centregr{M}}\bigr)\,$
in \cite[p.$\!\!$~397] {KZ82} by identifying $\gotht$ and $\gotht^*\!$.
\smallskip

{\bf (d)}
Consider the case where $G$ is the neutral component of the group of real points of a complex linear algebraic group defined over $\RR$ which is semisimple and simply connected.
This gives us a complex character $\,\xi_\rho\,$ of $\,\exp \gothh_\CC\,$ of differential~$\rho_{\gothg,\gothh,0}$.
Fix a unitary character $\eta$ of $T$ as in (b).
When $\nu$ is $(\gothg,\gotha)$-regular, the function $\theta(TA,\xi_\rho\,\eta,\nu)$ of \cite[p.$\!\!$~244] {Hr83}~is written:\\[1mm]
\hspace*{\fill}%
$\theta(TA,\xi_\rho\,\eta,\nu)
\,=\, \abs{C(\gothg(\lambda),\gothh)}^{-1}\;{\scriptstyle \times}\!
\Sum_{\Fp+ \in \, C(\gothg(\lambda),\gothh)_{reg}}\!\!
\sg_0 (\Fp+\!,\mu)\,{\scriptstyle \times}\,
\tr T_{\lambdatilde',\tau'}^G$
\hspace*{\fill}%
\\
where $\;\lambdatilde' = (\lambda,\Fp+)$,
$\lambda'_+$ is associated with $(\lambdatilde',\gotha^*)$
as in \ref{choix de racines positives} (b),
and $\tau' \in \Xirr_G(\lambdatilde')$ is determined by the equality
$\,\restriction{(\delta_{\lambda'_+}^{\gothg/\gothh} \tau')}{\centregr{M}}\!\!
=\, \restriction{\eta}{\centregr{M}}$.
\cqfr
\end{Remarque}
\medskip

The theorem will be proven in three steps.
The first two steps follow Bouaziz:
from $G$ to $G(\lambda_+) G_0$ (“beginning of the proof”), then from $G(\lambda_+) G_0$ to~$M'\!$ (“continuation of the proof”).
In the last step (“end of the proof”), a variant of Zuckerman's translation functor will enable me to relate the characters of the representations of $M'$ that interest me to the characters studied by Bouaziz.
\smallskip

\begin{Dem}{Beginning of the proof of the theorem}

Let $\;\widehat{e'_+} \in G(\lambda_+)^{\gothg/\gothh}\,$ and
$\;\widehat{e'} \in G(\lambdatilde)^{\gothg/\!\gothg(\lambda)(\mi \rho_\F+)}$
lie over the same elliptical element $e'$ of $G(\lambda_+)$.
The eigenvalues other than $-1$ and $1$ of the restriction of ${\Ad e'}^\CC$ to
$\;\gothw_\CC\egdef
\Sumpetit_{\alpha \in R^+(\gothg(\lambda)(\mi \rho_\F+)_\CC,\gothh_\CC)} \!\gothg_\CC^\alpha\;$
are nonreal pairwise conjugates of absolute value $1$.
Since $\call_{\lambdatilde,{\gotha^*}^+}/\gothh_\CC$
(see \ref{espace symplectique canonique} (a))
is a direct sum of the ${\Ad e'}^\CC\!$-invariants subspaces
$\call_{\lambdatilde}/\gothg(\lambda)(\mi \rho_\F+)_\CC\!$
(see \ref{espace symplectique canonique} (c)) and~$\gothw_\CC$, we have
$\;\,q_{\,\call_{\lambdatilde,{\gotha^*}^+}/\gothh_\CC}
((\Ad e')_{\gothg/\gothh})
=\, q_{\,\call_{\lambdatilde}/\gothg(\lambda)(\mi \rho_\F+)_\CC\!}
((\Ad e')_{\gothg / \gothg(\lambda)(\mi \rho_\F+)\!})\;\,$
relative to $B_{\lambda_+}\!$,
and then by \ref{espace symplectique canonique} and \ref{géométrie métaplectique} (d):\\
\hspace*{\fill}%
$\pm_{\widehat{e'_+}}\;
\rho_{\lambdatilde,{\gotha^*}^+}^{\gothg/\gothh}(\widehat{e'_+})
\;=\;
\pm_{\widehat{e'}}\;
\rho_\lambdatilde^{\gothg/\!\gothg(\lambda)(\mi \rho_\F+)}(\widehat{e'})\;
{\scriptstyle \times}\;
\mi^{-1/2 \dim (1-\Ad e')(\gothg(\lambda)(\mi\rho_\F+) / \gothh)}$
\hspace*{\fill}\\
where the signs $\;\pm_{\widehat{e'_+}}$ and $\pm_{\widehat{e'}}$ are defined as in the theorem.

Let $\lambdatilde'_0 \in G \cdot \lambdatilde \cap \gstregtilde(e)$,
$g_0 \in G\,$ satisfying $\,\lambdatilde'_0 = g_0\lambdatilde$, and
$\widehat{e'_0} \in G(\lambdatilde)^{\gothg/\!\gothg(\lambda)(\mi \rho_\F+)}\!$
lie over $e'_0 = g_0^{-1}eg_0$.
By Definition \ref{paramètres adaptes} (c) and the above, we obtain\\
\hspace*{\fill}%
$\pm_{\widehat{e'_0}}\;\,\mi^{-d_{e'_0,\lambdatilde}}
\tr\tau(\widehat{e'_0})
\;= \Sum_{\substack{\dot{x} \,\in\, G(\lambdatilde)/G(\lambda_+)\\
\mathrm{tel\ que\ } x^{-1}e'_0\,x \,\in\, G(\lambda_+)}}\!\!
\pm_{\widehat{e'_{+,x}}}\;\,\mi^{-d_e}\tr\tau_+(\widehat{e'_{+,x}})$
\hspace*{\fill}\\
independently of the choice of
$\;\widehat{e'_{+,x}} \in G(\lambda_+)^{\gothg/\gothh}\,$
that lies over $\,e'_x = x^{-1}e'_0\,x$.
\\[0.5mm]
Next, for every $\lambdatilde'_e \in \gestregtilde$ we obtain the equality of finite sums\\
\hspace*{\fill}%
$\Sum_{\substack
{\lambdatilde'_0 \in\, G \cdot \lambdatilde \, \cap \, \gstregtilde(e)\\
\textrm{such that } \,\lambdatilde'_0[e] \,=\, \lambdatilde'_e}}\!\!
\pm_{\widehat{e'_0}}\,\mi^{-d_{e'_0,\lambdatilde}}
\tr\tau(\widehat{e'_0})
\;\,=\,
\Sum_{\substack{\lambda_+' \in\, G \cdot \lambda_+ \cap\, \gothg(e)^*\\
\textrm{such that } \,\lambdatilde'\![e] \,=\, \lambdatilde'_e}}\!\!
\pm_{\widehat{e'_+}}\,\mi^{-d_e}
\tr\tau_+(\widehat{e'_+})$.
\hspace*{\fill}\\
The two expressions of character proposed in the theorem are therefore equal.
\smallskip

Let $\lambdatilde'\!$ (resp. $\lambda_+'$) be a term associated with a $\lambdatilde'_e$ as in the theorem. We~have:\\
$\{\; g\lambdatilde \in G \cdot \lambdatilde \cap \gstregtilde(e)\mid
(g\lambdatilde)[e] = \lambdatilde'_e \textrm{ and }
g\lambdatilde \in G(e)\cdot \lambdatilde' \,\}
= G(e)(\lambdatilde'_e)\cdot \lambdatilde'$\\
(resp.
$\;\{\; g\lambda_+ \in G \cdot \lambda_+ \cap \gothg(e)^* \mid
(g\lambdatilde)[e] = \lambdatilde'_e \textrm{ and }
g\lambda_+ \in G(e)\cdot \lambda_+' \,\}
= G(e)(\lambdatilde'_e)\cdot \lambda_+'\,$).\\[1mm]
The theorem is therefore equivalent to each of the following formulas $\!(\boldsymbol{F})\!$ and $\!(\boldsymbol{F\!}_+)$:\\[1mm]
$\bigl( \tr T_{\lambdatilde,\tau}^G \bigr)_{\!e}
\stackrel{\raisebox{.5ex}{$\scriptscriptstyle(\scriptstyle\boldsymbol{F}\scriptscriptstyle)$}}{=}\,
D_e^{-\frac{1}{2}}\!\!
\Sum_{\substack
{\classe{\lambdatilde'} \,\in\, G(e)\backslash\,
G \cdot \lambdatilde \, \cap \, \gstregtilde(e)\\
\textrm{such that } \,\lambdatilde'\![e] \,\in\, \gestregtilde}}\!\!
\pm_{\widehat{e'}}\; \mi^{-d_{e',\lambdatilde}}\,
\textstyle\frac
{\abs{G(e)(\lambdatilde'\![e])\cdot \lambdatilde'}}
{\abs{\{ {{\gotha'_e}^*}^+ \}}}\;
\tr\tau(\widehat{e'})\,
\restriction{~\widehat{\beta}_{G(e)\cdot \lambdatilde'\![e]}}{\calv_e}$\labind{F};\\
$\bigl( \tr T_{\lambdatilde,{\gotha^*}^+\!,\tau_+}^G \bigr)_{\!e}\!
\stackrel{\raisebox{.5ex}{$\scriptscriptstyle(\scriptstyle\boldsymbol{F\!}_+\scriptscriptstyle)$}}{=}\,
\mi^{-d_e} D_e^{-\frac{1}{2}}\!\!
\Sum_{\substack
{\classe{\lambda'_+} \,\in\, G(e)\backslash\,
G \cdot \lambda_+ \cap\, \gothg(e)^*\\
\textrm{such that }\,\lambdatilde'\![e] \,\in\, \gestregtilde}}\!\!
\pm_{\widehat{e'_+}}\,
\textstyle\frac
{\abs{G(e)(\lambdatilde'\![e])\cdot \lambda'_+}}
{\abs{\{ {{\gotha'_e}^*}^+ \}}}\,
\tr\tau_+(\widehat{e'_+})\,
\restriction{~\widehat{\beta}_{G(e)\cdot \lambdatilde'\![e]}}{\calv_e}$.\\[1mm]
In particular this makes it possible to specify the results set out in the introduction.%
\smallskip

According to Definition \ref{mesures de Liouville}, we have
$\;\widehat{\beta}_{\omega}
=\!\! \Sumpetit_{\omega_0 \in G(e)_0 \backslash \omega}\!\!\!
\widehat{\beta}_{\omega_0}\,$
for $\omega \in G(e)\backslash \gestregtilde$.
Thus, the theorem is equivalent to the statement obtained by replacing $G(e)$ with $G(e)_0$, and therefore also to the statement obtained by replacing $G(e)$ with $G(e)_0$ in $\!(\boldsymbol{F\!}_+)$.
Let $\lambda'_+$ and $\lambdatilde'$ be as in $\!(\boldsymbol{F\!}_+)$.
The group $G(e)_0(\lambdatilde'\![e])$ operates transitively on $\{ {{\gotha'_e}^*}^+ \}$ by the action of certain elements of $W(\gothg(e)_\CC,{\gothh'_e}_\CC)$ (theorem notation).
Since $\lambda'_+$ belongs to ${\gothh'_e}^* \cap \gestssreg$, the stabilizer of an element of $\{ {{\gotha'_e}^*}^+ \}$ is equal to $G(e)_0(\lambdatilde'\![e])(\lambda'_+)$ by Lemma \ref{lemme clef}.
Therefore
$\;\abs{G(e)_0(\lambdatilde'\![e])\cdot \lambda'_+}
= \abs{\{ {{\gotha'_e}^*}^+ \}}$.
%\medskip

Consequently, the formulas to be proven are equivalent to:%
\\[0.5mm]
$\!\bigl( \tr T_{\lambdatilde,{\gotha^*}^+\!,\tau_+}^G \bigr)_{\!e}\!
\stackrel{\raisebox{.5ex}{$\scriptscriptstyle(\scriptstyle\boldsymbol{F\!}_{+,0}\scriptscriptstyle)$}}{=}
\mi^{-d_e} \, D_e^{-\frac{1}{2}}\!\!\!\!\!\!\!\!
\Sum_{\substack
{\classe{\lambda'_+} \in\, G(e)_0 \backslash\,
G \cdot \lambda_+ \cap\, \gothg(e)^*\\
\textrm{such that }
(g\lambdatilde)[e] \,\in\, \gestregtilde}}\!\!\!\!\!\!\!
{\textstyle \frac
{\scalo (\widehat{e'_+})_{(1-\Ad e')(\gothg/\gothh)}}
{\scalo (B_{\lambda_+})_{(1-\Ad e')(\gothg/\gothh)}}}\;
\tr\tau_+(\widehat{e'_+})\,
\restriction{~\widehat{\beta}_{G(e)_0 \cdot (g\lambdatilde)[e]}}{\calv_e}$\labind{F_{+,0}}
\\
where
$\;\,g \in G\,$ satisfies $\;\lambda'_+ \!= g\lambda_+\;$
and
$\;\widehat{e'_+} \in G(\lambda_+)^{\gothg/\gothh}\,$ lies over $\,e' = g^{-1}eg$.
\\[0.5mm]
For later use (Section ~\ref{Rep11}) note that in the previous argument, we can replace $G(e)_0$ by any subgroup $L$ of $G(e)$ that contains $G(e)_0$ and whose adjoint group is included in $\interieur\gothg(e)_\CC$.
\smallskip

Since the orbits over which we sum are regular semisimple orbits, we can reproduce most of the arguments of \cite{Bo87}.
\smallskip

We have
$\;T_{\lambdatilde,{\gotha^*}^+\!,\tau_+}^G\!
\!=\, \Ind_{G(\lambda_+)\, G_0}^G
T_{\lambdatilde,{\gotha^*}^+\!,\tau_+}^{G(\lambda_+)\, G_0}\,$
by “induction in stages,” and then\\
\hspace*{\fill}%
$\displaystyle
\bigl( \tr T_{\lambdatilde,{\gotha^*}^+\!,\tau_+}^G \bigr)_{\!e}
\;= \Sum_{\dot{x} \,\in\, G / G(\lambda_+)\, G_0\,
\textrm{ such that } \,x^{-1}ex \,\in\, G(\lambda_+)\, G_0} \!\!\!
\Ad x \cdot \bigl(\tr T_{\lambdatilde,{\gotha^*}^+\!,\tau_+}^{G(\lambda_+)\, G_0}\bigr)_{\!x^{-1}ex}$.%
\hspace*{\fill}%

Since the set $\;\,G \cdot \lambda_+ \cap\, \gothg(e)^*\;$ is the union of the disjoint sets
$\,\Ad^* x \cdot \bigl( G(\lambda_+) G_0 \cdot \lambda_+ \cap\, \gothg(x^{-1}ex)^* \bigr)\,$
where $\,\dot{x} \in G / G(\lambda_+) G_0$ satisfies $\,x^{-1}ex \in G(\lambda_+) G_0$, we can assume that $\,G = G(\lambda_+)\, G_0$.
\cqfdpartiel
\vspace*{-5mm}
\end{Dem}

%-------------------------------------------------------------------------------
\section{Transition from \texorpdfstring{$\,T_{\lambdatilde,\tau_{M'}}^{M'}\!$}{T\_\{\~\lambda,\tau\_M'\}\^{}M'}
to \texorpdfstring{$\,T_{\lambdatilde,{\gotha^*}^+\!,\tau_+}^G\!$}{T\_\{\~\lambda,a*+,\tau+\}\^{}G}}\label{Rep11}
%------------------------------------------------------------------------------

%\noindent
In this section, we restrict ourselves to proving the theorem for $\,T_{\lambdatilde,\tau_{M'}}^{M'}$ repeating Bouaziz's proof for the regular semisimple case verbatim.%
\smallskip

The results stated in the following proposition are proved in \cite[Lem. 7.1.3 p.$\!\!$~65 and Lem. 4.2.1 p.$\!\!$~38]{Bo87}.

%------------------------------------------------------------------------------
\begin{Proposition}
\hspace{-15pt}{\bf{(Bouaziz)}}\quad %\proposition[Bouaziz]
\hlabel{caractères d'induites}
%------------------------------------------------------------------------------
%
{\bf (a)}
Let $\pi$ be a topologically irreducible unitary representation of $M'$.
Associate to it
$\;\Pi = \Ind_{M'U}^{M'G_0} (\, \pi \, \otimes \, \fonctioncar{U} \,)$.
The unitary representation $\Pi$ of $M'G_0$ is a Hilbertian sum of a finite number of topologically irreducible subrepresentations (see \cite[Prop. 5.4.13 (i) p.$\!\!$~110] {Di64} or
\cite[Prop. p.$\!\!$~25] {Wa88}).

Assume that $e \in M'G_0$.
For every $\,\gothj_e \in \Car \gothg(e)$ and every $\,X \in \gothj_e$ such that $\,e \exp X \in \MpGzerossreg$, we have%
\\[1mm]
\parbox{\textwidth}{
$\quad\; \abs{D_G(e\exp X)}^{1/2}\;\; (\tr\Pi) (e\exp X)$\smallskip

$\displaystyle
=\! \Sum_{ (M' \cap G_0)\,g \;\in\; M' \cap G_0 \,\backslash\,
\{ g_0 \in G_0 \,\mid\, g_0 \, e\exp \gothj_e\,g_0^{-1} \subseteq M' \} }
\!\!\!
|D_{M'}(g e \exp X g^{-1})|^{1/2}\; (\tr\pi) (g e \exp X g^{-1})$
}%
\\[1mm]
where the sum on the right relates to a finite set.
\smallskip

{\bf (b)}
Let $\nu' \in \gothg^*$ be a hyperbolic linear form, let $\gothh'$ be a fundamental Cartan subalgebra of $\gothg(\nu')$, and let $M'_1$ be a subgroup of $G_0(\nu')$ that contains the group $\centragr{G_0}{\gothh'} G(\nu')_0$.
For all $\lambda' \in {\gothh'}^* \cap \gstssreg$,
$\gothj \in \Car \gothg$ and $X \in \gothj \cap \gssreg$, we have\\[1mm]
\hspace*{\fill}%
$\abs{\cald_{\gothg}(X)}^{1/2}\;\;
\widehat{\beta}_{G_0 \cdot \lambda'} (X)\;
= \; \Sum_{ M'_1\,x \;\in\; M'_1 \,\backslash\,
\{ x_0 \in G_0 \,\mid\, x_0 \gothj \,\subseteq\, \gothg(\nu')\} }
|\cald_{\gothg(\nu')}(x\,X)|^{1/2}\;\;
\widehat{\beta}_{M'_1 \cdot \lambda'} (x\,X)$
\hspace*{\fill}\\
where the sum on the right relates to a finite set.
\vspace*{-1mm}
\end{Proposition}

\begin{Dem}{Continuation of the proof of the theorem}

We assume here the equality $G = G(\lambda_+)\, G_0$, to which we have reduced the problem.\\
We (provisionally) assume that the theorem holds for $\,T_{\lambdatilde,\tau_{M'}}^{M'}\!$.
\smallskip

To obtain the character formula for
$\,\tr T_{\lambdatilde,{\gotha^*}^+\!,\tau_+}^G\!$, it suffices to verify it on the open set formed by the $\;X \in \calv_e\,$ such that $\;e \exp X \in \Gssreg$, included in $\,\gessreg$ and of negligible complement (see \cite[before Lem. 1.4.1 p.$\!\!$~6] {Bo87}), on which the two sides of the equality are analytic functions.

Fix $\,\gothj_e \in \Car \gothg(e)$ and $\,X \in \gothj_e \cap \calv_e\,$ such that $\,e \exp X \in \Gssreg$.
We have:\\[1mm]
\hspace*{\fill}%
$k_e(X)
\,=\, |\det (1-\Ad e)_{\gothg/\gothg(e)}|^{-1/2} \,{\scriptstyle \times}\,
|\cald_{\gothg(e)}(X)|^{-1/2} \,{\scriptstyle \times}\,
|D_G(e \exp X)|^{1/2}$.
\hspace*{\fill}%
\smallskip

We paraphrase Bouaziz's calculations as follows.
We fix representatives $g_1,\dots,g_n$ of double classes in
$\;M' \cap G_0 \,\backslash\, \{ g_0 \in G_0 \mid g_0 e g_0^{-1} \in M' \} \,/\, G(e)_0$.,
and then set $\;e_1 =~g_1e g_1^{-1},\dots,e_n = g_ne g_n^{-1}$.
The disjoint union of the sets\\[0.5mm]
\hspace*{\fill}%
$M' \cap G(e_j)_0 \,\backslash\,
\{ x_0 \in G(e_j)_0 \mid x_0 g_j \, \gothj_e \subseteq \gothm'(e_j)\}\;$
where $\,1 \leq j \leq n$
\hspace*{\fill}\\
(we can write “$\subseteq \gothm'$” instead of “$\subseteq \gothm'(e_j)$”)
is in bijection with the set\\
\hspace*{\fill}%
$M' \cap G_0 \,\backslash\,
\{ g_0 \in G_0 \mid g_0 \, e \exp \gothj_e \, g_0^{-1} \subseteq M' \}$
\hspace*{\fill}\\
via the map that sends $(j,\dot{x})$ to $\dot{x g_j}$.
Write $\,\bigl( \tr T_{\lambdatilde,{\gotha^*}^+\!,\tau_+}^G \bigr)_{\!e} (X)\,$
using Proposition \ref{caractères d'induites} (a) by taking
$\,\pi = T_{\lambdatilde,\tau_{M'}}^{M'}$
and considering this bijection.
We then use the expression for
$\;\bigl( \tr T_{\lambdatilde,\tau_{M'}}^{M'} \bigr)_{\!e_j}\;$
deduced from formula $\!(\boldsymbol{F\!}_{+,0})\!$ page \pageref{F_{+,0}} with $\,L = M' \cap G(e_j)_0\,$ instead of $M'(e_j)_0$, as explained there.
This provides three summations.
We invert the last two and arrive at:\\[1mm]
\hspace*{\fill}%
$\Sum_{1 \leq j \leq n}\;\;\;
\Sum_{\substack
{\dot{g\lambda_+}\;\in\; M' \cap G(e_j)_0 \,\backslash\,
M' \cdot \lambda_+ \cap\, \gothm'(e_j)^*\\
\textrm{such that } (g\lambdatilde)[e_j] \,\in\, \gejstregtilde}}\;\;
\Sum_{\dot{x}\;\in\; M' \cap G(e_j)_0 \,\backslash\,
\{ x_0 \in G(e_j)_0 \,\mid\, x_0 \, g_j \gothj_e \,\subseteq\, \gothm' \}}
\;\;\cdots\;\;$.
\hspace*{\fill}%
%\smallskip

The last sum is calculated by applying Proposition \ref{caractères d'induites} (b) after replacing $G$ by $G(e_j)$, $\nu'$ by $\nu_+$,
$\gothh'$ by $(g\gothh)(e_j)$,
$M'_1$ by $M' \cap G(e_j)_0$,
$\lambda'$ by $g\lambda_t$ with $\,t \in \left]0,1\right]$ (see \ref{lemme clef}),
$\gothj$ by $g_j\gothj_e$,
and $X$ by $g_jX$,
then by going to the limit $t \to 0^+$ (see \ref{limites de mesures de Liouville}~(b)).
Furthermore, the disjoint union of the sets
$\;M' \cap G(e_j)_0 \,\backslash\, M' \cdot \lambda_+ \!\cap \gothm'(e_j)^*\;$
where $\,1 \leq j \leq n\,$ is in bijection with the set
$\;G(e)_0 \,\backslash\, G \cdot \lambda_+ \!\cap \gothg(e)^*$
by the map that sends $(j,\dot{l})$ to $\dot{g_j^{-1} l}$.
This makes it possible to combine the first two summations into a single one.
\smallskip

Let \smash{$\;\widehat{e'_+} \in G(\lambda_+)^{\gothg/\gothh}\;$}
and \smash{$\;\widehat{e'_{\gothm'}} \in M'(\lambdatilde)^{\gothm'/\gothh}\;$}
lie over an elliptic element $e'$ of $\,G(\lambda_+)$.
The eigenvalues other than $-1$ and $1$ of the restriction of ${\Ad e'}^\CC$ to $\gothu_\CC$ are nonreal pairwise conjugates of absolute value $1$.
Since $\,\call_{\lambda_t}/\gothh_\CC$ (see \ref{fonctions canoniques sur le revêtement double}) is a direct sum of the ${\Ad e'}^\CC$-invariant subspaces $\gothb_{M'} \!/ \gothh_\CC$ and $\gothu_\CC$, we have
$\,q_{\,\call_{\lambda_t} / \gothh_\CC}((\Ad e')_{\gothg/\gothh})
= q_{\,\gothb_{M'} \!/ \gothh_\CC}((\Ad e')_{\gothm'/\gothh})\,$
relative to $B_{\lambda_+}\!$,
then from \ref{géométrie métaplectique} (d):\\
\hspace*{\fill}%
$\frac
{\scalo (\widehat{e'_+})_{(1-\Ad e')(\gothg/\gothh)}}
{\scalo (B_{\lambda_+})_{(1-\Ad e')(\gothg/\gothh)}}\;
\rho_{\lambda_+}^{\gothg/\gothh} (\widehat{e'_+})
=
\frac
{\scalo (\widehat{e'_{\gothm'}})_{(1-\Ad e')(\gothm'/\gothh)}}
{\scalo (B_{\lambda_{+,\gothm'}})_{(1-\Ad e')(\gothm'/\gothh)}}\;
\rho_{\lambda_{+,\gothm'}}^{\gothm'\!/\gothh} (\widehat{e'_{\gothm'}})
\;{\scriptstyle \times}\;
\mi^{-1/2 \dim (1-\Ad e')(\gothg/\gothm')}$.%
\hspace*{\fill}%
\vskip1mm

Using the equality written in the proof of Lemma \ref{construction équivalente, cas connexe} (a), we deduce the character formula for $\tr T_{\lambdatilde,{\gotha^*}^+\!,\tau_+}^G\!$
(assuming this formula for $\tr T_{\lambdatilde,\tau_{M'}}^{M'}$).
\cqfdpartiel
\end{Dem}
\vspace*{-4mm}

%-------------------------------------------------------------------------------
\section{Translation in the sense of Zuckerman}\label{Rep12}
%------------------------------------------------------------------------------

%\noindent
In this section, the aim is to create a tool for the only tricky step: calculating the character of the representation defined after a passage to homology.
\smallskip

I begin by recalling a few classical results.
\vspace*{-1mm}

%------------------------------------------------------------------------------
\begin{Lemme}
\hlabel{composantes primaires}
%------------------------------------------------------------------------------

Let $\,\pi$ be a continuous linear representation of $G_0$ in a complex Hilbert space $V\!$ and let $K_0$ be a maximal compact subgroup of $G_0$.
Under these conditions, denote by $V_{K_0}$ the union of the finite-dimensional vector subspaces of~$V\!$ stable under $K_0$.

For each $N \!\in \NN$, the following properties (i) and (ii) are equivalent:

{\bf (i)}
there exists a finite sequence
$\;V_0=V \supseteq V_1 \supseteq \cdots \supseteq V_N=\{0\}\;$
of closed $G_0$-submodules of $V\!$ such that the quotients $V_{j-1}/V_j$
with $1 \leq j \leq N\,$ are topologically irreducible under $G_0$ and have infinitesimal character;

{\bf (ii)}
the multiplicities of the elements of $\widehat{K_0}$ in $V_{K_0}$ are finite, which implies that $\gothg$ operates on $V_{K_0}$ by derivation from the action of $G_0$ (see \cite[Th. 14 p.$\!\!$~313] {Va77}), and there is a finite sequence
$\;(V_{K_0})_0 \!=\! V_{K_0} \supseteq (V_{K_0})_1 \supseteq \!\cdots\! \supseteq (V_{K_0})_N \!=\! \{0\}\;$
of $(\gothg_\CC,K_0)$-submodules of $V_{K_0}$ such that the quotients
$\,(V_{K_0})_{j-1} \,/\, (V_{K_0})_j\,$ with $1 \leq j \leq N$
are irreducible $(\gothg_\CC,K_0)$-modules.

When these properties hold, $\pi$ is trace class, and the “$\centrealg{U \gothg_\CC}$-primary components”
$V^\chi$\labind{Vchi} of $V$ associated with the characters $\chi$ of the unital algebra $\centrealg{U \gothg_\CC}$, defined as the closures in $V\!$ of the vector subspaces\\
\hspace*{\fill}%
$(V^{\infty})^\chi
\egdef
\{\, v\in V^{\infty} \mid \forall u \in \centrealg{U \gothg_\CC}\;\;
(u-\chi(u))^n \! \cdot v = 0 \, \}\;$
for “$n$ large,”
\hspace*{\fill}\\
form a direct sum in a dense vector subspace of $V\!$;
furthermore, the set $P$ of~the characters $\chi\!$ of $\centrealg{U \gothg_\CC}\!$ such that $(V^{\infty})^\chi \!\!\not=\!\! \{0\}\!$ (“weights of $\centrealg{U \gothg_\CC}\!$ in $V^{\infty}\!$”) is~finite.%
\end{Lemme}
\vspace*{-3.5mm}

\begin{Dem}{Proof of the lemma}

First, note that for any closed $G_0$-submodule $W$ of~$V\!$, the canonical injection from $V_{K_0}/W_{K_0}$ to $(V/W)_{K_0}$ is surjective.
Indeed, if we denote by $\pr_{V/W}$ the canonical projection from $V$ onto $V/W\!$, for any finite-dimensional $K_0$-submodule $F$ of $V/W\!$ the dense vector subspace $\,\pr_{V/W}(V_{K_0} \cap \pr_{V/W}^{-1}(F))\,$ of $F$ is equal to $F$.
Furthermore, for each $\delta \in \widehat{K_0}$ we have a canonical $K_0$-module isomorphism from $V_\delta/W_\delta$ onto $(V_{K_0}/W_{K_0})_\delta\,$
(where $E_\delta$ denotes the isotypic component of type $\delta$ of a $K_0$-module $E$).
Thus for every $\delta \in \widehat{K_0}$, the $K_0$-modules $V_\delta/W_\delta$ and $(V/W)_\delta\,$ are isomorphic.

When (i) holds, the vector spaces $(V_{j-1}/V_j)_\delta$ with $1 \leq j \leq N$ and $\delta \in \widehat{K_0}$ are finite-dimensional according to \cite[Prop. 16 p.$\!\!$~314 and Th. 19 p.$\!\!$~316] {Va77};
hence the vector spaces $V_\delta$ with $\delta \in \widehat{K_0}$ are finite-dimensional;
the sequence $\;(V_0)_{K_0}=V_{K_0} \supseteq (V_1)_{K_0} \supseteq
\cdots \supseteq (V_N)_{K_0}=\{0\}\;$ satisfies condition (ii) according to \cite[end of Th. 14 p.$\!\!$~313] {Va77}.

When (ii) holds, the sequence
$\;\overline{(V_{K_0})_0\!}\,=V \supseteq \overline{(V_{K_0})_1\!}\, \supseteq
\cdots \supseteq \overline{(V_{K_0})_N\!}\,=\{0\}\;$
satisfies condition (i), again according to \cite[Th. 14 p.$\!\!$~313]{Va77}.
\smallskip

Now, assume (i) and (ii). The result “$\pi$ is trace class” is due to Harish-Chandra (see \cite[8.1.2 p.$\!\!$~292] {Wa88}).

By (ii) and \cite[Cor. 7.207 p.$\!\!$~530] {KV95}, we can apply \cite[Prop. 7.20 p.$\!\!$~446] {KV95} both to $V^{\infty}$ and to its $\centrealg{U \gothg_\CC}$-primary components (here \emph{a priori} in the sense of \cite{KV95}).
Consequently, there exists $N \in \NN$ such that the right-hand side of the equality defining $(V^{\infty})^\chi$ is the same for all $n \geq N$, and $V^{\infty}$ is the direct sum of a finite number of its $\centrealg{U \gothg_\CC}$-primary components (so $P$ is finite).\\
We fix $n \in \NN$ that is suitable for  the definitions of all $(V^{\infty})^\chi$, $\chi \in P$.

Let $Q \!\subseteq\! P\!$.
The characters of $\centrealg{U \gothg_\CC}$ being linearly independent, we have\\[1.25mm]
\hspace*{\fill}%
$\Sum_{\chi \in Q} (V^{\infty})^{\chi}
=
\Big\{\, v\in V^{\infty} \mid \forall u \in \centrealg{U \gothg_\CC}\;\;
\Prod_{\chi \in Q}(u-\chi(u))^n \cdot v = 0 \, \Big\}$.
\hspace*{\fill}\\
Consider a Haar measure $\diff_{G_0}$ on $G_0$ and a sequence $\,(\varphi_k)_{k \in \NN}\,$ of nonnegative elements of $C^\infty_c(G_0)$ of integral $1$ for $\diff_{G_0}$, whose supports “fit into any neighborhood of $1$ in~$G_0$.”
We have
$\;\,(\overline{\Sumpetit_{\chi \in Q} V^\chi\!}\,)^{\infty}\! = \Sumpetit_{\chi \in Q} (V^{\infty})^{\chi}\,$
because for all $u \in \centrealg{U \gothg_\CC}$ the continuous endomorphisms
$\,\pi\bigl( \varphi_k \diff_{G_0} \ast \Prodpetit_{\chi \in Q}(u-\chi(u))^n \bigr)\,$
of $V$ with $\,k \in \NN\,$ vanish everywhere on $\,\Sumpetit_{\chi \in Q} (V^{\infty})^{\chi}\,$
and then on $\,\overline{\Sumpetit_{\chi \in Q}\! V^\chi}$, and, on the other hand, converge pointwise on $(\overline{\Sumpetit_{\chi \in Q} V^\chi\!}\,)^{\infty}$
when $k \to +\infty$ to the operator associated with~$\Prodpetit_{\chi \in Q}(u-\chi(u))^n$.%

For each $\,\chi_0 \in P\!$, we therefore have
$\,\bigl( V^{\chi_0} \cap\, \overline{\!\Sumpetit_{\chi \in P \setminus \{\chi_0\}} V^\chi\!} \,\bigr)^{\!\infty}
= \{0\}$,
which shows that the sum of the vector spaces $V^\chi$ with $\chi \in P$ is a direct sum.
\cqfd
\end{Dem}

%------------------------------------------------------------------------------
\begin{Definition}\rm
%------------------------------------------------------------------------------

By “admissible \ $G_0$-module \ of \ finite length” \ (or \ “Harish-Chandra module of $G_0$”), we mean any complex Hilbert space $V\!$ equipped with a continuous linear representation $\pi$ of $G_0$ that satisfies the equivalent properties of the previous lemma.
\end{Definition}
\smallskip

The aim of the following proposition is to adapt Zuckerman's translation functor (see \cite{Zu77}) to the case of a representation of a disconnected group.%

%------------------------------------------------------------------------------
\begin{Proposition}
\hlabel{foncteur de translation}
%------------------------------------------------------------------------------

Fix a real Lie group $\underline{G}$ whose neutral component is equal to $G_0$ and an element $\underline{a}$ of $\underline{G}$ which operates on $G_0$ by the inner automorphism of some element $a$ of $G$.
For simplicity, we assume that the group $\underline{G}$ is generated by $\underline{a} G_0$.

Let $\pi$ be a continuous linear representation of $\underline{G}$ in a complex Hilbert space~$V$ whose restriction to $G_0$ is an admissible $G_0$-module of finite length.%
\smallskip

{\bf (a)}
Define a generalized function $(\tr\pi)_{a G_0}$ on $\,a G_0$, invariant under $\interieur{G}_0$, locally integrable on $\,a G_0$ and analytic on $\,a G_0 \cap \Gssreg$, by setting\\
\hspace*{\fill}%
$(\tr\pi)_{a G_0} (\varphi \, \diff_G)
= \tr\pi (\underline{\varphi} \, \diff_{\underline{G}})\;\;$
for all $\varphi \in C^\infty_c(a G_0)$,
\hspace*{\fill}\\
where $\;\underline{\varphi} \in C^\infty_c(\underline{a} G_0)\;$ is defined by
$\;\underline{\varphi}(\underline{a}x_0) = \varphi(ax_0)\;$
for $x_0 \in G_0$ and, $\diff_G$ and $\diff_{\underline{G}}$ are Haar measures on $G$ and $\underline{G}$ whose restrictions on $G_0$ coincide.

We have:
$\qquad (\tr\pi)_{a G_0}
= \!\Sumpetit_{\chi \in P \textrm{ such that } \underline{a} \, \chi = \chi}\!
(\tr\pi_{\dot{\chi}})_{a G_0}$,
\\[0.5mm]
where $\,P$ is the set of weights of $\centrealg{U \gothg_\CC}$ in $V^{\infty}\!$ and $\pi_{\dot{\chi}}$ denotes for each $\chi \in P$ such that $\underline{a} \, \chi = \chi$ the representation of $\underline{G}$ in $V^\chi$ resulting from $\pi$.

{\bf (b)}
It is assumed there exists a character $\chi$ of $\centrealg{U \gothg_\CC}$ such that $V\!=V^\chi$.

Let $x \in a G_0 \cap \Gssreg$.
Therefore $\,\gothj_x \egdef \gothg(x)\,$ is abelian (see \ref{notations pour les restrictions} (b)) and $\,\gothj \egdef \centraalg{\gothg}{\gothj_x}$ is a Cartan subalgebra of $\gothg$ (see \cite[Lem. 1.4.1 p.$\!\!$~6] {Bo87}).
Let $\,\gothj(x)_1$\labind{j (x)1} denote the connected component of $0$ in $\{ \, Y \in \gothj(x)\mid x\exp Y \in \Gssreg \}$.

There exist a family $(p_l)_l$ of elements of $S(\gothj(x)_\CC^*)$, indexed by those $l \!\in \gothj(x)_\CC^*$ that satisfy ${\chi \!=\! \chiinfty{l}{\gothg}}\!$ (notation of \ref{caractères canoniques} (a)), such that for every $\,Y \in \gothj(x)_1$ we have\\[1mm]
\hspace*{\fill}%
$\abs{D_G(x\exp Y)}^{1/2}\;\, (\tr\pi)_{a G_0}(x\exp Y)\;
= \Sum_{l \in \gothj(x)_\CC^* \textrm{ such that }
\chi = \chiinfty{l}{\gothg}}\!\!
p_l(Y)\;\, \me ^{l(Y)}$.%
\hspace*{\fill}%

Moreover, when $\chi$ is “regular” (i.e. associated with a regular semisimple orbit of $\interieur\gothg_\CC$ in $\gothg_\CC^*$) every $p_l$ is a scalar.
\smallskip

{\bf (c)}
Assume that $\pi$ has a regular infinitesimal character $\chi$.
Choose $\Lambda \in \gothg_\CC^*$ semisimple regular such that $\,\chi = \chiinfty{\Lambda}{\gothg}\!$.
Denote by $C$ the unique chamber of $R(\gothg_\CC,\gothg_\CC(\Lambda))$ in $\,\gothg_\CC(\Lambda)^* \cap \derivealg{\gothg_\CC}^*$ such that $\,\Lambda \in C + \centrealg{\gothg_\CC}^*\!$, and $\,\overline{C}$ its closure.

Let $\pi_F$ be a linear representation of $\underline{G}$ in a finite-dimensional complex Hilbert space $F\!$, whose restriction to $G_0$ is irreducible, and such that the lowest weight of $\gothg_\CC(\Lambda)$ in $F$ relative to $C\!$, denoted $\Lambda_F\!$, satisfies
$\Lambda + \Lambda_F \in \overline{C} + \centrealg{\gothg_\CC}^*\!$.

The restriction to $G_0$ of the representation $\pi \otimes \pi_F$ of $\underline{G}$ yields an admissible $G_0$-module of finite length and $\underline{G}$ operates in its $\,\centrealg{U \gothg_\CC}$-primary component $(V \otimes_\CC F)^{\chiinfty{_{\Lambda + \Lambda_F}}{\gothg}}$ by means of a representation $\pi_{Zuc}$ still having this property.

Let $x \in a G_0 \cap \Gssreg$.
Associate to it, as in (b), a Cartan subalgebra $\gothj$ of $\gothg$ and complex numbers $p_l$ (see end of (b)).
For all $Y \in \gothj(x)_1$, we have\\[1mm]
\hspace*{\fill}%
$\abs{D_G(x\exp Y)}^{1/2}\;\; (\tr\pi_{Zuc})_{a G_0}(x\exp Y)\;
= \Sum_{l \,\in\, \gothj(x)_\CC^* \textrm{ such that }
\chi = \chiinfty{l}{\gothg}}\!\!
\underline{x}_{~\!l_F}\;\, p_l \;\, \me ^{(l \,+\, l_F)(Y)}$,
\hspace*{\fill}\\
where $\;l_F \in \gothj_\CC^*$ is defined from $l$ as $\,\Lambda_F$ from $\Lambda$ and $\,\underline{x}_{~\!l_F}$ is the scalar by which the element $\;\underline{x} \egdef \underline{a}\,(a^{-1}x)\;$ of $\underline{G}$ (which fixes $\,l_F$) acts on the one-dimensional eigenspace of weight $\,l_F$ for the action of $\,\gothj(x)_\CC$ on $F\!$.
\vspace*{-4mm}
\end{Proposition}

\begin{Dem}{Proof of the proposition}

We will use the canonical action by left-invariant differential operators of the enveloping algebra of a real Lie group on the vector space of generalized functions on this Lie group.
\smallskip

{\bf (a)}
The representation $\pi$ is trace class, by Lemma \ref{composantes primaires}.
The generalized function $(\tr\pi)_{a G_0}$, equal to
$\delta_a \ast (\delta_{\underline{a}^{-1}} \ast \restriction{\tr\pi}{\underline{a} G_0})$,
is invariant under $\,\interieur{G}_0$.
\smallskip

Let $\,u_0 \in \centrealg{U \gothg_\CC}\,$ and $\,\underline{\varphi} \in C^\infty_c(\underline{a} G_0)$.
Set $\;u = \Prodpetit_{\smash{\chi \in P}} (u_0 - \chi(u_0))$.
The proof \ of \ \cite[Cor. 7.207 p.$\!\!$~530] {KV95} \ provides \ the \ existence \ of \ a \ finite \ sequence
$\;W_0=V \supseteq W_1 \supseteq \cdots \supseteq W_N=\{0\}\;$
of closed $\underline{G}$-submodules of $V$ such that the quotients $W_{k-1}/W_k$ with $1 \leq k \leq N'\,$ are topologically irreducible under $\underline{G}$.

Let $k \!\in\! \{1,\dots, N'\}$.
Denote by $\pi_k$ the representation of $\underline{G}$ in $W_{k-1}/W_k$.
Any subquotient of an admissible $G_0$-module of finite length is an admissible $G_0\text{-module}$ of finite length.
Therefore $W_{k-1}/W_k$ has a topologically irreducible closed $G_0$-submodule $E$ which has infinitesimal character.
The action of $u$ on $\;\Sumpetit_{n \in \ZZ}\underline{a}^n\,E^\infty\,$ is zero.
Therefore $\;\pi_k (\underline{\varphi} \diff_{\underline{G}} \ast u) = 0$.
\vspace*{-1mm}

Thus:
$\;\partial_u (\tr\pi)\,
(\underline{\varphi} \diff_{\underline{G}})
= \tr\pi_1 (\underline{\varphi} \diff_{\underline{G}} \ast u) + \cdots
+ \tr\pi_{N'} (\underline{\varphi} \diff_{\underline{G}} \ast u)
= 0 \;$
(see \cite[Lem. 8.1.3 p.$\!\!$~293] {Wa88}).
Furthermore, we have $\;\underline{a} \, u = u$.
It follows~that:\\
\hspace*{\fill}%
$\partial_u ((\tr\pi)_{a G_0})
= \delta_a \ast ((\delta_{\underline{a}^{-1}} \ast
\restriction{\tr\pi}{\underline{a} G_0})\ast
(\delta_{\underline{a}} \ast \check{u} \ast \delta_{\underline{a}^{-1}}))
= 0$.
\hspace*{\fill}\\
By \cite[Th. 2.1.1 p.$\!\!$~10] {Bo87} and \cite[Th. 4.95 p.$\!\!$~286 and Th. 7.30 (a) p.$\!\!$~450] {KV95},
$(\tr\pi)_{a G_0}$ is locally integrable on $a G_0$ and analytic on $a G_0 \cap \Gssreg$.
\smallskip

Let $\,\chi \in P$ such that $\underline{a} \, \chi \not= \chi$ and (again) $\,\underline{\varphi} \in C^\infty_c(\underline{a} G_0)$.
Denote by~$m$ the cardinality of the orbit $\dot{\chi}$ of $\chi$ under the action of the subgroup $\langle \underline{a} \rangle$ of $\underline{G}$ generated by $\underline{a}$, and denote by $\pi_{\dot{\chi}}$ the representation of $\underline{G}$ in $\,W \egdef \overline{\Sumpetit_{n \in \ZZ} \underline{a}^n \, V^\chi\!}$. The restriction of $\,\pi_{\dot{\chi}} (\underline{\varphi}\diff_{\underline{G}})\,$ to the direct sum
$\,V^\chi \oplus \underline{a}V^\chi \oplus \cdots \oplus \underline{a}^{m-1}V^\chi\,$
decomposes into blocks with zero diagonal blocks.
I exploit this property using a suggestion by G.~Skandalis to return to the finite-dimensional case.

The subalgebra of $\call (W)$ formed by endomorphisms that simultaneously stabilize $V^\chi$, $\underline{a}V^\chi,\dots,\underline{a}^{m-1}V^\chi$ is closed.
According to \cite[Rem. p.$\!\!$~55 and Prop. 7 p.$\!\!$~47] {B67}, the primary subspaces $W^z$ of the restriction of
$\,\pi_{\dot{\chi}} (\underline{\varphi} \diff_{\underline{G}})^m\,$
to $W$ associated with $\,z \in \CC\moins0$, which are finite-dimensional, are all sums of their intersections with $V^\chi$ and $\cdots$ and $\underline{a}^{m-1}V^\chi$.
Furthermore, each $W^z$ with $\,z \in \CC\moins0$ is equal to the sum of the primary subspaces of the restriction of
$\,\pi_{\dot{\chi}} (\underline{\varphi} \diff_{\underline{G}})\,$
to $W$ associated with the $m^{\textrm{th}}\!$ roots of $z$.
 Lidskij's formula (“the trace is the sum of the eigenvalues”) yields the equality
$\;\,\tr\pi_{\dot{\chi}} (\underline{\varphi} \, \diff_{\underline{G}})
= \Sumpetit_{z \in \CC\moins0}
\tr \big(
\restriction{\pi_{\dot{\chi}} (\underline{\varphi} \,
\diff_{\underline{G}})} {W^z} \big)
= 0$.

Hence, now using \cite{B67} relative to $\call (V)$ (see proof of \ref{composantes primaires}):\\
\hspace*{\fill}%
$(\tr\pi)_{a G_0}
=
\Sum_{\dot{\chi} \in \langle \underline{a} \rangle \backslash P}
(\tr\pi_{\dot{\chi}})_{a G_0}
= \Sum_{\chi \in P \textrm{ such that } \underline{a} \, \chi = \chi}
(\tr\pi_{\dot{\chi}})_{a G_0}$.
\hspace*{\fill}%
\smallskip

{\bf (b)}
I will now explain the calculations of \cite[p.$\!\!$~33 and p.$\!\!$~34] {Bo87}.
\smallskip

Set
$\;\,\theta(Y) = \abs{D_G(x\exp Y)}^{1/2}\; (\tr\pi)_{a G_0}(x\exp Y)\;\,$
for all $\,Y \!\in \gothj(x)_1$.\\
Let $W$ denote the group $W(\gothg_\CC,\gothj_\CC)$ and $\pr$ the canonical projection from $\,S\, \gothj_\CC\,$ onto $\,S\, \gothj(x)_\CC\,$ obtained as in Conventions \ref{conventions fréquentes} (b).
The proof of (a) shows that $(\tr\pi)_{a G_0}$ is an eigenvector of $\centrealg{U \gothg_\CC}$ associated with $\chi$.
Fix an element $l_0$ of $\gothj_\CC^*$ such that $\chi = \chiinfty{l_0}{\gothg}\!$.
According to \cite[line 4 of 2.5 p.$\!\!$~20] {Bo87},, we have:\\[1mm]
\hspace*{\fill}%
$q(l_0)\; \theta = \partial_{\pr(q)} \theta \;\,$
for all $\,q \in (S\, \gothj_\CC)^W\!$.
\hspace*{\fill}%
\medskip

In what follows, we repeat the arguments of \cite[p.$\!\!$~369 to p.$\!\!$~371] {Kn86}.

Let $X \in \gothj_\CC$.
Denote by $q_0,\dots,q_{\abs{W}-1}$ the elements of $(S\, \gothj_\CC)^W$ for which the following equality of polynomials in the indeterminate $T$ holds:\\
\hspace*{\fill}%
$\Prod_{w \in W} (T-wX)
= \; T^{\abs{W}} \,+\, q_{\abs{W}-1} T^{\abs{W}-1} \,+\, \cdots \,+\, q_0\;\;$
in $\,(S\, \gothj_\CC) [T]$.
\hspace*{\fill}\\
From this, first replace $T$ by $X$, then take the value of each member at $l_0$ and replace $T$ by $X$. This yields, in $\,S\, \gothj_\CC$:\\[1mm]
\hspace*{\fill}%
$X^{\abs{W}} \,+\, X^{\abs{W}-1} q_{\abs{W}-1} \,+\, \cdots \,+\, q_0 = 0$
\hspace*{\fill}\\
\hspace*{\fill}%
and
$\;\; X^{\abs{W}} \,+\, X^{\abs{W}-1} q_{\abs{W}-1}(l_0)\,+\, \cdots \,+\,
q_0(l_0)
= \Prod_{w \in W} (X-l_0(wX))$.
\hspace*{\fill}\\
Use the partial differential equation satisfied above by $\theta$ with $q$ successively equal to $\,q_{\abs{W}-1},\dots,q_0$.
This yields the following equality, which remains to be exploited:\\[1mm]
\hspace*{\fill}%
$\Prodpetit_{w \in W} (\partial_{\pr(X)} - l_0(wX))\cdot \theta = 0$.%
\hspace*{\fill}%

Let $X$ range over a base of $\gothj(x)_\CC$.
By \cite[Prop. 3 p.$\!\!$~58] {Va77}, there exists a unique family $(p_l)_{l \in \gothj(x)_\CC^*}$ of elements of $\,S(\gothj(x)_\CC^*)$ with finite support such that\\
\hspace*{\fill}%
$\theta(Y)\;= \Sum_{l \in \gothj(x)_\CC^*} p_l(Y)\; \me ^{l(Y)}\;\,$
for all $\,Y \!\in \gothj(x)_1$.
\hspace*{\fill}\\
More precise information will be obtained by substituting this expression of $\theta$ into the preceding equality (and extending analytically to $\gothj(x)$).
\smallskip

Let $l \in \gothj(x)_\CC^*$ such that $p_l \not= 0$.
Denote by $d$ the degree of $p_l$ in $S(\gothj(x)_\CC^*)$ and by $p_l^{[d]}$ the homogeneous component of $p_l$ of degree $d$.
Let $X \in \gothj_\CC$.
We have\\
\hspace*{\fill}%
$(\partial_{\pr(X)} - l_0(wX))\cdot p_l \; \me ^l \,
= (l(X) - l_0(wX))\, p_l \; \me ^l + (\partial_{\pr(X)} p_l)\; \me ^l\quad$
for all $w \in W$.
\hspace*{\fill}\\
The null polynomial
$\,\me ^{-l} {\scriptstyle \times}\! \Prodpetit_{w \in W} (\partial_{\pr(X)} - l_0(wX))\cdot p_l \, \me ^l\,$
admits
$\Prodpetit_{w \in W} (l - w^{-1} l_0) (X)\,{\scriptstyle \times}\, p_l^{[d]}$
as its homogeneous component of degree $d$.
Since the ring $\,S\, \gothj_\CC^*$ is an integral domain, this shows that $\;l \in W \!\cdot l_0$, hence $\;\chi = \chiinfty{l}{\gothg}\!$.

Now assume $p_l$ to be nonscalar.
Choose $X \in \gothj_\CC$ such that, for all $w \in W$ outside $W(l)$ we have $\,(w\,l)(X)\not= l(X)$, and $p_l^{[d]} (\pr(X))\not= 0$.
Let $(\pr(X)^*)\,$ denote the basis of $(\CC\pr(X))^*$ dual of $(\pr(X))$.
The restriction of the zero polynomial\\
\hspace*{\fill}%
$\me ^{-l} \,{\scriptstyle \times} \!\!
\Prodpetit_{w \in W(l)} \!\!
(\partial_{\pr(X)} - l(wX))\cdot \bigl( \Prodpetit_{w \in W \textrm{ and } w \notin W(l)} \!\!
(\partial_{\pr(X)} - l(wX))
\cdot p_l \; \me ^l \bigr)$
\hspace*{\fill}%
\\[-0.5mm]
to $\,\CC\pr(X)\,$ is the sum of a polynomial of the form
$\;\zeta \,{\scriptstyle\times}\, (\partial_{\pr(X)})^{\abs{W(l)}} \cdot (\pr(X)^*)^d\;$
with $\zeta \in \CC\moins0$ and monomials of (strictly) smaller degrees.
Hence $\,\abs{W(l)} > 1$.
\smallskip

{\bf (c)}
I follow \cite[p.$\!\!$~301] {Zu77}.
\smallskip

By \cite[Cor. 7.207 p.$\!\!$~530] {KV95}, $\pi \otimes \pi_F$ and the representations of $G_0$ on the $\centrealg{U \gothg_\CC}$-primary components of $V \otimes_\CC F$ provide admissible $G_0$-modules of finite length.

By hypothesis, we have $\,\underline{a} \chi = \chi$.
Hence, there exists
$\sigma \in (\Ad\underline{a})^\CC \interieur\gothg_\CC$
such that $\sigma(\Lambda) \!=\! \Lambda$.
We have: $\,\sigma(\gothg_\CC(\Lambda)) \!=\! \gothg_\CC(\Lambda)$ and $\sigma(C) = C$, then $\sigma(\Lambda_F) =~\Lambda_F$.
Therefore
$\;\underline{a} \, \chiinfty{_{\Lambda + \Lambda_F}}{\gothg}
= \sigma \, \chiinfty{_{\Lambda + \Lambda_F}}{\gothg}
= \chiinfty{_{\Lambda + \Lambda_F}}{\gothg}\!$
and the group $\underline{G}$ stabilizes $(V \otimes_\CC F)^{\chiinfty{_{\Lambda + \Lambda_F}}{\gothg}}\!$.

Denote by $P'$ the set of weights of $\centrealg{U \gothg_\CC}$ in $(V \otimes_\CC F)^{\infty}$, by $P(F,\gothj(x)_\CC)$ the set of weights of $\gothj(x)_\CC$ in $F$ and by $F^\gamma$ the eigenspace of $\gothj(x)_\CC$ associated with an element $\gamma$ of $P(F,\gothj(x)_\CC)$.
For all $Y \in \gothj(x)_1$, we obtain thanks to  (a):\\[1mm]
\hspace*{\fill}%
$(\tr\pi)_{a G_0}(x\exp Y)\,{\scriptstyle \times}\,
\tr\pi_F(\underline{x}\exp Y)\,
= \Sum_{\chi' \in P' \textrm{ such that } \underline{a} \, \chi' = \chi'}
(\tr(\pi \otimes \pi_F)_{\dot{\chi}'})_{a G_0}(x\exp Y)$
\hspace*{\fill}\\
\hspace*{\fill}%
and
$\;\; \tr\pi_F(\underline{x}\exp Y)\;
= \!\Sum_{\gamma \in P(F,\gothj(x)_\CC)}\!\!
\tr(\restriction{\pi_F(\underline{x})}{F^\gamma})\; \me ^{\gamma(Y)}$.
\hspace*{\fill}\\
Apply (b) to $\pi$, and to $(\pi \otimes \pi_F)_{\dot{\chi}'}$ with $\chi' \!\in\! P'$ such that $\underline{a} \, \chi'\! = \!\chi'$, among which is $\pi_{Zuc}$.
By identifying the polynomial factors written in front of the functions $\,\me ^l\,$ with $l \in \gothj(x)_\CC^*$, we find that for all $Y \in \gothj(x)_1$ we have\\[1mm]
\hspace*{\fill}%
$\abs{D_G(x\exp Y)}^{1/2}\;\; (\tr\pi_{Zuc})_{a G_0}(x\exp Y)
= \Sum_{(l,\gamma)\,\in\, \cala}
\tr(\restriction{\pi_F(\underline{x})}{F^\gamma})\;\,
p_l \;\, \me ^{(l \,+\, \gamma)(Y)}$,
\hspace*{\fill}\\
where we put $\;\,\cala
\,=\, \{ (l,\gamma)\in \gothj(x)_\CC^* \times P(F,\gothj(x)_\CC)\mid
\chi = \chiinfty{l}{\gothg} \textrm{ and }
\chiinfty{_{\Lambda + \Lambda_F}}{\gothg} = \chiinfty{l + \gamma}{\gothg} \}$.
\smallskip

Finally, we verify that
$\;\cala
= \{ (l,\gamma)\in \gothj(x)_\CC^* \times \gothj_\CC^* \mid
\chi = \chiinfty{l}{\gothg} \textrm{ and } \gamma = l_F \}$,
and $\dim F^\gamma = 1$ for $(l,\gamma)\in \cala$.
The inclusion $\,\supseteq\,$ is clear.
Consider $\,(l,\gamma)\in \cala$.
We modify the choice of $\Lambda$ by taking $\,\Lambda = l$.
We will show that $\,\gamma = \Lambda_F$, and calculate $\dim F^\gamma\!$.
Fix $\,w \in W(\gothg_\CC,\gothj_\CC)\,$ such that $\,l + \gamma = w (\Lambda + \Lambda_F)$.\\
We have
$\;\gamma - \Lambda_F = w (\Lambda + \Lambda_F) - (\Lambda + \Lambda_F)\,$
and $\,w (\Lambda + \Lambda_F) - (\Lambda + \Lambda_F)\leq 0\,$
for the order associated with $C$, because $\,\Lambda + \Lambda_F \in \overline{C} + \centrealg{\gothg_\CC}^* \!$.
On the other hand, there is a weight $\zeta$ of $\gothj_\CC$ in $F$ such that $\;\gamma = \restriction{\zeta}{\gothj(x)_\CC}$.
The inclusion $\,\gothj(x)_\CC^* \subseteq \gothj_\CC^*\;$ identifies $\gamma$ with the isobarycenter of the finite subset $\,\{ \Ad^*x^n \cdot \zeta \}_{n \in \ZZ}\,$ of the set of weights of $\gothj_\CC$ in $F\!$.
Therefore $\,\gamma - \Lambda_F \geq 0\,$ for the order associated with $C$.
Consequently, $\,\gamma = \Lambda_F$.
Since $\gamma$ is an extremal weight of $\gothj_\CC$ in $F$ which is isobarycenter of a finite set of weights of $\gothj_\CC$ in $F\!$, it is equal in particular to $\zeta$.
Thus the eigenspace of $F$ of weight $\gamma$ under the action of $\gothj(x)_\CC$ is of dimension $1$.
\cqfd
\end{Dem}

%-------------------------------------------------------------------------------
\section{Transition from \texorpdfstring{$\,T_{\lambda_r,\sigma_r}^{M'}\!$}{T\_\{\lambda\_r,\sigma\_r\}\^{}M'} to \texorpdfstring{$\,T_{\lambdatilde,\tau_{M'}}^{M'}\!$}{T\_\{\~\lambda,\tau\_M'\}\^{}M'}}\label{Rep13}
%------------------------------------------------------------------------------

%\noindent
To simplify the notation, in this section set:
$\MM'=\GG(\nu_+)$,
$\rho = \rho_{\gothm',\gothh}$, $\sigma = \tau_{M'}$,
and $\lambda_r = \lambda_{+,\gothm'}$ (see page \pageref{λ+m'}).
This yields $\lambda_r = (\mu-2\mi \rho) + \nu \,\in\, \mpstssreg$.
We'll see that, roughly speaking, the character $\tr T_{\lambdatilde,\sigma}^{M'}$ can be deduced from a character of the form $\tr T_{\lambda_r,\sigma_r}^{M'}$ by translation in the Zuckerman sense.
\medskip

The following proposition is obtained by proving a variant of Bouaziz's results in \cite[(i) and (ii) p.$\!\!$~550]{Bo84} (see also \cite[p.$\!\!$~547 and p.$\!\!$~548]{KV95}).
Its origin is Lemma 3.1 on p. 406 of \cite{KZ82}, which yields equality of restrictions to~$M'_0$ of the characters considered in (b) below.
\smallskip

%------------------------------------------------------------------------------
\begin{Proposition}
\hlabel{passage au cas régulier}
%------------------------------------------------------------------------------

Fix $\,\hat{a} \in M'(\lambdatilde)^{\gothm'/\gothh}\,$ lying \ over \ an \ elliptical \ element \ $a$ \ of $M'(\lambdatilde)$.
Denote by $\underline{M'\!}\,$ the semidirect product of $\ZZ$ by $\,M'_0$ such that $\,\underline{a} \egdef 1 \in \ZZ$ acts on $\,M'_0$ by $\interieur{a}$.
\smallskip

{\bf (a)}
There exists a continuous linear representation $\pi_F$ of $\underline{M'\!}\,$ in a finite-dimensional complex Hilbert space $F$, unique to within isomorphism, whose restriction to $M'_0$ is irreducible with lowest weight $-2\rho$ for the action of $\gothh_\CC$ and relative to the order deduced from $R^+(\gothm'_\CC,\gothh_\CC)$, such that $\underline{a}$ acts trivially on the eigenspace of weight $-2\rho$.
\smallskip

{\bf (b)}
Extend the unitary representations $T_{\lambdatilde}^{M'_0}\!$ and $T_{\lambda_r}^{M'_0}\!$ of $M'_0$ whose “spaces” are denoted $\calh$ and $V\!$, into continuous linear representations $\pi_{\lambdatilde}$ and $\pi_{\lambda_r}\!$ of $\underline{M'\!}\,$ by the conditions
$\,\pi_{\lambdatilde}(\underline{a}) = S_{\lambdatilde}(\hat{a})\,$
and
$\,\pi_{\lambda_r}(\underline{a}) = S_{\lambda_r}(\hat{a})$,
where $\,S_{\lambdatilde}\,$ and $\,S_{\lambda_r}$ are the representations of $M'(\lambdatilde)^{\gothm'/\gothh}$ (which is equal to $M'(\lambda_r)^{\gothm'/\gothh}$) attached to $\,\lambdatilde$ and $\lambda_r$ in Proposition \ref{construction de représentations, cas de $M'$} (b).

The representations of $\underline{M'\!}\,$ in $\calh$ and $(V \otimes_\CC F)^{\chiinfty{\mi \,\lambda}{\gothm'}}\!$ are trace class and have the same character.
\end{Proposition}
\vspace*{-3mm}

\begin{Dem}{Proof of the proposition}

{\bf (a)}
There exists --- and we fix --- an irreducible linear representation of $\gothm'_\CC$ in a finite-dimensional complex Hilbert space $F$, of lowest weight $-2\rho$ for the action of $\gothh_\CC$ relative to the order deduced from $R^+(\gothm'_\CC,\gothh_\CC)$.
It integrates into a representation of “the” universal covering of $M'_0$ with a trivial central character, and then into a representation of $M'_0$ denoted $\pi_{F,0}$.
Since $\underline{a}$ normalizes $\gothh$ and fixes $-2\rho$, there exists a unique intertwining operator $\Phi$ from $(F,\interieur\underline{a} \cdot \pi_{F,0})$ to $(F,\pi_{F,0})$ that acts trivially on the eigensubspace of $F$ of weight $-2\rho$ under the action of $\gothh_\CC$, common to both representations and of dimension~$1$.
Construct $\pi_F$ by extending $\pi_{F,0}$ and choosing $\pi_F(\underline{a})$ equal to $\Phi$.%
\smallskip

{\bf (b)}
First, note that $\underline{M'\!}\,$ stabilizes $\,\calh_0 \egdef (V \otimes_\CC F)^{\chiinfty{\mi \,\lambda}{\gothm'}}\!$.
By Lemma \ref{composantes primaires} and Proposition \ref{foncteur de translation} (c), the representations of $\underline{M'\!}\,$ in $\calh$ and $\calh_0$ are trace class.

Fix a maximal compact subgroup $K_{M'_0\!}$ of $M'_0$ stable under $\interieur{a}$ and whose Cartan involution normalizes $\gothh$.
For example, take $\exp (\gothc_{M'} \cap \gothm')$, where $\gothc_{M'}$ is the Lie algebra of a maximal compact subgroup of the semidirect product of $\{1,c\}$ ($c$:~conjugation of $\MM'(\CC)$) by $\MM'(\CC)$ that contains both $c$ and the subgroup generated by the product of the projection in $\MM'(\RR)$ of $a$ with $\exp_{\MM'(\CC)} (\mi \,\gothh_{(\RR)})$.
Denote by $\calh^f\!$, $\calh_0^f$ and $V^f$ the $(\gothm'_\CC,K_{M'_0})$-modules stable under $\underline{a}$ associated with $\calh$, $\calh_0$ and~$V\!$.
By \cite[Prop. 10.5 p.$\!\!$~336]{Kn86} with its proof (see \cite[Cor. 8.8 p.$\!\!$~211]{Kn86}), the representations of $\underline{M'\!}\,$ in $\calh$ and $\calh_0$ have the same character if there exists an isomorphism of $(\gothm'_\CC,K_{M'_0})$-modules from $\calh^f$ onto $\calh_0^f$ compatible with the action of~$\underline{a}$.
\smallskip

From the proof of Proposition \ref{construction de représentations, cas de $M'$} (a), we know that $\calh^f$ and $V^f$ are irreducible, isomorphic to $\,\calr_{M'_0}^{q'} (\CC_{\,\mi \,\lambda-\rho})\,$ and $\,\calr_{M'_0}^{q'} (\CC_{\,\mi \,\lambda_r-\rho})\,$ where $\,\calr_{M'_0}^{q'}$ is the cohomological induction functor relative to $\gothb_{M'}$, and $\CC_\Lambda$ is for each $\,\Lambda \in \mi \,\gothh^*\,$ the $(\gothh_\CC,T_0)$-module~$\CC$ (see \ref{conventions fréquentes} (a)) on which $\gothh$ acts by $\Lambda\id$.
From \cite[Th. 7.237 p.$\!\!$~544] {KV95} and using the notation $\,\psi_{l}^{l'}\,$ of \cite[(7.141) p.$\!\!$~493] {KV95}, we have\\
\hspace*{\fill}%
$\calr_{M'_0}^{q'} (\psi_{\,\mi \,\lambda_r-\rho}^{\,\mi \,\lambda-\rho}
(\CC_{\,\mi \,\lambda_r-\rho}))
\simeq
\psi_{\,\mi \,\lambda_r}^{\,\mi \,\lambda} (\calr_{M'_0}^{q'}
(\CC_{\,\mi \,\lambda_r-\rho}))$.
\hspace*{\fill}\\
Thus, the $(\gothm'_\CC,K_{M'_0})$-modules $\calh^f$ and $\calh_0^f$ are isomorphic.

Given the characterization of $\,S_{\lambdatilde}(\hat{a})\,$ in Proposition \ref{construction de représentations, cas de $M'$} (b), Schur's Lemma in \cite[Lem. 3.3.2 p.$\!\!$~80] {Wa88} and the identification $\,H_{q'} (\gothn_{M'},\calh^{\infty})^* = H_{q'} (\gothn_{M'},\calh^f)^*\,$ in \cite[Lem. 4 p.$\!\!$~165] {Df82a}, to obtain the compatibility with the action of $\,\underline{a}$ of a fixed $(\gothm'_\CC,K_{M'_0})$-module isomorphism from $\calh^f$ onto $\calh_0^f$, it suffices to show that $\,\underline{a}\,$ acts from $\,\pi_{\lambda_r}(\underline{a})\otimes \pi_F(\underline{a})\,$ in $\,(H_{q'} (\gothn_{M'},\calh_0^f)^*)_{-(\mi \,\lambda + \rho)}$ by $\rho_\lambdatilde^{\gothm'\!/\gothh} (\hat{a})^{-1} \id$, that is to say by $\rho_{\lambda_r}^{\gothm'\!/\gothh} (\hat{a})^{-1} \id$.

Bouaziz's idea basically consists in verifying the compatibility with the action of $\,\underline{a}$ of the isomorphism that would be obtained at the level of homology by taking $\gothq$ equal to the conjugate of $\gothb_{M'}$ and $X$ equal to $V^f$ in \cite[Th. 7.242 p.$\!\!$~546] {KV95}.
Unfortunately, the assumption “$\mi \lambda_r$ is at least as singular as $\mi \lambda$” does not hold.
Instead of going from $V^f$ to $\calh_0^f$, we will go from $\calh_0^f$ to $V^f$ using a property of the composite of two Zuckerman translation functors.
\smallskip

By \cite[Th. 7.220 p.$\!\!$~536] {KV95}, the homomorphism of $(\gothm'_\CC,K_{M'_0})$-modules compatible with the action of $\,\underline{a}$ from $(\calh_0^f \otimes_\CC F^*)^{\chiinfty{\mi \,\lambda_r}{\gothm'}}$ to $V^f\!$ obtained by restriction from the canonical map $\,V^f\! \otimes_\CC F \otimes_\CC F^* \to V^f$ is nonzero and therefore surjective.
The $\centrealg{U \gothh_\CC}$-primary components associated with $\chiinfty{\mi \,\lambda_r + \rho}{\gothh}\!$ in the exact sequence of homology $\gothh_\CC$-modules resulting from this surjection yield a homomorphism of $\gothh_\CC$--modules from
$\,H_{q'}(\gothn_{M'},(\calh_0^f \otimes_\CC F^*)^{\chiinfty{\mi \,\lambda_r}{\gothm'}})^{\chiinfty{\mi \,\lambda_r + \rho}{\gothh}}$
to $\,H_{q'} (\gothn_{M'},V^f)^{\chiinfty{\mi \,\lambda_r + \rho}{\gothh}}$
compatible with the action of $\,\underline{a}$, which is surjective according to \cite[Lem. 8.9 p.$\!\!$~553] {KV95}.

Let $v^*_0$ be a nonzero vector of $F^*$ of weight $2\rho$ for the action of $\gothh_\CC$.
We have $\;\gothn_{M'} \cdot v^*_0 = \{0\}$ and $\,\underline{a} \cdot v^*_0 = v^*_0\,$ according to (a).
It remains to prove that each of the two homomorphisms of $\gothh_\CC$-modules compatible with the action of $\,\underline{a}\,$ deduced from the canonical injections from
$(\calh_0^f \otimes_\CC F^*)^{\chiinfty{\mi \,\lambda_r}{\gothm'}}$
and $\calh_0^f \otimes_\CC\CC v^*_0$ into $\calh_0^f \otimes_\CC F^*\!$,
defined respectively on
$\,H_{q'} (\gothn_{M'},(\calh_0^f \otimes_\CC F^*)^{\chiinfty{\mi \,\lambda_r}{\gothm'}})^{\chiinfty{\mi \,\lambda_r + \rho}{\gothh}}$
and $\,H_{q'} (\gothn_{M'},\calh_0^f)^{\chiinfty{_{\mi \,\lambda + \rho}}{\gothh}} \otimes_\CC\CC v^*_0$,
and both with values in
$\,H_{q'} (\gothn_{M'},\calh_0^f \otimes_\CC F^*)^{\chiinfty{\mi \,\lambda_r + \rho}{\gothh}}$,
is bijective.\\
Indeed, in this case and with the notation $\CC_\Lambda$ introduced earlier in this proof, the line formed by homomorphisms of $\gothh_\CC$-modules from $\,H_{q'} (\gothn_{M'},V^f)\,$ to $\CC_{\mi \,\lambda_r + \rho}$ will be injected in a way compatible with the action of $\,\underline{a}$ into that formed by $\gothh_\CC$--module homomorphisms from $\,H_{q'} (\gothn_{M'},\calh_0^f)\,$ to $\,\CC_{\mi \,\lambda + \rho}$.
\smallskip

The first of the two preceding homomorphisms is bijective by \cite[(7.243) p.$\!\!$~547, see Prop. 7.166 p.$\!\!$~506] {KV95}.
For the second, set $\,W = F^*/\,\CC v^*_0$.

The exact sequence of $\!\gothb_{M'}$-modules
$0 \!\to\! \calh_0^f \!\otimes_\CC \CC v^*_0 \!\to\!\! \calh_0^f
\!\otimes_\CC F^* \!\!\to\!\! \calh_0^f \!\otimes_\CC W \!\to\!0$
 provides the following segment of the exact sequence of $\gothh_\CC$-modules in homology:\\
 %exact sequence of $\gothh_\CC$-homology moduli:\\
 \parbox{\textwidth}{
$H_{q'+1} (\gothn_{M'},\calh_0^f \otimes_\CC W)
\to H_{q'} (\gothn_{M'},\calh_0^f)\otimes_\CC\CC v^*_0\\
\hspace*{\fill}%
\to H_{q'} (\gothn_{M'},\calh_0^f \otimes_\CC F^*)
\to H_{q'} (\gothn_{M'},\calh_0^f \otimes_\CC W)$.
}
\\[0.5mm]
We complete the proof by showing that
$\;H_{\scriptscriptstyle\bullet}(\gothn_{M'},\calh_0^f \otimes_\CC W)^{\chiinfty{\mi \,\lambda_r + \rho}{\gothh}} = \{0\}$.

By Lie's theorem, there exists --- and we fix --- a sequence
$\;\{0\}
= W^{(-1)} \subseteq W^{(0)} \subseteq \cdots \subseteq W^{(N)}
= W\;$
of $\gothb_{M'}$-submodules of $W$ with $\dim W^{(j)}/\,W^{(j-1)} = 1$ for $0 \leq j \leq N$.
 By \cite[Prop. 3.12 p.$\!\!$~188] {KV95}, for any $(\gothb_{M'},T_0)$-module $E$ and any $n \in \NN$, the $\gothh_\CC$-module $H_n (\gothn_{M'},E)$ is isomorphic to $P_n(E)$, where $P_n$ is the $n^{\textrm{th}}\!$ derived functor from the right exact functor of the category of $(\gothb_{M'},T_0)$-modules to that of $(\gothh_\CC,T_0)$-modules denoted $\,P_{\gothb_{M'},T_0}^{\gothh_\CC,T_0}$ in
 \cite[(2.8) p.$\!\!$~104] {KV95}.
By \cite[Prop. D.57 (b) p.$\!\!$~887, see (2.123) p.$\!\!$~162] {KV95}, there exists --- and we fix --- a filtration
$\,\{0\}
\!=\! C_{\scriptscriptstyle\bullet}^{(-1)}
\subseteq C_{\scriptscriptstyle\bullet}^{(0)}
\subseteq \cdots
\subseteq C_{\scriptscriptstyle\bullet}^{(N)}
\!=\! C_{\scriptscriptstyle\bullet}\,$
of a chain complex $C_{\scriptscriptstyle\bullet}$ into $(\gothh_\CC,T_0)$-modules null in (strictly) negative degree, of which each $\gothh_\CC$-module of homology $H_n (C_{\scriptscriptstyle\bullet})$ is isomorphic to $H_n (\gothn_{M'},\calh_0^f\! \otimes_\CC W)$ for $n \in \NN$, and whose spectral sequence $(E^r)_{r \geq 0}$ provides $\gothh_\CC$-modules $E^1_{p,q}$ isomorphic to
$H_{p+q} (\gothn_{M'},\calh_0^f)\otimes_\CC(W^{(p)}/\,W^{(p-1)})\,$
for $0 \leq p \leq N$ and $p+q \geq 0$.
By \cite[Prop. 7.56 p.$\!\!$~460] {KV95}, the $\centrealg{U \gothh_\CC}$-primary components of each of these $\gothh_\CC$-modules $E^1_{p,q}$ are associated with characters of the form $\chiinfty{_{w\,\mi \,\lambda + \rho + \,\gamma}}{\gothh}$ with $w \in W(\gothm'_\CC,\gothh_\CC)$, where $\gamma$ denotes the weight of $\gothh_\CC$ in $\,W^{(p)}/\,W^{(p-1)}$.
Thus, given \cite[Prop. 7.166 p.$\!\!$~506] {KV95}, we have $\;(E^1_{p,q})^{\chiinfty{\mi \,\lambda_r + \rho}{\gothh}} = \{0\}$.
Furthermore, for $n \in \NN$ and $q \in \ZZ$, the $\gothh_\CC$-module $E^{\infty}_{n-q,q}$ is a subquotient of $E^1_{n-q,q}$ isomorphic to the quotient of $\,\pr_{H_n (C_{\scriptscriptstyle\bullet})} (H_n (C_{\scriptscriptstyle\bullet}^{(n-q)}))\,$ by $\,\pr_{H_n (C_{\scriptscriptstyle\bullet})} (H_n (C_{\scriptscriptstyle\bullet}^{(n-q-1)}))$, where we have used the same notation $\pr_{H_n(C_{\scriptscriptstyle\bullet})}$ for canonical maps with values in $H_n (C_{\scriptscriptstyle\bullet})$.
Hence for each $n \in \NN$, the following relations between $\gothh_\CC$-modules hold:\\
\parbox{\textwidth}{
$H_n(\gothn_{M'},\calh_0^f \otimes_\CC W)^{\chiinfty{\mi \,\lambda_r + \rho}{\gothh}}
\simeq
\pr_{H_n (C_{\scriptscriptstyle\bullet})}
(H_n (C_{\scriptscriptstyle\bullet}^{(N)}))^{\chiinfty{\mi \,\lambda_r + \rho}{\gothh}}\\[-1mm]
\hspace*{\fill}%
= \cdots
= \pr_{H_n (C_{\scriptscriptstyle\bullet})}
(H_n (C_{\scriptscriptstyle\bullet}^{(-1)}))^{\chiinfty{\mi \,\lambda_r + \rho}{\gothh}}
= \{0\}$.
}\\
This completes the proof.
\cqfd
\end{Dem}
\vspace*{-4mm}

%------------------------------------------------------------------------------
\begin{Remarque}\rm
%------------------------------------------------------------------------------

We retain the notation from the preceding proof.
In particular, $K_{M'_0\!}$ is stable under $\interieur{a}$.
Fix an intertwining operator $\Phi$ from $(\calh,\interieur{a} \cdot T_{\lambdatilde}^{M'_0})$ to $\smash{(\calh,T_{\lambdatilde}^{M'_0})}$ and an isomorphism $\Psi$ of $(\gothm'_\CC,K_{M'_0})$-modules from $\calh^f$ onto $\calh_0^f$.
Let $\Phi^f$ denote the restriction of $\Phi$ to $\calh^f$ and let $\pi_0(\underline{a})^f$ denote the operator by which $\underline{a}$ acts on~$\calh_0^f$.
According to Schur's Lemma, there exists $z\in\CC\moins0$ such that $\,\Phi^f = z\;\Psi^{-1}\!\circ\pi_0(\underline{a})^f\!\circ\Psi$.
The proof of Proposition \ref{passage au cas régulier} (b) shows that the representation $S$ of Proposition \ref{construction de représentations, cas de $M'$} (b) satisfies $S(\hat{a}) = z^{-1}\,\Phi$.
Furthermore, $K_{M'_0\!}$ is included in a maximal compact subgroup of $M'$ that contains $a$, and the maximal compact subgroups of $M'$ are conjugate under $M'_0$ (see \cite[Th. 3.1 p.$\!\!$~180 and lines before Th. 3.7 p.$\!\!$~186] {Ho65}).
It's easy to deduce that $z$ is independent of the choice of $K_{M'_0}$.
We could therefore have defined $S$ by returning to the regular semisimple case treated by Duflo (for values on elliptic elements) and making arrangements to ensure that
$\,S(\exp X) = \me ^{-\mi \lambda(X)} \, T_{\lambdatilde}^{M'_0}(\exp X)\,$
when $X\in\gothh$.
It would remain to use this new definition to prove that $S$ is a (unitary) representation of $M'(\lambdatilde)^{\gothm'/\gothh}$.
\cqfr
\end{Remarque}
\vspace*{-4mm}

\begin{Dem}{End of the proof of the theorem}

It remains to prove the theorem for the representation $\,T_{\lambdatilde,\sigma}^{M'}\!$.
The reference situation is that of \cite[Lem. 5.4 p.$\!\!$~304] {Zu77}, in which $M'$ is replaced by $M_0$ and $e$ is replaced by a regular semisimple element of $M_0$.

Assume $\,e \in M'$.
Let $\calv_e^{M'}\!$, $d_e^{M'}$ and $D_e^{M'}$ denote the objects analogous to $\calv_e$, $d_e$ and $D_e$ attached to $M'$.
Choose $\,\hat{a} \in M'(\lambdatilde)^{\gothm'/\gothh}$ lying over an elliptical $\,a \in M'(\lambdatilde)\,$ such that $\,e \in a\,M'_0$.
Associate it with $\underline{M'\!}\,$, $\pi_{\lambdatilde}$ and $\pi_{\lambda_r}\!$ as in Proposition \ref{passage au cas régulier}, so $\,\underline{a} \in \underline{M'\!}$.
We will use the notations of Proposition \ref{foncteur de translation} (a), with $M'$ and $\underline{M'\!}\,$ instead of $G$ and $\underline{G}$.
We introduce the element $\sigma_r$ of $\Xirr_{M'}(\lambda_r)$ such that:\\
\hspace*{\fill}%
$\sigma_r(\hat{u})
= \det(\Ad u^\CC)_{\gothn_{M'}}\; \sigma(\hat{u})\;\,$
for $\,\hat{u} \in M'(\lambdatilde)^{\gothm'/\gothh}$.
\hspace*{\fill}\\
\begin{tabular}{@{}rr@{$\,=\;$}l@{}}
It satisfies: &
$\restriction{\tr T_{\lambdatilde,\sigma}^{M'}}{e \, M'_0}$ &
$\tr \sigma (\hat{a})\,{\scriptstyle \times}\,
(\tr\pi_{\lambdatilde})_{a M'_0}$\\[-1mm]
and &
$\restriction{\tr T_{\lambda_r,\sigma_r}^{M'}}{e \, M'_0}$ &
$\tr \sigma_r (\hat{a})\,{\scriptstyle \times}\,
(\tr\pi_{\lambda_r})_{a M'_0}$.
\end{tabular}
\\[0.5mm]
To complete the proof, note that given Proposition \ref{passage au cas régulier}, Proposition \ref{foncteur de translation} will then allow us to link $\bigl( \tr T_{\lambdatilde,\sigma}^{M'} \bigr)_{\!e}$ to $\bigl( \tr T_{\lambda_r,\sigma_r}^{M'} \bigr)_{\!e}$.
%\smallskip

We verify the character formula for $\,\tr T_{\lambdatilde,\sigma}^{M'}$ (in the version proposed by formula $(\boldsymbol{F\!}_{+,0})$ page \pageref{F_{+,0}}) on the open set formed by $\;X \in \calv_e^{M'}\;$ such that $\;e \exp X \in \Mpssreg$, included in $\,\mpessreg$ and of negligible complement, on which the two members of the equality are analytic functions.

Let $\,X \in \calv_e^{M'}\!$ such that $\,x \egdef e \exp X \in \Mpssreg$.
Set $\,\gothj_x \egdef \gothm'(x)$.
Then $\gothj_x \!=\! \gothm'(e)(X)\in \Car \gothm'(e)$ (given the dimensions) and $\gothj \egdef \centraalg{\gothm'}{\gothj_x} \in \Car \gothm'$.
Let~$\Gamma$ denote the connected component of $X$ in $\gothj(x) \cap \mpessreg$.
To continue the cal\-culation, consider $\,Y \!\in \gothj(x)_1$ (see \ref{foncteur de translation} (b) for $G\!=\!M'$) such that $X\!+\!Y \!\in \Gamma \cap \calv_e^{M'}\!$.
\smallskip

According to \cite[Th. 5.5.3 p.$\!\!$~52]{Bo87} and by \ref{morphisme rho et fonction delta} (b), the character formula is acquired for $\,\tr T_{\lambda_r,\sigma_r}^{M'}\!$.
Therefore (see \ref{notations pour les restrictions} (a) (e) and $\!(\boldsymbol{F\!}_{+,0})\!$ page \pageref{F_{+,0}}), we have:

\noindent\parbox{\textwidth}{
$\quad\; \tr \sigma_r(\hat{a})
\;{\scriptstyle \times}\;
\abs{D_{M'}(x\exp Y)}^{1/2}\;\,
(\tr\pi_{\lambda_r})_{a M'_0} (x\exp Y)$\\
$\displaystyle
=\; |\det \;(1-\Ad e)_{_{\scriptstyle \gothm'/\gothm'(e)}}|^{1/2}\;\,
|\cald_{\gothm'(e)}(X\!+\!Y)|^{1/2}\;\,
\bigl( \tr T_{\lambda_r,\sigma_r}^{M'} \bigr)_{\!e}(X\!+\!Y)$
}\\
\parbox{\textwidth}{
$\stackrel
{\raisebox{.5ex}{\boldmath$
{\scriptscriptstyle (}\scriptstyle\star{\scriptscriptstyle )}$}}
{=}\;
\mi^{-d_e^{M'}}\;(D_e^{M'})^{-\frac{1}{2}}\;\,
|\det \;(1-\Ad e)_{_{\scriptstyle \gothm'/\gothm'(e)}}|^{1/2}
\Sum_{\classe{\lambda'_r} \, \in \, M'(e)_0 \backslash
M'_0 \cdot \lambda_r \cap\, \gothm'(e)^*}\\[-1mm]
\hspace*{\fill}%
{\textstyle \frac
{\scalo (\widehat{e'})_{(1-\Ad e')(\gothm'/\gothh)}}
{\scalo (B_{\lambda_r})_{(1-\Ad e')(\gothm'/\gothh)}}}\;
\tr\sigma_r(\widehat{e'})\;\,
|\cald_{\gothm'(e)}(X\!+\!Y)|^{1/2}\;\,
\widehat{\beta}_{M'(e)_0 \cdot \lambda'_r}(X\!+\!Y)$
}
\\
where $\;\,g \in M'_0\,$ satisfies $\;\lambda'_r = g\lambda_r\;$ and $\;\widehat{e'} \in M'(\lambdatilde)^{\gothm'/\gothh}\;$ lies over $\,e' = g^{-1}eg$.
\smallskip

Consider $\,g \in M'_0\,$ as above.
Denote by $\pr$ the canonical map from $M'$ to~$\MM'(\CC)$.
We associate with $g$ an element $y^g$ of $\MM'(\CC)(\pr(e))_0$ that satisfies $\gothj(x)_\CC= y^g \cdot (g \gothh)(e)_\CC$.
Fix a system of positive roots $R^+(\gothm'(e)_\CC,\gothj(x)_\CC)$ of $R(\gothm'(e)_\CC,\gothj(x)_\CC)$.
According to Proposition \ref{limites de mesures de Liouville} (b), there exist complex numbers $c^g_w$ independent of $Y$ and indexed by the $\;w \in W(\gothm'(e)_\CC,\gothj(x)_\CC)\,$ such that we have\\[-0.5mm]
\hspace*{\fill}%
$\Prod_{\alpha \in R^+(\gothm'(e)_\CC,\gothj(x)_\CC)}\!\! \alpha(X\!+\!Y)
\;{\scriptstyle \times}\;
\widehat{\beta}_{M'(e)_0 \cdot g\lambda_r}(X\!+\!Y)\,
\stackrel {\raisebox{.5ex}{\boldmath$
{\scriptscriptstyle (}\scriptstyle\star\star{\scriptscriptstyle )}$}}
{=}
\Sum_{w \in W(\gothm'(e)_\CC,\gothj(x)_\CC)} \!\! c^g_w \; \me ^{\mi \,wy^gg\lambda_r(X+Y)}$
\hspace*{\fill}%
\\
and
\hspace*{\fill}%
$\Prod_{\alpha \in R^+(\gothm'(e)_\CC,\gothj(x)_\CC)}\!\!\alpha(X)
\;{\scriptstyle \times}\;
\widehat{\beta}_{M'(e)_0 \cdot (g\lambdatilde)[e]}(X)\,
= \Sum_{w \in W(\gothm'(e)_\CC,\gothj(x)_\CC)} \!\! c^g_w \;\; \me ^{\mi \,wy^gg\lambda\,(X)}
\hspace*{\fill}$
\\
when $(g\lambdatilde)[e] \in \mpestregtilde$, \ and
$\Sum_{w \in W(\gothm'(e)_\CC,\gothj(x)_\CC)} \!\!\! c^g_w \; \me ^{\mi \,wy^gg\lambda\,(X)} = 0\;$
if $(g\lambdatilde)[e] \notin \mpestregtilde$.%

Furthermore, the function
$\,|\cald_{\gothm'(e)}|^{1/2} \,{\scriptstyle \times}\,
\bigl( \Prodpetit_{\alpha \in R^+(\gothm'(e)_\CC,\gothj(x)_\CC)} \! \alpha \bigr) ^{-1}$
on $\Gamma$ has values in the set of fourth roots of unity, and is therefore constant.
From the equalities $(\star)$ and $(\star\,\star)$, we derive a formula for
$\;\tr \sigma_r(\hat{a})\,{\scriptstyle \times}\, (\tr\pi_{\lambda_r})_{a M'_0} (x\exp Y)\;$
that involves this constant.

According to Proposition \ref{caractères canoniques} (a), the representation $\,\pi_{\lambda_r}\!$ of $\underline{M'\!}\,$ has the infinitesimal character $\,\chiinfty{\mi \, \lambda_r}{\gothm'}$.
Apply Proposition \ref{foncteur de translation} (c) to $\underline{M'\!}\,$ with $\pi = \pi_{\lambda_r}$, $\Lambda = \mi\,\lambda_r$ and $\Lambda_F = -2\rho$ ($\pi_F$ below).
We have a representation $(\pi_{\lambda_r})_{Zuc}$ of $\underline{M'\!}\,$ and, given the above, a formula for $\tr \sigma_r(\hat{a})\,{\scriptstyle \times}\, (\tr(\pi_{\lambda_r})_{Zuc})_{a M'_0} (x\exp Y)$.
By Proposition \ref{passage au cas régulier} (b), we also have $\,\tr(\pi_{\lambda_r})_{Zuc} = \tr\pi_{\lambdatilde}\,$ by choosing $\pi_F$.

We now choose $Y = 0$, and summarize.
We see that $\bigl( \tr T_{\lambdatilde,\sigma}^{M'} \bigr)_{\!e}(X)$ can be deduced from the expression for $\bigl( \tr T_{\lambda_r,\sigma_r}^{M'} \bigr)_{\!e}(X)$ obtained above, by replacing the coefficient
$\,c^g_w \, \me ^{\mi \,wy^gg\,\lambda_r\!(X)}\,$ by
$\;(\det(\Ad a^\CC)_{\gothn_{M'}})^{-1}
\,{\scriptstyle \times}\; \underline{x}_{~\!l_F}\; c^g_w \, \me ^{\mi \,wy^gg\,\lambda_r\!(X)}\;$
for each $w \in W(\gothm'(e)_\CC,\gothj(x)_\CC)$, where $\;l_F \egdef -2\,wy^gg\rho\;$
($\underline{x}_{~\!l_F}$ is defined in \ref{foncteur de translation} (c)).
%\smallskip

Let $\,w \in W(\gothm'(e)_\CC,\gothj(x)_\CC)$.
We are going to calculate the term $\underline{x}_{~\!l_F}$ associated with it.
Choose a representative $\wt{w}$ of $w$ in $\,\MM'(\CC)(\pr(e))_0$.
Let ${\underline{M'\!}\,}_\CC$ denote the semidirect product of $\ZZ$ by $\MM'(\CC)_0$ such that $1 \in \ZZ$ acts on $\MM'(\CC)_0$ by $\interieur\pr(a)$, and let $\,\underline{\pr}\,$ denote the canonical map from $\underline{M'\!}\,$ to ${\underline{M'\!}\,}_\CC$.
There is a unique holomorphic representation $\pi_{F,\CC}$ of ${\underline{M'\!}\,}_\CC$ in $F$ such that $\;\pi_F = \pi_{F,\CC} \circ \underline{\pr}$.
Fix a nonzero vector $v_0$ of $F$ of weight $\,-2\rho\,$ for the action of $\gothh_\CC$ on $F$.
We are interested in the nonzero vector $\,\wt{w}y^gg \cdot v_0\,$ of $F$ which has weight $\,l_F$ under the action of $\gothj(x)_\CC$.
Again, set $\,e' = g^{-1}eg$.
In the semidirect product of $\ZZ$ by $\MM'(\CC)$ by means of $\interieur\pr(a)$, we have
$\underline{\pr}(\underline{x}) = \underline{\pr}(\underline{a}\,(a^{-1}x))$ (see \ref{foncteur de translation} (c)).
We obtain:\\
$\underline{\pr}(\underline{x})\, \wt{w} y^g \pr(g)
\begin{array}[t]{cl}
\!=\!&
\pr(x)\, \wt{w} y^g \pr(g)\pr(a)^{-1} \, \underline{\pr}(\underline{a})\\
\!=\!&
\pr(\exp X)\, \wt{w} y^g \pr(g)\pr(e'a^{-1})\,
\underline{\pr}(\underline{a})
\end{array}$\\
then
$\;\,\underline{x} \cdot (\wt{w}y^gg \cdot v_0)
= \exp X \cdot (\wt{w}y^gg \cdot ((e'a^{-1})\cdot v_0))$.\\
The property $\;\,e'a^{-1} \!\in M'_0(\lambdatilde) = \exp {\gothh}\;$ allows us to deduce that:\\
\hspace*{\fill}%
$\underline{x}_{~\!l_F}
= \det(\Ad a^\CC)_{\gothn_{M'}}\; (\det({\Ad e'}^\CC)_{\gothn_{M'}})^{-1}
\; \me ^{-2\,wy^gg\,\rho(X)}$.
\hspace*{\fill}%
\smallskip

This leads to the conclusion.
\cqfd
\end{Dem}

%===============================================================================
%
\begin{index des notations}
%
%==============================================================================

\newcounter{foo}
\setcounter{foo}{0}
\newcommand{\myitem}{\refstepcounter{foo}[\thefoo]}
\hlabel{RepIndex}%

\item[] ${\scriptstyle \langle} ~,\!~ {\scriptstyle \rangle}$, \pageref{<,>ev}, \pageref{<>dual}

\indexspace

\item[] $\gotha$, \pageref{a}
\item[] $A$, \pageref{A}

\indexspace

\item[] $B_f$, \pageref{Bf}
\item[] $\gothb_{M'}$, \pageref{bM'}

\indexspace

\item[] $\Car \gothg$, \pageref{Carg}
\item[] $c_{\widehat{e'},\lambdatilde}$ cf. $(F)$ p.$\!$~\pageref{F}
\item[] $C(\gothg(\lambda ),\gothh)$, \pageref{Cgh}
\item[] $C(\gothg(\lambda ),\gothh)_{reg}$, \pageref{Cghreg}

\indexspace

\item[] $d_e$, \pageref{de}
\item[] $D_e$, \pageref{De}
\item[] $\cald_{\gothg}$, \pageref{Dg}
\item[] $D_G$, \pageref{DG}
\item[] $DL(V)$, \pageref{DL(V)}

\indexspace

\item[] $\F+[e]$, \pageref{F+[e]}
\item[] $\Fh+$, \pageref{Fh+}

\indexspace

\item[] $\gothg$, \pageref{g}
\item[] $G$, \pageref{G}
\item[] $\GG$, \pageref{GG}
\item[] $\widehat{G}$, \pageref{G^}
\item[] $\gothg(f)$, \pageref{g(f)}
\item[] $G(f)$, \pageref{G(f)}
\item[] $G(f)^{\gothg/\gothg(f)}$, \pageref{G(f)g/g(f)}
\item[] $\gothg(x)$, \pageref{g(x)}
\item[] $G(x)$, \pageref{G(x)}
\item[] $\gothg(X)$, \pageref{g(X)}
\item[] $G(X)$, \pageref{G(X)}
\item[] $\gothg(\lambdatilde)$, \pageref{g(λtilde)}
\item[] $G(\lambdatilde)$, \pageref{G(λtilde)}
\item[] $G(\lambda_+)^{\gothg / \gothh}$, \pageref{G(f)g/g(f)}
\item[] $G(\lambda_+)^{\gothg/\gothh}_0$, \pageref{G(f)g/g(f)0}
\item[] $G(\lambdatilde)^{\gothg/\!\gothg(\lambda)(\mi \rho_\F+)}\!\!$, \pageref{G(λtilde)g/glambdairho}
\item[] $G(\lambdatilde)^{\gothg/\gothh}$, \pageref{G(λtilde)g/h}
\item[] $G(\lambdatilde)^{\gothg/\gothh}_0$, \pageref{G(λtilde)g/h0}
\item[] $\gothg^*(x)$, \pageref{g*(x)}
\item[] $\gothg^*(X)$, \pageref{g*(X)}
\item[] $\gstfondtilde$, \pageref{g*fond tilde}
\item[] $\gstfondGtilde$, \pageref{g*fond,G tilde}
\item[] $\gstfondtilde(e)$, \pageref{g*fond tilde(e)}
\item[] $\gstItilde$, \pageref{g*I tilde}
\item[] $\gstItilde(e)$, \pageref{g*I tilde(e)}
\item[] $\gstIGtilde$, \pageref{g*I,G tilde}
\item[] $\gstInctilde$, \pageref{g*Inc tilde}
\item[] $\gstInctilde(e)$, \pageref{g*Inc tilde(e)}
\item[] $\gstIncGtilde$, \pageref{g*Inc,G tilde}
\item[] $\gstreg$, \pageref{g*reg}
\item[] $\gstregtilde$, \pageref{g*reg tilde}
\item[] $\gstreg(e)$, \pageref{g*reg(e)}
\item[] $\gstregtilde(e)$, \pageref{g*reg tilde(e)}
\item[] $\gstregGtilde$, \pageref{g*reg,G tilde}
\item[] $\gstss$, \pageref{g*ss}
\item[] $\gstssG$, \pageref{g*ss,G}
\item[] $\gstssI$, \pageref{g*ssI}
\item[] $\gstssIG$, \pageref{g*ssI,G}
\item[] $\gstssInc$, \pageref{g*ssInc}
\item[] $\gstssIncG$, \pageref{g*ssInc,G}
\item[] $\gssreg$, \pageref{g ssreg}
\item[] $\gstssreg$, \pageref{g*ssreg}
\item[] $\Gssreg$, \pageref{G ssreg}

\indexspace

\item[] $\gothh$, \pageref{lemme clef}
\item[] $\calh$, \pageref{H}
\item[] $\gothh_{(\RR)}$, \pageref{hr}
\item[] $H_{\alpha}$, \pageref{Hα}

\indexspace

\item[] $\interieur{\gothg}$, \pageref{int(a)}
\item[] $\interieur{g}$, \pageref{int(g)}
\item[] $\interieur{G}$, \pageref{int(G)}

\indexspace

\item[] $\gothj(x)_1$, \pageref{j (x)1}

\indexspace

\item[] $k_e$, \pageref{ke}

\indexspace

\item[] $\gothm$, \pageref{m}
\item[] $M$, \pageref{M}
\item[] $\gothm'$, \pageref{m'}
\item[] $M'$, \pageref{M'}
\item[] $Mp(V)$, \pageref{Mp(V)}
\item[] $Mp(V)_\call$, \pageref{Mp(V)L}

\indexspace

\item[] $n_\call$, \pageref{nL}
\item[] $\gothn_{M'}$, \pageref{nM'}

\indexspace

\item[] $\scalo (A)_{A \cdot V}$, \pageref{OA}
\item[] $\scalo (\hat{a})_{(1-a)\cdot V}$, \pageref{Oa}
\item[] $\scalo (B)_V$, \pageref{OB}

\indexspace

\item[] $q$, \pageref{q}
\item[] $q'$, \pageref{q'}
\item[] $q_\call$, \pageref{qL}

\indexspace

\item[] $R(\gothg_\CC,\gothh_\CC)$, \pageref{R(gC,hC)}
\item[] $R^+\!(\gothg_\CC,\!\gothh_\CC)$\!, \pageref{R+(gC,hC)}\!, \pageref{R+(gC,hC)noncan}
\item[] $R_G$, \pageref{RG}
\item[] $R^+_\lambdatilde$, \pageref{R+lambdatilde}
\item[] $R^+_{\lambdatilde,{\gotha^*}^+}$, \pageref{R+(gC,hC)can}

\indexspace

\item[] $S$, \pageref{S}
\item[] $\sg$, \pageref{sg}
\item[] $\Supp_{\gothg^*}\wt{\Omega}$, \pageref{SuppgstOmegatilde}

\indexspace

\item[] $\gotht$, \pageref{t}
\item[] $T_0$, \pageref{T0}
\item[] $T_{\lambdatilde,{\gotha^*}^+\!,\tau_+}^G\!$, \pageref{Tλtilde,a*,τ,connexe}, \pageref{Tλtilde,a*,τ}
\item[] $T_{\lambdatilde,\tau}^G$, \pageref{Tλtilde,τ}
\item[] $T_{\lambda_r,\sigma_r}^{M'}$, \pageref{Tλzr,sigmar,M'}
\item[] $T_{\lambda_r}^{M'_0}$, \pageref{Tλzr,M'0}
\item[] $T_{\lambdatilde}^{M'_0}$, \pageref{Tλztilde,M'0}

\indexspace

\item[] $\gothu$, \pageref{u}
\item[] $U$, \pageref{U}

\indexspace

\item[] $\calv_e$, \pageref{Ve}
\item[] $V^\chi$, \pageref{Vchi}

\indexspace

\item[] $W(G,\gothh)$, \pageref{W(G,h)}

\indexspace

\item[] $\Xfin_G\!$, \pageref{Xf}
\item[] $\Xfin_G\!(\lambdatilde,{\gotha^*}^+\!)\!$, \pageref{Xfinal}
\item[] $\XInd_G$, \pageref{XInd}
\item[] $\XInd_G(\lambdatilde)$, \pageref{XI}
\item[] $\Xirr_G(f)$, \pageref{Xirr}
\item[] $\Xirr_G(\lambdatilde)$, \pageref{Xirr tilde}
\item[] $\Xirrp_G(\lambdatilde,{\gotha^*}^+)$, \pageref{Xirr+}

\indexspace

\item[] $\beta_{\Omega}$, \pageref{β_Ω}
\item[] $\beta_{\wt{\Omega}}$, \pageref{Omegatilde}

\indexspace

\item[] $\delta$, \pageref{δ}
\item[] $\delta_{\lambda_+}^{\gothg/\gothh}$, \pageref{δ_λ+}

\indexspace

\item[] $\Theta_e$, \pageref{Θ_e}

\indexspace

\item[] $\iota$, \pageref{ι}

\indexspace

\item[] $\lambda$, \pageref{λ}
\item[] $\lambda_+=\lambda_{\gothg,\lambdatilde,{\gotha^*}^+\!,\epsilon}$, \pageref{λ+}
\item[] $\lambda_{+,\gothm'}$, \pageref{λ+m'}
\item[] $\lambdatilde$, \pageref{λtilde}
\item[] $\lambdatilde[e]$, \pageref{λtilde[e]}
\item[] $\lambda_{\textit{can}}$, \pageref{λcan}

\indexspace

\item[] $\mu$, \pageref{μ}
\item[] $\mu_+=\mu_{\gothg,\lambdatilde,{\gotha^*}^+}$, \pageref{μ+}
\item[] $\mu_{+,\gothm}$, \pageref{μ+m}
\item[] $\mutilde$, \pageref{μtilde}

\indexspace

\item[] $\nu$, \pageref{ν}
\item[] $\nu_+=\nu_{\gothg,\lambdatilde,{\gotha^*}^+\!,\epsilon}$, \pageref{ν+}

\indexspace

\item[] $\xi$, \pageref{ξ}

\indexspace

\item[] $\displaystyle \left |\Pi_{\gothg,\gothm}\right |$, \pageref{Π g,m}

\indexspace

\item[] $\rho_\F+$, \pageref{ρ F+}
\item[] $\rho_{\gothg,\gothh}$, \pageref{ρ g,h}
\item[] $\rho_\call$, \pageref{ρ L}
\item[] $\rho_{\lambda_+}^{\gothg/\gothh}$, \pageref{ρ λ+}
\item[] $\rho_\lambdatilde^{\gothg/\!\gothg(\lambda)(\mi \rho_\F+)}\!$, \pageref{ρ λtilde g,h}
\item[] $\rho_{\lambdatilde,{\gotha^*}^+}^{\gothg/\gothh}$, \pageref{ρ λtilde,a*+ g,h}

\indexspace

\item[] $\tau$, \pageref{τ}
\item[] $\tau_+$, \pageref{τ+}
\item[] $\tau_M$, \pageref{τ_M}
\item[] $\tau_{M'}$, \pageref{τ_M'}

\indexspace

\item[] $\Phi$, \pageref{Φ}

\indexspace

\item[] $\chiinfty{l}{\gothg}$, \pageref{χ lUgC}
\item[] $\chi_{\lambdatilde}^G$, \pageref{χ λ+G}

\end{index des notations}
%\vfill

%===============================================================================
%
\begin{thebibl}{AAA 00}{\slshape References}{0}
%
%==============================================================================

%\addcontentsline{toc}{part}{References}

\hlabel{RepRef}

\bibitem[ABV 92]{ABV92}
Adams,~J., D.~Barbasch, and D.~A.~Vogan,
\!“\href{https://doi.org/10.1007/978-1-4612-0383-4}{The Langlands Classification and Irreducible Characters for Real Reductive Groups},”
{Birkhäuser, Boston, 1992}.%
\vskip0.2mm

\bibitem[B.. 72]{B.72}
Bernat,~P., N.~Conze, M.~Duflo, M.~Lévy-Nahas, M.~Raïs, P.~Renouard, and M.~Vergne,
“Représentations des groupes de Lie résolubles,”
{Dunod, Paris, 1972}.%
\vskip0.2mm

\bibitem[BW 80]{BW80}
Borel,~A., and N.~Wallach,
“Continuous Cohomology, Discrete Subgroups, and Representations of Reductive Groups,”
{Annals of Math. Studies 94, Princeton University Press, 1980}.
Second edition:
\href{https://doi.org/10.1090/surv\%2F067}{Math. Surveys 67, Amer.\ Math.\ Soc.,\ 2000}.%
\vskip0.2mm

\bibitem[Bou 84]{Bo84}
Bouaziz,~A.,
\href{https://gdz.sub.uni-goettingen.de/download/pdf/PPN235181684_0268/PPN235181684_0268.pdf#page=544}{\it Sur les représentations des groupes de Lie réductifs non connexes},
{Math. Ann. {\bf268} (1984), 539--555}.%
\vskip0.2mm

\bibitem[Bou 87]{Bo87}
---,
\href{https://doi.org/10.1016/0022-1236%2887%2990122-4}{\it Sur les caractères des groupes de Lie réductifs non connexes},
{J. Funct. Anal. {\bf70} (1987), 1--79}.%
\vskip0.2mm

\bibitem[B~67]{B67}
Bourbaki,~N.,
“Théories spectrales, Chapitres 1 et 2,”
Hermann, Paris, 1967.
Second edition:
\href{https://doi.org/10.1007/978-3-030-14064-9}{Springer, Cham, 2019}.
\vskip0.2mm

\bibitem[Cha 96]{Ch96}
Charbonnel, J.-Y.,
\href{https://doi.org/10.1006/jfan.1996.0062}{{\it Orbites fermées et orbites tempérées, II}},
{J. Funct. Anal. {\bf138} (1996), 213--222}.%
\vskip0.2mm

\bibitem[Cow 88]{Co88}
Cowling,~M.,
\href{https://doi.org/10.1017/S1446788700032262}{\it On the characters of unitary representations},
{J. Austral. Math. Soc. Ser. A {\bf45} (1988), 62--65}.%
\vskip0.2mm

\bibitem[Dix 64]{Di64}
Dixmier,~J.,
“Les $C\sp{\ast}$-algèbres et leurs représentations,”
{Gauthier-Villars, Paris, 1964}.
English translation:
\href{https://archive.org/details/calgebras0000dixm/page/n9/mode/2up}{{$C\sp*$}-algebras,”
North-Holland Math. Library, Vol. 15}.
\vskip0.2mm

\bibitem[Dix 69]{Di69}
---,
\href{https://www.numdam.org/item/10.24033/asens.1180.pdf}{\it Sur la représentation régulière
d'un groupe localement compact connexe},
{Ann. Sci. École Norm. Sup.\ {\bf2} (1969), 423--436}.%
\vskip0.2mm

\bibitem[Duf 82a]{Df82a}
Duflo,~M.,
\href{https://doi.org/10.1007/978-3-642-11117-4_3}{\it Construction de représentations unitaires d'un groupe de Lie},
{In: “Harmonic analysis and group representations,”
Cours d'été du C.I.M.E., Cortona 1980,
Liguori, Naples, 1982, 129--221}%
\vskip0.2mm

\bibitem[Duf 82b]{Df82b}
---,
\href{https://doi.org/10.1007/BF02392353}{\it Théorie de Mackey pour les groupes de Lie algébriques},
{Acta Math. {\bf149} (1982), 153--213}.%
\vskip0.2mm

\bibitem[Duf 84]{Df84}
---,
\href{https://doi.org/10.1007/BFb0072338}{\it On the Plancherel formula for almost algebraic real Lie groups},
{In: “Lie group representations III,” Lect. Notes in Math. {\bf 1077},
Springer-Verlag, Berlin etc., 1984, 101--165}.%
\vskip0.2mm

\bibitem[DHV 84]{DHV84}
Duflo,~M., G.~Heckman, and M.~Vergne,
\href{https://numdam.org/item/10.24033/msmf.300.pdf}{\it Projection d'orbites, formule de Kirillov et formule de Blattner},
{Mém. Soc. Math. Fr. (N.S.){\bf15} (1984), 65--128}.%
\vskip0.2mm

\bibitem[DV 88]{DV88}
Duflo,~M., and M.~Vergne,
\href{https://projecteuclid.org/ebooks/advanced-studies-in-pure-mathematics/Representations-of-Lie-Groups-Kyoto-Hiroshima-1986/chapter/La-Formule-de-Plancherel-des-Groupes-de-Lie-Semi-Simples/10.2969/aspm/01410289.pdf}{\it La formule de Plancherel des groupes de Lie semi-simples réels},
{In: “Representations of Lie groups,” numéro~14 in Adv. Stud. Pure Math.,
Academic Press, 1988, 289--336}.%
\vskip0.2mm

\bibitem[DV 93]{DV93}
---,
\href{http://www.numdam.org/article/AST_1993__215__5_0.pdf}{\it Cohomologie équivariante et descente},
{Astérisque {\bf215} (1993), 5--108}.%
\vskip0.2mm

\bibitem[Har 65]{Ha65}
Harish-Chandra,
\href{http://dx.doi.org/10.1007/BF02391779}{\it Discrete series for semisimple Lie groups I. Construction of invariant eigendistributions},
{Acta Math. {\bf113} (1965), 241--318}.%
\vskip0.2mm

\bibitem[Har 76]{Ha76}
---,
\href{https://doi.org/10.2307/1971058}
{\it Harmonic analysis on real reductive groups III.
The Maass-Selberg relations
and the Plancherel formula},
{Ann. of Math. {\bf104} (1976), 117--201}.%
\vskip0.2mm

\bibitem[Her 83]{Hr83}
Herb,\!~R.\!~A.,
\!\href{https://www.ams.org/journals/tran/1983-277-01/S0002-9947-1983-0690050-7/S0002-9947-1983-0690050-7.pdf}
{\it Discrete series characters and Fourier inversion on semisimple real Lie groups},
{Trans. Amer. Math. Soc.\ {\bf277} (1983), 241--262}.%
\vskip0.2mm

\bibitem[Hoc 65]{Ho65}
Hochschild,~G.,
“The Structure of Lie Groups,”
{Holden-Day, San Francisco 1965}.%
\vskip0.2mm

\bibitem[Kna 86]{Kn86}
Knapp,~A.~W.,
“\href{https://www.jstor.org/stable/j.ctt1bpm9sn}{Representation Theory of Semisimple Groups},”
Princeton University Press, 1986.%
%https://www.projecteuclid.org/search?term=knapp+lie+groups
\vskip0.2mm

\bibitem[Kna 96]{Kn96}
---,
\!“\href{https://link.springer.com/book/10.1007/978-1-4757-2453-0}{Lie Groups beyond an Introduction},”
{Birkhäuser, Boston, 1996}.
\!\href{https://projecteuclid.org/ebooks/books-by-independent-authors/Lie-Groups-Beyond-an-Introduction/toc/10.3792/euclid/9798989504206}{Digital Second Edition: 2023}.%
%\!\href{https://www.math.stonybrook.edu/~aknapp/download/Beyond2-clickable.pdf}{Digital Second Edition: 2023}.%
\vskip0.2mm

\bibitem[KV 95]{KV95}
Knapp,~A.~W., and D.~A.~Vogan,
“\href{https://www.jstor.org/stable/j.ctt1bpm9qm}{Cohomological Induction and Unitary Representations},”
{Princeton University Press, 1995}.%
\vskip0.2mm

\bibitem[KZ 82]{KZ82}
Knapp,~A.~W., and G.~J.~Zuckerman,
\href{https://www.jstor.org/stable/2007066}{\it Classification of irreductible tempered representations of semisimple
groups},
Ann. of Math.\ {\bf116} (1982), 389--501.
Typesetter's correction: \href{https://www.jstor.org/stable/2007089}{Ann. of Math.\ {\bf119} (1984), 639}.%
\vskip0.2mm

\bibitem[Kos 59]{Ko59}
Kostant,~B.,
\href{https://www.jstor.org/stable/2372999}{\it The principal three-dimensional subgroup
and the Betti numbers of a complex
simple Lie group},
{Amer. J. Math., {\bf81} (1959), 973--1032}.%
\vskip0.2mm

\bibitem[Mac 58]{Ma58}
Mackey,~G.~W.,
\href{https://doi.org/10.1007/BF02392428}{\it Unitary representations of group extension I},
{Acta Math. {\bf99} (1958), 265--311}.%
\vskip0.2mm

\bibitem[Ran 72]{Ra72}
Ranga Rao,~R.,
\href{https://www.jstor.org/stable/1970822}{\it Orbital integrals in reductive groups},
{Ann. of Math., {\bf96} (1972), 505--510}.%
\vskip0.2mm

\bibitem[Ros 80]{Ro80}
Rossmann,~W.,
\href{https://link.springer.com/article/10.1007/BF01389894}{\it Limit characters of reductive Lie groups},
{Invent. Math.\ {\bf61} (1980), 53--66}.%
\vskip0.2mm

\bibitem[Ros 82]{Ro82}
---,
 \href{https://doi.org/10.1215/S0012-7094-82-04914-6}{\it Limit orbit in reductive Lie algebras},
{Duke Math. J.\ {\bf49} (1982), 215--229}.%
\vskip0.2mm

\bibitem[Spr 66]{Sp66}
Springer,~T.~A.,
\href{https://www.numdam.org/article/PMIHES_1966__30__115_0.pdf}{\it Some arithmetical results on semi-simple Lie algebras},
{Inst. Hautes Études Sci.\ {\bf30} (1966), 115--141}.%
\vskip0.2mm

\bibitem[Var 77]{Va77}
Varadarajan,~V.~S.,
\href{https://link.springer.com/book/10.1007/BFb0097814}{Harmonic analysis on real reductive groups},”
{Lect. Notes in Math.\ {\bf576}, Springer-Verlag, Berlin etc.,
1977}.%
\vskip0.2mm

\bibitem[Ver 94]{Ve94}
Vergne,~M.,
\href{https://doi.org/10.1007/978-3-0348-9110-3_8}{\it Geometric quantization and equivariant cohomology},
{In: “First European Congress of Mathematics I, Paris 1992,”
Progr. Math.\ {\bf119}, Birkhäuser, Boston, 1994, 249--295}.%
\vskip0.2mm

\bibitem[Wal 88]{Wa88}
Wallach,~N.,
“Real Reductive Groups I,”
{Academic Press, New York etc., 1988}.%
\vskip0.2mm

\bibitem[Zuc 77]{Zu77}
Zuckerman,~G.,
\href{https://www.jstor.org/stable/1971097}{\it Tensor products of finite and infinite dimensional representations of semisimple Lie groups},
{Ann. of Math.\ {\bf106} (1977), 295--308}.%

\end{thebibl}

%===============================================================================
%
%Last page
%
%==============================================================================

\font\eightrm=ecrm0800 scaled\magstep1

\par\vspace*{1cm}\noindent
\begin{minipage}[t]{.4\hsize}
\parindent0pt\noindent
\fontsize{9}{8}\selectfont
\href{http://www.imj-prg.fr/~jean-yves.ducloux/}{Jean-Yves Ducloux}\\
Université Paris 7\\
UFR de Mathématiques, UMR 7586\\
2 place Jussieu,\\
F-75251 Paris Cedex 05, France\\
\href{mailto:jean-yves.ducloux@u-pariscite.fr}{ducloux{@}math.jussieu.fr}\hfil
\end{minipage}
\par\vspace*{0.8cm}\noindent
\begin{minipage}[t]{.4\hsize}
{\obeylines\eightrm
Received \fontsize{9}{8}\selectfont September 22, 2000\vskip-2mm
\noindent and in final form June 28, 2001}
\end{minipage}

\end{document}